%% file: consolidated.tex
\documentclass[pdflatex,default,iicol]{sn-jnl}

\usepackage{csquotes}
\usepackage{lipsum}
\usepackage{graphicx}
\usepackage{dblfloatfix}
\usepackage{transparent}
\usepackage{amsmath,amssymb,amsfonts}
\usepackage{mathtools}
\usepackage{stmaryrd}
\usepackage{natbib}
\usepackage{cleveref}
\usepackage{xcolor}
\usepackage{multirow}%
\usepackage{amsthm}%
\usepackage{mathrsfs}%
\usepackage[title]{appendix}%
\usepackage{textcomp}%
\usepackage{manyfoot}%
\usepackage{booktabs}%
\usepackage{algorithm}%
\usepackage{algorithmicx}%
\usepackage{algpseudocode}
\usepackage{listings}%
\usepackage{lmodern}%
\usepackage{hyperref}
\usepackage{microtype}
\usepackage{tabulary}
\usepackage{tabularx}
\usepackage{multicol}
\usepackage{array}
\usepackage{bm}
\usepackage{siunitx}

\MakeOuterQuote{"}
\graphicspath{{./figures/}}
\AtBeginEnvironment{appendices}{\crefalias{section}{appendix}}

\newcommand{\closure}[1]{\overline{#1}}
\newcommand{\SetDef}[1]{\left\{ #1 \right\}}
\newcommand{\SetOfMaterials}{\mathcal{M}}

\renewcommand{\vec}[1]{\bm{\mathrm{#1}}}
\newcommand{\mat}[1]{\bm{\mathrm{#1}}}

\newcommand{\pfrac}[2]{\frac{\partial #1}{\partial #2}}

\newcommand{\abs}[1]{  \left\vert #1 \right\vert }
\newcommand{\jump}[1]{ \left[ \! \left[ #1 \right] \! \right] }
\newcommand{\Jump}[1]{ \left[ \kern-0.26em \left[ #1 \right] \kern-0.26em \right] }
\newcommand{\WeightedJump}[1]{ \left\{ #1 \right\} }
\newcommand{\norm}[2]{ \left\vert \! \left\vert #1 \right\vert \! \right\vert_{#2} }
\newcommand{\NormalDeriv}[1]{\partial^{#1}_{\normal}}

\newcommand{\AmbientDomain}{\mathcal{A}}
\newcommand{\Domain}{\Omega}

\newcommand{\InteriorDomain}{\mathring{\Domain}}
\newcommand{\MaterialDomain}{\Domain_{\MaterialIndex}}

\newcommand{\Boundary}[1]{\partial#1}
\newcommand{\BoundarySegment}{\Gamma}
\newcommand{\DirichletBoundary}{\BoundarySegment^D}
\newcommand{\NeumannBoundary}{\BoundarySegment^N}
\newcommand{\InterfaceBoundary}{\BoundarySegment^I}
\newcommand{\MaterialInterface}[2]{\InterfaceBoundary_{#1,#2}}
\newcommand{\InteriorInterfaces}{\mathring{\BoundarySegment}^I}

\newcommand{\normal}{\vec{n}}
\newcommand{\pos}{\vec{x}}
\newcommand{\displ}{\vec{u}}
\newcommand{\TestDispl}{\vec{v}}
\newcommand{\strain}{\vec{\epsilon}}

\newcommand{\stress}{\mat{\sigma}}

\newcommand{\traction}{\vec{t}}

\newcommand{\MaterialIndex}{m}

\newcommand{\Residual}{\mathcal{R}}

\newcommand{\PolyOrder}{p}

\newcommand{\BasisFunction}{N}
\newcommand{\BF}{\BasisFunction}
\newcommand{\BasisFunctionIndex}{B}
\newcommand{\BfIndex}{\BasisFunctionIndex}
\newcommand{\LocalBfIndex}{b}
\newcommand{\EnrBfIndex}{\ell}
\newcommand{\NumBFs}{n_{\BfIndex}}
\newcommand{\SetOfBfs}{\mathcal{B}}

\newcommand{\EnrichedBasisFunction}{\widetilde{\BF}}
\newcommand{\EnrBF}{\EnrichedBasisFunction}
\newcommand{\SetOfEnrBfs}{\widetilde{\SetOfBfs}}

\newcommand{\dof}{d}

\newcommand{\EnrichmentLevel}{\varepsilon}
\newcommand{\EnrLvl}{\EnrichmentLevel}
\newcommand{\NumEnrLvls}{n_{\EnrLvl}}
\newcommand{\IndicatorFunction}{\psi}
\newcommand{\EnrichmentFunction}{\IndicatorFunction^{\EnrLvl}_{\BfIndex}}
\newcommand{\EnrFnct}{\EnrichmentFunction}

\newcommand{\Unzipping}{u}
\newcommand{\NumUnzippings}{n_{\Unzipping}}

\newcommand{\BfSubPhase}[1]{R_{\BasisFunctionIndex}^{#1}}
\newcommand{\SubPhase}{S}
\newcommand{\ElemSubPhase}[1]{\SubPhase_{\BgElemInd}^{#1}}

\newcommand{\BackgroundElementIndex}{E}
\newcommand{\BgElemInd}{\BackgroundElementIndex}
\newcommand{\NumBgElems}{n_{\BgElemInd}}
\newcommand{\BackgroundElementDomain}{\Domain^{\BgElemInd}}
\newcommand{\BgElemDom}{\BackgroundElementDomain}

\newcommand{\Facet}{F}
\newcommand{\SetOfFacets}{\mathcal{\Facet}}
\newcommand{\SetOfGhostFacets}{\SetOfFacets_{G}}

\newcommand{\DiscreteSpace}{\mathcal{V}^{h}}

\newcommand{\DirichletPenalty}{\gamma_D}
\newcommand{\InterfacePenalty}{\gamma_I}
\newcommand{\GhostPenalty}{\gamma_G}

\newcommand{\DiscDispl}{\displ^{h}}
\newcommand{\DiscTestDispl}{\TestDispl^{h}}
\newcommand{\DiscDisplExt}{\tilde{\displ}^{h}}
\newcommand{\DiscTestDisplExt}{\tilde{\TestDispl}^{h}}
\newcommand{\DiscTrialStress}{\stress(\DiscDispl)}
\newcommand{\DiscTestStress}{\stress(\DiscTestDispl)}

\newcommand{\DiscTestStrain}{\strain(\DiscTestDispl)}

\newcommand{\Ancestor}{a}
\newcommand{\Cell}{\mathrm{c}}
\newcommand{\ChildMesh}{\mathrm{CM}}
\newcommand{\ListOfCMs}{\mathrm{CMS}}
\newcommand{\Cluster}{D}
\newcommand{\BgElem}{\BgElemInd}

\newcommand{\Entity}{e}
\newcommand{\FacetConn}{\mathcal{F}}
\newcommand{\Geometry}{G}
\newcommand{\ID}{I}
\newcommand{\Material}{\MaterialIndex}
\newcommand{\Ordinal}{o}
\newcommand{\Processor}{p}
\newcommand{\Proc}{\Processor}
\newcommand{\LocalCommTable}{\mathrm{P_{loc}}}

\newcommand{\Proximity}{P}
\newcommand{\Queue}{\mathrm{Q}}
\newcommand{\Rank}{r}
\newcommand{\Subphase}{\SubPhase}
\newcommand{\Vertex}{\mathrm{v}}
\newcommand{\IEN}{\mathrm{IEN}}

\newcommand{\BackgroundMesh}{\mathcal{H}}
\newcommand{\BgMesh}{\BackgroundMesh}
\newcommand{\ForegroundMesh}{\mathcal{T}}
\newcommand{\FgMesh}{\ForegroundMesh}

\newcommand*{\Interface}{\Gamma}

\newcommand*{\SubphaseGraph}{\mathcal{G}}
\newcommand*{\SpGraph}{\SubphaseGraph}

\begin{document}


\title[Enriched Immersed Finite Element and Isogeometric Analysis -- Algorithms and Data Structures]{Enriched Immersed Finite Element and Isogeometric Analysis -- Algorithms and Data Structures}

\author[1]{\fnm{Nils} \sur{Wunsch}}\email{nils.wunsch@colorado.edu}
\equalcont{These authors contributed equally to this work.}
\author[1]{\fnm{Keenan} \sur{Doble}}\email{keenan.doble@colorado.edu}
\equalcont{These authors contributed equally to this work.}
\author[2]{\fnm{Mathias R.} \sur{Schmidt}}\email{schmidt43@llnl.gov}
\author[3]{\fnm{Lise} \sur{No\"el}}\email{l.f.p.noel@tudelft.nl}
\author[1]{\fnm{John A.} \sur{Evans}}\email{john.a.evans@colorado.edu}
\author*[1]{\fnm{Kurt} \sur{Maute}}\email{kurt.maute@colorado.edu}

\affil[1]{ \orgdiv{Smead Aerospace Engineering Sciences}, \orgname{University of Colorado Boulder}, \\ \orgaddress{ \street{3775 Discovery Dr.}, \city{Boulder}, \postcode{80303}, \state{CO}, \country{USA}}}

\affil[2]{ \orgdiv{Computational Engineering Division}, \orgname{Lawrence Livermore National Laboratory}, \\ \orgaddress{ \street{7000 East Ave.}, \city{Livermore}, \postcode{94550}, \state{CA}, \country{USA}}}

\affil[3]{ \orgdiv{Department of Precision and Microsystems Engineering}, \orgname{Delft University of Technology}, \orgaddress{ \street{Mekelweg 2}, \city{Delft}, \postcode{2628 CD}, \country{The Netherlands}}}


\abstract{
Immersed finite element methods provide a convenient analysis framework for problems involving geometrically complex domains, such as those found in topology optimization and microstructures for engineered materials. However, their implementation remains a major challenge due to, among other things, the need to apply nontrivial stabilization schemes and generate custom quadrature rules. 
This article introduces the robust and computationally efficient algorithms and data structures comprising an immersed finite element preprocessing framework.
The input to the preprocessor consists of a background mesh and one or more geometries defined on its domain. 
The output is structured into groups of elements with custom quadrature rules formatted such that common finite element assembly routines may be used without or with only minimal modifications.
The key to the preprocessing framework is the construction of material topology information, concurrently with the generation of a quadrature rule,
which is then used to perform enrichment and generate stabilization rules.  
While the algorithmic framework applies to a wide range of immersed finite element methods using different types of meshes, integration, and stabilization schemes, the preprocessor is presented within the context of the extended isogeometric analysis. This method utilizes a structured B-spline mesh, a generalized Heaviside enrichment strategy considering the material layout within individual basis functions' supports, and face-oriented ghost stabilization.
Using a set of examples, the effectiveness of the enrichment and stabilization strategies is demonstrated alongside the preprocessor's robustness in geometric edge cases. Additionally, the performance and parallel scalability of the implementation are evaluated.
}

\maketitle


\onecolumn
\newpage
\tableofcontents
\twocolumn


\section{Introduction}
\label{sec:introduction}

\begin{figure}[btp]
    \centering
    \includegraphics[width = 7.5cm]{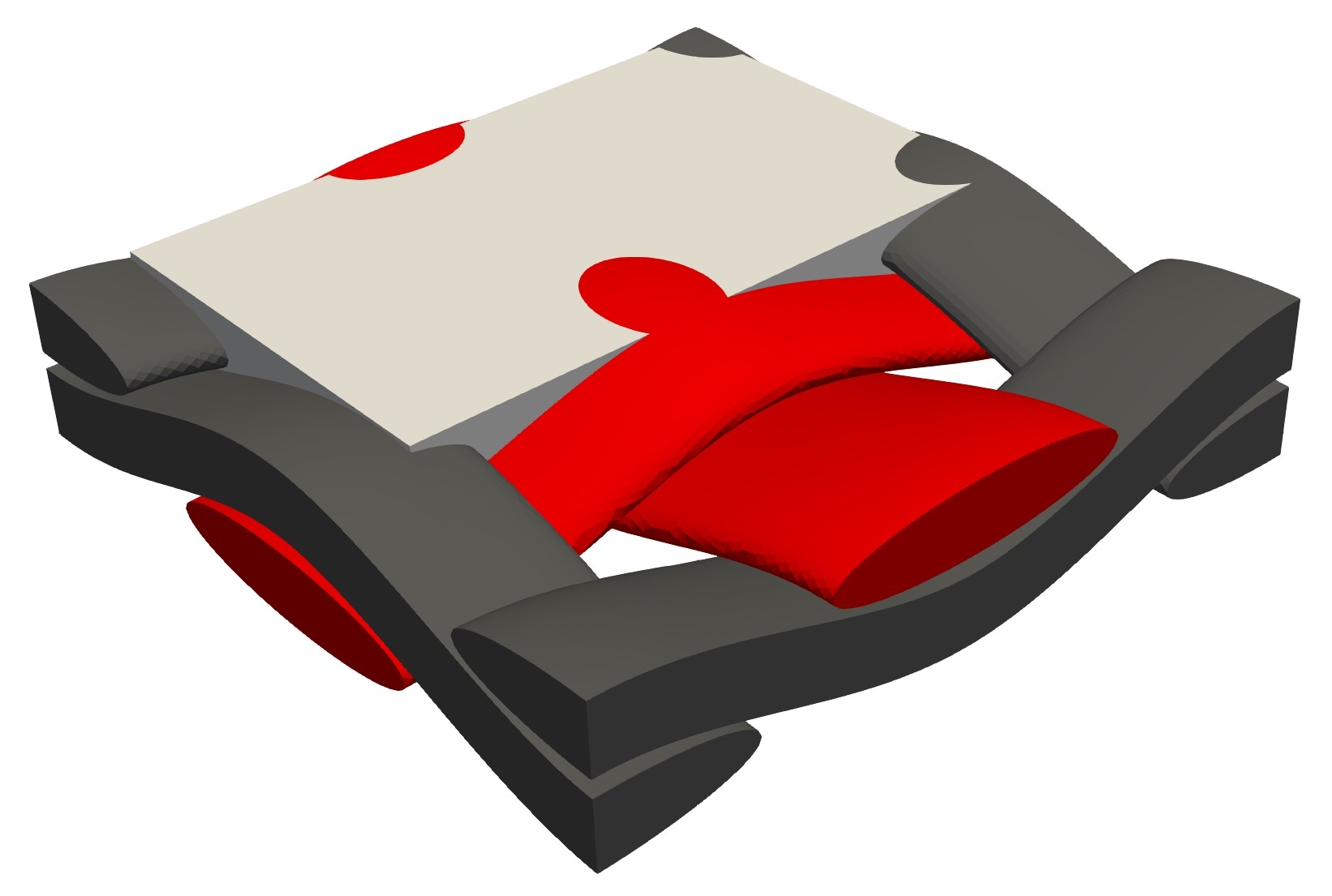}
    \caption{Fiber composite weave with two types of fibers and a binding matrix material. The matrix material is shown in off-white and partially removed for better visibility.}
    \label{fig:fiber_patch_3D}
\end{figure}

The finite element method (FEM) \cite{hughes2012finite} has been adopted as a standard tool to solve continuum mechanics problems in engineering analysis. In its traditional form, the FEM requires boundary-fitted meshes consisting of elements that define both the approximation of the state variable fields, e.g. temperature or displacement, and the geometric domain over which the weak statement of the governing partial differential equations (PDEs) is integrated. To obtain high-fidelity solutions, the finite element (FE) function spaces need to be of sufficient quality which may result in nontrivial geometric constraints on the elements comprising the mesh. For complex domain shapes, such as the fiber composite shown in \Cref{fig:fiber_patch_3D}, boundary-fitted meshes satisfying these quality measures can generally not be generated reliably without manual intervention. This renders mesh generation a major bottleneck in computational engineering analysis \cite{cottrell2009isogeometric,bazilevs2010isogeometric}. 

"Unfitted" \cite{hansbo2002unfitted} or "immersed" \cite{zhang2004immersed} FEMs have gained significant traction in recent years as an elegant solution to automate finite element analysis. The commonly used methods for immersed finite element analysis include the CutFEM \cite{burman2015cutfem} (which originated from partition of unity (PU) based methods \cite{claus2015cutfem}, namely, the PU-FEM \cite{melenk1996partition, babuvska1997partition}), the generalized FEM \cite{duarte2000generalized, duarte2001generalized, soghrati2012interface, soghrati20123d}, and the extended FEM \cite{moes1999finite, belytschko1999elastic}. Other lineages of immersed methods include the finite cell method \cite{parvizian2007finite, duster2008finite, schillinger2015finite} (which originated from fictitious domain methods) and the shifted boundary method \cite{main2018shifted1, main2018shifted2}.
These methods are specifically designed to accommodate arbitrary strong or weak discontinuities within elements, representing features such as cracks, domain boundaries, or material interfaces. As a result, a mesh, often referred to as the "background mesh," can be constructed on a geometrically simple domain into which a geometrically complex domain is immersed.
The immersed approach is of particular interest and has been successfully applied to problems where the meshing effort is exacerbated by complex or evolving interface geometries such as in scan-based analysis \cite{ruess2012finite, verhoosel2015image, elhaddad2018multi}, topology optimization \cite{van2013level, groen2017higher, noel2017shape, villanueva2017cutfem, noel2023xfem}, or physical phenomena such as multi-phase flows \cite{sauerland2013stable} and fluid-structure interaction \cite{hsu2014fluid,hsu2015dynamic,kamensky2015immersogeometric}.

\begin{figure*}[h]
    \centering 
    \def\svgwidth{13.0cm}
    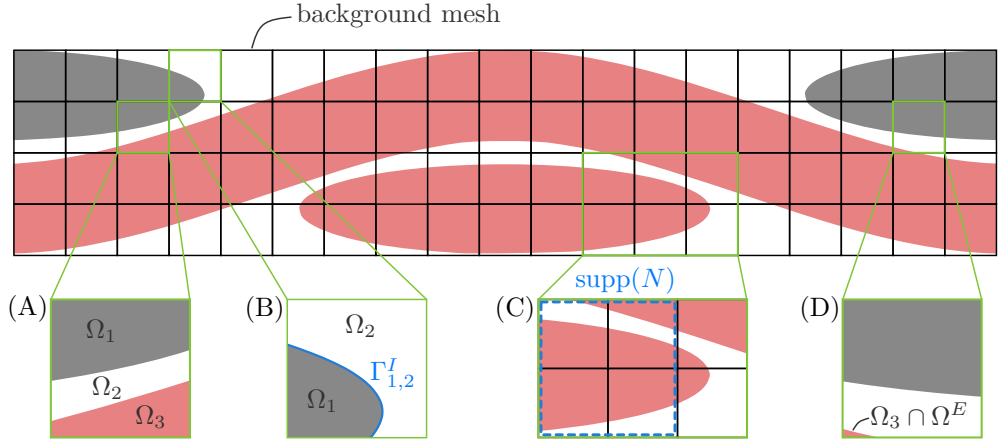
    \caption{Cross-section of the fiber weave shown in \Cref{fig:fiber_patch_3D} with a coarse background mesh. Details: (A) Cut background element occupied by three different materials. (B) Non-trivially intersected background element with highlighted material interface $\Gamma^I_{1,2}$. (C) Sub-element scale material features in neighboring elements. Highlighted is a single basis function's support $\mathrm{supp}(N)$ which spans multiple features. (D) Cut background element where one material subdomain occupies only a small volume fraction.}
    \label{fig:fiber_patch_cross_section}
\end{figure*}

While immersed methods provide a convenient and time-saving design-through-analysis workflow, the approach departs from traditional FEMs in four key aspects. To illustrate these non-standard aspects, consider the cross-section of a fiber composite immersed into a background mesh, as shown in \Cref{fig:fiber_patch_cross_section}.
\begin{enumerate}
\item[(i)] Multiple materials, components, or other features like cracks may be present within a cut element, as shown in \Cref{fig:fiber_patch_cross_section} detail (A). The finite element basis needs to resolve the associated discontinuities inside the background elements. This is commonly done through basis function enrichment \cite{fries2010extended,schweitzer2012generalizations}. 
\item[(ii)] Essential boundary and interface conditions cannot be built into the test and trial function spaces. Instead, these conditions require weak enforcement through the penalty method \cite{zhu1998modified}, Lagrange-multiplier methods \cite{burman2010fictitious}, or Nitsche's method \cite{burman2012fictitious}. For an overview of the issue, the reader is referred to \cite{fernandez2004imposing}.
\item[(iii)] The governing equations must be satisfied on given material domains $\MaterialDomain$. However, the background elements 
do not coincide with these material regions. As a result, these elements cannot be reused to generate quadrature rules to evaluate the integrals of the associated weak form. Instead, custom quadrature rules need to be constructed for each background element intersected by a material interface $\Interface^I_{m,l}$ and the material subdomains within them. 
Furthermore, quadrature rules must be constructed for segments of the interface $\Interface^I_{m,l}$ within background elements, such as the one in \Cref{fig:fiber_patch_cross_section} detail (B), to evaluate the integrals enforcing the boundary conditions.
\item[(iv)] Lastly, the arbitrary shape and location of an interface within a background element may lead to small volume fractions $\mathrm{vol}(\MaterialDomain \cap \BgElemDom) \ll \mathrm{vol}(\BgElemDom)$ within that background element's domain $\BgElemDom$, as shown in \Cref{fig:fiber_patch_cross_section} detail (D). Such configurations can lead to stability and conditioning issues and require stabilization. For an overview of this issue, the reader is referred to \cite{de2023stability}.
\end{enumerate}

The generation of quadrature rules and the construction of the enriched finite element space, as well as the application of appropriate stabilization schemes, require extensive supplemental geometric information that needs to be generated during preprocessing. These requirements lead to a challenging implementation of immersed methods.
Multiple approaches have been developed to individually address each of the four aspects. However, the literature mostly presents these methods from a theoretical standpoint or covers the implementation of individual methods independently of how other aspects are addressed.
The implementation of immersed finite elements and the associated methods is indirectly addressed by publically available source code. Examples of published open-source FE packages with immersed FE capabilities are \texttt{ngsxfem} \cite{lehrenfeld2021ngsxfem}, \texttt{deal.II} \cite{dealII95}, \texttt{MOOSE-XFEM} \cite{zhang2018modified, jiang2020ceramic}, and the discontinued libraries \texttt{LibCutFEM} \cite{burman2015cutfem} and \texttt{MultiMesh} \cite{johansson2019multimesh}.
An exception is the work by Zhang et al. \cite{zhang2022object}, which provides an overview of the algorithms and data structures underlying an enriched/immersed FE preprocessor and discusses the interaction of various components of the framework. However, it does not cover the stabilization aspect or the realization of a parallel implementation.
The authors believe that there is significant value in the discussion of the implementation of an immersed FE framework as not only is the translation from theory into code challenging, but so is doing so efficiently.

This article presents a comprehensive implementation of a preprocessor that generates all data necessary to perform immersed finite element analysis on multi-material problems, including the construction of enriched finite element spaces, as well as the generation of custom quadrature and stabilization rules. We discuss the algorithmic design choices made and their impact on robustness and performance characteristics, specifically runtime, memory consumption, and parallel scalability. A further focus is placed on how the four key aspects described previously are addressed cohesively.
Although the implementation is presented in the context of the extended isogeometric analysis (XIGA) introduced by \cite{noel2022xiga} and \cite{schmidt2023extended}, the algorithms and data structures are general and adaptable to suit other immersed analysis frameworks that may utilize different types of background meshes, geometry descriptions, or use alternative methods to generate quadrature and stabilization rules. The preprocessor can also be used for interpolation-based immersed analysis \cite{fromm2023interpolation, fromm2024interpolation}.
The implementation, as presented in this paper, is part of an immersed FE and optimization package \texttt{Moris}, whose \texttt{C++} source code is also \href{https://github.com/kkmaute/moris/}{publicly available} \cite{moris}. The framework has been developed with a focus on solving multi-material level-set (LS) topology optimization problems, where the shape of objects and their material composition are defined implicitly by a set of LS functions. Such functions can be created, e.g., from images or 3D scans or via geometry modeling tools using function representations as done by \texttt{nTop} \cite{vlahinos2020unlocking} and \texttt{openVCAD} \cite{wade2023openvcad}.

The XIGA employs a generalized Heaviside enrichment strategy in combination with Nitsche's method \cite{nitsche1971variationsprinzip} to resolve material interfaces and the associated discontinuities inside background elements.
Most Heaviside enrichment strategies consider either the material domains or their globally connected subdomains \cite{terada2003finite, hansbo2004finite}, i.e., material components, to construct the enriched function space. However, configurations with features on a length scale smaller than a background element may arise, e.g., in the context of topology optimization or engineered materials. This leads to numerical artifacts due to basis functions that span features of the same globally connected material domain but disconnected within that function's support, as shown in \Cref{fig:fiber_patch_cross_section} detail (C). This issue is exacerbated when basis functions with larger supports are employed, such as high-continuity splines. To this end, the enrichment strategy employed herein considers the material layout on the level of the individual basis function supports \cite{noel2022xiga}. 

The information about the material layout necessary to apply the generalized Heaviside enrichment strategy is obtained by building connectivity graphs after the tessellation of the background mesh. The tessellation procedure employed relies on the application of subdivision templates as done by, e.g., \cite{soghrati2014hierarchical}. Custom quadrature rules on the various material subdomains are generated by mapping a Gauss-Legendre quadrature rule from the resulting boundary-fitted "foreground" mesh onto the background elements. Alternative, potentially higher-order accurate, quadrature methods and tessellation procedures may be employed, such as, e.g., the ones presented in \cite{badia2022geometrical} and \cite{martorell2024high} for explicit geometry descriptions or the one in \cite{fries2016higher} for implicit geometry descriptions, so long as the necessary arrangements to build a material connectivity graph are made.

Lastly, the stabilization aspect can be addressed using methods such as basis function removal \cite{elfverson2018cutiga}, basis function extension approaches \cite{badia2018aggregated, burman2022cutfem, burman2023extension}, or ghost stabilization \cite{burman2010ghost,burman2014fictitious}. With the information about the material layout within intersected background elements previously obtained, either of these stabilization methods can be applied with relative ease. For this work, a face-oriented ghost stabilization \cite{schott2014new} is chosen as, unlike the other methods, the ghost stabilization does not rely on non-intersected elements in the vicinity of an interface, resulting in better robustness in pathological geometric configurations.

The remainder of this article is structured as follows. First, \Cref{sec:framework} lays out the general challenges of an enriched immersed method and the mathematical formulations of the enrichment and ghost stabilization strategies. Using this knowledge,  the output information, which the preprocessor needs to generate, is formalized.
\Cref{sec:main} covers the algorithms and data structures underlying the implementation of the preprocessing framework.
Lastly, \Cref{sec:examples} demonstrates robustness, efficiency, and parallel scalability using a suite of examples.

\clearpage
\section{Numerical Framework}
\label{sec:framework}

This section discusses the theoretical aspects of the computational framework and focuses on each of the framework's components, which were first introduced in \cite{noel2022xiga}, with the intent of providing the requirements for its implementation. Some details are omitted for brevity; the reader is referred to \cite{noel2022xiga} for a detailed exposure to the underlying theory.
We start the discussion by first presenting the discretized weak form for the fiber patch shown in \Cref{fig:fiber_patch_3D} as a linear elastic model problem and then address each of the aspects necessary for constructing the information to assemble the system of equations within the given framework.

\subsection{Discretization}
\label{sec:framework-discretization}

\paragraph*{Model Problem and Weak Form}

The fiber patch, shown in \Cref{fig:fiber_patch_3D}, is composed of three materials $m \in \SetOfMaterials = \SetDef{1,2,3}$ (red fibers, gray fibers, and matrix material), all modeled as linear elastic, but with different material parameters. 
The displacements at the interfaces are assumed to be matching.
The patch is subjected to external loads and displacement constraints. For convenience, we define the union $\InteriorDomain$ of domains occupied by each material $\MaterialDomain$, i.e., the physical domain without the internal interfaces, as
\begin{equation}
    \label{eqn:interior_domain}
    \InteriorDomain = \bigcup_{\MaterialIndex \in \SetOfMaterials } \MaterialDomain.
\end{equation}
The discretized weak form for the problem with discrete trial and test displacement $\DiscDispl$ and $\DiscTestDispl$ states as: find $\DiscDispl \in \DiscreteSpace \subset H^1(\InteriorDomain)$ such that, $\forall \DiscTestDispl \in \DiscreteSpace \subset H^1(\InteriorDomain)$: 
\begin{equation}
    \begin{split}
    \label{eqn:weak_form_discrete}
    \Residual_{\Domain}(\DiscDispl, \DiscTestDispl) + 
    \Residual_{N}(\DiscDispl, \DiscTestDispl) + 
    \Residual_{D}(\DiscDispl, \DiscTestDispl) + 
    \cdots \\ \cdots +
    \Residual_{I}(\DiscDispl, \DiscTestDispl) + 
    \Residual_{G}(\DiscDispl, \DiscTestDispl) = 0,
    \end{split} 
\end{equation}
where
\begin{align}
    \label{eqn:weak_form_bulk}
    \Residual_{\Domain}(\DiscDispl, \DiscTestDispl) &= 
    \sum_{m \in \SetOfMaterials}
    \int_{\MaterialDomain} \DiscTestStrain : \DiscTrialStress \, d\Domain, \\
    \label{eqn:weak_form_neumann}
    \Residual_{N}(\DiscDispl, \DiscTestDispl) &=
    \sum_{m \in \SetOfMaterials}
    \int_{\NeumannBoundary_{\MaterialIndex}} - \DiscTestDispl \cdot \traction \, d\BoundarySegment,
\end{align}
\begin{align}
    \label{eqn:weak_form_dirichlet_Nitsche}
    & \Residual_{D}(\DiscDispl, \DiscTestDispl) = \nonumber
    \\
    & \sum_{m \in \SetOfMaterials} \left[
    - \int_{\DirichletBoundary_{\MaterialIndex}} \DiscTestDispl \cdot \left( \DiscTrialStress \cdot \normal_m \right) \, d\BoundarySegment \right. \nonumber
    \\
    & \qquad \cdots \mp \int_{\DirichletBoundary_{\MaterialIndex}} \left( \DiscTestStress \cdot \normal_m \right) \cdot (\DiscDispl - \displ_0 )\ d\BoundarySegment \nonumber 
    \\
    & \qquad \left. \cdots + \int_{\DirichletBoundary} \DirichletPenalty \, \DiscTestDispl \cdot (\DiscDispl - \displ_0 ) \, d\BoundarySegment
    \right], \\
    \label{eqn:weak_form_interface_Nitsche}
    & \Residual_{I}(\DiscDispl, \DiscTestDispl) = \nonumber
    \\
    & \sum_{ \substack{ m,l \in \SetOfMaterials \\ m \neq m } } \left[
    - \int_{\MaterialInterface{m}{l}} \Jump{ \DiscTestDispl } \cdot \WeightedJump{ \DiscTrialStress \cdot \normal_m} \, d\BoundarySegment \right. \nonumber
    \\
    & \qquad \cdots \mp \int_{\MaterialInterface{m}{l}} \WeightedJump{\DiscTestStress \cdot \normal_i} \cdot \Jump{\DiscDispl} \, d\BoundarySegment \nonumber 
    \\
    & \qquad \left. \cdots + \int_{\MaterialInterface{m}{l}} \InterfacePenalty \Jump{\DiscTestDispl} \cdot \Jump{\DiscDispl} \, d\BoundarySegment
    \right].
\end{align}
\eqref{eqn:weak_form_bulk} and \eqref{eqn:weak_form_neumann} are the bulk and Neumann terms, originating from the strong form and integration by parts. Here $\strain$, $\stress$, and $\traction$ denote the strain, stress, and prescribed traction, respectively; $\normal_m$ denotes the outward unit normal to the interface associated with material domain $\MaterialDomain$. The terms used to weakly enforce the Dirichlet boundary conditions via Nitsche's method \cite{nitsche1971variationsprinzip, burman2012fictitious} are collected in \eqref{eqn:weak_form_dirichlet_Nitsche}. 
The contribution \eqref{eqn:weak_form_interface_Nitsche} weakly enforces matching displacements and tractions at interfaces between materials 
\begin{equation}
    \label{eqn:interface_boundary}
    \MaterialInterface{m}{l} = 
    \Boundary{\Domain_{m}} \cap \Boundary{\Domain_{l}} 
    \quad \text{for} \quad m, l \in \SetOfMaterials, m \neq l,
\end{equation}
where $\Boundary{(\cdot)}$ denotes the boundary of a domain.
The braces and double brackets are jump terms defined as
\begin{align}
    \label{eqn:weighted_jump}
    \WeightedJump{f(\pos)} &= 
    \lim_{\delta \to 0^{+}}
    w_{m} \, f(\pos - \delta \, \normal_m) + w_{l} \, f(\pos + \delta \, \normal_m),
    \\
    \label{eqn:jump_operator}
    \jump{f(\pos)} &= 
    \lim_{\delta \to 0^{+}}
    f(\pos - \delta \, \normal_m) - f(\pos + \delta \, \normal_m),
\end{align}
with $w_{m} + w_{l} = 1$ and $\pos \in \MaterialInterface{m}{l}$.
The choice of weights $w_{m}$ and $w_{l}$, as well as the penalty parameters $\DirichletPenalty$ and $\InterfacePenalty$ are adopted from \cite{annavarapu2012robust}. Further, the signs of the symmetry terms in \eqref{eqn:weak_form_dirichlet_Nitsche} and \eqref{eqn:weak_form_interface_Nitsche} are chosen to yield either a symmetric or a non-symmetric formulation. We briefly discuss the consequences of this choice in \Cref{sec:framework-ghost} alongside the ghost term $\Residual_{G}$.

Note that every material subdomain $\MaterialDomain$ is treated separately to allow for, e.g., the use of different constitutive laws or physics for each material, or allow for a wider range of interface conditions to be enforced via the term $\Residual_{I}$. 
Materials may be chosen to be void, thus not contributing to the residual, and some interfaces $\MaterialInterface{m}{l}$ become part of the Dirichlet or Neumann boundaries.

\paragraph*{Approximation Spaces and Enrichment}

The finite element space $\DiscreteSpace = \mathrm{span}(\SetOfBfs)$ is provided by a set of basis functions $\SetOfBfs = \SetDef{ N_{\BfIndex}(\pos) }_{\BfIndex = 1}^{\NumBFs}$ with
\begin{equation}
    \label{eqn:interpolation}
    \DiscDispl(\pos) = \sum_{\BfIndex = 1}^{\NumBFs} \BasisFunction_{\BfIndex}(\pos) \, \dof_{\BfIndex},
\end{equation}
where $\dof_{\BfIndex}$ are the weights associated with each basis function $\BasisFunction_{\BfIndex}$. Throughout the remainder of this paper, we are using $B$ to index basis functions. Further, $n_{(\cdot)}$ denotes the number of whatever the subscript index is associated with.
The basis functions are defined on the elements of the background mesh $\BackgroundMesh = \SetDef{ \BgElemDom }_{\BgElemInd = 1}^{\NumBgElems}$ which does not conform to interfaces, and we assume the basis functions to be compactly supported by a set of background elements, without making any further assumptions on the type of basis function used at this time.

\begin{figure}[b]
    \centering 
    \def\svgwidth{4.8cm}
    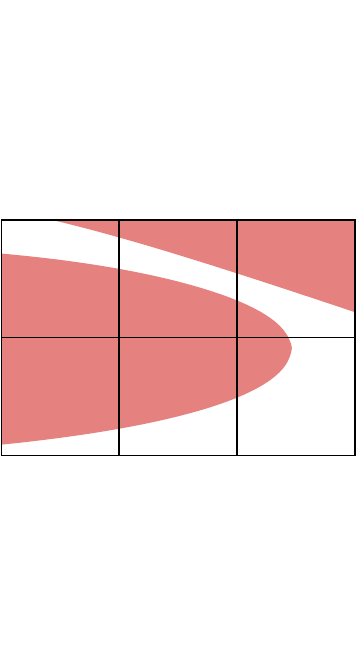
    \caption{Enrichment strategy for a single basis function for the material layout shown in \Cref{fig:fiber_patch_cross_section}, detail (C). The basis function's support is shown in dashed blue lines. The support contains four disconnected material subdomains $\BfSubPhase{\EnrLvl}$ for each of which an enriched basis function is defined, where $\EnrBF_{\BfIndex}^{\EnrLvl} = \EnrFnct \BF_{\BfIndex}$ with $\mathrm{supp}(\EnrBF_{\BfIndex}^{\EnrLvl}) = \BfSubPhase{\EnrLvl}$ }
    \label{fig:enrichment_strategy}
\end{figure}

To restrict the support of basis functions to the material domains and allow for the enforcement of arbitrary interface conditions, a Heaviside enrichment strategy is employed. The specific strategy used in this framework considers different material subdomains within a single basis function's support, as illustrated in \Cref{fig:enrichment_strategy}. This prevents artificial numerical coupling between small geometry features 
without the need for local refinement, as demonstrated in \Cref{sec:examples-multi_beam}.

With this enrichment strategy, the approximation of the state variables reads as follows:
\begin{equation}
    \label{eqn:enriched_interpolation}
    \DiscDispl(\pos) = 
    \sum_{\BfIndex = 1}^{\NumBFs} 
    \sum_{\EnrLvl = 1}^{\NumEnrLvls(\BfIndex)}
    \EnrichmentFunction(\pos) \, 
    \BasisFunction_{\BfIndex}(\pos) \, 
    \dof_{\BfIndex}^{\EnrLvl},
\end{equation}
where we define the indicator, or "enrichment", function $\EnrichmentFunction$ as
\begin{equation}
    \label{eqn:enrichment_function}
    \EnrichmentFunction(\pos) = 
    \begin{cases}
        1, & \pos \in \BfSubPhase{\EnrLvl}, \\
        0, & \text{otherwise}.
    \end{cases}
\end{equation}
Here $\BfSubPhase{\EnrLvl} \subseteq \mathrm{supp}(\BF_{\BfIndex})$ with 
$ \bigcup_{\EnrLvl=1}^{\NumEnrLvls(\BfIndex)} \BfSubPhase{\EnrLvl} = \mathrm{supp}(\BF_{\BfIndex}) \cap \InteriorDomain$ 
correspond to the disconnected subdomains within the basis function's support, as is shown in \Cref{fig:enrichment_strategy}.

\begin{figure*}[btp]
    \vspace*{0.3cm}
    \centering 
    \def\svgwidth{13.0cm}
    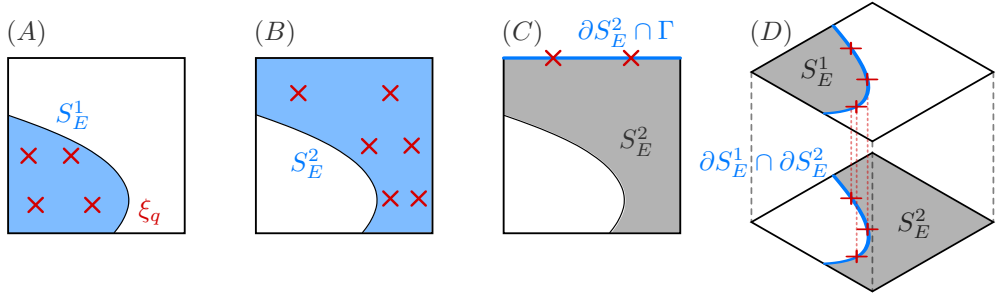
    \caption{Concept of a cluster: a background element with a single active subdomain $\ElemSubPhase{\Unzipping}$, a set of enriched basis functions supported inside that subdomain, and a quadrature rule associated with either that subdomain (A, B) or any part of the subdomain's boundary $\Boundary{\ElemSubPhase{\Unzipping}}$ (C). Green indicates the volume or surface integrated over; blue crosses represent exemplary quadrature points. Additionally, linked pairs of clusters are constructed for interfaces (D).}
    \label{fig:cluster_concept}
\end{figure*}

We refer to the product $\EnrichedBasisFunction_{\BfIndex}^{\EnrLvl} = \EnrichmentFunction \cdot \BasisFunction_{\BfIndex}$ as an "enriched basis function". The set of all enriched basis functions is defined as
\begin{align}
    \SetOfEnrBfs := \{ \, \EnrichmentFunction(\pos) \BF_{\BfIndex}(\pos) : ~
    & \BfIndex \in \SetDef{1, \dots, \NumBFs}, \nonumber 
    \\
    & \EnrichmentLevel \in \SetDef{1, \dots, \NumEnrLvls(\BfIndex)} \, \} 
    \label{eqn:enriched_basis_function_set}
\end{align}

From this arises the requirement on the implementation to identify the number of enrichment levels $\NumEnrLvls(\BfIndex)$ for each basis function $\BasisFunction_{\BfIndex}$ and construct their enriched counterparts $\EnrichedBasisFunction_{\BfIndex}^{\EnrLvl}$. As is detailed in \Cref{sec:enrichment}, a flood-fill algorithm is used to identify the disconnected subdomains. This approach subsequently requires that the material connectivity is available.

\subsection{Element Formation and Integration}
\label{sec:framework-ig}

In a body-fitted finite element approach, each finite element has a set of non-zero basis functions supported within it and belongs to part of a material domain. The integrals \eqref{eqn:weak_form_bulk} -\eqref{eqn:weak_form_interface_Nitsche} can hence be written as a sum of integrals over each element. The evaluation of these elemental integrals is referred to as "element formation". In contrast, for an immersed approach with multiple material domains, the elemental integrals need to be further split into the different material or interface subdomains within a given background element. Taking into account the enrichment strategy \eqref{eqn:enriched_interpolation}, we specifically consider materially disconnected subdomains within a background element -- rather than material subdomains -- such that it can be ensured that each subdomain supports a set of non-zero enriched basis functions. 
The generalized form of the bulk, boundary, and interface integrals \eqref{eqn:weak_form_bulk} -\eqref{eqn:weak_form_interface_Nitsche} can then be stated as
\begin{align}
    \label{eqn:elemental_bulk_integral}
    \sum_{\MaterialIndex \in \SetOfMaterials}
    \int_{\MaterialDomain} (\cdot) \, d\Domain &= 
    \sum_{\BgElemInd = 1}^{\NumBgElems} \
    \sum_{\Unzipping = 1}^{\NumUnzippings(\BgElemInd)}
    \int_{\ElemSubPhase{\Unzipping} } (\cdot) \, d\Domain,
    \\
    \label{eqn:elemental_boundary_integral}
    \sum_{\MaterialIndex \in \SetOfMaterials}
    \int_{\BoundarySegment_{\MaterialIndex}} (\cdot) \, d\Domain &= 
    \sum_{\BgElemInd = 1}^{\NumBgElems} \
    \sum_{\Unzipping = 1}^{\NumUnzippings(\BgElemInd)}
    \int_{\Boundary{\ElemSubPhase{\Unzipping}} \cap \BoundarySegment} (\cdot) \, d\BoundarySegment, 
    \\
    \label{eqn:elemental_interface_integral}
    \sum_{ \substack{ m,l \in \SetOfMaterials \\ m \neq l } }
    \int_{\InterfaceBoundary_{m,l}} (\cdot) \, d\Domain &= 
    \sum_{\BgElemInd = 1}^{\NumBgElems} \
    \sum_{ \substack{ \Unzipping,v = 1 \\ u \neq v } }^{\NumUnzippings(\BgElemInd)}
    \int_{\Boundary{\ElemSubPhase{\Unzipping}} \cap \Boundary{\ElemSubPhase{v}}} (\cdot) \, d\BoundarySegment.
\end{align}
The background elements are indexed by $\BgElemInd$, and we correspondingly denote the number of background elements by $\NumBgElems$. The number of disconnected subdomains $\ElemSubPhase{\Unzipping}$ within a given background element's domain $\BgElemDom$ is denoted by $\NumUnzippings(\BgElemInd)$. 
From an implementation standpoint, the challenge of evaluating the integrals is twofold. The first challenge is to identify the disconnected subdomains $\SetDef{\ElemSubPhase{\Unzipping}}_{\Unzipping = 1}^{\NumUnzippings(\BgElemInd)}$ within each background element $E$ and the enriched basis functions supported within each subphase $\ElemSubPhase{\Unzipping}$. 
The second challenge is to construct custom quadrature rules for each subdomain to evaluate the integrals in \eqref{eqn:elemental_bulk_integral} - \eqref{eqn:elemental_interface_integral}. 

\begin{figure*}[btp]
    \centering 
    \def\svgwidth{13.0cm}
    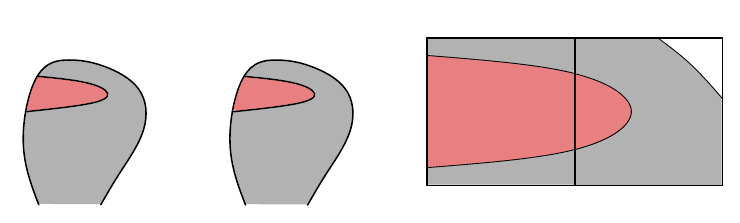
    \caption{Ghost faces for an example material layout. (A) Material layout of a two-material problem with materials $\SetOfMaterials = \SetDef{1,2}$ occupying a domain smaller than the ambient domain $\AmbientDomain$ filled by the background mesh. The white domain $\Domain_{0}$ is assumed to be void. (B) The set of interior element faces $\SetOfGhostFacets$ which the ghost penalty is applied to. (C) Exemplary ghost facet $\Facet$ and the material layout in its adjacent background elements.}
    \label{fig:ghost_for_enriched}
\end{figure*}

Element formation in the immersed setting is generalized for any type of integral by introducing the concept of a cluster, as shown in \Cref{fig:cluster_concept}. We define a cluster to be a background element with domain $\BgElemDom$, a single subdomain $\ElemSubPhase{\Unzipping}$, the set of enriched basis functions supported inside that subdomain, and a quadrature rule associated with either that subdomain's volume or any part of its boundary $\Boundary{\ElemSubPhase{\Unzipping}}$. The quadrature rules themselves consist of a set of quadrature points and associated weights. In the case of interface integrals, a linked pair of clusters is needed to evaluate the jump terms \eqref{eqn:weighted_jump} and \eqref{eqn:jump_operator}. The step of element formation is then reduced to evaluating integrands using the basis defined on the cluster (pair)'s background element(s) at the quadrature points and summing up the contributions. 

Within the presented framework, the disconnected subdomains are obtained by a tessellation of intersected background elements using a fast, templated subdivision scheme detailed in \Cref{sec:fg_mesh_generation}. 

\subsection{Face-Oriented Ghost Stabilization}
\label{sec:framework-ghost}

Small material slivers lead to ill-conditioning due to the vanishing contributions of sparsely supported basis functions to the global system of equations \cite{de2023stability}. Further, stability issues may arise due to the penalty terms in \eqref{eqn:weak_form_dirichlet_Nitsche} and \eqref{eqn:weak_form_interface_Nitsche} not being any longer bounded for arbitrarily small material fractions\footnote{The stability parameters $\gamma$ in \eqref{eqn:weak_form_dirichlet_Nitsche} and \eqref{eqn:weak_form_interface_Nitsche} need to scale inversely with elemental material volumes (which can be arbitrarily small) to retain coercivity and, hence, well-posedness for the symmetric Nitsche formulation \cite{de2018note}. An asymmetric version of Nitsche's method avoids this issue at the cost of losing symmetry in the linear system of equations and adjoint consistency of the bilinear form resulting in potentially reduced $L^2$-convergence \cite{schillinger2016non}.}. 

The framework employs a version of the face-oriented ghost stabilization proposed by \cite{burman2014fictitious} adapted to the enrichment strategy used. 
The basic idea underlying face-oriented ghost stabilization is to penalize the jumps in derivatives across elemental faces adjacent to intersected background elements, thereby ensuring a sufficient contribution of sparsely supported basis functions. The penalties are only applied to those faces across which connected material regions transition, to not introduce spurious numerical coupling. 

To illustrate the adaptions necessary for the enrichment strategy, consider the material layout shown in \Cref{fig:ghost_for_enriched} (C). The background elements $\BgElemInd^+$ and $\BgElemInd^-$ contain multiple material subdomains $\SetDef{\SubPhase_{\BgElemInd^+}^u}_{u = 1}^{\NumUnzippings(\BgElemInd^+)}$ and $\SetDef{\SubPhase_{\BgElemInd^-}^u}_{u = 1}^{\NumUnzippings(\BgElemInd^-)}$ with different sets of enriched basis functions supported within each subdomain.
Additionally, the subdomains $\SubPhase_{\BgElemInd^+}^{1}$ and $\SubPhase_{\BgElemInd^+}^{3}$ support different sets of enriched basis functions but are part of the same material domain $\Domain_{1}$ connected in the neighboring element $\BgElemInd^-$. Hence, the penalty needs to be applied to two different pairs of element bases on the same facet to involve all sparsely supported basis functions and appropriately stabilize the problem.

\begin{figure*}[h]
    \centering 
    \def\svgwidth{13.0cm}
    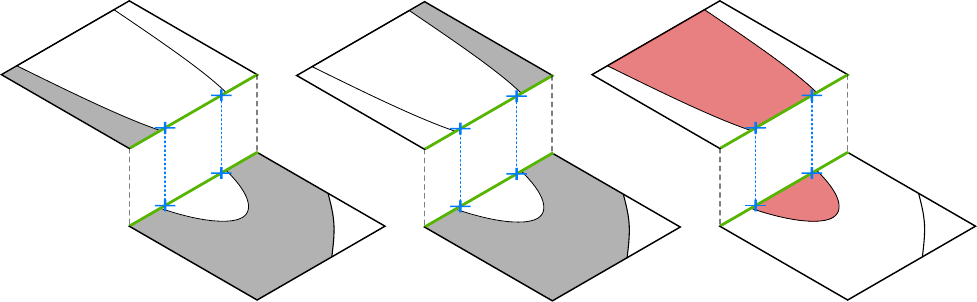
    \caption{Ghost cluster pairs constructed for the example ghost facet and material layout shown in \Cref{fig:ghost_for_enriched} (C).}
    \label{fig:ghost_clusters}
\end{figure*}

The procedure is generalized and formalized as follows. Let's define the union $\InteriorInterfaces$ of all material interfaces and boundaries inside the ambient domain $\AmbientDomain$ as 
\begin{equation}
    \label{eqn:material_interfaces}
    \InteriorInterfaces = \bigcup_{\substack{ m,l \in \SetOfMaterials \cup \SetDef{0} \\ m \neq l } } \MaterialInterface{m}{l}.
\end{equation}
The set of all ghost faces is defined as the set of faces $\Facet$ between neighboring background elements $\Omega^{E^+}$ and $\Omega^{E^+}$, at least one of which is intersected by an interior interface $\InteriorInterfaces$
\begin{align}
    \label{eqn:set_of_ghost_facets}
    \SetOfGhostFacets = \Bigl\{ 
        \Facet & : \Domain^{\BgElemInd^+} \cap \InteriorInterfaces \neq \emptyset
        \ \vee \ 
        \Domain^{\BgElemInd^-} \cap \InteriorInterfaces \neq \emptyset \nonumber
        \\
        & \overline{\Facet} = \overline{\Omega^{E^+}} \cap \overline{\Omega^{E^-}}, ~ E^+ \neq E^- \Bigl\},
\end{align}
where $\overline{(\cdot)}$ denotes the completion of a domain.

The contribution to the weak form \eqref{eqn:weak_form_discrete} then reads as:
\begin{align}
    \Residual_{G}(\DiscDispl, \DiscTestDispl) & =
    \sum_{\Facet \in \SetOfGhostFacets}
    \sum_{\Unzipping=1}^{\NumUnzippings(\BgElemInd^{+})} 
    \sum_{v \in \mathcal{U}_{u,\BgElemInd^{-}}} \cdots \nonumber
    \\
    \cdots & \sum_{k=1}^{\PolyOrder} \int_{F}  
    \gamma_{G}^k 
    \Big\llbracket \NormalDeriv{k} \, \DiscTestDisplExt \Big\rrbracket 
    \cdot
    \Big\llbracket \NormalDeriv{k} \, \DiscDisplExt \Big\rrbracket 
    \, 
    d\BoundarySegment
    \label{eqn:ghost_residual}
\end{align}
where $\mathcal{U}_{u,\BgElemInd^{-}}$ is the set of all subdomains in $\Domain^{\BgElemInd^{-}}$ that are connected to $\SubPhase_{\BgElemInd^+}^{u}$ across the facet $\Facet$ defined as follows:
\begin{align}
    \label{eqn:ghost_connected_subdomains}
    & \mathcal{U}_{u,\BgElemInd^{-}} = \left\{ \,
    v \in \SetDef{1, \dots, \NumUnzippings(\BgElemInd^{-})} : \right. \nonumber
    \\
    & \cdots ~ \left. \MaterialIndex(\SubPhase_{\BgElemInd^+}^{u}) = \MaterialIndex(\SubPhase_{\BgElemInd^-}^{v}) 
    \, \wedge \,
    \left| \closure{\SubPhase_{\BgElemInd^+}^{u}} \cap \closure{\SubPhase_{\BgElemInd^-}^{v}} \right| >0
    \, \right\},
\end{align}
where $\MaterialIndex(\ElemSubPhase{u})$ indicates the materials the respective subdomains are part of, and the term $|\cdot|>0$ indicates that the subphases need to connect via an interface of non-zero length/area. The operator $\NormalDeriv{k}$ denotes the $k$-th normal derivative. 
The jump in normal derivatives operates on the extension of the test and trial functions $\DiscTestDisplExt$ and $\DiscDisplExt$ from their respective subdomains $\SubPhase_{E^+}^u$ and $\SubPhase_{E^-}^u$ to the background elements $\Omega^{E^+}$ and $\Omega^{E^-}$ which are obtained by omitting the enrichment function $\EnrichmentFunction$ in \eqref{eqn:enriched_interpolation}. 
For details on the scaling of the penalty parameter $\gamma_{G}^k$, the reader is referred to \cite{noel2022xiga}.

Using the definition of clusters introduced in \Cref{sec:framework-ig}, the goal of the implementation of the ghost stabilization in the preprocessor consists of finding and constructing the cluster pairs corresponding to each ghost facet and valid subdomain combination. For the example shown in \Cref{fig:ghost_for_enriched} (C), we need to construct the cluster pairs shown in \Cref{fig:ghost_clusters}.

\clearpage
\section{Mesh Generation}
\label{sec:main}

\begin{figure}[btp]
    \centering 
    \def\svgwidth{4.4cm}
    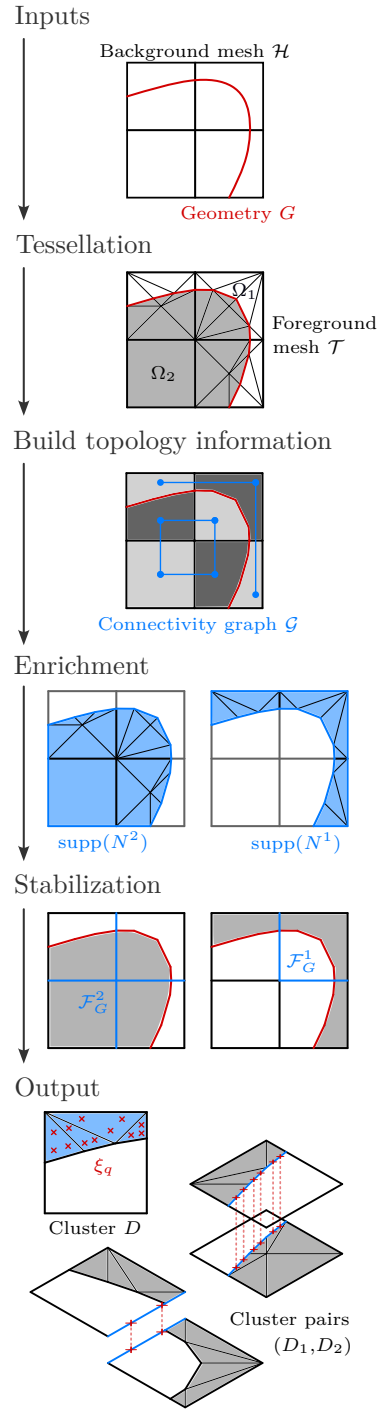
    \vspace*{0.2cm}
    \caption{Sketch of the workflow used in the implementation of the preprocessor.}
    \label{fig:workflow}
\end{figure}

This section outlines the implementation of the preprocessor with its associated algorithms and data structures. Furthermore, it discusses the concepts and strategies used to achieve robustness and parallel scalability. 
The input to the preprocessor consists of one or multiple geometries and a background mesh in which they are immersed. The output are clusters, as introduced in \Cref{sec:framework-ig}. These inputs and outputs are formalized in \Cref{sec:preliminaries}, alongside some basic concepts and notation used throughout this section.
The preprocessor's workflow follows a series of four main steps as visualized in \Cref{fig:workflow} and later summarized in \Cref{alg:driver_algorithm}.
First, a grid conforming to the geometric boundaries referred to as the "foreground mesh" $\ForegroundMesh$, is generated. This is achieved by triangulating the background mesh in two sub-steps: regular and templated subdivision, both outlined in \Cref{sec:fg_mesh_generation}. In the second step, topological information needed for the enrichment and ghost stabilization is computed and stored for the foreground mesh $\ForegroundMesh$, see \Cref{sec:subphase_generation}. 
Using the topology information, enrichment is performed in the third step outlined in \Cref{sec:enrichment}. The generation of ghost clusters is discussed as the last step in \Cref{sec:ghost}. The parallel implementation of each of these steps is discussed separately at the end in \Cref{sec:parallel}, with the intent of not providing too much information at once.

\subsection{Preliminaries}
\label{sec:preliminaries}
  
\begin{figure}[h!]
    \centering 
    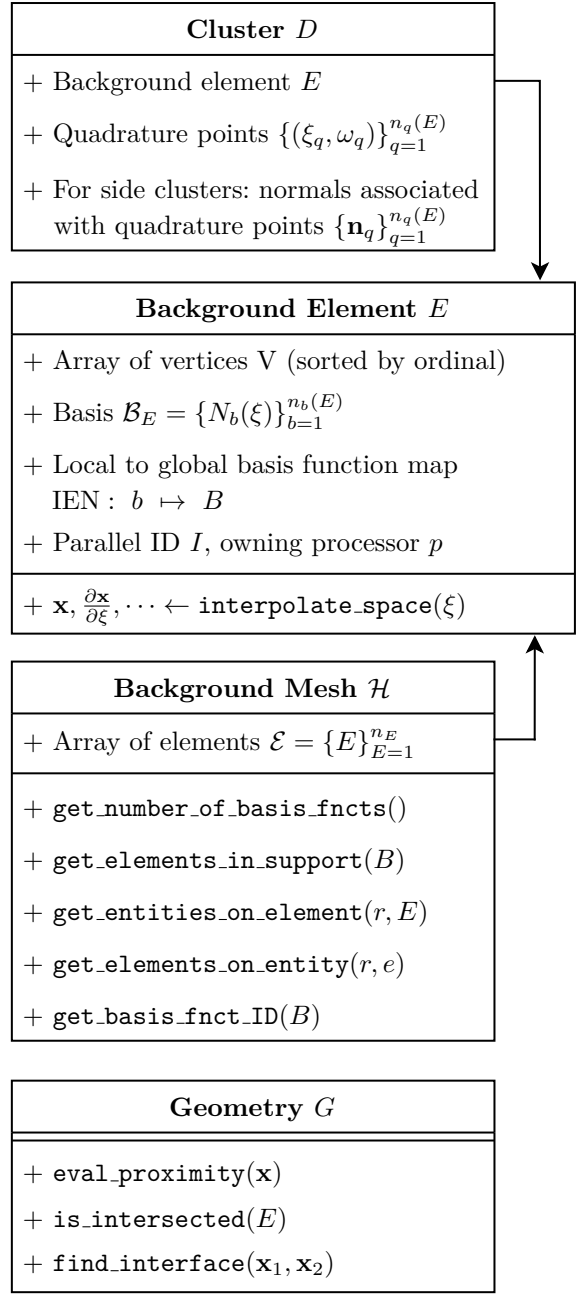
    \vspace*{0.2cm}
    \caption{Unified modeling language (UML) chart outlining the data structures for the input and output data for the preprocessor.}
    \label{fig:IO_data_structures}
\end{figure}

\paragraph*{Inputs and Outputs}
The goal of the immersed FE preprocessor is to generate sets of clusters associated with the various subdomains, namely, material domains $\Omega_m$, boundaries $\Gamma_m$, interfaces $\MaterialInterface{m}{l}$, and ghost faces $\mathcal{F}_G$, to support immersed FE analysis.
This information is formalized in the data structure shown in \Cref{fig:IO_data_structures} containing the background element, which itself consists of geometry (vertices) information and a finite element basis $\mathcal{B}_E$, and the custom quadrature rules needed in the immersed setting. For internal interfaces $\MaterialInterface{m}{l}$ and ghost faces $\mathcal{F}_G$, clusters are stored in pairs $(D_1, D_2)$, with matching local ordering of the quadrature points to compute the integrals \eqref{eqn:elemental_interface_integral}.
The local-to-global basis function map on each background element, or the "$\mathrm{IEN}$"\footnote{The naming of the IEN in \cite{hughes2012finite} stems from the fact that it is used to "index elemental nodes". However, its usage is equally valid for non-nodal finite element bases.}, adopting the language in \cite{hughes2012finite}, provides the information necessary to assemble the elemental residual and tangent matrices into a global system. 

The inputs needed for constructing clusters are a background mesh $\BackgroundMesh$ and a set of geometries $\SetDef{G_g}_{g=1}^{n_g}$ defining the various material domains $\Omega_m$. 
The background mesh is a finite element mesh that provides: (i) a list of elements that make up the support of a basis function with global index $B$ and (ii) the connectivity information of its mesh entities as discussed in the following paragraph. Without loss of generality, we assume a rectilinear background mesh of Quad or Hex elements, but the algorithm also applies to unstructured meshes made of different element types. 
Non-conformal meshes, such as those found in locally refined hierarchical meshes, require additional modifications which we may address in future work.
The representation of each geometry needs to support the following queries:
\begin{itemize}
    \item Evaluate the proximity of a point in space, i.e., whether a point lies outside, inside, or on the boundary of the geometry within some tolerance. Implicit geometry definitions naturally lend themselves to this query. Explicit geometries can be translated into implicit geometries \cite{stanford2019higher}, or approaches such as ray-tracing \cite{arvo1989survey} may be used.
    \item Check whether a given background element $E$ is intersected by a geometry interface. Note that, for robustness, the query may not return any false negatives, though false positives are acceptable and only affect computational efficiency. For implicit geometries, the range evaluation procedure presented by Saye \cite{saye2022high} may be employed. 
    \item Find the location of an interface between two points with opposite proximity within a background element. 
    For implicit geometries, this can be done by interpolating the level-set field along the edge or on the background mesh and then performing a polynomial root-finding operation \cite{saye2022high,reuter2008solving,mourrain2004bernstein}. For geometries defined through boundary representations, approaches such as the one in \cite{badia2022geometrical} may be utilized.
\end{itemize}

\begin{figure}[btp]
    \centering 
    \def\svgwidth{6.9cm}
    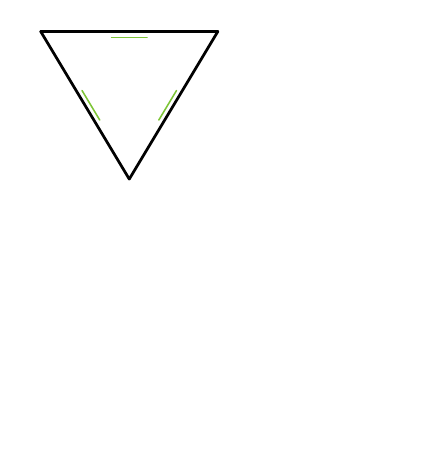
    \vspace*{0.2cm}
    \caption{(A) Vertex and edge ordinals for a triangular cell. (B) Vertex, edge, and cell indices for an example mesh. (C) Cell-vertex connectivity $\mathcal{C}$ defining the mesh in (B) and the resulting edge connectivity.}
    \label{fig:entity_connectivity}
\end{figure}

\paragraph*{Mesh Entities and Connectivity}
\label{pg:mesh_entities}
Understanding the relationships between mesh entities -- vertices, edges, faces (2D cells), and 3D cells -- is essential for the algorithms presented in this paper. The following mesh entity and connectivity terminology is used throughout this paper:
\begin{itemize}
    \item \textbf{Rank $r$}: The dimensionality of a mesh entity, e.g., vertices have rank 0, edges rank 1, faces rank 2, and cells rank 3.
    \item \textbf{Facet $F$}: A $(d-1)$-dimensional entity of a $d$-dimensional mesh, e.g., edges in 2D and faces in 3D.
    \item \textbf{Ordinal $o$}: The local index of an entity with respect to another entity, as shown in \Cref{fig:entity_connectivity} (A).
    \item \textbf{Cells $\mathrm{c}$ and Vertices $\mathrm{v}$}: Purely geometric entities, distinct from elements and nodes, which also include finite element basis information. Vertices specifically refer to the corner points of a cell.
    \item \textbf{Entity Index:} A (processor-local) identifier for mesh entities of a specific type, ranging from 1 to the total number of such entities on the processor. 
\end{itemize}
Background meshes generated with commonly used mesh generators do not provide this information directly, but it can be computed using \Cref{alg:compute_entity_connectivity} provided in the appendix.

\begin{figure}[btp]
    \centering 
    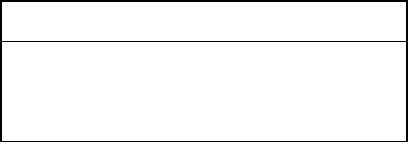
    \caption{Entity connectivity consisting of two maps.} 
    \label{fig:entity_connectivity_data_structure}
\end{figure}
Using these definitions, we introduce two maps, implemented as nested arrays, to express the connectivity of mesh entities: (i) an entity-to-cell map ($\mathrm{EtC}$), which contains the cells attached to an entity of a given index, and (ii) a cell-to-entity map ($\mathrm{CtE}$), which contains the indices of entities attached to a cell of a given index. The list of entities attached to a cell is ordered by the ordinal of the entity with respect to the cell, as illustrated in \Cref{fig:entity_connectivity_data_structure}. An example of these maps is shown in \Cref{fig:entity_connectivity} (C) for the edge connectivity of the small four-element mesh shown in \Cref{fig:entity_connectivity} (B). 

\begin{algorithm*}[btp]
    \caption{
        Driver algorithm. \\
        Input: Background mesh $\BackgroundMesh$, Geometries $\SetDef{G_g}_{g=1}^{n_g}$ \\
        Output: Sets of bulk $\mathcal{D}_B$ and side clusters $\mathcal{D}_S$; interface $\mathcal{D}_I$ and ghost side cluster pairs $\mathcal{D}_G$ }
    \label{alg:driver_algorithm}
    \begin{algorithmic}[1]
        \State \textit{\# Foreground mesh generation} \Comment{see \Cref{sec:fg_mesh_generation}}
        \State $\ForegroundMesh \leftarrow \mathtt{initialize\_foreground\_mesh}(\BackgroundMesh)$
        \State $\ForegroundMesh.\mathtt{regular\_subdivision}(\SetDef{G_g}_{g=1}^{n_g})$ \Comment{see \Cref{alg:regular_subdivision} in \Cref{sec:fg_mesh_generation-reg_sub}}
        \For{each geometry $G \in \SetDef{G_g}_{g=1}^{n_g}$}
            \State $\ForegroundMesh.\mathtt{templated\_subdivision}(G)$ \Comment{see \Cref{alg:templated_subdivision} in \Cref{sec:fg_mesh_generation-templ_sub}}
        \EndFor
        \State {\color{gray}assign parallel IDs to the foreground cells $\ForegroundMesh.\mathtt{communicate\_cell\_IDs}()$}
        \For{each fg. cell $c \in \ForegroundMesh.\mathcal{C}$}
            \State $c.m \leftarrow \mathtt{material\_map}(c.m)$
        \EndFor
        \State
        \State \textit{\# Generate topology information} \Comment{see \Cref{sec:subphase_generation}}
        \State $\ForegroundMesh.\mathcal{F} \leftarrow \mathtt{compute\_entity\_connectivity}(\ForegroundMesh, d-1)$ \Comment{see \Cref{alg:compute_entity_connectivity}}
        \State fg. facets attached to each bg. facet $\mathrm{DF} \leftarrow \ForegroundMesh.\mathtt{compute\_bg\_facet\_descendants}()$ \Comment{see \Cref{alg:generate_bg_facet_descendants}}
        \State $\ForegroundMesh.\mathtt{generate\_subphases}()$ \Comment{see \Cref{alg:generate_subphases}}
        \State $\ForegroundMesh.\mathtt{generate\_subphase\_graphs}(\mathrm{DF})$ \Comment{see \Cref{alg:generate_subphase_graphs}}
        \State
        \State \textit{\# Enrichment - create unzipped background elements} \Comment{see \Cref{sec:enrichment}}
        \State $\Xi \leftarrow \mathtt{unzip\_interpolation\_mesh}(\BackgroundMesh, \ForegroundMesh)$ \Comment{see \Cref{alg:unzipping}}
        \State
        \State \textit{\# Cluster generation}
        \State $\mathcal{D}_B \leftarrow \mathtt{create\_bulk\_clusters}(\ForegroundMesh, \Xi)$ \Comment{see \Cref{alg:bulk_cluster_generation}}
        \State $\mathcal{D}_S, \mathcal{D}_I  \leftarrow \mathtt{create\_side\_clusters}(\ForegroundMesh, \Xi)$ \Comment{see \Cref{alg:side_cluster_generation}}
        \State $\mathcal{D}_G \leftarrow \mathtt{generate\_ghost\_clusters}(\ForegroundMesh, \Xi)$ \Comment{see \Cref{alg:ghost} in \Cref{sec:ghost}}
        \State
        \State \Return{$\mathcal{D}_B, \, \mathcal{D}_S, \, \mathcal{D}_I, \, \mathcal{D}_G$ }
    \end{algorithmic}
\end{algorithm*}

\paragraph*{Notation and Conventions}
The following notation and conventions are used in the presented algorithms.
\begin{itemize}
    \item The "." operator denotes access to a member function or variable of an object.
    \item Arrays are represented using braces "$\{\}$" and tuples using parentheses "$()$".
    \item The presented algorithms use $1$-based indexing unless otherwise specified, for the sake of readability.
    \item Access to an element of an array or a map is denoted using the "[]" operator.
    \item "[-1]" is shorthand for accessing the last element of an array.
    \item Loops over each element of an array are condensed using the shorthand "[:]".
    \item All relevant objects have an index. The symbols representing an object are used interchangeably as the index or the object itself; e.g., in "$E \in \BackgroundMesh.\mathcal{E}$", $E$ may represent the element or its index in the global list $\BackgroundMesh.\mathcal{E}$.
    \item Lines in italics and starting with a "\#", are comments.
    \item Components of algorithms are colored in gray to indicate they are only relevant to parallel implementation. We only discuss these additions at the end of \Cref{sec:parallel} and recommend ignoring them on the first read-through.
    \item The words foreground and background are abbreviated as "fg." and "bg.", respectively.
\end{itemize}
Using these conventions, the workflow is summarized in \Cref{alg:driver_algorithm}.

\begin{figure}[h!]
    \centering 
    \def\svgwidth{6.2cm}
    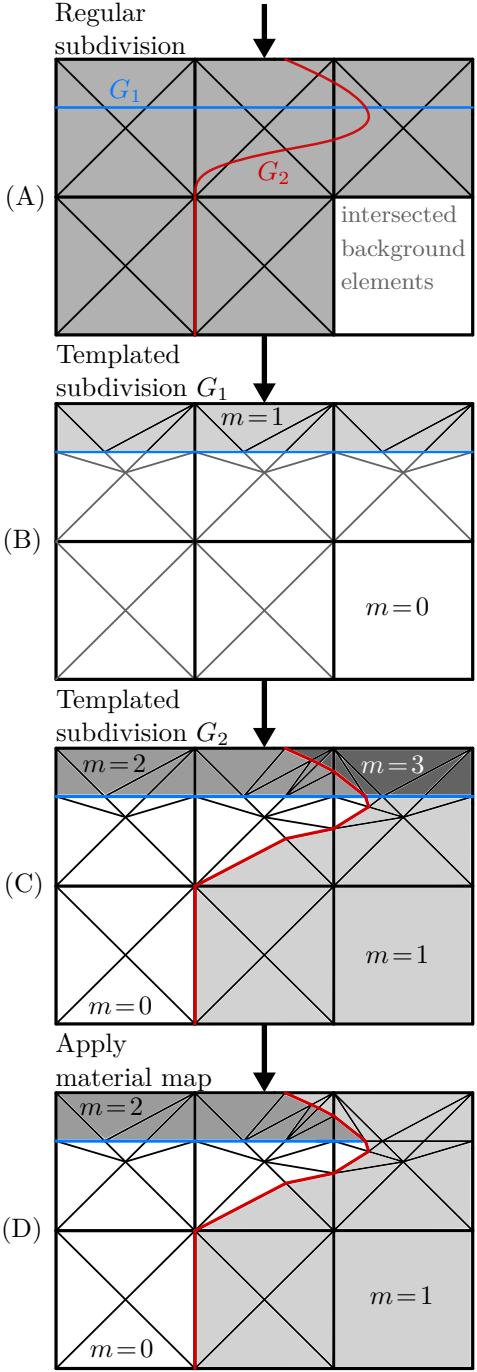
    \vspace*{0.2cm}
    \caption{Subdivision process overview. (A) Regular subdivision: application of templates to bg. elements intersected by geometric interfaces, shown with their original accuracy. (B) and (C) Templated subdivision with recursive material assignment for the first $G_1$ and second geometry $G_2$, respectively. (D) Application of a material map and the resulting foreground mesh.}
    \label{fig:subdivision_process}
\end{figure}

\subsection{Conformal Foreground Mesh Generation}
\label{sec:fg_mesh_generation}

The goal of this part of the preprocessor algorithm is to subdivide the background mesh's elements intersected by geometric interfaces into groups of triangular or tetrahedral cells with facets conforming to those interfaces. As a result, the geometric interfaces are approximated by sets of straight edges in 2D or flat triangles in 3D, while the volumetric domains for each of the materials $\Omega_m$ are approximated by groups of the generated triangular or tetrahedral cells. The focus for this step, besides enabling the generation of topology information for enrichment and ghost stabilization, is on robustness and speed. To this end, a combination of templated subdivision procedures is employed. 

A drawback of the presented approach is that the quadrature rules, which are subsequently constructed from the triangulation in \Cref{sec:enrichment}, only capture the geometry and its interfaces with linear accuracy. However, modifications using isoparametric projection \cite{lehrenfeld2016high} or manipulation \cite{scholz2019numerical, scholz2020high} of the resulting triangulated foreground mesh are possible as part of the foreground mesh generation; although, they are not further addressed in this work. 

To obtain a Tri/Tet mesh in the vicinity of the interfaces, background elements cut by any geometric interfaces are first subdivided into a pre-defined set of triangular or tetrahedral cells. This step is referred to as the "regular subdivision". On the resulting cells, the edges are probed for intersections with geometric interfaces and further subdivided into triangular or tetrahedral cells conforming to the interface using a set of subdivision templates. This "templated subdivision" approach has been used commonly in other works, see \cite{soghrati20123d}. The main difference in the presented approach is the addition of background ancestry information which enables the computation of topology information without floating point operations at a later stage. 
Lastly, a recursive procedure to determine the material membership of each generated cell is incorporated into the subdivision processes. The complete process is illustrated in \Cref{fig:subdivision_process}.

In what follows, we first introduce the concept of background ancestry and the initialization of the foreground mesh, before discussing the regular subdivision, the templated subdivision, and material assignment.

\begin{figure}[btp]
    \centering
    \def\svgwidth{7.2cm}
    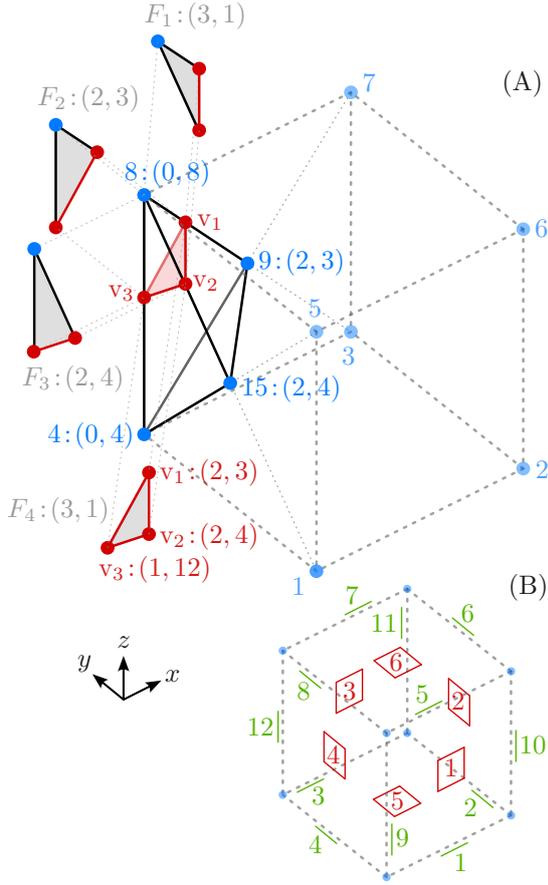
    \vspace*{0.1cm}
    \caption{(A) Concept of background ancestry demonstrated using a Tet which gets subdivided by an interface. The "parent" background element is indicated in dashed lines. The ancestry (entity rank $\Rank$, entity index $\Ancestor$) of the existing (blue) and added (red) vertices are shown alongside some of the newly formed faces (gray). (B) Indices of the edges (green) and faces (red) attached to the parent background element.}
    \label{fig:ancestry_deduction}
\end{figure}

\paragraph{Background Ancestry}
\label{pg:background_ancestry}
The background ancestor of an entity of the foreground mesh is the background entity with index $\Ancestor$ of the lowest rank $\Rank$, i.e., dimensionality, which fully contains the foreground entity. This information, stored in the form of index pairs $(\Rank,\Ancestor)$, is used to unambiguously identify vertices during the subdivision and to generate topology information without the use of floating point arithmetic. Examples of the background ancestry of various entities formed during the conformal mesh generation process are shown in \Cref{fig:ancestry_deduction}. 

Consider the tetrahedral foreground cell $\SetDef{4,8,9,15}$ in \Cref{fig:ancestry_deduction} that was formed during the subdivision process. Vertices $4$ and $8$ were part of the original background mesh and are descendants of those vertices; hence, their background ancestry is denoted as $(0,4)$ and $(0,8)$ for background vertices (rank $\Rank=0$) with indices $a=4$ and $a=8$, respectively. Vertices $9$ and $15$ were formed on a face of the background element; hence, their ancestors are denoted as $(2,4)$ and $(2,3)$, for being descendants for background faces (dimensionality, i.e. rank $\Rank=2$) with indices $a=4$ and $a=3$. A new vertex $\Vertex_3$ is formed on the edge connecting vertices $4$ and $8$. It can be checked using entity connectivity information whether those vertices are attached to the same edge or face, or are part of the same background element. In this case, they are part of the same background edge $(\Rank=1,a=12)$, and therefore any new vertex formed by them would be a descendant of that edge. Considering entities of one dimension higher, the face formed between the vertices $\Vertex_1$, $\Vertex_2$, and $\Vertex_3$ would be a descendant of the background element (dimensionality $\Rank=3$) itself, as two of the vertices descend from different faces. 

\begin{figure}[btp]
    \centering 
    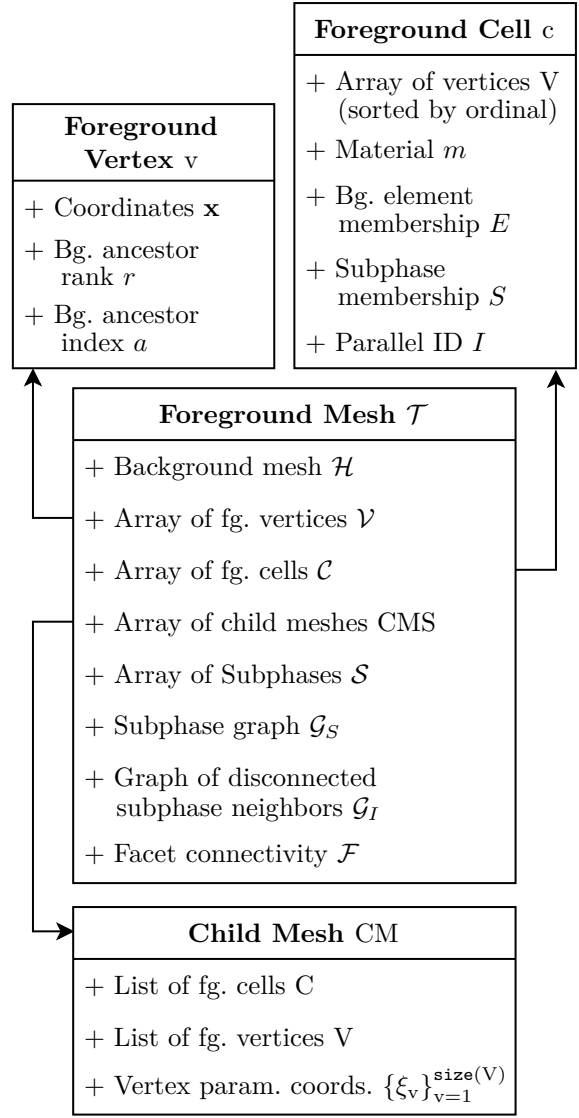
    \vspace*{0.2cm}
    \caption{UML chart showing the data stored in the foreground mesh. The first four members listed for the foreground mesh are populated during the subdivision process in \Cref{sec:fg_mesh_generation}; the remaining four are filled during the generation of the topology information in \Cref{sec:subphase_generation}.}
    \label{fig:fg_mesh_data_structure}
\end{figure}

\paragraph*{Initialization}

To initialize the foreground mesh data, shown in \Cref{fig:fg_mesh_data_structure}, the geometric information, i.e., the vertex coordinates and vertex connectivity of the background elements, are copied from the background to the foreground mesh. The ancestry $(r, a)$ for every foreground vertex and cell $E$ is initialized trivially, given that each foreground vertex or cell is a descendant of the background vertex or element with the same index. The material index $m=0$ is assigned to all cells.

To access foreground cell information from the perspective of individual background elements, child mesh containers $\mathrm{CM}$ are created containing the vertices and cells formed within each background element. 
This includes the parametric coordinates of the vertices relative to the parametric coordinate space of the background element, which are essential for geometry interpolation operations.
During initialization, each child mesh container is populated with the vertices from the single cell that originally constitutes the background element.

The remaining information in \Cref{fig:fg_mesh_data_structure} is related to the material topology and is populated after the subdivision process.

\subsubsection{Regular Subdivision}
\label{sec:fg_mesh_generation-reg_sub}

\begin{figure}[btp]
    \centering 
    \def\svgwidth{6.0cm}
    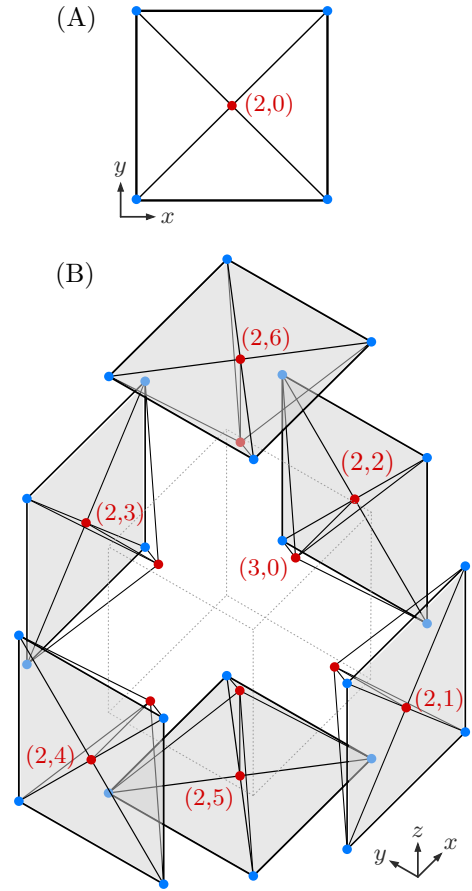
    \vspace*{0.2cm}
    \caption{Regular subdivision templates in 2D (A) and 3D (B). New vertices generated by the templates are shown with their ancestry (ancestor entity rank, ordinal) in red and existing vertices in blue.}
    \label{fig:reg_sub_template}
\end{figure}

The regular subdivision serves the purpose of generating Tri or Tet cells near the geometric interfaces as these allow for a robust conformal mesh generation procedure using templates. If the algorithm is to be applied to a pure Tri/Tet background mesh, this step is skipped. For the procedure outlined in \Cref{alg:regular_subdivision}, the templates shown in \Cref{fig:reg_sub_template} are applied to every background element intersected by a geometric interface. 

The templates employed here do not minimize the number of Tris or Tets required to subdivide a Quad or Hex element. However, they offer two distinct advantages over templates generating a minimal number of Tri/Tet cells.
First, the cells generated by these templates do not further subdivide existing background edges, ensuring that each generated foreground cell includes at most one complete background edge in both 2D and 3D. This simplifies the determination of the background ancestry of newly formed vertices and facets during the templated subdivision, reducing the process to a series of straightforward checks outlined in \Cref{sec:apx-ancestry_reg_sub_template}.
Second, the templates exhibit symmetry about all axes, which ensures that, for higher-order geometry descriptions, the resulting lower-order approximation does not introduce a bias depending on the orientation of the applied template.

\begin{algorithm*}[btp]
    \caption{
        Regular Subdivision. \\
        Input: Initialized foreground mesh $\ForegroundMesh$ (self), Geometries $\SetDef{G_g}_{g=1}^{n_g}$ \\
        Output: Foreground mesh $\ForegroundMesh$ with triangulated cut bg. elements}
    \label{alg:regular_subdivision}
    \begin{algorithmic}[1]
        \State initialize new vertex $\mathrm{Q_v}$ and cell queues $\mathrm{Q_c}$, get template $T \leftarrow \mathtt{regular\_subdivison\_template}(d)$
        \State get location of new vertices in bg. element $\SetDef{(\xi,\Rank,\Ordinal)}_{i=1}^{\mathtt{NumNewVerts}} \leftarrow T.\mathtt{get\_new\_vertex\_locations}()$
        \For{each bg. element $\BgElemInd \in \BgMesh.\mathcal{\BgElemInd}$ }
            \If{$G.\mathtt{is\_element\_intersected}(\BgElemInd)$ for any $G \in \SetDef{G_g}_{g=1}^{n_g}$ }
                \For{each new vertex in template $(\xi,r,o) \in \SetDef{(\xi,\Rank,o)}_{i=1}^{\mathtt{NumNewVerts}}$}
                    \State get the entity index $\Entity = \BackgroundMesh.\mathtt{get\_entities\_on\_element}(\BgElemInd,\Rank)[\Ordinal]$
                    \If{this is a new vertex $!\mathrm{Q_v}.\mathtt{request\_exists}((\Rank,\Entity))$ }
                        \State compute the physical coordinates $\pos \leftarrow \BgElemInd.\mathtt{interpolate\_space}(\xi)$
                        \State request new vertex and remember its index $\mathrm{NVI}.\mathtt{append}(\mathrm{Q_v}.\mathtt{queue}(\pos,\Rank,\Entity))$
                    \Else{ \textit{\# another bg. element has already requested this vertex}}
                        \State store the index for the new vertex $\mathrm{NVI}.\mathtt{append}(\mathrm{Q_v}.\mathtt{get\_index}(\Rank,\Entity))$
                    \EndIf
                    \State store the new vertex in the child mesh $\ListOfCMs[\BgElemInd].\mathtt{add\_vertex}(\mathrm{NVI}[-1],\xi)$
                \EndFor 
                \State get the cell-vertex connectivity for the new cells $\mathrm{C}_{new} \leftarrow T.\mathtt{generate\_cells}(\BgElemInd.\mathrm{V},\mathrm{NVI})$
                \For{each new cell in template except the first $\Cell_{new} \in \mathrm{C}_{new}[-1:2]$}
                    \State queue cell and add its index to the child mesh $\ListOfCMs[\BgElemInd].\mathrm{C}.\mathtt{append}(\mathrm{Q_c}.\mathtt{queue}(\Cell_{new},\BgElemInd,0))$
                \EndFor 
                \State replace the subdivided element with the first new cell $\mathcal{C}[\BgElemInd].\mathrm{V} \leftarrow \mathrm{C}_{new}[1]$
            \EndIf 
        \EndFor 
        \State $\mathtt{create\_vertices}(\mathrm{Q_v})$
        \State $\mathtt{create\_cells}(\mathrm{Q_c})$
    \end{algorithmic}
\end{algorithm*} 

To perform the regular subdivision, as outlined in \Cref{alg:regular_subdivision}, the templates contain the following information: 
\begin{enumerate}
    \item A list of new vertices to be generated with their location $\xi$ in the parametric space of the background element in addition to the rank $r$ and ordinal $o$ of the background entity they are located on. This information is provided through the function \texttt{get\_new\_vertex\_locations}().
    \item A pre-defined cell-vertex connectivity for the new cells in the format shown in \Cref{fig:entity_connectivity}. The template expresses the connectivity using the indices of both the existing and new vertices. This information is provided through the function $\mathtt{generate\_cells}(E.\mathrm{V},\mathrm{NVI})$, where $E.\mathrm{V}$ is the list of vertices of the existing cell and $\mathrm{NVI}$ is the list of indices of the new vertices generated by the template.
\end{enumerate}

New vertices and cells are temporarily stored in a queue, rather than directly added into the foreground mesh data structure. A queue consists of an array of requested entities with their initialization data and a (hash) map, which relates identifying information with the position of the entity in the array of requests. This approach
\begin{enumerate}
    \item avoids repeated resizing operations of the global lists of vertices $\mathcal{V}$ and cells $\mathcal{C}$ on the foreground mesh which can be computationally costly.\footnote{Dynamic arrays, such as C++'s \texttt{std::vector}, are equipped with smart resizing strategies when repeatedly appending elements to the array. Such strategies may lead to over-allocation of memory and repeated copying of the existing array elements, resulting in significant computational cost and memory fragmentation in the context of the large arrays being edited in this case.} 
    \item facilitates the generation of unique vertices. The regular subdivision templates, shown in \Cref{fig:reg_sub_template}, generate at most a single vertex on a given background entity. The rank and index of the background entity are used to identify identical vertices that may have already been requested by a neighboring background element if the new vertex is located on a shared background facet. Otherwise, newly created vertices would need to be merged using their physical coordinates, as is done by \cite{zhang2022object}, requiring floating point arithmetic.
\end{enumerate}
Lastly, to eliminate the subdivided cell from the foreground mesh, its vertex list is overwritten with the vertex list of one of the new cells. This avoids a delete operation and the need for re-indexing.

\subsubsection{Templated Subdivision}
\label{sec:fg_mesh_generation-templ_sub}

The templated subdivision completes the generation of the foreground mesh by subdividing the cells generated in the regular subdivision into cells conforming to the geometric interfaces. This procedure is outlined in \Cref{alg:templated_subdivision}.

\begin{figure}[btp]
    \centering 
    \def\svgwidth{5.5cm}
    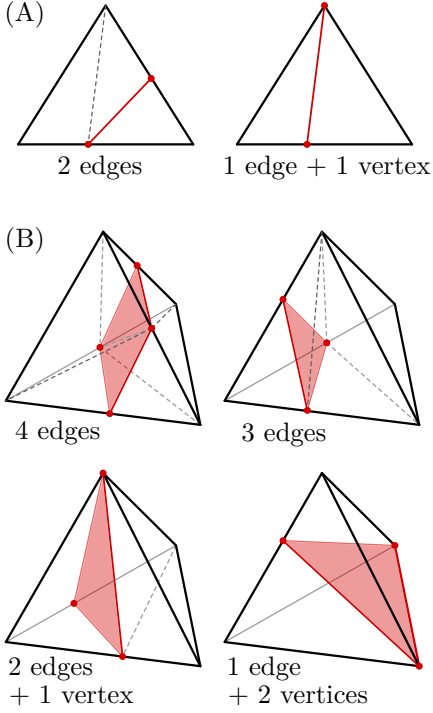
    \vspace*{0.3cm}
    \caption{Possible intersection configurations for foreground cells during the templated subdivision in 2D (A) and 3D (B). The templates are permutations of the configurations shown.}
    \label{fig:subdivision_templates}
\end{figure}

The procedure assesses foreground edges generated during the regular subdivision if they are intersected by a geometric interface and then generates additional vertices and cells according to one of the templates shown in \Cref{fig:subdivision_templates} or their permutations. Depending on the proximity of vertices, a new material index $m$ is assigned to the elements. The process is repeated for every geometry applied.

First, the proximity of each foreground vertex, i.e., whether a vertex is located inside, outside, or on the boundary of a geometry, is determined and stored to avoid error-prone re-computation. Subsequently, the edges of every potentially intersected cell are checked for differences in the proximity of the adjacent two vertices. If so, it is considered intersected and the intersection point is determined. 
The approach does not allow for geometry intersections on edges that contain a vertex that lies on the interface. It also limits the number of intersections to one per edge. This ensures robustness toward unwanted intersection configurations but leads to a dependency of the final foreground mesh on the order in which the geometries are applied. Consider $G_2$ in \Cref{fig:subdivision_process} (A): the protrusion on the right would be cut off if $G_2$ were to be applied first, as it intersects the background edge twice. The process would not consider the edge intersected in this case as both its vertices, i.e., its end-points, lie on the same side of the interface. The issue of insufficiently resolved geometry can be overcome by local mesh refinement, either by refining the background mesh or by applying additional, e.g., quadtree/octree subdivision templates before the regular subdivision. 

Computing the interface location along the edge, both the physical and parametric coordinates of the new vertex are found. The background ancestry of the new vertex is determined by finding a common ancestor of the two edge vertices. For the regular subdivision template presented here, determining the ancestor is simplified as outlined in Appendix \ref{sec:apx-ancestry_reg_sub_template}.
If a different regular subdivision template or an unstructured Tri or Tet input background mesh were to be used, search operations would need to be performed on the entity connectivity of the background mesh to determine the background ancestors of the new vertices. 
Requests for the new vertices on the intersected edges are queued in the same manner as for the regular subdivision. A sorted pair of vertex indices serves as the identifier for the edge. 
The cell-vertex connectivity for the new cells is generated by providing the indices of the new and existing vertices to the template through the method $\mathtt{generate\_cells}(c.\mathrm{V},\mathrm{V}_{new})$ in \Cref{alg:templated_subdivision}, and, again, one of the new cells is used to replace the information of the subdivided cell to forego an explicit deletion. 

\begin{algorithm*}[btp]
    \caption{
        Templated subdivision. \\
        Input: Foreground mesh $\ForegroundMesh$ (self), Geometry $G$ \\
        Output: Fully subdivided foreground mesh $\ForegroundMesh$ (self) }
    \label{alg:templated_subdivision}
    \begin{algorithmic}[1]

        \State evaluate each vertex's proximity $\mathrm{P}[:] \leftarrow G.\mathtt{compute\_proximity}(\mathcal{V}[:].\pos)$ \Comment{possible values $(+,0,-)$}

        \State initialize new vertex $\mathrm{Q_v}$ and cell queues $\mathrm{Q_c}$
        \For{each intersected bg. element $\BgElem \in \SetDef{1,\cdots,n_{\BgElem}}$ }
            \If{bg. element is intersected $G.\mathtt{is\_intersected}(\BgElem)$}

                \For{each child cell in bg. element $\Cell \in \ListOfCMs[\BgElem].\mathrm{C}$}
                    \State initialize array with empty entry for every edge on cell $\mathrm{V}_{new} \leftarrow \SetDef{0}_{\Entity=1}^{\mathtt{NumEdges}}$ 
                    \For{each edge on cell $(\Vertex_1,\Vertex_2)_{\Entity} \in \SetDef{(\Vertex_1,\Vertex_2)_{\Entity}}_{\Entity=1}^{\mathtt{NumEdges}}$}
                        \If{edge is intersected $\mathrm{P}[\Vertex_1] \neq \mathrm{P}[\Vertex_2] \neq 0$}
                            \State get the location of the interface on edge $\xi_e \leftarrow G.\mathtt{find\_interface}(\Vertex_1.\pos,\Vertex_2.\pos)$
                            \State $\xi \leftarrow \mathtt{interpolate\_on\_edge}(\xi_e, \ListOfCMs[\BgElem].\mathtt{get\_coords}(\Vertex_1), \ListOfCMs[\BgElem].\mathtt{get\_coords}(\Vertex_2))$
                            \If{existing vertex $\mathrm{Q_v}.\mathtt{request\_exists}((\Vertex_1,\Vertex_2))$}
                                \State $\mathrm{V}_{new}[\Entity] \leftarrow \mathrm{Q_v}.\mathtt{get\_index}((\Vertex_1,\Vertex_2))$
                            \Else{} 
                                \State get the physical coordinates $\pos \leftarrow \BgElem.\mathtt{interpolate\_space}(\xi)$
                                \State compute the ancestry $(\Rank,\Ancestor) \leftarrow \mathtt{find\_common\_ancestor}(\Vertex_1,\Vertex_2)$
                                \State $\mathrm{V}_{new}[\Entity] \leftarrow \mathrm{Q_v}.\mathtt{queue}((\Vertex_1,\Vertex_2),\pos,\Rank,\Ancestor)$
                            \EndIf 
                            \State $\ListOfCMs[\BgElem].\mathtt{add\_vertex\_if\_not\_found}(\mathrm{V}_{new}[\Entity])$
                        \EndIf 
                    \EndFor 
                    \State copy the proximity of vertices on $\Cell$ into an array $\mathrm{P_c}[:] = \mathrm{P}[\Cell.\mathrm{V}[:]]$
                    \State select template using intersected edges \& proximity $T \leftarrow \mathtt{select\_template}(\mathrm{V}_{new} \neq 0, \mathrm{P_c}, \Cell)$ 
                    \For{cell-vertex connectivity for each new cell $\Cell_{new} \in T.\mathtt{generate\_cells}(\Cell.\mathrm{V},\mathrm{V}_{new})$}
                        \State determine new cell's proximity and material $\Material \leftarrow 2 \cdot \Cell.\Material + \mathtt{is\_positive}(\mathtt{vote}(\mathrm{P}[\mathrm{\Cell}_{new}[:]]))$
                        \State \ \ \textit{(\# for vote function, ignore new vertices which are outside the index range of $\mathrm{P}$)}
                        \If{iteration is not the last in loop}
                            \State $\ListOfCMs[\BgElem]\mathrm{V}.\mathtt{append}(\mathrm{Q_c}.\mathtt{queue}(\Cell_{new},\BgElem,\Material))$
                        \Else{ \textit{\# replace information for subdivided cell with one of the new cells}}
                            \State $\Cell.\mathrm{V} \leftarrow \Cell_{new}$, $\Cell.\Material \leftarrow \Material$
                        \EndIf
                    \EndFor 
                \EndFor 

            \Else{ \textit{\# non-cut bg. elements}}
                \State determine proximity of the cell $\Proximity \leftarrow \mathtt{vote}(\mathrm{P}[\BgElem.\mathrm{V}[:]])$
                \For{each of the child cells $\Cell \in \ListOfCMs[\BgElem].\mathrm{C}$}
                    \State increment the material index $\Cell.\Material \leftarrow 2 \cdot \Cell.\Material+ \mathtt{is\_positive}(\Proximity)$
                \EndFor
            \EndIf 
        \EndFor 

        \State $\mathtt{create\_vertices}(\mathrm{Q_v})$
        \State $\mathtt{create\_cells}(\mathrm{Q_c})$
            
    \end{algorithmic}
\end{algorithm*} 

\begin{figure}[btp]
    \centering 
    \def\svgwidth{5.5cm}
    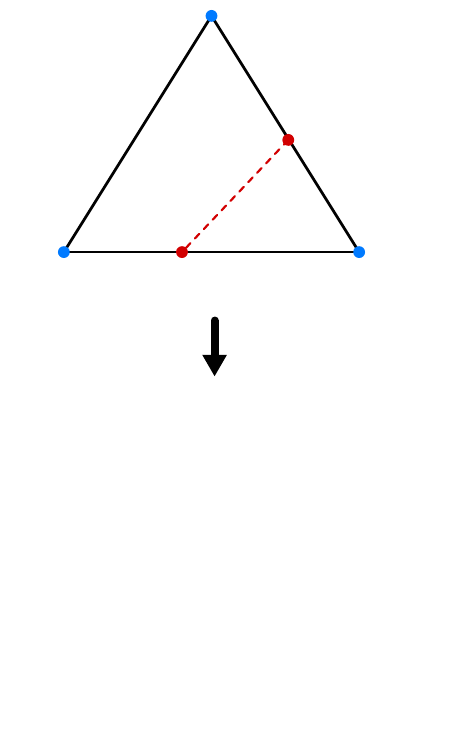
    \vspace*{0.2cm}
    \caption{Recursive material assignment rule \eqref{eqn:phase_assignment_logic} applied during the templated subdivision. The proximity $P$ of each cell after subdivision is obtained by applying a $\mathtt{vote}$ to their attached vertices. The angle brackets $\langle \cdot,\cdots \rangle$ indicate arrays of Booleans which are interpreted as binary numbers to obtain material index $m$.}
    \label{fig:proximity_deduction}
\end{figure}

\paragraph*{Material Assignment}
\label{pg:material_assignment}
The material membership of all cells is determined concurrently with the templated subdivision in a recursive process. We have found that performing material assignment after the subdivision process leads to robustness issues, due to the mismatch of discretized interfaces and the geometries of the LS isocontours.

The sign of the proximity of each cell is stored for every geometry applied. Interpreting the resulting Boolean array as an integer, as shown in \Cref{fig:proximity_deduction}, yields a material index.
Alternatively, this can be stated as a recursive rule applied during the subdivision process:
\begin{equation}
    \label{eqn:phase_assignment_logic}
    m = 2 \cdot m + 
    \begin{cases}
        1, & \text{if} ~ P > 0,\\
        0, & \text{otherwise},
      \end{cases}  
\end{equation}
where $P$ represents the proximity of the cell to the geometry.

Given the restrictions on the possible intersection configurations for the subdivision templates, each of the cells generated by a template can necessarily only have vertices attached to it with non-conflicting proximity values, i.e., $0$ and $+$, or $0$ and $-$, but never $+$ and $-$ or only $0$'s. 
Ignoring vertices with 0, the proximity values of the remaining vertices uniquely determine the proximity of the cell.
\Cref{fig:proximity_deduction} demonstrates this logic for a single subdivided cell in 2D.
The rule \eqref{eqn:phase_assignment_logic} is also applied to the cells without edges intersected by the geometric interface within a given subdivision step. 

In the last step, a user-defined material map is applied to the foreground mesh to merge different material regions, as demonstrated in \Cref{fig:subdivision_process} (D). This avoids the need to construct interface conditions between regions of the same material and physics but separated by an interface of one of the input geometries. It also allows for the exact representation of sharp corners and edges in the interior of background elements, e.g., when defining a square domain using four planes and merging the outside regions.

\subsection{Generation of Topological Information} 
\label{sec:subphase_generation}

\begin{figure}[b]
    \centering 
    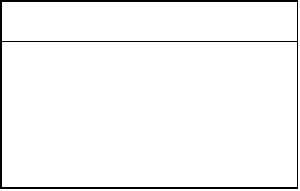
    \caption{data structure for the subphase.}
    \label{fig:subphase_data_structure}
\end{figure}

As discussed in \Cref{sec:framework-ig}, the goal is to evaluate the integrals \eqref{eqn:elemental_bulk_integral}-\eqref{eqn:elemental_interface_integral}.  To this end, the disconnected subdomains $\SetDef{S_E^u}_{u=1}^{n_u(E)}$ within a given background element $E$ need to be identified.
From here on the term "subphase" is used to refer to these disconnected subdomains. In code, a subphase is represented by a list of foreground cells that constitute a connected material subdomain within a background element. 
The associated data structure is shown in \Cref{fig:subphase_data_structure}.

\begin{algorithm*}[btp]
    \caption{
        Generate descendants for the interior bg. facets \\
        Input: fg. mesh $\ForegroundMesh$ (self) with facet connectivity $\mathcal{F}$ \\
        Output: fg. mesh $\ForegroundMesh$ (self), descendant fg. facets of each bg. facet $\mathrm{DF}$}
    \label{alg:generate_bg_facet_descendants}
    \begin{algorithmic}[1]
        \State initialize list of descendant facets for each bg. facet $\mathrm{DF} = \SetDef{\SetDef{}}_{f=1}^{\BackgroundMesh.\mathtt{get\_num\_entities}(d-1)}$
        \For{each fg. facet $\Facet_{fg} \in \FacetConn.\mathrm{CtE}[\Cell]$ on each fg. cell $\Cell \in \FgMesh.\mathcal{C}$}
                \If{$\Facet_{fg}$ has not already been processed}
                    \State get the vertices attached to the current facet $\mathrm{V}_{\Facet} \leftarrow \Cell.\mathtt{get\_entity}(d-1,\Facet_{fg})$ 
                    \State bg. elements vertices are contained in $\mathrm{A} \leftarrow \bigcap_{i=1}^{\mathtt{size}(\mathrm{V}_{\Facet})}\BackgroundMesh.\mathtt{get\_entities\_on\_cell}(\mathrm{V_f}[i].\Rank,\mathrm{V_f}[i].\Ancestor)$
                    \If{$\mathtt{size}(\mathrm{A}) = 2$, i.e., the fg. facet is contained on a bg. facet}
                        \State get the bg. facet $\Facet_{bg} \leftarrow \bigcap_{i=1}^2 \BackgroundMesh.\mathtt{get\_entities\_on\_cell}(d-1,A[i])$
                        \State store that the current fg. facet coincides with bg. facet $\mathrm{DF}[\Facet_{bg}].\mathtt{append}(\Facet_{fg})$
                    \EndIf
                \EndIf
        \EndFor 
        \State \Return{$\mathrm{DF}$}
    \end{algorithmic}
\end{algorithm*} 

Following the discussion in \Cref{sec:framework-discretization} and \Cref{sec:framework-ghost}, connected material regions within each basis function's support need to be identified for the enrichment and ghost stabilization strategies. Additionally, an easy identification of interfaces between each of the material regions is useful.
The preprocessing steps in this subsection generate, besides the subphases themselves, two graphs: one to represent the connection between the subphases $\mathcal{G}_S$ and one to represent adjacent but disconnected material regions $\mathcal{G}_I$. Both graphs are stored in the form of adjacency lists. 

\begin{figure}[btp]
    \centering 
    \def\svgwidth{\columnwidth}
    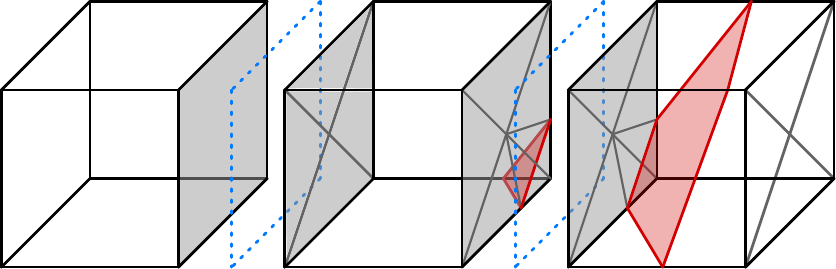
    \vspace*{0.2cm}
    \caption{Example for the foreground facets descending from two background facets in 3D. On the left: non-conformal transition between a non-cut and a cut element due to the regular subdivision template. On the right: conformal transition with an intersecting interface.}
    \label{fig:bg_facet_descendants}
\end{figure}

\paragraph*{Generation of Facet Connectivity Information}
A connectivity for the facets on the foreground mesh is generated using \Cref{alg:compute_entity_connectivity}. The regular subdivision templates in 3D, see \Cref{fig:reg_sub_template} (B), resulting in a foreground mesh with hanging vertices, where subdivided and non-subdivided background elements meet. The left background facet in \Cref{fig:bg_facet_descendants} depicts such a configuration. Due to the mismatching facets, the connectivity information between the adjacent foreground cells needs to be established separately.
The background ancestry of the foreground vertices and the entity connectivity of the background mesh are used to find the missing connections in \Cref{alg:generate_bg_facet_descendants}. The procedure is based on the fact that any foreground facet coinciding with an interior background facet is attached to exactly two background elements. First, the background ancestors for each vertex on each foreground facet are found. If the intersection between the lists of cells attached to the respective ancestors returns exactly two background element indices, the foreground facet must be a subset of the background facet connecting the two background elements. The resulting list of descendant foreground facets $\mathrm{DF}$ for each interior background facet is used to establish the subphase connectivity across background elements.

\begin{figure}[btp]
    \centering 
    \def\svgwidth{6.5cm}
    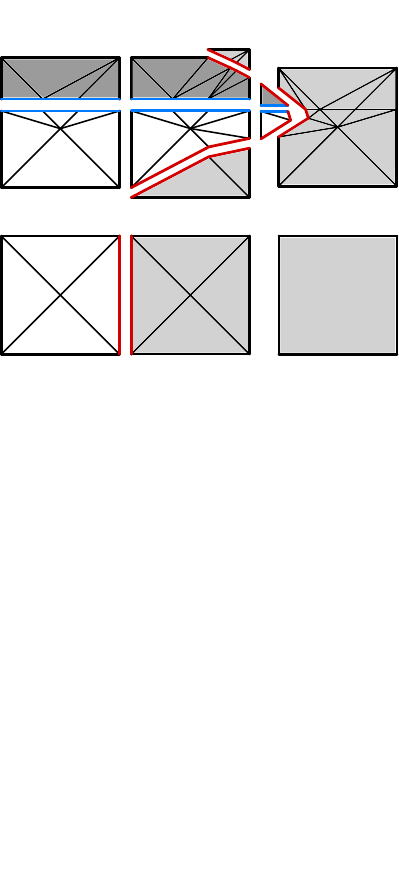
    \vspace*{0.2cm}
    \caption{Generation of subphases and the subphase graph. (A) Identification of disconnected groups of fg. cells in each bg. element using a flood-fill. (B) Subphase interface graph $G_I$ containing neighboring subphases within the background elements. (C) Complete subphase neighbor and interface graphs $\mathcal{G}_S$ and $\mathcal{G}_I$.}
    \label{fig:subphase_generation}
\end{figure}

\paragraph*{Identification and Generation of Subphases}
\label{pg:subphase_generation}
The subphases themselves are created by component analysis. A cell connectivity graph $\mathcal{G}_{local}$ is generated on each child mesh and followed by a flood-fill operation considering their material membership, as stated in \Cref{alg:generate_subphases}. The previously generated facet connectivity $\mathcal{F}$ lends itself to constructing the cell connectivity, as it represents its dual-graph.  
Subphases are initialized from each of the containers of foreground cells returned by the flood-fill. In addition, the index of the subphases is stored with the cells. The subphases resulting from performing these operations on the foreground mesh in \Cref{fig:subdivision_process} (D) are visualized in \Cref{fig:subphase_generation} (A).
Further, the local connectivity information already constructed is used to identify neighboring subphases and store them in the interface graph $\mathcal{G}_I$. As the subphases are disconnected by definition, subphase neighbors inside a background element are never part of the graph $\mathcal{G}_S$. The initialized graph $\mathcal{G}_I$ is shown in \Cref{fig:subphase_generation} (B).

\begin{algorithm*}[btp]
    \caption{
        Generate the subphases. \\
        Input: subdivided fg. mesh $\ForegroundMesh$ (self) \\
        Output: fg. mesh $\ForegroundMesh$ (self) with subphases $\mathcal{S}$ }
    \label{alg:generate_subphases}
    \begin{algorithmic}[1]
        \State initialize interface subphase graph $\mathcal{G}_I \leftarrow \SetDef{\SetDef{}}$ and list of subphases $\mathcal{S} \leftarrow \SetDef{}$ 
        \For{each bg. element $\BgElem \in \BgMesh.\mathcal{\BgElem}$}
            \State generate graph of elements in child mesh $\mathcal{G}_{local} \leftarrow \mathtt{generate\_element\_graph}(\ListOfCMs[\BgElem].\mathrm{C},\mathcal{F})$
            \State get the materials for each of the fg. cells $\mathrm{M}[:] \leftarrow \ListOfCMs[\BgElem].\mathrm{C}[:].m$
            \State get groups of connected fg. cells within bg. element $\SetDef{\mathrm{C}_{\Unzipping}}_{\Unzipping=1}^{n_{\Unzipping}(\BgElem)} \leftarrow \mathtt{flood\_fill}(\mathcal{G}_{local}, \mathrm{M})$
            \For{each connected group of cells $\mathrm{C}_{\Unzipping} \in \SetDef{\mathrm{C}_{\Unzipping}}_{\Unzipping=1}^{n_{\Unzipping}(\BgElem)}$}
                \State create a new subphase $\mathcal{S}.\mathtt{append}(\mathtt{create\_subphase}(\mathrm{C}_{\Unzipping}, \BgElem, \Unzipping))$
                \For{each fg. cell in group $\Cell \in \mathrm{C}_{\Unzipping}$}
                    \State assign the subphase index to the cell $\Cell.\Subphase \leftarrow \mathtt{size}(\mathcal{S})$
                \EndFor
            \EndFor 
            \For{each subphase $\mathrm{C}_{\Unzipping} \in \SetDef{\mathrm{C}_{\Unzipping}}_{\Unzipping=1}^{n_{\Unzipping}(\BgElem)}$}
                \State find subphase indices of adjacent elements $\mathrm{S}_a \leftarrow \mathtt{collect\_subphase\_neighbors}(\mathrm{C}_{\Unzipping}, \mathcal{G}_{local},\mathcal{F})$
                \State add the neighbors to the subphase interface graph $\mathcal{G}_I.\mathtt{append}(\mathrm{S}_a)$
            \EndFor
        \EndFor 
        \State {\color{gray} $\mathtt{communicate\_subphase\_IDs}()$}
    \end{algorithmic}
\end{algorithm*} 

\begin{algorithm*}[btp]
    \caption{
        Generate the subphase graphs. \\
        Input: fg. mesh $\ForegroundMesh$ (self) with subphases $\mathcal{S}$, descendant facets of each bg. facet $\mathrm{DF}$ \\
        Output: finalized fg. mesh $\ForegroundMesh$ (self) with subphase graphs $\mathcal{G}_S$ and $\mathcal{G}_I$ }
    \label{alg:generate_subphase_graphs}
    \begin{algorithmic}[1]
        \State initialize empty subphase graph $\mathcal{G}_S \leftarrow \SetDef{\SetDef{}}_{\Subphase=1}^{\mathtt{size}(\mathcal{S})}$
        \For{each bg. facet $\Facet_{bg} \in \SetDef{1,\cdots,\mathtt{size}(\mathrm{DF})}$}
            \State get the bg. elements connected to bg. facet $\SetDef{\BgElem_1,\BgElem_2} \leftarrow \BackgroundMesh.\mathtt{get\_cell\_on\_entity}(d-1,\Facet_{bg})$
            \If{facet is on exterior $\BgElem_2 = \mathrm{NULL}$}
                \State skip it 
            \State \kern-1.5em \textit{\# if one of the cells is not cut, the material of the bg. facet adjacent cells must be the same}
            \ElsIf{one of the elements is not cut $\mathtt{size}(\ListOfCMs[\BgElem_1].\mathrm{C}) = 1$ OR $\mathtt{size}(\ListOfCMs[\BgElem_2].\mathrm{C}) = 1$}
                \State get the subphases on either side of the bg. facet $(\Subphase_1,\Subphase_2)$
                \State mark the two subphases as neighbors $\mathcal{G}_S[S_1].\mathtt{append}(\Subphase_2)$, $\mathcal{G}_S[\Subphase_2].\mathtt{append}(\Subphase_1)$
            \Else{ \textit{\# both adjacent bg. elements are cut}}
                \For{each fg. facet on bg. facet $\Facet_{fg} \in \mathrm{DF}[\Facet_{bg}]$}
                    \State get the subphases on either side $\SetDef{\Subphase_1,\Subphase_2} \leftarrow \FacetConn.\mathrm{EtC}[:].\Subphase$
                    \If{the neighboring subphases are of the same material $\Subphase_1.\Material = \Subphase_2.\Material$}
                        \State $\mathcal{G}_S[\Subphase_i].\mathtt{append\_if\_not\_found}(\Subphase_j)$, for $(i,j) \in \SetDef{(1,2),(2,1)}$
                    \Else{}
                        \State $\mathcal{G}_I[\Subphase_i].\mathtt{append\_if\_not\_found}(\Subphase_j)$, for $(i,j) \in \SetDef{(1,2),(2,1)}$
                    \EndIf
                \EndFor
            \EndIf 


        \EndFor 
    \end{algorithmic}
\end{algorithm*} 

\paragraph*{Generation of the Subphase Graphs}
For the subsequent construction of $\mathcal{G}_S$ and completion of $\mathcal{G}_I$, only subphase connections across background elements need to be considered. Hence, only foreground facets fully contained on the interior background facets, i.e.,  those in $\mathrm{DF}$, are needed to find the connections. The procedure outlined in \Cref{alg:generate_subphase_graphs} iterates over these facets and stores the subphase connections in the adjacency lists $\mathcal{G}_S$ and $\mathcal{G}_I$.

Two cases need to be considered for each background facet:
\begin{enumerate}
    \item One or both background elements are non-cut, as shown on the left of \Cref{fig:bg_facet_descendants}. In this case, any geometric interfaces are necessarily some non-zero distance away from the background facet, and only one subphase is attached to either side of the background facet. The two attached subphases must be of the same material and, hence, are connected. In \Cref{alg:generate_subphase_graphs}, the two subphases are collected from the cells attached to the facets.
    \item Both background elements are intersected, as shown on the right of \Cref{fig:bg_facet_descendants}. In this case, the foreground mesh is conformal and, by iterating over the foreground facets, the neighboring subphases are collected and stored. For the special case where a geometric interface coincides with the background facet, the connection is stored in $\mathcal{G}_I$ rather than $\mathcal{G}_S$.
\end{enumerate}

\subsection{Enrichment}
\label{sec:enrichment}

\begin{figure}[btp]
    \centering 
    \def\svgwidth{6.2cm}
    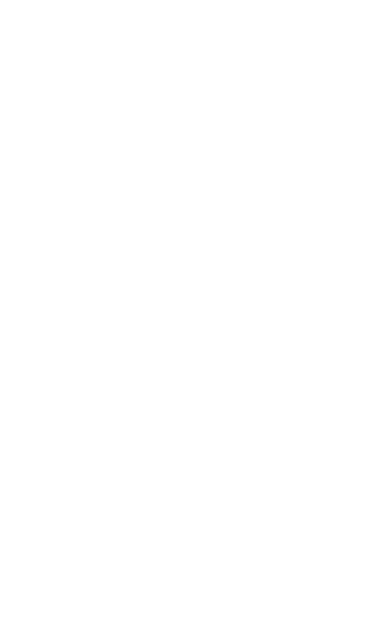
    \vspace*{0.2cm}
    \caption{Basis function enrichment using flood-fill. (A) Restriction of the subphase graph $\mathcal{G}_S$ shown in \Cref{fig:subphase_generation} to the support of a single background basis function $N_{\BfIndex}$. (B) Lists of connected subphases within $\mathrm{supp}(N_{\BfIndex})$ resulting from a flood-fill operation. (C) Support of an example enriched basis function whose support consists of a list of connected subphases.}
    \label{fig:flood_fill}
\end{figure}

To enrich the background basis and construct the enriched basis \eqref{eqn:enriched_basis_function_set}, the supports of the enriched basis functions $\EnrBF_{\BfIndex}^{\EnrLvl}$ are identified. This requires finding where the enrichment functions $\EnrichmentFunction(\pos)$ are non-zero and the number of "enrichment levels" $n_{\EnrLvl}(\BfIndex)$ for a given background basis function $N_{\BfIndex}$. The number of enrichment levels $n_{\EnrLvl}(\BfIndex)$ is computed by performing a flood-fill on the subphase graph within each basis function's support, as outlined in \Cref{fig:flood_fill}. A list of subphases $\mathrm{S}_{\BfIndex}$ within the elements comprising $\mathrm{supp}(N_{\BfIndex})$ is collected by pruning the global adjacency list $\mathcal{G}_S$ to $\mathrm{R}_{\BfIndex} = \mathrm{supp}(N_{\BfIndex})$. 
Lists of subphases connected within $\mathrm{supp}(N_{\BfIndex})$ are then obtained by performing a flood-fill algorithm on the pruned graph $\mathcal{G}_P$.
The number of lists is $n_{\EnrLvl}(\BfIndex)$ and each list corresponds to where $\EnrBF_{\BfIndex}^{\EnrLvl}(\pos) \neq 0$.

\begin{figure}[btp]
    \centering 
    \def\svgwidth{6.0cm}
    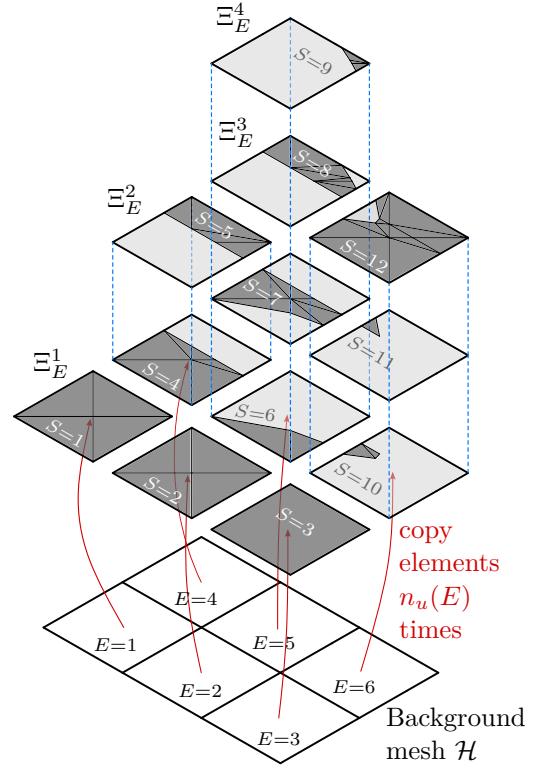
    \vspace*{0.2cm}
    \caption{Unzipping of background mesh. Copies of background elements are created for each of the subphases contained in them, see \Cref{fig:subphase_generation}. The $\mathrm{IEN}$-arrays of the resulting "unzipped" background elements $\Xi_E^u$ are edited to encode the enrichment information.}
    \label{fig:unzipping}
\end{figure}

\paragraph*{Unzipping}
\label{pg:unzipping}
After identifying the enriched basis functions $\EnrBF_{\BfIndex}^{\EnrLvl}(\pos)$, the challenge remains to represent and store the resulting enriched basis on the background elements associated with the final cluster output. 
The strategy chosen here is to represent the enrichment function $\EnrichmentFunction(\pos)$ implicitly by ensuring that a given enriched basis function $\EnrBF_{\BfIndex}^{\EnrLvl}(\pos)$ is only evaluated at quadrature points $\pos \in \mathrm{supp}(\EnrBF_{\BfIndex}^{\EnrLvl})$, where $\BF_{\BfIndex}(\pos) =\EnrBF_{\BfIndex}^{\EnrLvl}(\pos)$.
To achieve this, copies of each of the background elements are created for every subphase within that background element, as shown in \Cref{fig:unzipping}. We refer to this process, outlined in \Cref{alg:unzipping}, as "unzipping". The basis functions of the copied background elements $\SetDef{\Xi_E^u}_{\Unzipping=1}^{\NumUnzippings(\BgElemInd)}$ will, during the element formation, only be evaluated at points within the associated subphase $S_E^u$. Hence, the basis functions $\SetDef{N_{\LocalBfIndex}(\pos)}_{\LocalBfIndex=1}^{n_{\LocalBfIndex}}$ on the copied elements do not need to be altered by adding the associated indicator functions $\psi_{b}^{\EnrLvl}$. 

In a next step, a single index $\EnrBfIndex$ is created to replace the double-index $(\BfIndex,{\EnrLvl})$ used for the functions $\EnrBF_{\BfIndex}^{\EnrLvl}(\pos)$. This allows for \eqref{eqn:enriched_interpolation} to be stated as 
\begin{equation}
    \label{eqn:enriched_interpolation_simplified}
    \DiscDispl(\pos) = 
    \sum_{\EnrBfIndex = 1}^{n_{\EnrBfIndex}} 
    \EnrBF_{\EnrBfIndex}(\pos) \, 
    \dof_{\EnrBfIndex},
\end{equation}
where $n_{\EnrBfIndex}$ is the total number of enriched basis functions $\EnrBF_{\EnrBfIndex}(\pos) = \EnrBF_{\BfIndex}^{\EnrLvl}(\pos)$.
Subsequently, each of the basis function indices $\BfIndex$ listed in the $\mathrm{IEN}$ of each of the copied background elements $\Xi_E^u$ is replaced with the index of the enriched counterpart $\EnrBfIndex \leftarrow (\BfIndex,{\EnrLvl})$ supported within the subphase $S_E^u$. As a result, the enrichment information is fully encoded in the basis function indices in the $\mathrm{IEN}$-arrays and not further exposed to the outside, so, existing element assembly routines apply to the clusters.

The procedure combining both the enrichment of the basis functions and the unzipping of the background elements is presented in \Cref{alg:unzipping}.

\begin{algorithm*}[h!]
    \caption{Enrichment: Unzip background elements \\
        Input: Background mesh $\BgMesh$, Foreground mesh $\FgMesh$ \\
        Output: Unzipped background elements $\Xi$ }
    \label{alg:unzipping}
    \begin{algorithmic}[1]
        \State initialize empty container for unzipped background elements $\Xi \leftarrow \SetDef{}_{\BgElem=1}^{n_{\BgElem}}$
        \For{each background element $\BgElem \in \SetDef{1, \cdots, n_{\BgElem}} $ } 
            \State size inner container $\Xi[\BgElem] \leftarrow \SetDef{}_{\Unzipping=1}^{n_{\Unzipping}(\BgElem)}$
            \For{each subphase in the element $\Unzipping \in \SetDef{1, \cdots, n_{\Unzipping}(\BgElem)} $ }
                \State copy initialize from background element $\Xi[\BgElem][\Unzipping] = \mathtt{create\_element}(\BgMesh.\mathtt{get\_element}(\BgElem))$
            \EndFor
        \EndFor
        \State initialize enriched basis function counter $\EnrBfIndex \leftarrow 0$
        \For{each basis function $\BfIndex \in \SetDef{1, \cdots, n_{\BfIndex}}$ }  
            \State initialize list of subphases in $\BfIndex$'s support $\mathrm{S}_{\BfIndex} \leftarrow \SetDef{}$
            \For{each background element in the support $\BgElem \in \BgMesh.\mathtt{get\_elements\_in\_support}(\BfIndex) $ } 
                \State collect subphases in support $\mathrm{S}_{\BfIndex}.\mathtt{append}( \FgMesh.\mathtt{get\_subphases\_in\_element}(\BgElem) )$
            \EndFor
            \State $\mathcal{G}_{p} \leftarrow \mathtt{prune}(\FgMesh.\mathtt{get\_subphase\_graph}(), \, \mathrm{S}_{\BfIndex} ) $
            \State groups of connected subphases $\SetDef{ R_{\BfIndex}^{\EnrLvl} }_{\EnrLvl=1}^{n_{\EnrLvl}(\BfIndex)} \leftarrow \mathtt{flood\_fill}(\mathcal{G}_{p}) $
            \For{each enrichment level ${\EnrLvl} \in \SetDef{1, \cdots, n_{\EnrLvl}(\BfIndex)} $ }
                \State increment the basis function counter $\EnrBfIndex \leftarrow \EnrBfIndex + 1$
                \For{each subphase in the connected group $\Subphase \in R_{\BfIndex}^{\EnrLvl} $ }
                    \State replace basis function index $\Xi[\Subphase.\BgElem][\Subphase.\Unzipping].\mathrm{IEN}.\mathtt{find\_and\_replace}(\BfIndex,\EnrBfIndex)$ 
                    \State \ \ \ \textit{(\# operation includes safeguard preventing IEN entries from being replaced multiple times)}
                \EndFor
            \EndFor
        \EndFor
        \State {\color{gray}assign basis function IDs and get index to ID map $I_{\EnrBfIndex} \leftarrow \FgMesh.\mathtt{communicate\_basis\_function\_IDs}(\Xi)$}
        \State {\color{gray}replace the basis function indices with IDs in IEN arrays $\mathtt{replace\_IEN\_with\_IDs}(\Xi,I_{\EnrBfIndex})$}
        \State \Return $\Xi$ 
    \end{algorithmic}
\end{algorithm*} 

\paragraph*{Construction of Clusters}

\begin{figure}[btp]
    \centering 
    \def\svgwidth{6.0cm}
    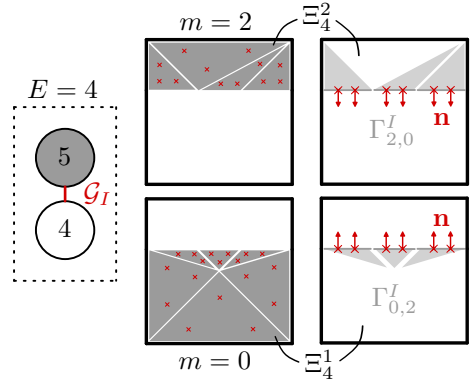
    \vspace*{0.2cm}
    \caption{Construction of bulk and side clusters in background element $E=4$ from \Cref{fig:subphase_generation}. Left: Mapping Gauss points from the foreground cells to the parametric space of the background elements yields quadrature rules for the volume of every subphase. Right: A connection in the interface subphase graph $\mathcal{G}_I$ corresponds to a pair of side clusters constructed with matching Gauss points on every foreground facet.}
    \label{fig:cluster_generation}
\end{figure}

The remaining task is to generate the quadrature points associated with each subphase and package the information up in the cluster data structure introduced in \Cref{sec:preliminaries}. The creation of volumetric \textit{bulk} cluster is discussed first, before the creation of  \textit{side} clusters. An example of clusters generated on a single background element is shown in \Cref{fig:cluster_generation}.

For the bulk clusters, the creation is summed up in \Cref{alg:bulk_cluster_generation}. For every subphase in the foreground mesh, the foreground cells are collected. Using the parametric coordinates of their vertices stored on the corresponding child meshes, a standard Gauss-Legendre quadrature rule is mapped into the parametric coordinates of the background element. Pairing the quadrature rule with the unzipped background element associated with the subphase yields a complete cluster. These clusters are grouped by their respective material domains, such that, during assembly, these groups are iterated over considering the same set of equations associated with the given material domain.

\begin{algorithm*}[btp]
    \caption{Create bulk clusters. \\
        Input: Foreground mesh $\FgMesh$, unzipped background elements $\Xi$ \\
        Output: Sets of bulk clusters for every material $\mathcal{D}_B$ }
    \label{alg:bulk_cluster_generation}
    \begin{algorithmic}[1]
        \State initialize containers of bulk clusters for every material $\mathcal{D}_B \leftarrow \SetDef{\SetDef{}}_{\Material=0}^{n_{\Material}-1}$
        \For{each subphase $\Subphase \in \FgMesh.\mathcal{S}${\color{gray}, for which $\Subphase.\BgElem.\mathtt{is\_owned}()$}}
            \State generate quadrature rule $\SetDef{\xi_q,\omega_q}_{q=1}^{n_q} \leftarrow \FgMesh.\mathtt{map\_param\_gauss\_points}(\Subphase.\mathrm{C}[:])$ 
            \State create bulk cluster $\mathcal{D}_B[S.\mathtt{get\_material}()].\mathtt{append}(\mathtt{create\_cluster}(\Xi[\Subphase.\BgElem][\Subphase.\Unzipping],\SetDef{\xi_q,\omega_q}_{q=1}^{n_q}))$
        \EndFor 
        \State \Return $\mathcal{D}_B$
    \end{algorithmic}
\end{algorithm*} 

\begin{algorithm*}[btp]
    \caption{Create side clusters. \\
        Input: Foreground mesh $\FgMesh$, unzipped background elements $\Xi$ \\
        Output: Sets of side clusters and side cluster pairs for every interface $\mathcal{D}_S$ and $\mathcal{D}_I$}
    \label{alg:side_cluster_generation}
    \begin{algorithmic}[1]
        \State initialize containers of side clusters for every material combination $\mathcal{D}_S, \, \mathcal{D}_I \leftarrow \SetDef{\SetDef{}}_{k=1}^{n_{\Material}^2}$
        \For{each subphase pair $(\Subphase_1,\Subphase_2) \in \FgMesh.\mathcal{G}_I$, for which $\Subphase_1{\color{gray}.I} > \Subphase_2{\color{gray}.I}$ {\color{gray} AND $\Subphase_1.\BgElemInd.\mathtt{is\_owned}()$} } 
            \State find facets between subphases $\mathrm{F}_S \leftarrow \FgMesh.\mathtt{find\_connecting\_facets}(\Subphase_1,\Subphase_2)$
            \For{each side of the interface $(i,j) \in \SetDef{(1,2),(2,1)}$}
                \State get the cells attached to the facets on this side of the interface with the facet ordinals $\mathrm{C},\mathrm{\Ordinal}$
                \State generated quadrature rule $\SetDef{\xi_q,\omega_q,\normal_q}_{q=1}^{n_q} \leftarrow \FgMesh.\mathtt{map\_facet\_gauss\_points}(\mathrm{C}[:],\mathrm{\Ordinal}[:],i)$
                \State \ \ \ (\textit{\# the last input indicates whether (counter-)clockwise ordering is used})
                \State ravel the multi-index $k_i \leftarrow \Subphase_i.\mathtt{get\_material}() \cdot n_{\Material} + \Subphase_j.\mathtt{get\_material}() + 1$
                \State create side cluster $\Cluster_i \leftarrow \mathtt{create\_cluster}(\Xi[\Subphase_i.\BgElem][\Subphase_i.\Unzipping],\SetDef{\xi_q,\omega_q,\normal_q}_{q=1}^{n_q})$
            \EndFor
            \State store side clusters $\mathcal{D}_S[k_1].\mathtt{append}(D_1)$, \ $\mathcal{D}_S[k_2].\mathtt{append}(D_2)$
            \State store side cluster pairs $\mathcal{D}_I[k_1].\mathtt{append}((\Cluster_1,\Cluster_2))$, \ $\mathcal{D}_I[k_2].\mathtt{append}((\Cluster_2,\Cluster_1))$
            \State {\color{gray} collect non-owned subphases \textbf{if} !$\Subphase_2.\BgElem.\mathtt{is\_owned}(): \mathrm{S_{comm}}.\mathtt{append}(\Subphase_2)$}
        \EndFor 
        \State {\color{gray} communicate IENs for non-owned unzipped bg. elements $\mathtt{communicate\_IEN}(\Xi,\mathrm{S_{comm}})$}
        \State \Return $\mathcal{D}_S$, $\mathcal{D}_I$
    \end{algorithmic}
\end{algorithm*} 

To construct the side clusters and cluster pairs to evaluate integral equations of the form \eqref{eqn:elemental_boundary_integral} and \eqref{eqn:elemental_interface_integral}, the steps outlined in \Cref{alg:side_cluster_generation} are necessary. The specific foreground facets $\mathrm{F}_c$ between each neighboring but (materially) disconnected subphase pair $(S_1, S_2)$ are determined. Information about the facet connectivity on the foreground mesh $\mathcal{F}$ and the foreground cells $\mathrm{c}$ that comprise the subphases $S_i$ are employed. The ordinals $o$ of these facets with respect to the foreground cells attached to either side are obtained concurrently. Knowing the facet ordinal and foreground cell on either side, Gauss points on the facets are mapped into the coordinate space of the background elements. In addition, the facet has an outward pointing normal which is stored with the quadrature points. 

With a matched ordering of quadrature points on either side of the interface, the jump terms \eqref{eqn:weighted_jump} and \eqref{eqn:jump_operator} are computed without storing additional information. This is achieved by iterating over the attached foreground facets in the same order from either side and then requesting quadrature points from either side in opposite orientation, i.e., clockwise or counter-clockwise. The quadrature rules generated from either side are matched up to the unzipped background elements associated with the subphases on either side to form the side clusters $D_1$, $D_2$ and pairs $(D_1, D_2)$, $(D_2, D_1)$, before being stored on the list of side cluster (pairs) associated with the interfaces $\Gamma^I_{m,l}$. The bases $\mathcal{B}_E$ belonging to the different unzipped background elements allow the terms \eqref{eqn:weighted_jump} and \eqref{eqn:jump_operator} to be evaluated. 

\begin{figure*}[b]
    \centering 
    \def\svgwidth{14.5cm}
    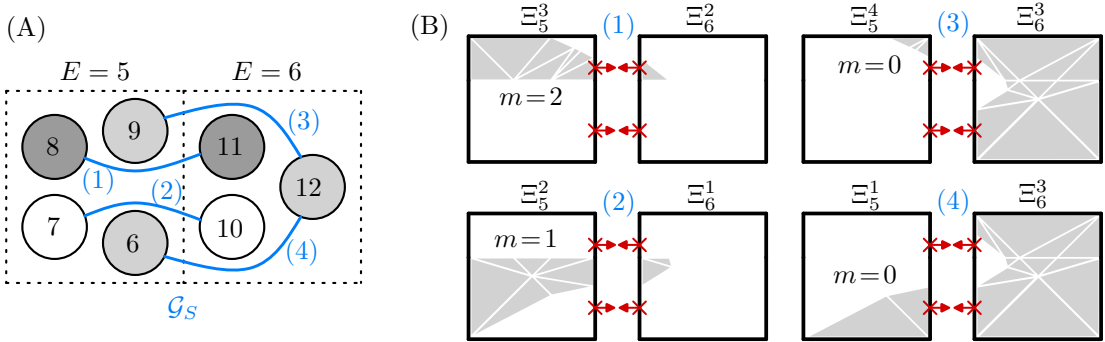
    \vspace*{0.2cm}
    \caption{Ghost side cluster pairs generated from graph. (A) Section of $\mathcal{G}_S$ for a single background facet. (B) Each of the connections in $\mathcal{G}_S$ results in a ghost cluster pair being constructed, since at least one of the adjacent background elements is cut.}
    \label{fig:ghost_cluster_generation}
\end{figure*}

\begin{algorithm*}[b]
    \caption{Construct side cluster pairs for ghost stabilization\\ 
        Input: Background mesh $\BgMesh$, Foreground mesh $\FgMesh$, Unzipped background elements $\Xi$ \\
        Output: Ghost side cluster pairs for every material $\mathcal{D}_G$ } 
    \label{alg:ghost}
    \begin{algorithmic}[1]
        \State initialize empty container for ghost side cluster pairs $\mathcal{D}_G \leftarrow \SetDef{}_{\Material=0}^{n_{\Material}-1}$
        \For{each subphase $\Subphase_1 \in \SetDef{1,\cdots,\FgMesh.\SpGraph.\mathtt{size()}}$}
            \For{each neighbor subphase $\Subphase_2 \in \FgMesh.\SpGraph[\Subphase_1]$}
                \State \textit{\# check if at least one of the adjacent bg. elements is cut and each pair is only constructed once}
                \If{$(\FgMesh.n_{\Unzipping}[\Subphase_1.\BgElem] > 1 \ \text{OR} \ \FgMesh.n_{\Unzipping}[\Subphase_2.\BgElem] > 1) \ \text{AND} \ \Subphase_1{\color{gray}.I} > \Subphase_2{\color{gray}.I} {\color{gray} \ \text{AND} \ \Subphase_1.\BgElem.\mathtt{is\_owned}()}$}
                    \State get the connecting facet ordinals $\SetDef{\Ordinal_{1},\Ordinal_{2}} \leftarrow \BgMesh.\mathtt{get\_connecting\_facets}(\Subphase_1.\BgElem, \Subphase_2.\BgElem)$
                    \For{primary and secondary sides $i \in \SetDef{1,2}$}
                        \State get quadrature points on facet $\SetDef{(\xi_{q},\omega_{q},\normal_q)}_{q=1}^{n_q} \leftarrow \mathtt{get\_quadrature\_points}(\Subphase_i.\BgElem,\Ordinal_{i},i)$
                        \State \ \ \ \textit{(\# the last input indicates whether to output clockwise or counter-clockwise ordering)}
                        \State create side cluster $D_i \leftarrow \mathtt{create\_cluster}(\Xi[\Subphase_i.\BgElem][\Subphase_i.\Unzipping], \SetDef{(\xi_{q},\omega_{q},\normal_q)}_{q=1}^{n_q})$
                    \EndFor
                    \State store cluster pair for material $\mathcal{D}_G[\Subphase_1.\mathtt{get\_material}()].\mathtt{append}((\Cluster_1,\Cluster_2))$
                    \State {\color{gray} collect non-owned subphases \textbf{if} !$\Subphase_2.\BgElem.\mathtt{is\_owned}(): \mathrm{S_{comm}}.\mathtt{append}(\Subphase_2)$}
                \EndIf
            \EndFor
        \EndFor
        \State {\color{gray} communicate IENs for non-owned unzipped bg. elements $\mathtt{communicate\_IEN}(\Xi,\mathrm{S_{comm}})$}
        \State \Return $\mathcal{D}_G$
    \end{algorithmic}
\end{algorithm*}

\paragraph*{Note on Quadrature Rules}
We chose Gauss-Legendre quadrature rules generated on the foreground mesh cells for their simplicity. Other, potentially more efficient, quadrature rules could be generated using the existing subphase information. This could include Stokes' theorem-based methods \cite{gunderman2021spectral, gunderman2021high} or moment-fitting \cite{sudhakar2013quadrature, jiang2020ceramic, muller2013highly}, both on the volume of the subphases directly or on their boundaries using Stokes' theorem.

\subsection{Face-Oriented Ghost Stabilization}
\label{sec:ghost}

This subsection details the implementation of the residual term \eqref{eqn:ghost_residual} for face-oriented ghost stabilization. From a theoretical point of view, the challenge here is twofold. The first is to identify the ghost facets $\mathcal{F}_G$ and the adjacent material regions $(u,v): u \in \SetDef{1,\cdots,n_u(E^+)}, v \in \mathcal{U}_{u,E^-}$ for each ghost facet $F \in \mathcal{F}_G$, that is finding the sets of the first three sums in \eqref{eqn:ghost_residual}. The second is to construct the polynomial extensions $\DiscTestDisplExt$ and $\DiscDisplExt$ and evaluate the jumps in their derivatives at the facets $F \in \mathcal{F}_G$.

The triple sum in \Cref{eqn:ghost_residual} is restated as a sum over all subphase neighbors in $\mathcal{G}_S$, where at least one of the background elements contains two or more subphases, i.e., it contains a material interface. \Cref{fig:ghost_cluster_generation} (A) shows an example of this for a single facet $F \in \mathcal{F}_G$ in the mesh in \Cref{fig:subphase_generation}. On this facet, there are four subphase pairs which the second and third sum in \eqref{eqn:ghost_residual} amount to.
During the unzipping process, outlined in \Cref{sec:enrichment}, the basis functions on the background elements $\mathcal{B}_E$ are not changed to include the enrichment function $\EnrichmentFunction$.  
Hence, the basis already represents the polynomial extension $\DiscTestDisplExt$ and $\DiscDisplExt$ without further modification. Pairing up the unzipped background elements corresponding to each of the subphases and adding matching integration points to either side, as is illustrated in \Cref{fig:ghost_cluster_generation} (B), allows for convenient evaluation of the jump terms. 

\Cref{alg:ghost} outlines the compact procedure to generate the cluster pairs. Analogous to the side cluster construction, a standard Gauss-Legendre quadrature rule is generated on the relevant background element facets with mirrored spatial ordering on either side of each facet. 
Note, that this is done using the background element's facet directly, as opposed to collecting the attached foreground cells, reducing the number of quadrature points.

\subsection{Parallelization}
\label{sec:parallel}

This subsection addresses the parallel implementation of the presented preprocessing framework. First, the general approach taken to parallelization is discussed before covering the additional considerations for executing each of the previously presented steps in parallel.

The preprocessor uses a domain decomposition approach.
The processors independently execute the presented algorithms on parts of the domain assigned to them, allowing for efficient use of distributed memory systems. 

\begin{figure}[btp]
    \centering 
    \def\svgwidth{6.5cm}
    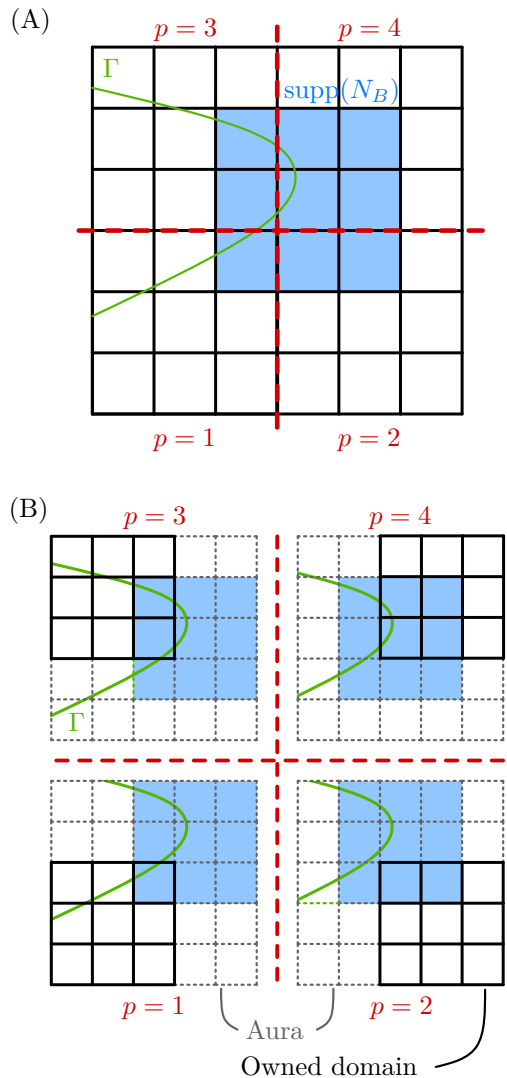
    \vspace*{0.2cm}
    \caption{Domain decomposition in parallel. (A) Global background mesh being split across four processors with the support of one basis function $N_{\BfIndex}$ highlighted in blue. (B) Each processor has access to a layer of background elements, the "aura", beyond the owned domain such that the support of any basis function in the owned domain can be analyzed for enrichment.}
    \label{fig:aura}
\end{figure}

An extra layer of elements, referred to as the "aura", is added beyond the part of the background mesh assigned to, or "owned" by, a given processor. This is to keep parallel communication on such systems to a minimum. Due to the enrichment strategy requiring the topological analysis of the full support of each basis function, the aura's width is such that the support of any basis function, partially supported within a given processor's owned domain, is fully contained in the owned domain and the aura. 
The resulting parallel decomposition scheme is illustrated in \Cref{fig:aura} using a domain split into four parts and a single basis function. Assuming all basis functions have a support size of $3 \times 3$ background elements, like the one shown in the figure, the aura needs to extend by two elements beyond the domain. 

An advantage of this approach is that the amount of data communicated between each of the independent processors is minimal. As such, only indexing information within the aura is communicated, as discussed in the next paragraph. While, from a hardware perspective, this allows for efficient scaling across multiple compute nodes, a weakness becomes apparent from the illustration in \Cref{fig:aura}. If the size of the processors' domains is very small and the aura constitutes a considerable fraction of a mesh handled by a given processor, the computational overhead of processing the aura is large. We are quantifying this overhead 
in the scaling study presented in \Cref{sec:examples-nTop_sandwich}.

\begin{algorithm}[btp]
    \caption{Communicate IDs. \\ 
        Inputs: My proc. $p_c$, list of adjacent procs. $\mathrm{P}_{loc}$ }
    \label{alg:comm}
    \begin{algorithmic}[1]
        \State \textit{\# Assign IDs to owned entities}
        \State initialize $\mathrm{O}, \mathrm{N} = \SetDef{}$
        \For{each entity $e$}
            \If{$\mathtt{is\_owned}(e)$}
                \State $\mathrm{O}.\mathtt{append}(e)$
            \Else{ \textit{\# entity is not owned}}
                \State $\mathrm{N}.\mathtt{append}(e)$
            \EndIf
        \EndFor
        \State $\SetDef{n_{\Entity}^i}_{i=1}^{\mathtt{NumProcs}} \leftarrow \mathtt{MPI\_Allgather}(\mathtt{size}(\mathrm{O}))$
        \State first ID $I \leftarrow 1 + \sum_{i=1}^{\Proc_c} n_{\Entity}^i$ 
        \For{each owned entity $e \in \mathrm{O}$}
            \State $\Entity.\ID \leftarrow \ID++$
        \EndFor
        \State 
        \State \textit{\# Send and receive ID requests}
        \State initialize $\mathrm{R},\mathrm{E},\mathrm{A} \leftarrow \SetDef{\SetDef{}}_{\Proc \in \LocalCommTable}$
        \For{each non-owned entity $\Entity \in \mathrm{N}$}
            \State find owning proc. $\Proc_o$
            \State get identifying information $i_d$
            \State store request $\mathrm{R}[\Proc_o].\mathtt{append}(i_d)$
            \State ... and entity $\mathrm{E}[\Proc_o].\mathtt{append}(\Entity)$
        \EndFor
        \State send requests $\mathtt{MPI\_Isend}(\SetDef{\mathrm{R}_{\Proc}}_{\Proc \in \LocalCommTable}, \LocalCommTable)$
        \State receive $\SetDef{\mathrm{R}_{\Proc}}_{\Proc \in \LocalCommTable} \leftarrow \mathtt{MPI\_Irecv}(\LocalCommTable)$
        \State 
        \State \textit{\# Answer ID requests}
        \For{each proc. neighbor $\Proc \in \LocalCommTable$}
            \For{each request $r \in \mathrm{R}_{\Proc}$}   
                \State find entity $\Entity \leftarrow \mathtt{find\_entity}(r)$
                \State store ID answer $\mathrm{A}[\Proc].\mathtt{append}(\Entity.\ID)$
            \EndFor
        \EndFor
        \State send answers $\mathtt{MPI\_Isend}(\SetDef{\mathrm{A}_{\Proc}}_{\Proc \in \LocalCommTable}, \LocalCommTable)$
        \State receive $\SetDef{\mathrm{A}_{\Proc}}_{\Proc \in \LocalCommTable} \leftarrow \mathtt{MPI\_Irecv}(\LocalCommTable)$
        \For{each proc. neighbor $\Proc \in \LocalCommTable$}
            \For{each answer $i \in \SetDef{1,\cdots,\mathtt{size}(\mathrm{A}_{\Proc})}$}
                \State store the ID $\mathrm{E}[\Proc][i].I \leftarrow \mathrm{A}_{\Proc}[i]$
            \EndFor
        \EndFor
    \end{algorithmic}
\end{algorithm}

\paragraph*{Concept of Parallel IDs}
\label{pg:parallel_IDs}
Entities, that are present on two or more processor domains simultaneously, need to be handled in a parallel consistent manner. As the indices of entities differ between processors, a parallel consistent, and globally unique ID is assigned to each entity. Further, it is assumed that any entity is owned by exactly one processor.
Upon creation of entities of any type, \Cref{alg:comm} is used to assign IDs to the owned entities and communicate IDs with neighboring processors that have access to those entities but do not own them. The first step of the procedure, i.e., the assignment of IDs to owned entities, only requires the communication of the number of owned entities by each processor. In the second step, information is collected for non-owned entities with which they are identified unambiguously by the owning processor. The content of this information is discussed in the following. In the last step, using the identifying information received from adjacent processor domains, the IDs of the associated entities are collected and sent back to the requesting processors. The communication calls in \Cref{alg:comm} are adapted from the \texttt{OpenMPI} library \cite{gabriel2004open}. 

\paragraph*{Background Mesh}
For the presented algorithms to work in parallel, the data structure of the input mesh, shown in \Cref{fig:IO_data_structures}, needs to meet the following requirements:
\begin{itemize}
    \item The list of elements $\mathcal{E}$ provided to a given processor $p$ includes both the owned and the aura elements. Each element $E \in \mathcal{E}$ stores its owning processor $p$ and has a parallel ID assigned to it.
    \item The basis functions provided by the mesh and their indexing are restricted to those basis functions that are at least partially supported within the owned domain. Basis functions only supported within the aura are ignored. The $\mathrm{IEN}$ arrays of the aura elements are left empty.
    \item All basis functions are assigned an owner and a parallel consistent ID.
\end{itemize}

\paragraph*{Foreground Mesh}
The foreground cells and the subphases need parallel consistent IDs. 
The parallel consistent subdivision in each processor's aura is ensured by the order in which the $\mathtt{for}$ loops in \Cref{alg:templated_subdivision} are executed.
After the regular subdivision templates are applied in identical orientation by each processor, the foreground cells generated during the subdivision are stored in identical order and with identical orientation on the child meshes across the processors. Knowing the ID of the background element a child mesh is located in and the local index of a given cell inside the child mesh, cells are identified across processors and their IDs communicated in the function $\FgMesh.\mathtt{communicate\_cell\_IDs}()$ in \Cref{alg:driver_algorithm}.
Subsequently, subphases are identified by providing the ID of one of the foreground cells as part of the function $\FgMesh.\mathtt{communicate\_subphase\_IDs}()$ in \Cref{alg:generate_subphases}.

\paragraph*{Enrichment}
The enriched basis functions require parallel consistent IDs to assemble the global system of equations. As part of the enrichment process in \Cref{alg:unzipping}, one of the subphases comprising the support of an enriched basis function is collected to identify enriched basis functions across processors. The basis function indices listed in the IEN arrays of the unzipped background elements are replaced with their respective parallel IDs.

\paragraph*{Clusters}
Generally, only clusters and cluster pairs contained within a processor's domain are constructed to ensure each cluster or cluster pair's contribution is only added once to the global linear system during assembly. For side cluster pairs, particularly for those constructed for the purpose of ghost stabilization, situations may arise, where the two adjacent unzipped background elements are part of different processors' domains. In this case, the adjacent subphases are collected in Algorithms \ref{alg:side_cluster_generation} and \ref{alg:ghost}, and the content of the IENs of the associated unzipped background elements is communicated. Recall, that the IEN of background elements within the aura was initially left empty.

\newpage
\vspace*{8.0cm}
\begin{figure*}[btp]
    \centering 
    \def\svgwidth{13.6cm}
    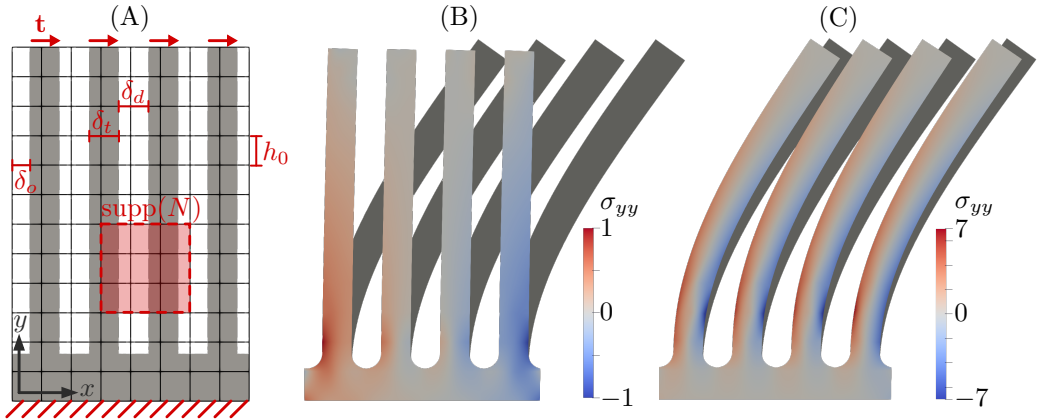
    \caption{(A) Setup of the example with multiple connected beams including the background mesh and the support of a single basis function. (B) Normal stresses in the solution to the non-enriched problem and (C) the problem using the enrichment strategy \eqref{eqn:enriched_interpolation} with fillets added to the interior corners. The reference displacement solution generated using a highly refined mesh is shown in gray in the background.}
    \label{fig:multi_beam_setup}
\end{figure*}
\newpage

\section{Numerical Examples}
\label{sec:examples}

This section presents a series of numerical examples which study the features and performance of the preprocessing framework and its implementation specifically. For a detailed analysis of the numerical framework's performance, e.g., its convergence behavior and stability properties, the reader is referred to \cite{noel2022xiga}.
The first example in \Cref{sec:examples-multi_beam} is used to validate and discuss the impacts of the enrichment strategy and the ghost stabilization for a geometry with small features. The "Brick Wall" example in \Cref{sec:examples-brick_wall} demonstrates the preprocessing framework's robustness for pathological geometric cases with multiple geometries, sharp features, vanishing material domains, and coinciding interfaces. Lastly, we demonstrate the performance and parallel scalability with a "scaffold sandwich" structure in \Cref{sec:examples-nTop_sandwich}.

The condition numbers presented were estimated using either Matlab's $\mathtt{cond}()$ function for small systems or by solving for the smallest and largest eigenvalues using a Krylov-Schur Method \cite{stewart2002krylov}, as implemented in the $\mathtt{SLEPc}$ package \cite{hernandez2007krylov}.
To solve linear systems, the direct solver implemented in $\mathtt{PARDISO}$ \cite{schenk2004solving,schenk2006fast} and a flexible generalized minimal residual solver (FGMRES) \cite{saad1993flexible} as implemented in the $\mathtt{PETSc}$ library \cite{balay1997efficient} were used.

\subsection{Multiple Thin Beams}
\label{sec:examples-multi_beam}

We consider the setup shown in \Cref{fig:multi_beam_setup} (A) with multiple connected beams under a bending load. The beams are immersed in a coarse background mesh with element size $h_0=1$. The material has a Poisson ratio of $\nu=0.3$ and Young's modulus of $E=1000$. A traction load $\traction=(1,0)$ is applied in $x$-direction. For the initial setup, the beams have a width of $\delta_t = h_0$, are spaced apart by $\delta_d=h_0$, and, to avoid generating a body-fitted mesh, are offset by $\delta_o=0.6\,h_0$ relative to the background mesh. The problem is discretized using quadratic $C^1$-continuous B-spline basis functions. 

In this configuration, without basis function enrichment\footnote{The numerical framework without basis function enrichment is equivalent to CutFEM \cite{burman2015cutfem} with B-spline basis functions, sometimes also referred to as CutIGA \cite{elfverson2018cutiga}.}, the support of individual basis functions spans across multiple features, as shown in \Cref{fig:multi_beam_setup} (A). This results in severe artificial stiffening, or "cross-talk", as the basis functions incorrectly couple the individual beams. This results in the normal stress, $\sigma_{yy}$, being underpredicted by a factor of approximately 7, see \Cref{fig:multi_beam_setup} (B). The enrichment strategy \eqref{eqn:enriched_interpolation} remedies this issue by effectively creating separate enriched basis functions for each beam. Even with the excessively coarse background mesh, the enriched solution follows the overall stress response of the reference solution created using a highly refined mesh, as shown in \Cref{fig:multi_beam_setup} (C).

\begin{figure}[h!]
    \vspace*{0.4cm}
    \centering
    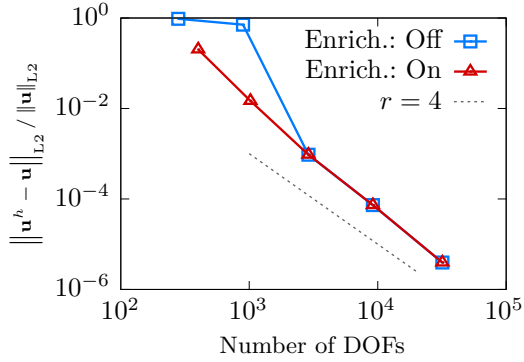
    \caption{Relative L2-error with increasing mesh resolution with and without the generalized enrichment strategy using quadratic $p=2$ B-splines. Note that the convergence rate $r$ is double the slope when plotted against the number of DOFs in 2D.}
    \label{fig:multi_beam_L2_Error_over_nDof}
\end{figure}

\begin{figure}[h!]
    \centering
    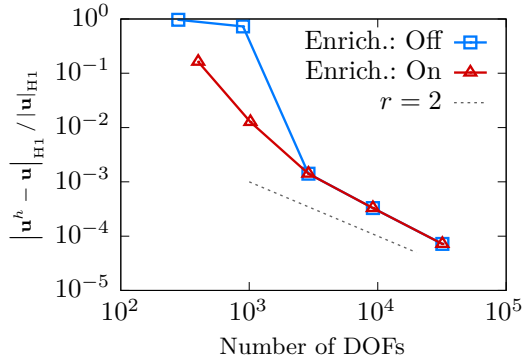
    \caption{Relative H1-semi-error with increasing mesh resolution with and without the generalized enrichment strategy using quadratic $p=2$ B-splines.}
    \label{fig:multi_beam_H1s_Error_over_nDof}
\end{figure}

\begin{figure}[h!]
    \centering 
    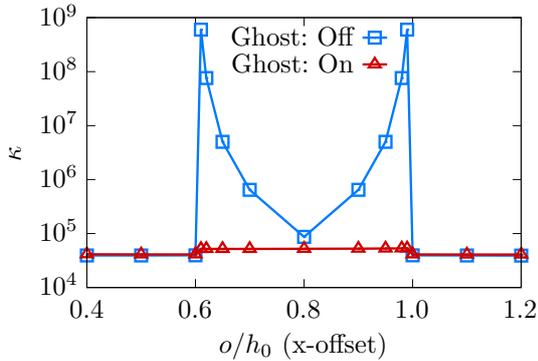
    \caption{Condition number $\kappa$ of the linear system with and without the face-oriented ghost stabilization applied as the structure is translated in x-direction through the background mesh. A penalty of $\GhostPenalty^k = 0.01 \cdot E h_0^{2k-1}$ is used in the stabilized case.}
    \label{fig:Ghost_cond_over_offset}
\end{figure}

The impacts of the enrichment strategy can further be quantified by comparing the error in the solution with and without enrichment under quad-tree, i.e., isotropic, background mesh refinement. Omitting enrichment, both the errors in the displacements, see \Cref{fig:multi_beam_L2_Error_over_nDof}, and in the stresses, see \Cref{fig:multi_beam_H1s_Error_over_nDof}, drastically decrease as the size of the support of basis functions shrinks to below the distance between the beams $\delta_d$; only once the mesh is sufficiently refined the errors converge optimally\footnote{To mitigate stress singularities and obtain optimal convergence rates, fillets are added to the interior corners, as shown in \Cref{fig:multi_beam_setup} (B),(C). Additionally, the foreground mesh is refined to limit the geometric error.}. 
Employing the proposed enrichment strategy, we do not observe the initial region of high error due to under-resolved features. Note that the enriched problem has a greater number of degrees of freedom compared to the non-enriched problem on an equal background mesh with element size $h$. For higher mesh resolutions, the enrichment has no impact on the number of degrees of freedom and the error in the solution.
The enrichment, in essence, adds anisotropic refinement in regions where, and only where, the effects of the material topology would otherwise not be captured by the background mesh. Hence, the generalized Heaviside enrichment strategy loosens the requirements on the mesh resolution relative to the feature size in a given problem.

\paragraph*{Ghost Stabilization}
We investigate the impact of the ghost stabilization on the enriched problem by considering the setup shown in \Cref{fig:multi_beam_setup} (A). Compared to the previous study, the beam spacing is decreased to a size smaller than the background element size $\delta_d = 0.6 \, h_0$, and the beam width is increased to $\delta_t = 1.4 \, h_0$, such that multiple beams connect in the same background element. The configuration requires the construction of multiple ghost penalty terms on the same facet, as shown in \Cref{fig:ghost_clusters}. The position of the beams $\delta_o$ relative to the background mesh is then shifted to obtain various cut configurations. The impact on the condition number is shown in \Cref{fig:Ghost_cond_over_offset}. The condition number for the non-stabilized system spikes as the beams' edges transition through the background mesh lines, while the stabilized system only experiences small changes in the condition number due to the double penalization of a single facet leading to an increased contribution of the ghost penalty term. 
The conditioning of the linear system remains unaffected by the location of intersections besides this small perturbation. This aligns with the expected behavior for ghost stabilization and indicates that the adapted ghost stabilization scheme successfully targets any insufficiently supported enriched basis functions.


\begin{figure}[btp]
    \centering 
    \def\svgwidth{7.0cm}
    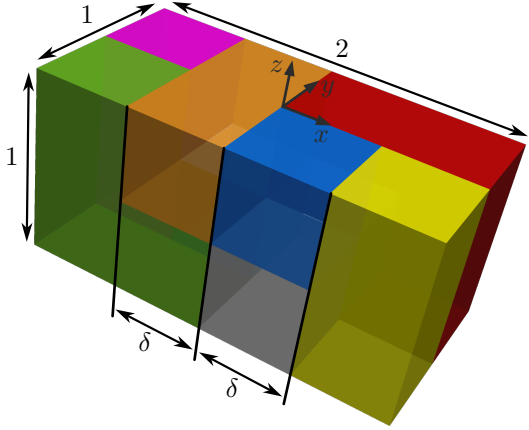
    \vspace*{0.2cm}
    \caption{Geometry of the "Brick Wall" example.}
    \label{fig:brick_wall_setup}
\end{figure}

\subsection{Brick Wall}
\label{sec:examples-brick_wall}

This example aims to demonstrate the robustness of the preprocessing framework in a multi-material setting, specifically in terms of handling geometric edge cases and in terms of retaining good conditioning of the linear system. To this end, a multi-material geometry involving straight edges, vanishing material domains, and coinciding interfaces, is considered and randomly placed within a background mesh.

\begin{figure}[btp]
    \centering
    \def\svgwidth{4.3cm}
    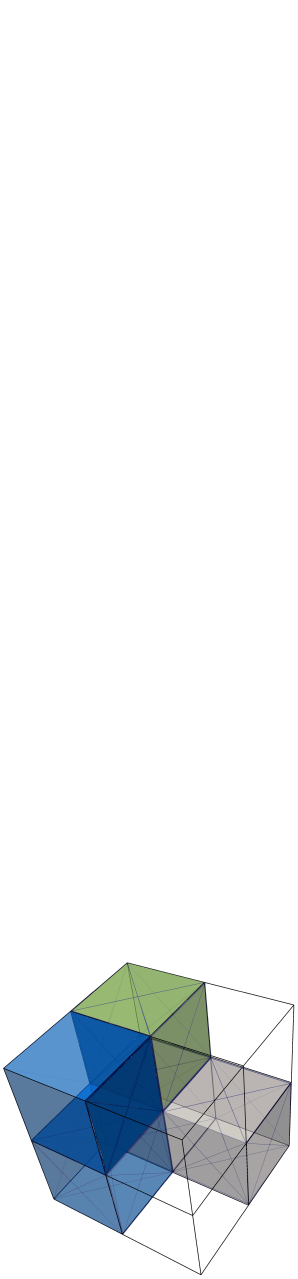
    \caption{Edge cases for cut configurations and resulting foreground mesh on a $2 \times 2 \times 2$ element section of the bg. mesh, shown in black. 
    (A) Intersections align with the centers of bg. elements. (B) Intersections align with the edges of bg. elements. (C) Rotation of the geometry in the bg. mesh with interfaces coinciding at a background vertex. (D) Small offset of interfaces to the bg. mesh resulting in slim material slivers. }
    \label{fig:brick_wall_cut_configs}
\end{figure}

The geometric setup consists of a $2 \times 1 \times 1$ cuboid cut in half in each spatial dimension by three axis-aligned planes. Additionally, two planes are placed symmetrically at varying positions $x=-\delta$ and $x=+\delta$ from the center of the cuboid. Some of the resulting subdomains are then joined together to generate material blocks, or "bricks", of various shapes, as shown in \Cref{fig:brick_wall_setup}.
Examples of various edge cases in the alignment of the geometry with the background mesh and the resulting foreground meshes generated by the preprocessor are shown in \Cref{fig:brick_wall_cut_configs}. 

\begin{figure}[btp]
    \centering
    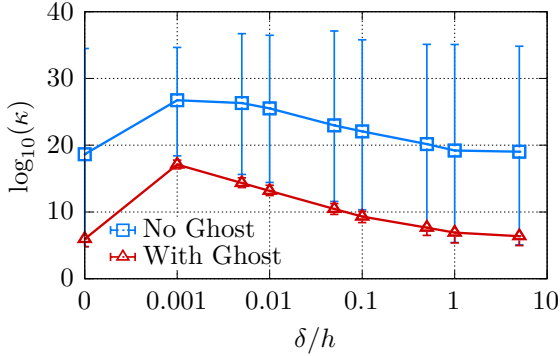
    \caption{Mean and range of $\mathrm{log}_{10}(\kappa)$ of the linear systems resulting from the "Brick Wall" example for varying cut distances $r$ and random translations of the geometry inside the background mesh. The stabilization penalty is set to $\GhostPenalty^k = 0.01 \cdot E h^{2k-1}$. Note, the left-most entry violates the scale and is $\delta/h = 0$.}
    \label{fig:brick_wall_cond}
\end{figure}

The impact of the geometry's position in the background mesh on the conditioning of the linear system is determined by immersing the geometry in a Cartesian grid background mesh with element size $h = 0.1$, such that for $\delta = 5h$ all interfaces align with the background mesh lines. From this starting configuration, the distance $\delta$ is reduced to $0$ to generate increasingly smaller material domains. In addition, the geometry is translated randomly in all coordinate directions $100$ times for each value of $\delta$, where the translation in each spatial dimension is given by
\begin{equation}
    \label{eqn:offset_distribution}
    \Delta x_i = \mathrm{sign}(Z) \cdot 10^{-|Z|} \cdot h,
\end{equation}
with $Z$ being a random number drawn from a normal distribution with mean $0$ and standard deviation $1$.
Quadratic $C^1$-continuous B-spline basis functions are used for discretization, and random, uniformly distributed, Young's moduli $E \in [1,5]$ and Poisson ratios $\mu \in [0.25,0.35]$ are assigned to the different material blocks. The whole structure is clamped at $x=-1$ and a traction load of $\traction = (1,1,1)$ is applied at $x=-1$. The condition number ranges of the linear system resulting from the different cut configurations, with and without ghost stabilization, are shown in \Cref{fig:brick_wall_cond}.

The results demonstrate that the problem exhibits a conditioning behavior expected for the face-oriented ghost stabilization: the stabilization method keeps the conditioning close to the best-case cut configuration for the given geometry, independent of the geometry's position in the background mesh. For small values of $\delta$, thin slices of material are created which are only supported by cut background elements. Hence, the ghost stabilization is not able to stabilize cut basis functions against non-cut basis functions and prevent a slight increase in the condition number $\kappa$. 
The system remains solvable for the considered material sliver sizes with the linear solvers employed in this study. However, this may not remain true for smaller sliver sizes. This observation informs the choice of numerical tolerances that should be used in the $\mathtt{compute\_proximity}$ and $\mathtt{find\_interface}$ functions in \Cref{alg:templated_subdivision}. Tolerances for the proximity of a vertex to be considered "on the interface" ($P(v) = 0$) should be chosen large enough to ensure numerical stability in geometric edge cases, like the one presented, as opposed to only being large enough to avoid errors due to machine precision arithmetic.


\begin{figure*}[btp]
    \centering 
    \def\svgwidth{14.0cm}
    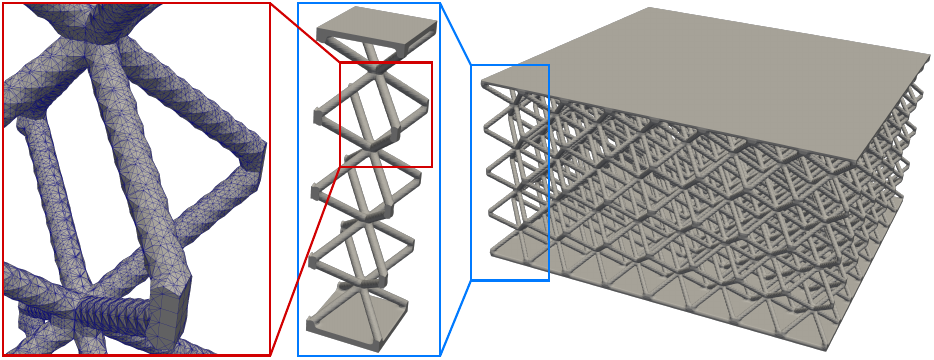
    \vspace*{0.2cm}
    \caption{Geometry for the scaffold sandwich. (A) Largest problem size constructed from an $8 \times 8$ grid of unit blocks. (B) Single unit block which is periodically repeated to construct the example. (C) Detail of the tessellation generated during the decomposition. }
    \label{fig:nTop_sandwich_setup}
\end{figure*}

\subsection{Scaffold Sandwich}
\label{sec:examples-nTop_sandwich}

We consider the scaffold structure with thin geometric features, shown in \Cref{fig:nTop_sandwich_setup} (A), to investigate the computational performance and scalability of the preprocessing framework. The structure consists of a periodic grid of unit blocks shown in subfigure (B). This problem setup is chosen as the size of the numerical model can be altered by increasing the number of unit blocks, whereas the global geometric complexity -- the portion of background elements intersected by an interface, the number of enriched basis functions generated relative to the number of basis functions on the background mesh, the density of the skeleton mesh of ghost facets, etc. -- stays constant. For the following parallel scaling studies, the domain is only decomposed in $x$- and $y$-directions to ensure that the geometry near processor interfaces remains similar. Further, the effect of the parallel domain decomposition on the scaling behavior due to unequally distributed geometric complexity is mitigated by ensuring that each processor's domain contains approximately the same geometry.

Each unit block is discretized using a part of the background mesh with $24 \times 24 \times 96 $ elements and quadratic $C^1$-continuous B-spline basis functions. Each unit block contains $\sim \! 7400$ intersected background elements and $\sim \! 530,000$ cells in the foreground mesh\footnote{A voxel-based geometry definition for the entire structure is used which does not align exactly with the background mesh; hence, the geometry and number of cut elements in each unit block may differ slightly.}. A detailed view of the tessellation generated is shown in \Cref{fig:nTop_sandwich_setup} (C). 

Three types of scaling are investigated. 
\begin{enumerate}
    \item \textbf{Time complexity:} The problem size is increased by increasing the number of unit blocks while keeping the number of processors constant. 
    \item \textbf{Weak parallel scaling:} The problem size and the number of processors are increased simultaneously. The weak scalability is an indicator of how well the framework performs in scenarios where both the computational resources available and the problem size are large.
    \item \textbf{Strong parallel scaling:} The number of processors is increased while keeping the problem size constant. The strong scalability is a crucial factor in how viable the framework is in use cases where the solution time is critical, and a large number of processors is available. Examples of such scenarios are optimization or design studies during which many designs are analyzed in rapid succession.
\end{enumerate}

\paragraph*{Quantifying Scalability}
The generalized efficiency measure $\eta$ is used to quantify the scalability, i.e., how efficiently the framework uses computational resources in the three scaling scenarios. 
\begin{equation}
    \label{eqn:efficiency_measure}
    \eta = \frac{\tfrac{s}{s_0}}{\tfrac{n_p}{n_{p,0}} \cdot \tfrac{t}{t_0}},
\end{equation}
where $s$, $n_p$, and $t$ denote the problem size, the number of processors used, and the time to perform the various preprocessing steps, respectively. The subscript $(\cdot)_0$ denotes the baseline values relative to which the efficiency is evaluated. An efficiency value of $\eta = 1$ indicates perfect scaling. The measure is generalized in the sense that it is applicable to evaluate strong and weak parallel scalability, as well as problem size scaling assuming that the algorithms in question scale linear in time, $t = \mathcal{O}(s)$. Additionally, the memory efficiency is defined as 
\begin{equation}
    \label{eqn:memory_efficiency}
    \mu = \frac{\tfrac{s}{s_0}}{\tfrac{M}{M_0}},
\end{equation}
where $M$ denotes the memory consumption.

We define two mesh overhead values to better understand the scaling results and how they are impacted by the chosen parallel strategy. 
The processor local mesh overhead $\lambda_{loc}$ is defined as the maximum ratio across all processors of background elements on a processor including the aura over the number of elements owned by it 
\begin{equation}
    \label{eqn:local_mesh_overhead}
    \lambda_{loc} =  \max_{ p \in \mathrm{P}_{glob} }
        \left( 
            \frac{ \text{\# bg. elems. on } p }{ \text{\# bg. elems. owned by } p } 
        \right).
\end{equation}
Meanwhile, the global mesh overhead $\lambda_{glob}$ is the total number of background elements processed across all processors -- counting aura elements multiple times as they are shared between different processors -- over the actual number of background elements.
\begin{equation}
    \label{eqn:global_mesh_overhead}
    \lambda_{glob} = \frac{ 
            \sum_{ p \in \mathrm{P}_{glob} } \text{\# bg. elems. on } p 
        }{ 
            \sum_{ p \in \mathrm{P}_{glob} } \text{\# bg. elems. owned by } p 
        }.
\end{equation}
The efficiency measures $\eta$ and $\mu$ are compared to the local or global mesh overheads relative to the baseline configurations $\lambda' = \lambda / \lambda_0$.

\paragraph*{Hardware and Measurements}
All results were run on a single server node using a single AMD EPYC\textsuperscript{\texttrademark} 7713 64-core processor with an eight-channel DDR4 memory configuration. Dynamic frequency scaling was disabled to prevent any parallel efficiency losses due to clock-speed variations.

The timing and memory consumption of the presented preprocessor is split into three sections: the \textit{decomposition}, \textit{enrichment}, and \textit{ghost facet generation} consisting of the steps described in \Cref{sec:fg_mesh_generation} and \Cref{sec:subphase_generation}, \Cref{sec:enrichment}, and \Cref{sec:ghost}, respectively.

The memory consumption is measured using the memory allocation reported by the operating system. This measurement slightly overestimates the actual amount of memory used, though this error is negligible when measuring memory consumption in the decomposition and enrichment steps. However, the creation of the ghost stabilization facets consumes comparatively negligible amounts of memory, on the order of $1e0 \! \sim 1e1$ \! megabytes for problems consuming $1e1 \! \sim \! 1e2$ gigabytes, which cannot be measured accurately. Given their insignificance, these measurements are omitted from the results.
Running the large-scale problems presented with a memory profiler to obtain precise memory consumption would be cost-prohibitive.
For runtime measurements, the maximum core time across all processes is used.

All measurements are averaged over three independent runs of each problem setup. The baseline values for timing and memory consumption, $t_0$ and $M_0$, for the results presented in the following, are summarized in \Cref{tab:timing_and_memory_baseline}.


\paragraph*{Time \& Memory Complexities}

\begin{figure}[btp]
    \centering
    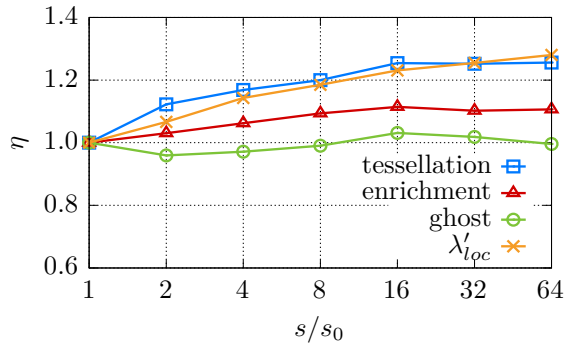
    \caption{Runtime efficiency $\eta$ when scaling up the problem size from one to 64 unit blocks using four processors.}
    \label{fig:size_scaling_eta}
\end{figure}

\begin{figure}[h]
    \centering
    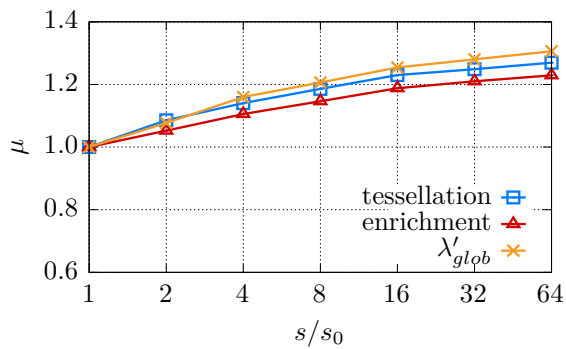
    \caption{Memory efficiency $\mu$ when scaling up the problem size from one to 64 unit blocks using four processors.}
    \label{fig:size_scaling_nu}
\end{figure}

The time and memory complexities of the framework are analyzed by scaling the problem size from a single unit block to a grid of $8\times8=64$ unit blocks. The problems were run in parallel on four processors with the domain split in half in $x$- and $y$-directions in each case.
The runtime and memory efficiencies $\eta$ and $\mu$ are shown in \Cref{fig:size_scaling_eta} and \Cref{fig:size_scaling_nu}, respectively. The relative mesh overhead decreases as the problem size increases, given that the number of processors is constant. Hence, efficiency measures $\eta,\mu > 1$ are observed.

The runtime for the decomposition closely follows $\lambda'_{loc}$ which indicates that the algorithms employed have an approximately linear time complexity $t \approx \mathcal{O}(s)$. For the ghost facet generation, $\eta \approx 1$ throughout. Given that the individual processors do not generate ghost facets in the aura, these operations are not expected to be affected by the decreasing relative mesh overhead. For the enrichment, only part of the basis functions supported in the aura is processed. As the runtime for the enrichment scales between the mesh overhead and $1.0$, the operations have an approximately linear time complexity $t \approx \mathcal{O}(s)$.

The memory consumption of the decomposition step closely follows the expected memory consumption due to the mesh overhead. The enrichment is less affected by the mesh overhead, and therefore its memory efficiency remains slightly below $\lambda'_{glob}$. These results suggest that the memory consumption of the algorithms scales linearly with the problem size $M = \mathcal{O}(s)$.


\begin{figure}[btp]
    \centering
    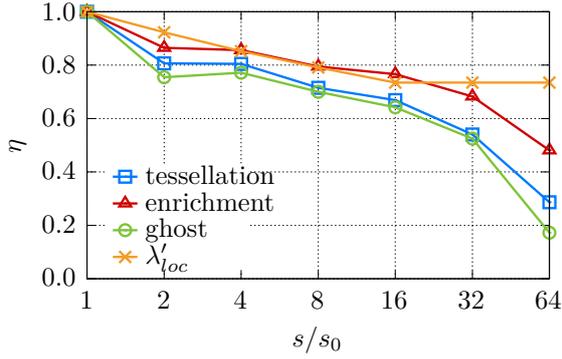
    \caption{Weak scaling efficiency $\eta$ when increasing the number of both the processors and unit blocks in the problem from one to 64.}
    \label{fig:weak_scaling_eta}
\end{figure}

\begin{figure}[btp]
    \centering
    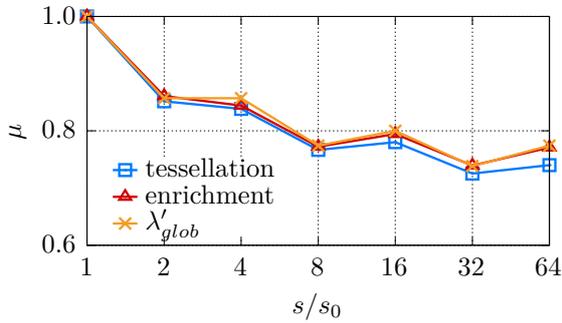
    \caption{Memory efficiency $\mu$ when increasing the number of both the processors and unit blocks in the problem from one to 64.}
    \label{fig:weak_scaling_nu}
\end{figure}

\paragraph*{Weak Parallel Scaling}

The weak parallel scalability is investigated by increasing both the number of parallel processors and the number of unit blocks in the problem from one to 64, yielding a domain size of one unit block per processor.

The weak parallel efficiency is shown in \Cref{fig:weak_scaling_eta}. For the decomposition, $\eta$ mostly follows the relative mesh overhead with a few percentage points offset, which can be explained by the approximately constant time it requires to perform the communication procedure in \Cref{alg:comm}. The mesh overhead increases as new processor boundaries and, therefore, aura elements are created when adding unit blocks to the initial 1-processor case with no aura. 
Notably, though, the runtime efficiency decreases significantly when scaling to above 16 processors. This is caused by a saturation of the available memory bandwidth. This assumption was confirmed by repeating the study on a server node with half the number of memory channels, and therefore memory bandwidth, but otherwise identical hardware. Its parallel scaling behavior showed the same significant performance breakdown scaling above 8 processors. 
The enrichment is not as strictly affected by the mesh overhead but suffers from the same problem. 
It is a known problem that the improvements to raw compute performance have outpaced improvements to memory performance for multiple decades of high-performance hardware development \cite{mccalpin1995memory,shalf2020future}. The result is that for modern hardware, such as the one used here, the memory performance is generally the bottleneck. Additionally, the presented algorithms heavily rely on storing and accessing data at different points in the workflow, which results in cache misses and frequent random memory access operations. Hence, the algorithms presented end up being memory-bound.

Memory consumption can, again, accurately be predicted by the mesh overhead as shown in \Cref{fig:weak_scaling_nu}.

\begin{figure}[b]
    \centering
    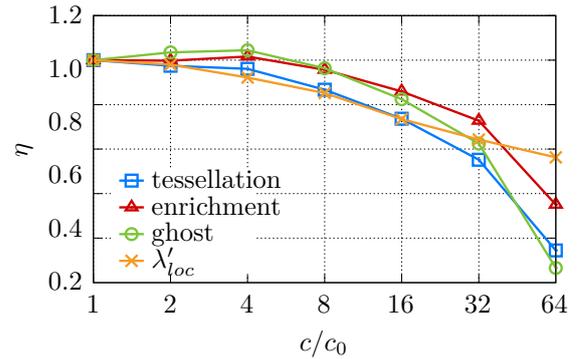
    \caption{Strong scaling efficiency $\eta$ when increasing the number of processors from one to 64 on a \textbf{large} problem consisting of 64 unit blocks.}
    \label{fig:strong_scaling_eta}
\end{figure}

\begin{figure}[btp]
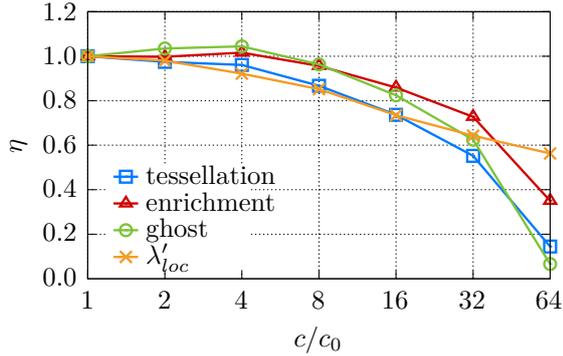

    \centering
    
    \caption{Strong scaling efficiency $\eta$ when increasing the number of processors from one to 64 on a \textbf{small} problem consisting of $4 \times 4 = 16$ unit blocks.}
    \label{fig:strong_scaling_eta_small}
\end{figure}

\begin{table*}[!b]
    \caption{Runtime $t_0$ and memory consumption $M_0$ for the baseline cases in Figures \ref{fig:size_scaling_nu} - \ref{fig:strong_scaling_eta_small}.}
    \label{tab:timing_and_memory_baseline}
    \centering
    \begin{tabular}{c|c|c|c|c}
    scaling type           & problem size   & weak      &  strong (large) &  strong (small) \\ \hline \hline
    tessellation $t_0$     &  $1.740s$      & $4.801s$  &  $349.5s$       &  $84.05s$       \\ \hline
    enrichment $t_0$       &  $1.895s$      & $6.097s$  &  $439.8s$       &  $110.6s$       \\ \hline
    ghost facet gen. $t_0$ &  $0.084s$      & $0.265s$  &  $24.71s$       &  $5.41s$        \\ \hline
    $\lambda_{loc,0}$      &  $1.361$       & $1.0$     &  $1.0$          &  $1.0$          \\ \hline \hline
    tessellation $M_0$     &  $1250$ MB     & $919$ MB  &  $60775$ MB     &  $13311$ MB     \\ \hline
    enrichment $M_0$       &  $1621$ MB     & $1239$ MB &  $81519$ MB     &  $29395$ MB     \\ \hline
    $\lambda_{glob,0}$     &  $1.333$       & $1.0$     &  $1.0$          &  $1.0$          \\ 
    \end{tabular}
\end{table*}

\paragraph*{Strong Parallel Scaling}

We investigate the strong parallel scalability by increasing the number of processors from one to 64 on two problems, a larger one with 64 unit blocks and a smaller one with 16 unit blocks. Timing results for these two cases are shown in \Cref{fig:strong_scaling_eta} and \Cref{fig:strong_scaling_eta_small}, respectively.
Scaling to more than 16 processors, again, results in a significant performance breakdown due to memory bandwidth saturation, regardless of the problem size. Ignoring this, the runtime efficiency for the decomposition follows the mesh overhead, while the enrichment is not as affected by the mesh overhead. Compared to the weak scaling scenario, the overhead created by communication is negligible. Comparing the two size cases, the decreasing performance ceiling due to the mesh aura for smaller problem sizes presents itself. E.g., for the smaller problem, the overall parallel efficiency for running the smaller scale problem on 16 processors is limited to only $\eta \approx 80\%$, as opposed to $\eta \approx 91\%$ for the larger problem. 

\paragraph*{Summary}
Both the compute time and memory consumption for the presented preprocessing framework scale approximately linear with the problem size 
if the geometric complexity is kept constant. In weak and strong scaling scenarios, the compute time and memory consumption scale close to optimally, once the mesh overhead is accounted for, with only a couple of percentage points decrease due to communication in the weak scaling scenarios.

The memory-bound design of the algorithms and the mesh overhead created by the aura present arguably the biggest shortfalls of the preprocessing framework. 
Re-designing the algorithms to more smartly (re-)compute data in place of access would mitigate the former bottleneck. However, it remains an open question to what extent this is possible.

The mesh overhead can constrain the performance severely when moving to very small background meshes on each processor's domain. As an example, consider a problem using quadratic $C^1$-continuous B-splines and running on a large number of processors such that the processors' domains are surrounded by those of adjacent processors to all sides. Reducing the processor-local mesh size to, e.g., $10 \times 10 \times 10$ elements results in a mesh overhead limiting the performance to $\eta \approx \lambda^{\prime}_{loc} = 2.6\%$, whereas a local mesh size of $25 \times 25 \times 25$ elements leads to a more modest loss $\eta \approx \lambda^{\prime}_{loc} = 64\%$. Given the overall speed of the procedures presented, it is likely unnecessary to move to such small processor-local domains during preprocessing. Re-distributing the generated clusters among a larger number of processors after preprocessing would be a solution if small processor-local domains are necessary to rapidly perform element assembly.

\clearpage
\section{Conclusion}
\label{sec:conclusion}

In this article, we discussed the implementation of an efficient, parallel, and robust immersed finite element preprocessing framework, including the underlying algorithms and data structures. 

A tessellation procedure creates a grid of cells fitted to the interfaces and generates quadrature rules on the different material regions and their interfaces. Material connectivity data is built on this grid and provides the information necessary to subsequently perform basis function enrichment and apply a stabilization scheme.
This approach of building information about the material connectivity from a geometric decomposition of the background mesh forms the backbone of the pre-processing framework. It allows some of the challenges associated with immersed finite element methods -- the generation of custom quadrature rules, enrichment if needed, and stabilization -- to be addressed cohesively.

The tessellation procedure applies subdivision templates repeatedly for an arbitrary number of geometries, enabling multi-material domains to be defined. Further, the presented algorithms for building topological information are designed without floating-point arithmetic, making them highly robust. 
The preprocessor's output is organized into elements paired with custom quadrature rules. Any information about the enriched finite element basis and stabilization methods is implicitly contained in them and not further exposed to the outside as such. Hence, assembly may be performed with existing, standard finite element routines without or with only minimal modifications.

The presented examples demonstrate the robustness and performance characteristics of the preprocessor. 
The generalized enrichment strategy prevents spurious coupling between close but disconnected features and thereby eliminates the case where the lack of mesh resolution leads not only to a large error but also to the simulation of a geometrically different, de-featured domain. Robustness is demonstrated in multi-material cases, both in terms of handling geometric edge cases and in terms of retaining good conditioning.
The implementation's time complexity and memory consumption are found to scale approximately linear with the problem size $\eta, \mu \approx \mathcal{O}(s)$, and the overhead created by parallel communication is shown to be minimal.

However, as summarized at the end of \Cref{sec:examples-nTop_sandwich}, the outlined parallel strategy creates a performance overhead for very small processor domains. Further, the scaling study highlighted that the implementation is limited by memory performance, rather than compute performance due to the algorithms' reliance on frequently accessing data stored in memory. 
Lastly, future work substituting the tessellation procedure and quadrature method, which only yield a linear geometry approximation, with ones that capture the geometry with higher-order accuracy, would likely improve the framework's numerical efficiency. This especially applies to problems involving non-linear interface phenomena, such as those encountered in contact mechanics, where highly refined meshes are currently needed.

\section*{Statements and Declarations}
\label{sec:declarations} 

\paragraph*{Acknowledgements}
N. Wunsch, J. A. Evans, and K. Maute were supported by the National Science Foundation (NSF), United States award OAC-2104106. N. Wunsch and K. Maute were additionally supported by the Air Force Office of Scientific Research (AFOSR), United States grant FA9550-20-1-0306. K. Doble was supported by Sandia National Laboratories through the contract PO2120843. M. R. Schmidt acknowledges the partial auspice of the U.S. Department of Energy by Lawrence Livermore National Laboratory under Contract DE-AC52-07NA27344 (LLNL-JRNL-87153).
L. No\"el, M. Schmidt, J.A. Evans, and K. Maute received support from the Defense Advanced Research Projects Agency (DARPA) under the TRADES program (agreement HR0011-17-2-0022).
The geometry data for the example presented in \Cref{sec:examples-nTop_sandwich} was provided by nTopology, Inc. free of charge.
The opinions and conclusions presented in this paper are those of the authors and do not necessarily reflect the views of the sponsoring organization.

\paragraph*{Conflicts of Interest}
J. A. Evans is an editor of Engineering with Computers.
The other authors have no competing interests to declare that are relevant to the content of this article.

\paragraph*{Replication of Results}
The numerical studies presented in this paper used the open-source software \texttt{Moris}, available at \href{https://github.com/kkmaute/moris/}{www.github.com/kkmaute/moris}. Further, the input decks for the presented example problems are archived in the same repository under 
\href{https://github.com/kkmaute/moris/tree/main/share/papers/XTK\_paper}{/share/papers/XTK\_paper}.

\paragraph*{Author Contributions}
The conceptual design and software development of the presented preprocessor were led by K. Doble and N. Wunsch. M. Schmidt and L. No\"el have made significant contributions to both the design and development of the preprocessor. J. A. Evans and K. Maute supervised the presented work and contributed to the conceptual design.
The numerical studies presented were performed by N. Wunsch. The first draft of the manuscript was written by N. Wunsch. All authors revised and commented on previous versions of the manuscript.
All authors read and approved the final manuscript.

\clearpage
\begin{appendices}

\section{Nomenclature}
\label{sec:nomenclature}

\subsection{Numerical Symbols}
\label{sec:nomenclature-numerical}

\renewcommand{\arraystretch}{1.4}

\paragraph*{Operators \& Conventions}
\begin{tabular}{ll}
    $(\cdot)^h$ & discrete quantity \\
    $\widetilde{(\cdot)}$ & enriched function (space) \\
    $\tilde{(\cdot)}$ & polynomial extension \\
    $\overline{(\cdot)}$ & closure of a domain \\
    $\partial(\cdot)$ & boundary of a domain \\
    $\WeightedJump{\cdot}$ & weighted interface operator, see \eqref{eqn:weighted_jump} \\
    $\jump{\cdot}$ & jump operator, see \eqref{eqn:jump_operator} \\
    $\norm{\cdot}{L2}$ & L2 norm \\
    $\abs{\cdot}_{H1}$ & H1 semi-norm \\
    $\NormalDeriv{(i)}$ & $i$-th derivative in normal direction $\normal$
\end{tabular}

\renewcommand{\arraystretch}{1.15}
\paragraph*{Numerical Quantities}
\begin{tabular}{ll}
    
    $\mathcal{B},\mathcal{B}_E$ & finite element basis, \\ 
                                & local to the bg. element $E$ \\
    $\dof,\vec{\dof}$ & degree of freedom, solution vector \\
    $E$ & Young's modulus \\
    $h,h_0$ & (initial) background element size \\
    $\normal,\normal_m$ & normal, pointing out from $\MaterialDomain$ \\
    $\BF$ & basis function \\
    $\PolyOrder$ & polynomial degree \\
    $\Residual_{(\cdot)}$ & residual for (B-bulk, D-Dirichlet, \\
                          & N-Neumann, I-interface, G-ghost) \\
    $\traction$ & traction vector \\
    $\displ$ & state variable/displacement field \\
    $\displ_0$ & prescribed displacement at boundary \\
    $\TestDispl$ & test function/virtual displacement \\
    $\mathcal{V}$ & space of test/trial functions \\
    $w$ & interface weight, see \eqref{eqn:weighted_jump} \\
    $\pos$ & physical coordinates \\

    $\gamma_{(\cdot)}$ & penalty parameters \\
                       & (D-Dirichlet, I-interface, G-ghost) \\
    $\strain$ & strain vector/tensor \\
    $\kappa$ & condition number \\
    $\stress$ & stress tensor \\
    $\nu$ & Poisson's ratio \\
    $\psi$ & enrichment function, see \eqref{eqn:enrichment_function} \\
    $\vec{\xi}$ & parametric coordinates \\
                          
\end{tabular}

\paragraph*{Domains}
\begin{tabular}{ll}
    $\AmbientDomain$ & ambient domain \\
                     & (domain of the background mesh) \\
    $\mathcal{F}_G$ & skeleton domain for face-oriented \\ 
                    & ghost stabilization, see \eqref{eqn:set_of_ghost_facets} \\
    $R,R_B^{\varepsilon}$ & connected domain within B.F. support, \\
                          & support of the enriched B.F. $\EnrichedBasisFunction_{\BfIndex}^{\EnrLvl}$, \\
                          & see \Cref{fig:enrichment_strategy} \\
    $S$ & connected domain within bg. element \\
    $\DirichletBoundary,\NeumannBoundary$ & Dirichlet \& Neumann boundaries \\
    $\MaterialInterface{m}{l}$ & material interface between \\
                               & materials with indices $m$ and $l$  \\
    $\InteriorInterfaces$ & union of material interfaces \\ 
                          & in the interior of $\AmbientDomain$, see \eqref{eqn:material_interfaces} \\
    $\MaterialDomain$ & material domain \\
    $\BgElemDom$ & domain of a background element \\
    $\InteriorDomain$ & union of material domains, see \eqref{eqn:interior_domain}
\end{tabular}

\subsection{Implementation}
\label{sec:nomenclature-implementation}

\paragraph*{Indices \& Entities/Objects}
\begin{tabular}{ll}
    $a$ & ancestor entity (index) \\
    $b,B$ & local \& global basis function indices \\
    $\mathrm{c}$ & cell (index) \\
    $\mathrm{CM}$ & child mesh \\
    $d$ & dimensionality \\
    $D$ & cluster \\
    $e$ & entity (index) \\
    $\varepsilon$ & enrichment level \\
    $E$ & background element (index) \\
    $F$ & facet (index) \\
    $G$ & geometry \\
    $g$ & geometry index \\
    $\mathcal{H}$ & background mesh \\
    $i,j,k,l$ & general indices \\
    $I$ & parallel ID \\
    $\EnrBfIndex$ & enriched B.F. multi-index $(\BfIndex,\EnrLvl)$ \\
    $m$ & material (index) \\
    $n_{(\cdot)}$ & number of $(\cdot)$ \\
    $o$ & ordinal \\
    $P$ & proximity \\
    $p$ & (current) processor ID \\
    $q$ & quadrature point index \\
    $\mathrm{Q_v},\mathrm{Q_c}$ & vertex \& cell queues
\end{tabular}

\paragraph*{Indices \& Entities/Objects (continued)}
\begin{tabular}{ll}
    $r$ & (entity) rank \\
    $S$ & subphase \\
    $\mathcal{T}$ & foreground mesh \\
    $u,v$ & unzipping index \\
    $\mathrm{v}$ & vertex (index)
\end{tabular}

\paragraph*{Sets \& Graphs}
\begin{tabular}{ll}
    $\mathcal{C}$ & set of all fg. cells \\
    $\mathcal{D}_{(\cdot)}$ & set of all clusters of type \\
                            & (B-bulk, S-side, I-interface, G-ghost) \\
    $\mathcal{E}$ & set of all bg. elements \\
    $\mathcal{F}$ & facet connectivity (graph) \\
    $\mathcal{G}_S,\mathcal{G}_I$ & graph of connected, \\
                                  & disconnected subphases \\
    $\mathcal{M}$ & set of (non-void) material indices \\
    $\mathcal{S}$ & set of all subphases \\
    $\mathcal{U}_{u,E}$ & set of subphases connected to $S_E^u$, \\
                        & see \eqref{eqn:ghost_connected_subdomains} \\
    $\mathrm{CtE}$ & cell-to-entity connectivity, \\
                   & see \Cref{fig:entity_connectivity_data_structure} \\
    $\mathrm{DF}$ & sets of descendant facets, \\
                  & see \Cref{alg:generate_bg_facet_descendants} \\
    $\mathrm{EtC}$ & entity-to-cell connectivity, \\
                   & see \Cref{fig:entity_connectivity_data_structure} \\
    $\mathrm{NVI}$ & set of new vertex indices \\
    $\mathrm{P}_{glob}$ & global communication table \\ 
                   & (set of all processors) \\
    $\mathrm{P}_{loc}$ & local communication table \\ 
                   & (set of neighboring processors) \\
    
    $\Xi$ & set of unzipped bg. elements    
\end{tabular}

\subsection{Compute Performance}
\label{sec:nomenclature-compute}

\begin{tabular}{ll}
    $M$ & memory consumption \\
    $n_p$ & number of processors \\
    $s$ & problem size \\
    $t$ & compute time \\
    $\eta$ & compute efficiency, see \eqref{eqn:efficiency_measure} \\
    $\mu$ & memory efficiency, see \eqref{eqn:memory_efficiency} \\ 
    $\lambda_{loc}$ & local mesh overhead, see \eqref{eqn:local_mesh_overhead} \\
    $\lambda_{glob}$ & global mesh overhead, see \eqref{eqn:global_mesh_overhead} \\
    $\mathcal{O}(\cdot)$ & order of limiting behavior \\
    $(\cdot)_0$ & value in reference run \\
    $(\cdot)^{\prime}$ & relative quantity $=(\cdot)/(\cdot)_0$
\end{tabular}

\subsection{Dictionary}
\label{sec:nomenclature-dictionary}

\ \ \textbf{Adjacency list} $\bullet$ Data structure to represent the graphs $\mathcal{G}_{(\cdot)}$ as an array of arrays. In our case, the outer array entries correspond to each subphase and the inner (unordered) array contains the adjacent subphase neighbors, in the sense of the particular graph.
\vspace{0.15cm}

\textbf{(Background) ancestry/ancestor} $\bullet$ Symbol: $\Ancestor$ $\bullet$ The background ancestor of an entity of the foreground mesh is the background entity of the lowest rank, i.e., dimensionality, which fully contains the foreground entity. See: \Cref{sec:fg_mesh_generation}, par. \hyperref[pg:background_ancestry]{Background Ancestry}.
\vspace{0.15cm}

\textbf{Background mesh} $\bullet$ Symbol: $\BgMesh$ $\bullet$ Input mesh constructed on a geometrically simple domain into which the geometry to be analyzed is immersed. It defines the (non-enriched) finite element basis. The associated data structure is shown in \Cref{fig:IO_data_structures}.
\vspace{0.15cm}

\textbf{Cell} $\bullet$ Symbol: $\Cell$ $\bullet$ Purely geometric, $d$-dimensional entity on a $d$-dimensional mesh. They are, e.g., the triangles and quadrilateral in 2D, or tetrahedrons and hexahedrons in 3D comprising the foreground and background meshes. Note, that unlike elements, cells do not have a basis associated with them. See: \Cref{sec:preliminaries}, par. \hyperref[pg:mesh_entities]{Mesh Entities and Connectivity}.
\vspace{0.15cm}

\textbf{Child mesh}: $\bullet$ Symbol: $\ChildMesh$ $\bullet$ Group of foreground cells generated from a single background element during the tessellation procedure outlined in \Cref{sec:fg_mesh_generation}. The associated data structure is shown in \Cref{fig:fg_mesh_data_structure}.
\vspace{0.15cm}

\textbf{Cluster} $\bullet$ Symbol: $\Cluster$ $\bullet$ Data structure that combines a (background) element with a (custom) quadrature rule defined on its domain, or a subdomain thereof. The output of the preprocessors are groups of clusters, intended for standard finite element assembly routines. The associated data structure is shown in \Cref{fig:IO_data_structures}. The concept of a cluster is discussed from a theoretical standpoint in \Cref{sec:framework-ig}.
\vspace{0.15cm}

\textbf{(Parallel, mesh) decomposition} $\bullet$ Process of dividing the background mesh into subdomains to be handled by processors independently. See: \Cref{sec:parallel}.
\vspace{0.15cm}

\textbf{(Background, finite) element} $\bullet$ Symbol: $\BgElemInd$ $\bullet$ A finite element consists of a domain $\BgElemDom$ and a set of basis functions $\mathcal{B}_E$ with their associated degrees of freedom $\IEN$. In the article, the only finite elements used are those comprising the background mesh. Hence, they are referred to as background elements. The associated data structure is shown in \Cref{fig:IO_data_structures}.
\vspace{0.15cm}

\textbf{Enrichment function} $\bullet$ Symbol: $\IndicatorFunction$ $\bullet$ The indicator function associated with an enriched basis function with multi-index $\EnrBfIndex=(\BfIndex,\EnrLvl)$. For a definition, see \eqref{eqn:enrichment_function}.
\vspace{0.15cm}

\textbf{(Mesh) entity} $\bullet$ Symbol: $\Entity$ $\bullet$ Geometric objects comprising a mesh or grid, namely: vertices, edges, faces, 2D/3D cells. See: \Cref{sec:preliminaries}, par. \hyperref[pg:mesh_entities]{Mesh Entities and Connectivity}.
\vspace{0.15cm}

\textbf{Facet} $\bullet$ Symbol: $\Facet$ $\bullet$ The $(d-1)$-dimensional entity of a $d$-dimensional mesh, e.g.,  edges in 2D and faces in 3D. See: \Cref{sec:preliminaries}, par. \hyperref[pg:mesh_entities]{Mesh Entities and Connectivity}.
\vspace{0.15cm}

\textbf{Foreground Mesh} $\bullet$ Symbol: $\FgMesh$ $\bullet$ Body-fitted mesh 
resulting from the subdivision process outlined in \Cref{sec:fg_mesh_generation}. Additionally, the associated data structure, shown in \Cref{fig:fg_mesh_data_structure}, carries information about the material topology constructed in \Cref{sec:subphase_generation}.
\vspace{0.15cm}

\textbf{Geometry} $\bullet$ Symbol: $\Geometry$ $\bullet$ Object defining a material interface or volumetric domain. The associated data structure is shown in \Cref{fig:IO_data_structures}.
\vspace{0.15cm}

\textbf{(Parallel) ID} $\bullet$ Symbol: $\ID$ $\bullet$ An entity ID is a parallel-globally unique identifier for a mesh entity that ensures consistency and facilitates communication across processors in parallel computations. See: \Cref{sec:parallel}, par. \hyperref[pg:parallel_IDs]{Concept of Parallel IDs}.
\vspace{0.15cm}

\textbf{Index} $\bullet$ A (processor-local) identifier for entities or other programming objects of a specific type, ranging from 1 to the total number of such entities or objects on a processor. Indices are used to access entities or objects that are stored in arrays. See: \Cref{sec:preliminaries}, par. \hyperref[pg:mesh_entities]{Mesh Entities and Connectivity}.
\vspace{0.15cm}

\textbf{Initializer list} $\bullet$ A list of values, entities, and/or objects that are used to create/initialize another entity or object.
\vspace{0.15cm}

\textbf{Ordinal} $\bullet$ Symbol: $\Ordinal$ $\bullet$ The local index of an entity with respect to another entity, as shown in \Cref{fig:entity_connectivity} (A). See: \Cref{sec:preliminaries}, par. \hyperref[pg:mesh_entities]{Mesh Entities and Connectivity}.
\vspace{0.15cm}

\textbf{Proximity} $\bullet$ Symbol: $\Proximity$ $\bullet$ Whether an entity lies outside ($+$), inside ($-$), or on the boundary ($0$) of the geometry within some tolerance. The concept is first introduced in \Cref{sec:preliminaries} and used extensively in the step of \hyperref[pg:material_assignment]{material assignment}.
\vspace{0.15cm}

\textbf{Queue} $\bullet$ Symbol: $\Queue$ $\bullet$ Data structure collecting new entities/objects to be created, also ensuring their uniqueness. The concept is discussed in \Cref{sec:fg_mesh_generation-reg_sub}. The associated data structure is shown in \Cref{fig:fg_mesh_full_data_structure}.
\vspace{0.15cm}

\textbf{Rank} $\bullet$ Symbol: $\Rank$ $\bullet$ Dimensionality of an entity in the foreground mesh. See: \Cref{sec:preliminaries}, par. \hyperref[pg:mesh_entities]{Mesh Entities and Connectivity}.
\vspace{0.15cm}

\textbf{Subphase} $\bullet$ Symbol: $\Subphase$ $\bullet$ Topologically connected material region within a single background element's domain $\BgElemDom$. Their generation is discussed in \Cref{sec:subphase_generation}, par. \hyperref[pg:subphase_generation]{Identification and Generation of Subphases}. data structure: \Cref{fig:subphase_data_structure}.
\vspace{0.15cm}

\textbf{Tessellation} $\bullet$ Process of dividing a domain into geometric primitives. In the context of this article, it is used to refer to the subdivision process generating the foreground mesh.
\vspace{0.15cm}

\textbf{Unzipping} $\bullet$ Index: $\Unzipping$ $\bullet$ Process of duplicating background elements intersected by material interfaces and re-indexing their basis functions to represent the enriched basis. The process is described in \Cref{sec:enrichment}, par. \hyperref[pg:unzipping]{Unzipping} and illustrated in \Cref{fig:unzipping}.
\vspace{0.15cm}

\textbf{Vertex} $\bullet$ Symbol: $\Vertex$  $\bullet$ Purely geometric $0$-dimensional entity on a mesh. They are the corners (or endpoints) of cells (or edges). Note, that unlike nodes, vertices do not have the notion of a basis function associated with them. See: \Cref{sec:preliminaries}, par. \hyperref[pg:mesh_entities]{Mesh Entities and Connectivity}.
\vspace{0.15cm}

\onecolumn
\newpage
\section{Additional Algorithms and Data Structures}
\label{sec:apx-algs_and_dss}

\subsection{Determining Ancestry of New Vertices}
\label{sec:apx-ancestry_reg_sub_template}

Using the particular regular subdivision templates shown in \Cref{fig:reg_sub_template}, the background ancestry of a new vertex added at any stage during the hierarchical templated subdivision process outlined in \Cref{sec:fg_mesh_generation-templ_sub} can be determined by the ancestry of the two vertices attached to that edge. One only needs to check which of the following conditions is true:
\begin{enumerate}
    \item One of the vertices is inside the bg. element: the new vertex is inside the element.
    \item The edge vertices descend from the same edge or face: the new vertex is a descendant of that edge or face.
    \item Both edge vertices each descend from a background vertex: these bg. vertices are necessarily part of the same background edge, and the new vertex is a descendant of that edge.
    \item The edge vertices descend from different background faces: the new vertex is inside the background element.
    \item One edge vertex descends from a vertex and the other from a face or edge: the new vertex is a descendant of that face or edge.
\end{enumerate}

\subsection{Entity Connectivity}
\label{sec:apx-entity_connectivity}

The algorithms presented require information about which lower-order entities, vertices, edges, and faces, are connected to cells and vice versa. This information is stored in the data structure shown in \Cref{fig:entity_connectivity}, and can be generated on any conformal mesh using \Cref{alg:compute_entity_connectivity}.

\begin{algorithm}[h!]
    \caption{
        Generate entity connectivity.\\
        Input: mesh $\ForegroundMesh$ containing a cell-vertex connectivity $\mathcal{C}$, rank of entities $r$ \\
        Output: entity connectivity $\mathcal{F}$ }
    \label{alg:compute_entity_connectivity}
    \begin{algorithmic}[1]
        \State initialize empty entity connectivity $\mathcal{F}$ and hash map $\mathrm{VtE}:(\Vertex_1,\Vertex_2,\cdots) \mapsto \Entity$, and running index $i \leftarrow 0$
        \State \textit{\# the number of cells is known and the cell-to-entity array in $\mathcal{F}$ can be initialized accordingly}
        \For{each cell $\Cell \in \ForegroundMesh.\mathcal{C}$}
            \For{each entity ordinal on cell $o \in \SetDef{1,\cdots,\Cell.\mathtt{num\_entities}(r)}$}
                \State get the entities expressed as a list of vertices $(\Vertex_1,\Vertex_2,\cdots) \leftarrow \mathtt{sort}(\Cell.\mathtt{get\_entity}(r,o))$
                \State find entity in the map $\Entity \leftarrow \mathrm{VtE}.\mathtt{find}((\Vertex_1,\Vertex_2,\cdots))$
                \If{the entity does not exist in the map already $!j$}
                    \State store the entity in the map $\mathrm{VtE}.\mathtt{insert}((\Vertex_1,\Vertex_2,\cdots),++i)$
                    \State add information $\mathcal{F}.\mathrm{EtC}.\mathtt{append}(\SetDef{\Cell})$, $\mathcal{F}.\mathrm{CtE}[\Cell].\mathtt{append}(i)$ 
                \Else{} 
                    \State add information $\mathcal{F}.\mathrm{EtC}[\Entity].\mathtt{append}(\Cell)$, $\mathcal{F}.\mathrm{CtE}[\Cell].\mathtt{append}(\Entity)$
                \EndIf 
            \EndFor
        \EndFor 
    \end{algorithmic}
\end{algorithm} 

\subsection{Foreground Mesh -- Complete Data Structure}
\label{sec:apx-fg_mesh_data_structure}

\begin{figure}[h!]
    \centering 
    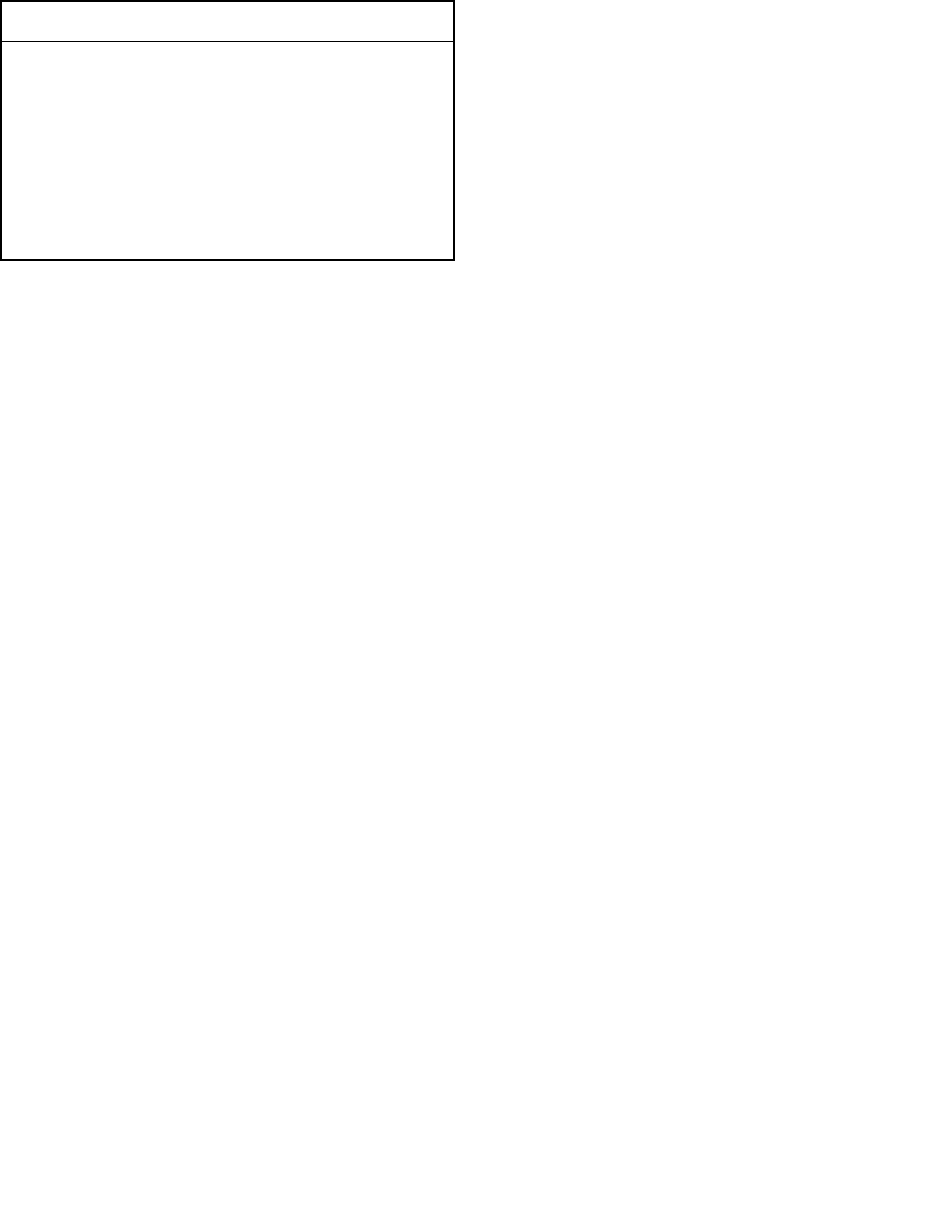
    \vspace*{0.2cm}
    \caption{UML chart of the data structures used in the presented preprocessor.}
    \label{fig:fg_mesh_full_data_structure}
\end{figure}

The presented algorithms make use of data structures introduced in Figures \ref{fig:IO_data_structures},  \ref{fig:entity_connectivity_data_structure}, \ref{fig:fg_mesh_data_structure}, and \ref{fig:subphase_data_structure} throughout \Cref{sec:main}. The unified modeling language (UML) chart in \Cref{fig:fg_mesh_full_data_structure} illustrates their interaction with each other.

\twocolumn

\end{appendices}

\bibliographystyle{sn-basic}
\bibliography{consolidated} 


\end{document}

%% file: figures/cross_section_with_details.pdf_tex
\begingroup%
  \makeatletter%
  \providecommand\color[2][]{%
    \errmessage{(Inkscape) Color is used for the text in Inkscape, but the package 'color.sty' is not loaded}%
    \renewcommand\color[2][]{}%
  }%
  \providecommand\transparent[1]{%
    \errmessage{(Inkscape) Transparency is used (non-zero) for the text in Inkscape, but the package 'transparent.sty' is not loaded}%
    \renewcommand\transparent[1]{}%
  }%
  \providecommand\rotatebox[2]{#2}%
  \newcommand*\fsize{\dimexpr\f@size pt\relax}%
  \newcommand*\lineheight[1]{\fontsize{\fsize}{#1\fsize}\selectfont}%
  \ifx\svgwidth\undefined%
    \setlength{\unitlength}{501.53803307bp}%
    \ifx\svgscale\undefined%
      \relax%
    \else%
      \setlength{\unitlength}{\unitlength * \real{\svgscale}}%
    \fi%
  \else%
    \setlength{\unitlength}{\svgwidth}%
  \fi%
  \global\let\svgwidth\undefined%
  \global\let\svgscale\undefined%
  \makeatother%
  \begin{picture}(1,0.44103799)%
    \lineheight{1}%
    \setlength\tabcolsep{0pt}%
    \put(0,0){\includegraphics[width=\unitlength,page=1]{cross_section_with_details.pdf}}%
    \put(0.2908575,0.4243744){\color[rgb]{0.19607843,0.19607843,0.19607843}\makebox(0,0)[lt]{\lineheight{1.25}\smash{\begin{tabular}[t]{l}background mesh\end{tabular}}}}%
    \put(-0.0019384,0.12565772){\color[rgb]{0,0,0}\makebox(0,0)[lt]{\lineheight{1.25}\smash{\begin{tabular}[t]{l}(A)\end{tabular}}}}%
    \put(0.2385715,0.12453068){\color[rgb]{0,0,0}\makebox(0,0)[lt]{\lineheight{1.25}\smash{\begin{tabular}[t]{l}(B)\end{tabular}}}}%
    \put(0.49167117,0.12513894){\color[rgb]{0,0,0}\makebox(0,0)[lt]{\lineheight{1.25}\smash{\begin{tabular}[t]{l}(C)\end{tabular}}}}%
    \put(0,0){\includegraphics[width=\unitlength,page=2]{cross_section_with_details.pdf}}%
    \put(0.57306301,0.1550154){\color[rgb]{0,0.47843137,1}\makebox(0,0)[lt]{\lineheight{1.25}\smash{\begin{tabular}[t]{l}$\mathrm{supp}(N)$\end{tabular}}}}%
    \put(0.36523246,0.05948917){\color[rgb]{0,0.47843137,1}\makebox(0,0)[lt]{\lineheight{1.25}\smash{\begin{tabular}[t]{l}$\Gamma^I_{1,2}$\end{tabular}}}}%
    \put(0.30131847,0.03108122){\color[rgb]{0.19607843,0.19607843,0.19607843}\makebox(0,0)[lt]{\lineheight{1.25}\smash{\begin{tabular}[t]{l}$\Omega_{1}$\end{tabular}}}}%
    \put(0.34003924,0.107732){\color[rgb]{0.19607843,0.19607843,0.19607843}\makebox(0,0)[lt]{\lineheight{1.25}\smash{\begin{tabular}[t]{l}$\Omega_{2}$\end{tabular}}}}%
    \put(0.07766338,0.10251268){\color[rgb]{0.19607843,0.19607843,0.19607843}\makebox(0,0)[lt]{\lineheight{1.25}\smash{\begin{tabular}[t]{l}$\Omega_{1}$\end{tabular}}}}%
    \put(0.08357643,0.04430239){\color[rgb]{0.19607843,0.19607843,0.19607843}\makebox(0,0)[lt]{\lineheight{1.25}\smash{\begin{tabular}[t]{l}$\Omega_{2}$\end{tabular}}}}%
    \put(0.12744018,0.01485735){\color[rgb]{0.19607843,0.19607843,0.19607843}\makebox(0,0)[lt]{\lineheight{1.25}\smash{\begin{tabular}[t]{l}$\Omega_{3}$\end{tabular}}}}%
    \put(0.79984496,0.12508046){\color[rgb]{0,0,0}\makebox(0,0)[lt]{\lineheight{1.25}\smash{\begin{tabular}[t]{l}(D)\end{tabular}}}}%
    \put(0,0){\includegraphics[width=\unitlength,page=3]{cross_section_with_details.pdf}}%
    \put(0.87137689,0.01429657){\color[rgb]{0.19607843,0.19607843,0.19607843}\makebox(0,0)[lt]{\lineheight{1.25}\smash{\begin{tabular}[t]{l}$\Omega_{3} \cap \Omega^{E}$\end{tabular}}}}%
    \put(0,0){\includegraphics[width=\unitlength,page=4]{cross_section_with_details.pdf}}%
  \end{picture}%
\endgroup%

%% file: figures/enrichment_strategy.pdf_tex
\begingroup%
  \makeatletter%
  \providecommand\color[2][]{%
    \errmessage{(Inkscape) Color is used for the text in Inkscape, but the package 'color.sty' is not loaded}%
    \renewcommand\color[2][]{}%
  }%
  \providecommand\transparent[1]{%
    \errmessage{(Inkscape) Transparency is used (non-zero) for the text in Inkscape, but the package 'transparent.sty' is not loaded}%
    \renewcommand\transparent[1]{}%
  }%
  \providecommand\rotatebox[2]{#2}%
  \newcommand*\fsize{\dimexpr\f@size pt\relax}%
  \newcommand*\lineheight[1]{\fontsize{\fsize}{#1\fsize}\selectfont}%
  \ifx\svgwidth\undefined%
    \setlength{\unitlength}{171.09682459bp}%
    \ifx\svgscale\undefined%
      \relax%
    \else%
      \setlength{\unitlength}{\unitlength * \real{\svgscale}}%
    \fi%
  \else%
    \setlength{\unitlength}{\svgwidth}%
  \fi%
  \global\let\svgwidth\undefined%
  \global\let\svgscale\undefined%
  \makeatother%
  \begin{picture}(1,1.87877562)%
    \lineheight{1}%
    \setlength\tabcolsep{0pt}%
    \put(0,0){\includegraphics[width=\unitlength,page=1]{enrichment_strategy.pdf}}%
    \put(0.8028673,0.67890178){\color[rgb]{0.19607843,0.19607843,0.19607843}\makebox(0,0)[lt]{\lineheight{1.25}\smash{\begin{tabular}[t]{l}$\Omega_{2}$\end{tabular}}}}%
    \put(0.81580854,1.16044116){\color[rgb]{0.19607843,0.19607843,0.19607843}\makebox(0,0)[lt]{\lineheight{1.25}\smash{\begin{tabular}[t]{l}$\Omega_{3}$\end{tabular}}}}%
    \put(0.14363104,0.78653266){\color[rgb]{0.19607843,0.19607843,0.19607843}\makebox(0,0)[lt]{\lineheight{1.25}\smash{\begin{tabular}[t]{l}$\Omega_{3}$\end{tabular}}}}%
    \put(0.12804542,1.29182329){\color[rgb]{0,0.47843137,1}\makebox(0,0)[lt]{\lineheight{1.25}\smash{\begin{tabular}[t]{l}$\mathrm{supp}(N_{\BfIndex})$\end{tabular}}}}%
    \put(0,0){\includegraphics[width=\unitlength,page=2]{enrichment_strategy.pdf}}%
    \put(0.71031213,1.38773168){\color[rgb]{0.19607843,0.19607843,0.19607843}\makebox(0,0)[lt]{\lineheight{1.25}\smash{\begin{tabular}[t]{l}$\EnrLvl=1$\end{tabular}}}}%
    \put(0.1303071,1.38936524){\color[rgb]{0.19607843,0.19607843,0.19607843}\makebox(0,0)[lt]{\lineheight{1.25}\smash{\begin{tabular}[t]{l}$\EnrLvl=2$\end{tabular}}}}%
    \put(0.70985007,0.01412413){\color[rgb]{0.19607843,0.19607843,0.19607843}\makebox(0,0)[lt]{\lineheight{1.25}\smash{\begin{tabular}[t]{l}$\EnrLvl=4$\end{tabular}}}}%
    \put(0.13030613,0.01173144){\color[rgb]{0.19607843,0.19607843,0.19607843}\makebox(0,0)[lt]{\lineheight{1.25}\smash{\begin{tabular}[t]{l}$\EnrLvl=3$\end{tabular}}}}%
    \put(0.83130947,0.18777989){\color[rgb]{0,0,0}\makebox(0,0)[lt]{\lineheight{1.25}\smash{\begin{tabular}[t]{l}$\BfSubPhase{4}$\end{tabular}}}}%
    \put(0.82723661,1.73359681){\color[rgb]{0,0,0}\makebox(0,0)[lt]{\lineheight{1.25}\smash{\begin{tabular}[t]{l}$\BfSubPhase{1}$\end{tabular}}}}%
    \put(0,0){\includegraphics[width=\unitlength,page=3]{enrichment_strategy.pdf}}%
    \put(0.0832472,1.72043225){\color[rgb]{0,0,0}\makebox(0,0)[lt]{\lineheight{1.25}\smash{\begin{tabular}[t]{l}$\BfSubPhase{2}$\end{tabular}}}}%
    \put(0.06242457,0.20144726){\color[rgb]{0,0,0}\makebox(0,0)[lt]{\lineheight{1.25}\smash{\begin{tabular}[t]{l}$\BfSubPhase{3}$\end{tabular}}}}%
    \put(0,0){\includegraphics[width=\unitlength,page=4]{enrichment_strategy.pdf}}%
  \end{picture}%
\endgroup%

%% file: figures/clusters.pdf_tex
\begingroup%
  \makeatletter%
  \providecommand\color[2][]{%
    \errmessage{(Inkscape) Color is used for the text in Inkscape, but the package 'color.sty' is not loaded}%
    \renewcommand\color[2][]{}%
  }%
  \providecommand\transparent[1]{%
    \errmessage{(Inkscape) Transparency is used (non-zero) for the text in Inkscape, but the package 'transparent.sty' is not loaded}%
    \renewcommand\transparent[1]{}%
  }%
  \providecommand\rotatebox[2]{#2}%
  \newcommand*\fsize{\dimexpr\f@size pt\relax}%
  \newcommand*\lineheight[1]{\fontsize{\fsize}{#1\fsize}\selectfont}%
  \ifx\svgwidth\undefined%
    \setlength{\unitlength}{395.6311708bp}%
    \ifx\svgscale\undefined%
      \relax%
    \else%
      \setlength{\unitlength}{\unitlength * \real{\svgscale}}%
    \fi%
  \else%
    \setlength{\unitlength}{\svgwidth}%
  \fi%
  \global\let\svgwidth\undefined%
  \global\let\svgscale\undefined%
  \makeatother%
  \begin{picture}(1,0.29654451)%
    \lineheight{1}%
    \setlength\tabcolsep{0pt}%
    \put(0,0){\includegraphics[width=\unitlength,page=1]{clusters.pdf}}%
    \put(0.61938499,0.14715153){\color[rgb]{0.19607843,0.19607843,0.19607843}\makebox(0,0)[lt]{\lineheight{1.25}\smash{\begin{tabular}[t]{l}$\ElemSubPhase{2}$\end{tabular}}}}%
    \put(0.90082826,0.06283312){\color[rgb]{0.19607843,0.19607843,0.19607843}\makebox(0,0)[lt]{\lineheight{1.25}\smash{\begin{tabular}[t]{l}$\ElemSubPhase{2}$\end{tabular}}}}%
    \put(0.80254833,0.21851466){\color[rgb]{0.19607843,0.19607843,0.19607843}\makebox(0,0)[lt]{\lineheight{1.25}\smash{\begin{tabular}[t]{l}$\ElemSubPhase{1}$\end{tabular}}}}%
    \put(-0.00110171,0.25746314){\color[rgb]{0.19607843,0.19607843,0.19607843}\makebox(0,0)[lt]{\lineheight{1.25}\smash{\begin{tabular}[t]{l}$(A)$\end{tabular}}}}%
    \put(0.24966872,0.25780905){\color[rgb]{0.19607843,0.19607843,0.19607843}\makebox(0,0)[lt]{\lineheight{1.25}\smash{\begin{tabular}[t]{l}$(B)$\end{tabular}}}}%
    \put(0.50043915,0.25780905){\color[rgb]{0.19607843,0.19607843,0.19607843}\makebox(0,0)[lt]{\lineheight{1.25}\smash{\begin{tabular}[t]{l}$(C)$\end{tabular}}}}%
    \put(0.75120954,0.25780926){\color[rgb]{0.19607843,0.19607843,0.19607843}\makebox(0,0)[lt]{\lineheight{1.25}\smash{\begin{tabular}[t]{l}$(D)$\end{tabular}}}}%
    \put(0,0){\includegraphics[width=\unitlength,page=2]{clusters.pdf}}%
    \put(0.13295556,0.07742787){\color[rgb]{0.82352941,0,0.00392157}\makebox(0,0)[lt]{\lineheight{1.25}\smash{\begin{tabular}[t]{l}$\xi_q$\end{tabular}}}}%
    \put(0.28649078,0.1238762){\color[rgb]{0,0.47843137,1}\makebox(0,0)[lt]{\lineheight{1.25}\smash{\begin{tabular}[t]{l}$\ElemSubPhase{2}$\end{tabular}}}}%
    \put(0.57760165,0.25809068){\color[rgb]{0,0.47843137,1}\makebox(0,0)[lt]{\lineheight{1.25}\smash{\begin{tabular}[t]{l}$\Boundary{\ElemSubPhase{2}} \cap \BoundarySegment$\end{tabular}}}}%
    \put(0.04905237,0.17367956){\color[rgb]{0,0.47843137,1}\makebox(0,0)[lt]{\lineheight{1.25}\smash{\begin{tabular}[t]{l}$\ElemSubPhase{1}$\end{tabular}}}}%
    \put(0.70028582,0.1250495){\color[rgb]{0,0.47843137,1}\makebox(0,0)[lt]{\lineheight{1.25}\smash{\begin{tabular}[t]{l}$\Boundary{\ElemSubPhase{1}} \cap \Boundary{\ElemSubPhase{2}}$\end{tabular}}}}%
    \put(0,0){\includegraphics[width=\unitlength,page=3]{clusters.pdf}}%
  \end{picture}%
\endgroup%

%% file: figures/ghost_for_enriched.pdf_tex
\begingroup%
  \makeatletter%
  \providecommand\color[2][]{%
    \errmessage{(Inkscape) Color is used for the text in Inkscape, but the package 'color.sty' is not loaded}%
    \renewcommand\color[2][]{}%
  }%
  \providecommand\transparent[1]{%
    \errmessage{(Inkscape) Transparency is used (non-zero) for the text in Inkscape, but the package 'transparent.sty' is not loaded}%
    \renewcommand\transparent[1]{}%
  }%
  \providecommand\rotatebox[2]{#2}%
  \newcommand*\fsize{\dimexpr\f@size pt\relax}%
  \newcommand*\lineheight[1]{\fontsize{\fsize}{#1\fsize}\selectfont}%
  \ifx\svgwidth\undefined%
    \setlength{\unitlength}{354.91214061bp}%
    \ifx\svgscale\undefined%
      \relax%
    \else%
      \setlength{\unitlength}{\unitlength * \real{\svgscale}}%
    \fi%
  \else%
    \setlength{\unitlength}{\svgwidth}%
  \fi%
  \global\let\svgwidth\undefined%
  \global\let\svgscale\undefined%
  \makeatother%
  \begin{picture}(1,0.30049687)%
    \lineheight{1}%
    \setlength\tabcolsep{0pt}%
    \put(0,0){\includegraphics[width=\unitlength,page=1]{ghost_for_enriched.pdf}}%
    \put(0.6575297,0.25624055){\color[rgb]{0.19607843,0.19607843,0.19607843}\makebox(0,0)[lt]{\lineheight{1.25}\smash{\begin{tabular}[t]{l}$\Omega^{E^+}$\end{tabular}}}}%
    \put(0.81914561,0.21770174){\color[rgb]{0,0.47843137,1}\makebox(0,0)[lt]{\lineheight{1.25}\smash{\begin{tabular}[t]{l}$\normal$\end{tabular}}}}%
    \put(0.74117516,0.01830311){\color[rgb]{0,0.47843137,1}\makebox(0,0)[lt]{\lineheight{1.25}\smash{\begin{tabular}[t]{l}$F \in \SetOfGhostFacets$\end{tabular}}}}%
    \put(0.46470671,0.22872784){\color[rgb]{0,0.47843137,1}\makebox(0,0)[lt]{\lineheight{1.25}\smash{\begin{tabular}[t]{l}$\SetOfGhostFacets$\end{tabular}}}}%
    \put(0.72574248,0.22102076){\color[rgb]{0.19607843,0.19607843,0.19607843}\makebox(0,0)[lt]{\lineheight{1.25}\smash{\begin{tabular}[t]{l}$S_{E^+}^{1}$\end{tabular}}}}%
    \put(0.01910722,0.22456646){\color[rgb]{0.19607843,0.19607843,0.19607843}\makebox(0,0)[lt]{\lineheight{1.25}\smash{\begin{tabular}[t]{l}$\Omega_0$\end{tabular}}}}%
    \put(0.07951859,0.10380929){\color[rgb]{0.19607843,0.19607843,0.19607843}\makebox(0,0)[lt]{\lineheight{1.25}\smash{\begin{tabular}[t]{l}$\Omega_1$\end{tabular}}}}%
    \put(0.06801147,0.16609916){\color[rgb]{0.19607843,0.19607843,0.19607843}\makebox(0,0)[lt]{\lineheight{1.25}\smash{\begin{tabular}[t]{l}$\Omega_2$\end{tabular}}}}%
    \put(0.19316901,0.27139635){\color[rgb]{0.19607843,0.19607843,0.19607843}\makebox(0,0)[lt]{\lineheight{1.25}\smash{\begin{tabular}[t]{l}$\AmbientDomain$\end{tabular}}}}%
    \put(0.00009191,0.27658023){\color[rgb]{0.19607843,0.19607843,0.19607843}\makebox(0,0)[lt]{\lineheight{1.25}\smash{\begin{tabular}[t]{l}$(A)$\end{tabular}}}}%
    \put(0.28061914,0.27664628){\color[rgb]{0.19607843,0.19607843,0.19607843}\makebox(0,0)[lt]{\lineheight{1.25}\smash{\begin{tabular}[t]{l}$(B)$\end{tabular}}}}%
    \put(0.57745773,0.27542953){\color[rgb]{0.19607843,0.19607843,0.19607843}\makebox(0,0)[lt]{\lineheight{1.25}\smash{\begin{tabular}[t]{l}$(C)$\end{tabular}}}}%
    \put(0.72776163,0.0643779){\color[rgb]{0.19607843,0.19607843,0.19607843}\makebox(0,0)[lt]{\lineheight{1.25}\smash{\begin{tabular}[t]{l}$S_{E^+}^{3}$\end{tabular}}}}%
    \put(0.63712739,0.14512089){\color[rgb]{0.19607843,0.19607843,0.19607843}\makebox(0,0)[lt]{\lineheight{1.25}\smash{\begin{tabular}[t]{l}$S_{E^+}^{2}$\end{tabular}}}}%
    \put(0.79350506,0.14408938){\color[rgb]{0.19607843,0.19607843,0.19607843}\makebox(0,0)[lt]{\lineheight{1.25}\smash{\begin{tabular}[t]{l}$S_{E^-}^{1}$\end{tabular}}}}%
    \put(0.88380269,0.10474915){\color[rgb]{0.19607843,0.19607843,0.19607843}\makebox(0,0)[lt]{\lineheight{1.25}\smash{\begin{tabular}[t]{l}$S_{E^-}^{2}$\end{tabular}}}}%
    \put(0.92883329,0.22319688){\color[rgb]{0.19607843,0.19607843,0.19607843}\makebox(0,0)[lt]{\lineheight{1.25}\smash{\begin{tabular}[t]{l}$S_{E^-}^{3}$\end{tabular}}}}%
    \put(0.85627023,0.25624055){\color[rgb]{0.19607843,0.19607843,0.19607843}\makebox(0,0)[lt]{\lineheight{1.25}\smash{\begin{tabular}[t]{l}$\Omega^{E^-}$\end{tabular}}}}%
    \put(0,0){\includegraphics[width=\unitlength,page=2]{ghost_for_enriched.pdf}}%
  \end{picture}%
\endgroup%

%% file: figures/ghost_clusters.pdf_tex
\begingroup%
  \makeatletter%
  \providecommand\color[2][]{%
    \errmessage{(Inkscape) Color is used for the text in Inkscape, but the package 'color.sty' is not loaded}%
    \renewcommand\color[2][]{}%
  }%
  \providecommand\transparent[1]{%
    \errmessage{(Inkscape) Transparency is used (non-zero) for the text in Inkscape, but the package 'transparent.sty' is not loaded}%
    \renewcommand\transparent[1]{}%
  }%
  \providecommand\rotatebox[2]{#2}%
  \newcommand*\fsize{\dimexpr\f@size pt\relax}%
  \newcommand*\lineheight[1]{\fontsize{\fsize}{#1\fsize}\selectfont}%
  \ifx\svgwidth\undefined%
    \setlength{\unitlength}{468.95154578bp}%
    \ifx\svgscale\undefined%
      \relax%
    \else%
      \setlength{\unitlength}{\unitlength * \real{\svgscale}}%
    \fi%
  \else%
    \setlength{\unitlength}{\svgwidth}%
  \fi%
  \global\let\svgwidth\undefined%
  \global\let\svgscale\undefined%
  \makeatother%
  \begin{picture}(1,0.30880001)%
    \lineheight{1}%
    \setlength\tabcolsep{0pt}%
    \put(0,0){\includegraphics[width=\unitlength,page=1]{ghost_clusters.pdf}}%
    \put(0.69537247,0.23762125){\color[rgb]{0.19607843,0.19607843,0.19607843}\makebox(0,0)[lt]{\lineheight{1.25}\smash{\begin{tabular}[t]{l}$S_{E^+}^{2}$\end{tabular}}}}%
    \put(0.80780918,0.0969663){\color[rgb]{0.19607843,0.19607843,0.19607843}\makebox(0,0)[lt]{\lineheight{1.25}\smash{\begin{tabular}[t]{l}$S_{E^-}^{1}$\end{tabular}}}}%
    \put(0.44671321,0.2353997){\color[rgb]{0.19607843,0.19607843,0.19607843}\makebox(0,0)[lt]{\lineheight{1.25}\smash{\begin{tabular}[t]{l}$S_{E^+}^{1}$\end{tabular}}}}%
    \put(0.56557051,0.05808267){\color[rgb]{0.19607843,0.19607843,0.19607843}\makebox(0,0)[lt]{\lineheight{1.25}\smash{\begin{tabular}[t]{l}$S_{E^-}^{2}$\end{tabular}}}}%
    \put(0.11076492,0.20256238){\color[rgb]{0.19607843,0.19607843,0.19607843}\makebox(0,0)[lt]{\lineheight{1.25}\smash{\begin{tabular}[t]{l}$S_{E^+}^{3}$\end{tabular}}}}%
    \put(0.26333843,0.05894358){\color[rgb]{0.19607843,0.19607843,0.19607843}\makebox(0,0)[lt]{\lineheight{1.25}\smash{\begin{tabular}[t]{l}$S_{E^-}^{2}$\end{tabular}}}}%
    \put(0,0){\includegraphics[width=\unitlength,page=2]{ghost_clusters.pdf}}%
  \end{picture}%
\endgroup%

%% file: figures/workflow.pdf_tex
\begingroup%
  \makeatletter%
  \providecommand\color[2][]{%
    \errmessage{(Inkscape) Color is used for the text in Inkscape, but the package 'color.sty' is not loaded}%
    \renewcommand\color[2][]{}%
  }%
  \providecommand\transparent[1]{%
    \errmessage{(Inkscape) Transparency is used (non-zero) for the text in Inkscape, but the package 'transparent.sty' is not loaded}%
    \renewcommand\transparent[1]{}%
  }%
  \providecommand\rotatebox[2]{#2}%
  \newcommand*\fsize{\dimexpr\f@size pt\relax}%
  \newcommand*\lineheight[1]{\fontsize{\fsize}{#1\fsize}\selectfont}%
  \ifx\svgwidth\undefined%
    \setlength{\unitlength}{139.3823122bp}%
    \ifx\svgscale\undefined%
      \relax%
    \else%
      \setlength{\unitlength}{\unitlength * \real{\svgscale}}%
    \fi%
  \else%
    \setlength{\unitlength}{\svgwidth}%
  \fi%
  \global\let\svgwidth\undefined%
  \global\let\svgscale\undefined%
  \makeatother%
  \begin{picture}(1,4.22752519)%
    \lineheight{1}%
    \setlength\tabcolsep{0pt}%
    \put(0,0){\includegraphics[width=\unitlength,page=1]{workflow.pdf}}%
    \put(0.24543672,2.3448707){\color[rgb]{0,0.47843137,1}\makebox(0,0)[lt]{\lineheight{1.25}\smash{\begin{tabular}[t]{l}\footnotesize Connectivity graph $\mathcal{G}$\end{tabular}}}}%
    \put(0,0){\includegraphics[width=\unitlength,page=2]{workflow.pdf}}%
    \put(0.00009104,3.48058603){\color[rgb]{0.19607843,0.19607843,0.19607843}\makebox(0,0)[lt]{\lineheight{1.25}\smash{\begin{tabular}[t]{l}Tessellation\end{tabular}}}}%
    \put(-0.00622555,2.21247783){\color[rgb]{0.19607843,0.19607843,0.19607843}\makebox(0,0)[lt]{\lineheight{1.25}\smash{\begin{tabular}[t]{l}Enrichment\end{tabular}}}}%
    \put(-0.00516363,1.5532989){\color[rgb]{0.19607843,0.19607843,0.19607843}\makebox(0,0)[lt]{\lineheight{1.25}\smash{\begin{tabular}[t]{l}Stabilization\end{tabular}}}}%
    \put(-0.0083559,2.89080729){\color[rgb]{0.19607843,0.19607843,0.19607843}\makebox(0,0)[lt]{\lineheight{1.25}\smash{\begin{tabular}[t]{l}Build topology information\end{tabular}}}}%
    \put(0,0){\includegraphics[width=\unitlength,page=3]{workflow.pdf}}%
    \put(0.49212921,3.58307389){\color[rgb]{0.82352941,0,0.00392157}\makebox(0,0)[lt]{\lineheight{1.25}\smash{\begin{tabular}[t]{l}\footnotesize Geometry $G$\end{tabular}}}}%
    \put(0,0){\includegraphics[width=\unitlength,page=4]{workflow.pdf}}%
    \put(-0.00343887,4.17014663){\color[rgb]{0.19607843,0.19607843,0.19607843}\makebox(0,0)[lt]{\lineheight{1.25}\smash{\begin{tabular}[t]{l}Inputs\end{tabular}}}}%
    \put(0.24977129,4.0734444){\color[rgb]{0,0,0}\makebox(0,0)[lt]{\lineheight{1.25}\smash{\begin{tabular}[t]{l}\footnotesize Background mesh $\mathcal{H}$\end{tabular}}}}%
    \put(0.76529642,3.25294322){\color[rgb]{0,0,0}\makebox(0,0)[lt]{\lineheight{1.25}\smash{\begin{tabular}[t]{l}\footnotesize Foreground\end{tabular}}}}%
    \put(0.76845238,3.1804783){\color[rgb]{0,0,0}\makebox(0,0)[lt]{\lineheight{1.25}\smash{\begin{tabular}[t]{l}\footnotesize mesh $\mathcal{T}$\end{tabular}}}}%
    \put(0,0){\includegraphics[width=\unitlength,page=5]{workflow.pdf}}%
    \put(0.70394161,1.68633848){\color[rgb]{0,0.47843137,1}\makebox(0,0)[lt]{\lineheight{1.25}\smash{\begin{tabular}[t]{l}$\scriptstyle \mathrm{supp}(N^1)$\end{tabular}}}}%
    \put(0,0){\includegraphics[width=\unitlength,page=6]{workflow.pdf}}%
    \put(0.12637152,1.6875537){\color[rgb]{0,0.47843137,1}\makebox(0,0)[lt]{\lineheight{1.25}\smash{\begin{tabular}[t]{l}$\scriptstyle \mathrm{supp}(N^2)$\end{tabular}}}}%
    \put(0,0){\includegraphics[width=\unitlength,page=7]{workflow.pdf}}%
    \put(0.40388043,3.10278633){\color[rgb]{0,0,0}\makebox(0,0)[lt]{\lineheight{1.25}\smash{\begin{tabular}[t]{l}$\scriptstyle \Omega_2$\end{tabular}}}}%
    \put(0.64580359,3.35217519){\color[rgb]{0,0,0.05490196}\makebox(0,0)[lt]{\lineheight{1.25}\smash{\begin{tabular}[t]{l}$\scriptstyle \Omega_1$\end{tabular}}}}%
    \put(0,0){\includegraphics[width=\unitlength,page=8]{workflow.pdf}}%
    \put(0.18668396,1.20756653){\color[rgb]{0,0.47843137,1}\makebox(0,0)[lt]{\lineheight{1.25}\smash{\begin{tabular}[t]{l}$\scriptstyle \mathcal{F}_G^2$\end{tabular}}}}%
    \put(0.80874756,1.32958972){\color[rgb]{0,0.47843137,1}\makebox(0,0)[lt]{\lineheight{1.25}\smash{\begin{tabular}[t]{l}$\scriptstyle \mathcal{F}_G^1$\end{tabular}}}}%
    \put(0,0){\includegraphics[width=\unitlength,page=9]{workflow.pdf}}%
    \put(-0.00428977,0.94197892){\color[rgb]{0.19607843,0.19607843,0.19607843}\makebox(0,0)[lt]{\lineheight{1.25}\smash{\begin{tabular}[t]{l}Output\end{tabular}}}}%
    \put(0,0){\includegraphics[width=\unitlength,page=10]{workflow.pdf}}%
    \put(0.0979579,0.52435386){\color[rgb]{0,0,0}\makebox(0,0)[lt]{\lineheight{1.25}\smash{\begin{tabular}[t]{l}\footnotesize Cluster $D$\end{tabular}}}}%
    \put(0.63871058,0.25159038){\color[rgb]{0,0,0}\makebox(0,0)[lt]{\lineheight{1.25}\smash{\begin{tabular}[t]{l}\footnotesize Cluster pairs\end{tabular}}}}%
    \put(0.76958567,0.17329515){\color[rgb]{0,0,0}\makebox(0,0)[lt]{\lineheight{1.25}\smash{\begin{tabular}[t]{l}$\scriptstyle (D_1,D_2)$\end{tabular}}}}%
    \put(0.23296058,0.71596716){\color[rgb]{0.82352941,0,0.00392157}\makebox(0,0)[lt]{\lineheight{1.25}\smash{\begin{tabular}[t]{l}$\scriptstyle \xi_q$\end{tabular}}}}%
  \end{picture}%
\endgroup%

%% file: figures/inputs_outputs.pdf_tex
\begingroup%
  \makeatletter%
  \providecommand\color[2][]{%
    \errmessage{(Inkscape) Color is used for the text in Inkscape, but the package 'color.sty' is not loaded}%
    \renewcommand\color[2][]{}%
  }%
  \providecommand\transparent[1]{%
    \errmessage{(Inkscape) Transparency is used (non-zero) for the text in Inkscape, but the package 'transparent.sty' is not loaded}%
    \renewcommand\transparent[1]{}%
  }%
  \providecommand\rotatebox[2]{#2}%
  \newcommand*\fsize{\dimexpr\f@size pt\relax}%
  \newcommand*\lineheight[1]{\fontsize{\fsize}{#1\fsize}\selectfont}%
  \ifx\svgwidth\undefined%
    \setlength{\unitlength}{211.5bp}%
    \ifx\svgscale\undefined%
      \relax%
    \else%
      \setlength{\unitlength}{\unitlength * \real{\svgscale}}%
    \fi%
  \else%
    \setlength{\unitlength}{\svgwidth}%
  \fi%
  \global\let\svgwidth\undefined%
  \global\let\svgscale\undefined%
  \makeatother%
  \begin{picture}(1,2.29787234)%
    \lineheight{1}%
    \setlength\tabcolsep{0pt}%
    \put(0,0){\includegraphics[width=\unitlength,page=1]{inputs_outputs.pdf}}%
    \put(0.42553191,0.31914894){\makebox(0,0)[t]{\lineheight{1.25}\smash{\begin{tabular}[t]{c}\textbf{Geometry }$\Geometry$\end{tabular}}}}%
    \put(0,0){\includegraphics[width=\unitlength,page=2]{inputs_outputs.pdf}}%
    \put(0.0212766,0.20212766){\makebox(0,0)[lt]{\lineheight{1.25}\smash{\begin{tabular}[t]{l}+ $\mathtt{eval\_proximity}(\pos)$\end{tabular}}}}%
    \put(0,0){\includegraphics[width=\unitlength,page=3]{inputs_outputs.pdf}}%
    \put(0.0212766,0.12056738){\makebox(0,0)[lt]{\lineheight{1.25}\smash{\begin{tabular}[t]{l}+ $\mathtt{is\_intersected}(\BgElemInd)$\end{tabular}}}}%
    \put(0,0){\includegraphics[width=\unitlength,page=4]{inputs_outputs.pdf}}%
    \put(0.0212766,0.04255319){\makebox(0,0)[lt]{\lineheight{1.25}\smash{\begin{tabular}[t]{l}+ $\mathtt{find\_interface}(\pos_1,\pos_2)$\end{tabular}}}}%
    \put(0,0){\includegraphics[width=\unitlength,page=5]{inputs_outputs.pdf}}%
    \put(0.4964539,1.73758865){\makebox(0,0)[t]{\lineheight{1.25}\smash{\begin{tabular}[t]{c}\textbf{Background Element }$\BgElemInd$\end{tabular}}}}%
    \put(0,0){\includegraphics[width=\unitlength,page=6]{inputs_outputs.pdf}}%
    \put(0.0212766,1.64539007){\makebox(0,0)[lt]{\lineheight{1.25}\smash{\begin{tabular}[t]{l}+ Array of vertices $\mathrm{V}$ (sorted by ordinal)\end{tabular}}}}%
    \put(0.0212766,1.55319149){\color[rgb]{0,0,0}\makebox(0,0)[lt]{\lineheight{1.25}\smash{\begin{tabular}[t]{l}+ Basis $\mathcal{B}_{\BgElemInd} = \SetDef{N_{\LocalBfIndex}(\xi)}_{\LocalBfIndex=1}^{n_{\LocalBfIndex}(\BgElemInd)}$\end{tabular}}}}%
    \put(0.0212766,1.46099291){\color[rgb]{0,0,0}\makebox(0,0)[lt]{\lineheight{1.25}\smash{\begin{tabular}[t]{l}+ Local to global basis function map\end{tabular}}}}%
    \put(0.0099911,1.39714965){\color[rgb]{0,0,0}\makebox(0,0)[lt]{\lineheight{0}\smash{\begin{tabular}[t]{l}    $\mathrm{IEN}: \LocalBfIndex \mapsto \BfIndex$\end{tabular}}}}%
    \put(0,0){\includegraphics[width=\unitlength,page=7]{inputs_outputs.pdf}}%
    \put(0.0212766,1.3221224){\makebox(0,0)[lt]{\lineheight{1.25}\smash{\begin{tabular}[t]{l}+ Parallel ID $I$, owning processor $p$\end{tabular}}}}%
    \put(0,0){\includegraphics[width=\unitlength,page=8]{inputs_outputs.pdf}}%
    \put(0.0212766,1.21276596){\makebox(0,0)[lt]{\lineheight{1.25}\smash{\begin{tabular}[t]{l}+ $\pos,\pfrac{\pos}{\xi},\cdots \leftarrow \mathtt{interpolate\_space}(\xi)$\end{tabular}}}}%
    \put(0,0){\includegraphics[width=\unitlength,page=9]{inputs_outputs.pdf}}%
    \put(0.42553191,2.23404255){\makebox(0,0)[t]{\lineheight{1.25}\smash{\begin{tabular}[t]{c}\textbf{Cluster }$\Cluster$\end{tabular}}}}%
    \put(0,0){\includegraphics[width=\unitlength,page=10]{inputs_outputs.pdf}}%
    \put(0.0212766,2.14184397){\makebox(0,0)[lt]{\lineheight{1.25}\smash{\begin{tabular}[t]{l}+ Background element $\BgElemInd$\end{tabular}}}}%
    \put(0.0212766,2.04964539){\color[rgb]{0,0,0}\makebox(0,0)[lt]{\lineheight{1.25}\smash{\begin{tabular}[t]{l}+ Quadrature points $\SetDef{(\xi_q,\omega_q)}_{q=1}^{n_q(\BgElemInd)}$\end{tabular}}}}%
    \put(0.0212766,1.95035461){\color[rgb]{0,0,0}\makebox(0,0)[lt]{\lineheight{1.25}\smash{\begin{tabular}[t]{l}+ For side clusters: normals associated \end{tabular}}}}%
    \put(0.01252428,1.88577879){\color[rgb]{0,0,0}\makebox(0,0)[lt]{\lineheight{1.25}\smash{\begin{tabular}[t]{l}    with quadrature points $\SetDef{\normal_q}_{q=1}^{n_q(\BgElemInd)}$\end{tabular}}}}%
    \put(0,0){\includegraphics[width=\unitlength,page=11]{inputs_outputs.pdf}}%
    \put(0.42553191,1.06382979){\makebox(0,0)[t]{\lineheight{1.25}\smash{\begin{tabular}[t]{c}\textbf{Background Mesh }$\BgMesh$\end{tabular}}}}%
    \put(0,0){\includegraphics[width=\unitlength,page=12]{inputs_outputs.pdf}}%
    \put(0.0212766,0.96808511){\makebox(0,0)[lt]{\lineheight{1.25}\smash{\begin{tabular}[t]{l}+ Array of elements $\mathcal{\BgElemInd}= \SetDef{\BgElemInd}_{\BgElemInd=1}^{n_{\BgElemInd}}$ \end{tabular}}}}%
    \put(0,0){\includegraphics[width=\unitlength,page=13]{inputs_outputs.pdf}}%
    \put(0.0212766,0.85106383){\makebox(0,0)[lt]{\lineheight{1.25}\smash{\begin{tabular}[t]{l}+ $\mathtt{get\_number\_of\_basis\_fncts}()$\end{tabular}}}}%
    \put(0,0){\includegraphics[width=\unitlength,page=14]{inputs_outputs.pdf}}%
    \put(0.0212766,0.75886525){\makebox(0,0)[lt]{\lineheight{1.25}\smash{\begin{tabular}[t]{l}+ $\mathtt{get\_elements\_in\_support}(\BfIndex)$\end{tabular}}}}%
    \put(0,0){\includegraphics[width=\unitlength,page=15]{inputs_outputs.pdf}}%
    \put(0.0212766,0.66666667){\makebox(0,0)[lt]{\lineheight{1.25}\smash{\begin{tabular}[t]{l}+ $\mathtt{get\_entities\_on\_element}(\Rank,\BgElemInd)$\end{tabular}}}}%
    \put(0,0){\includegraphics[width=\unitlength,page=16]{inputs_outputs.pdf}}%
    \put(0.0212766,0.57446809){\makebox(0,0)[lt]{\lineheight{1.25}\smash{\begin{tabular}[t]{l}+ $\mathtt{get\_elements\_on\_entity}(\Rank,\Entity)$\end{tabular}}}}%
    \put(0,0){\includegraphics[width=\unitlength,page=17]{inputs_outputs.pdf}}%
    \put(0.0212766,0.4822695){\makebox(0,0)[lt]{\lineheight{1.25}\smash{\begin{tabular}[t]{l}+ $\mathtt{get\_basis\_fnct\_ID}(\BfIndex)$\end{tabular}}}}%
    \put(0,0){\includegraphics[width=\unitlength,page=18]{inputs_outputs.pdf}}%
  \end{picture}%
\endgroup%

%% file: figures/entity_connectivity_data_structure.pdf_tex
\begingroup%
  \makeatletter%
  \providecommand\color[2][]{%
    \errmessage{(Inkscape) Color is used for the text in Inkscape, but the package 'color.sty' is not loaded}%
    \renewcommand\color[2][]{}%
  }%
  \providecommand\transparent[1]{%
    \errmessage{(Inkscape) Transparency is used (non-zero) for the text in Inkscape, but the package 'transparent.sty' is not loaded}%
    \renewcommand\transparent[1]{}%
  }%
  \providecommand\rotatebox[2]{#2}%
  \newcommand*\fsize{\dimexpr\f@size pt\relax}%
  \newcommand*\lineheight[1]{\fontsize{\fsize}{#1\fsize}\selectfont}%
  \ifx\svgwidth\undefined%
    \setlength{\unitlength}{195.75bp}%
    \ifx\svgscale\undefined%
      \relax%
    \else%
      \setlength{\unitlength}{\unitlength * \real{\svgscale}}%
    \fi%
  \else%
    \setlength{\unitlength}{\svgwidth}%
  \fi%
  \global\let\svgwidth\undefined%
  \global\let\svgscale\undefined%
  \makeatother%
  \begin{picture}(1,0.36015326)%
    \lineheight{1}%
    \setlength\tabcolsep{0pt}%
    \put(0,0){\includegraphics[width=\unitlength,page=1]{entity_connectivity_data_structure.pdf}}%
    \put(0.49808429,0.29118774){\makebox(0,0)[t]{\lineheight{1.25}\smash{\begin{tabular}[t]{c}\textbf{Entity Connectivity}\end{tabular}}}}%
    \put(0,0){\includegraphics[width=\unitlength,page=2]{entity_connectivity_data_structure.pdf}}%
    \put(0.02298851,0.19157088){\makebox(0,0)[lt]{\lineheight{1.25}\smash{\begin{tabular}[t]{l}+ Entity to Cell $\mathrm{EtC} = \SetDef{\SetDef{\cdots}}_{e=1}^{\mathtt{NumEntities}}$\end{tabular}}}}%
    \put(0,0){\includegraphics[width=\unitlength,page=3]{entity_connectivity_data_structure.pdf}}%
    \put(0.02298851,0.09195402){\makebox(0,0)[lt]{\lineheight{1.25}\smash{\begin{tabular}[t]{l}+ Cell to Entity $\mathrm{CtE} = \SetDef{\SetDef{\cdots}}_{c=1}^{\mathtt{NumCells}}$\end{tabular}}}}%
    \put(0.02298851,0.03831418){\makebox(0,0)[lt]{\lineheight{1.25}\smash{\begin{tabular}[t]{l}   (entities sorted by ordinal)\end{tabular}}}}%
  \end{picture}%
\endgroup%

%% file: figures/entity_connectivity.pdf_tex
\begingroup%
  \makeatletter%
  \providecommand\color[2][]{%
    \errmessage{(Inkscape) Color is used for the text in Inkscape, but the package 'color.sty' is not loaded}%
    \renewcommand\color[2][]{}%
  }%
  \providecommand\transparent[1]{%
    \errmessage{(Inkscape) Transparency is used (non-zero) for the text in Inkscape, but the package 'transparent.sty' is not loaded}%
    \renewcommand\transparent[1]{}%
  }%
  \providecommand\rotatebox[2]{#2}%
  \newcommand*\fsize{\dimexpr\f@size pt\relax}%
  \newcommand*\lineheight[1]{\fontsize{\fsize}{#1\fsize}\selectfont}%
  \ifx\svgwidth\undefined%
    \setlength{\unitlength}{207.45820774bp}%
    \ifx\svgscale\undefined%
      \relax%
    \else%
      \setlength{\unitlength}{\unitlength * \real{\svgscale}}%
    \fi%
  \else%
    \setlength{\unitlength}{\svgwidth}%
  \fi%
  \global\let\svgwidth\undefined%
  \global\let\svgscale\undefined%
  \makeatother%
  \begin{picture}(1,1.09926548)%
    \lineheight{1}%
    \setlength\tabcolsep{0pt}%
    \put(0,0){\includegraphics[width=\unitlength,page=1]{entity_connectivity.pdf}}%
    \put(0.07672694,0.95905809){\color[rgb]{0,0.47843137,1}\makebox(0,0)[lt]{\lineheight{1.25}\smash{\begin{tabular}[t]{l}$1$\end{tabular}}}}%
    \put(0.23618468,0.66199386){\color[rgb]{0,0.47843137,1}\makebox(0,0)[lt]{\lineheight{1.25}\smash{\begin{tabular}[t]{l}$2$\end{tabular}}}}%
    \put(0.50281599,0.9631929){\color[rgb]{0,0.47843137,1}\makebox(0,0)[lt]{\lineheight{1.25}\smash{\begin{tabular}[t]{l}$3$\end{tabular}}}}%
    \put(0.28126774,0.96499926){\color[rgb]{0.49411765,0.77254902,0.21568627}\makebox(0,0)[lt]{\lineheight{1.25}\smash{\begin{tabular}[t]{l}$3$\end{tabular}}}}%
    \put(0.21887257,0.8613957){\color[rgb]{0.49411765,0.77254902,0.21568627}\makebox(0,0)[lt]{\lineheight{1.25}\smash{\begin{tabular}[t]{l}$1$\end{tabular}}}}%
    \put(0.34834596,0.85834067){\color[rgb]{0.49411765,0.77254902,0.21568627}\makebox(0,0)[lt]{\lineheight{1.25}\smash{\begin{tabular}[t]{l}$2$\end{tabular}}}}%
    \put(0,0){\includegraphics[width=\unitlength,page=2]{entity_connectivity.pdf}}%
    \put(0.22573823,0.45731594){\color[rgb]{0,0.47843137,1}\makebox(0,0)[lt]{\lineheight{1.25}\smash{\begin{tabular}[t]{l}$1$\end{tabular}}}}%
    \put(0.90254993,0.45307168){\color[rgb]{0,0.47843137,1}\makebox(0,0)[lt]{\lineheight{1.25}\smash{\begin{tabular}[t]{l}$2$\end{tabular}}}}%
    \put(0.66931311,0.89253223){\color[rgb]{0,0.47843137,1}\makebox(0,0)[lt]{\lineheight{1.25}\smash{\begin{tabular}[t]{l}$4$\end{tabular}}}}%
    \put(0.2632004,1.05898053){\color[rgb]{0,0,0}\makebox(0,0)[lt]{\lineheight{1.25}\smash{\begin{tabular}[t]{l}(A)\end{tabular}}}}%
    \put(0.73377515,0.79759211){\color[rgb]{0,0,0}\makebox(0,0)[lt]{\lineheight{1.25}\smash{\begin{tabular}[t]{l}(B)\end{tabular}}}}%
    \put(0.50660237,0.68588072){\color[rgb]{0.49411765,0.77254902,0.21568627}\makebox(0,0)[lt]{\lineheight{1.25}\smash{\begin{tabular}[t]{l}$6$\end{tabular}}}}%
    \put(0.6342478,0.6765328){\color[rgb]{0.49411765,0.77254902,0.21568627}\makebox(0,0)[lt]{\lineheight{1.25}\smash{\begin{tabular}[t]{l}$4$\end{tabular}}}}%
    \put(0.67206775,0.45068299){\color[rgb]{0.49411765,0.77254902,0.21568627}\makebox(0,0)[lt]{\lineheight{1.25}\smash{\begin{tabular}[t]{l}$1$\end{tabular}}}}%
    \put(0.76074414,0.51577756){\color[rgb]{0.49411765,0.77254902,0.21568627}\makebox(0,0)[lt]{\lineheight{1.25}\smash{\begin{tabular}[t]{l}$2$\end{tabular}}}}%
    \put(0.43008968,0.449377){\color[rgb]{0.49411765,0.77254902,0.21568627}\makebox(0,0)[lt]{\lineheight{1.25}\smash{\begin{tabular}[t]{l}$3$\end{tabular}}}}%
    \put(0.67254345,0.73074329){\color[rgb]{0.49411765,0.77254902,0.21568627}\makebox(0,0)[lt]{\lineheight{1.25}\smash{\begin{tabular}[t]{l}$5$\end{tabular}}}}%
    \put(0,0){\includegraphics[width=\unitlength,page=3]{entity_connectivity.pdf}}%
    \put(0.61708077,0.58365174){\color[rgb]{0,0.47843137,1}\makebox(0,0)[lt]{\lineheight{1.25}\smash{\begin{tabular}[t]{l}$3$\end{tabular}}}}%
    \put(0,0){\includegraphics[width=\unitlength,page=4]{entity_connectivity.pdf}}%
    \put(0.48474485,0.59286498){\color[rgb]{0.82352941,0,0.00392157}\makebox(0,0)[lt]{\lineheight{1.25}\smash{\begin{tabular}[t]{l}$3$\end{tabular}}}}%
    \put(0.69892153,0.6063083){\color[rgb]{0.82352941,0,0.00392157}\makebox(0,0)[lt]{\lineheight{1.25}\smash{\begin{tabular}[t]{l}$2$\end{tabular}}}}%
    \put(0.58338549,0.48247293){\color[rgb]{0.82352941,0,0.00392157}\makebox(0,0)[lt]{\lineheight{1.25}\smash{\begin{tabular}[t]{l}$1$\end{tabular}}}}%
    \put(0,0){\includegraphics[width=\unitlength,page=5]{entity_connectivity.pdf}}%
    \put(-0.00168656,0.24395893){\color[rgb]{0,0,0}\makebox(0,0)[lt]{\lineheight{1.25}\smash{\begin{tabular}[t]{l}$\mathcal{C}=\SetDef{\SetDef{1,2,3},\SetDef{4,3,2},\SetDef{3,4,1}}$\end{tabular}}}}%
    \put(-0.00340612,0.0805399){\color[rgb]{0,0,0}\makebox(0,0)[lt]{\lineheight{1.25}\smash{\begin{tabular}[t]{l}$\mathrm{EtC}=\SetDef{\SetDef{1},\SetDef{1,2},\SetDef{1,3},\SetDef{2,3},\SetDef{2},\SetDef{3}}$\end{tabular}}}}%
    \put(-0.00181231,0.16243305){\color[rgb]{0,0,0}\makebox(0,0)[lt]{\lineheight{1.25}\smash{\begin{tabular}[t]{l}\underline{Edge connectivity}:\end{tabular}}}}%
    \put(-0.00331459,0.01248653){\color[rgb]{0,0,0}\makebox(0,0)[lt]{\lineheight{1.25}\smash{\begin{tabular}[t]{l}$\mathrm{CtE}=\SetDef{\SetDef{1,2,3},\SetDef{4,2,5},\SetDef{4,6,3}}$\end{tabular}}}}%
    \put(-0.00359137,0.32213322){\color[rgb]{0,0,0}\makebox(0,0)[lt]{\lineheight{1.25}\smash{\begin{tabular}[t]{l}(C)\end{tabular}}}}%
    \put(0,0){\includegraphics[width=\unitlength,page=6]{entity_connectivity.pdf}}%
  \end{picture}%
\endgroup%

%% file: figures/subdivision_process.pdf_tex
\begingroup%
  \makeatletter%
  \providecommand\color[2][]{%
    \errmessage{(Inkscape) Color is used for the text in Inkscape, but the package 'color.sty' is not loaded}%
    \renewcommand\color[2][]{}%
  }%
  \providecommand\transparent[1]{%
    \errmessage{(Inkscape) Transparency is used (non-zero) for the text in Inkscape, but the package 'transparent.sty' is not loaded}%
    \renewcommand\transparent[1]{}%
  }%
  \providecommand\rotatebox[2]{#2}%
  \newcommand*\fsize{\dimexpr\f@size pt\relax}%
  \newcommand*\lineheight[1]{\fontsize{\fsize}{#1\fsize}\selectfont}%
  \ifx\svgwidth\undefined%
    \setlength{\unitlength}{192.35458014bp}%
    \ifx\svgscale\undefined%
      \relax%
    \else%
      \setlength{\unitlength}{\unitlength * \real{\svgscale}}%
    \fi%
  \else%
    \setlength{\unitlength}{\svgwidth}%
  \fi%
  \global\let\svgwidth\undefined%
  \global\let\svgscale\undefined%
  \makeatother%
  \begin{picture}(1,2.92126576)%
    \lineheight{1}%
    \setlength\tabcolsep{0pt}%
    \put(0,0){\includegraphics[width=\unitlength,page=1]{subdivision_process.pdf}}%
    \put(0.4623159,2.02034094){\color[rgb]{0,0,0}\makebox(0,0)[lt]{\lineheight{1.25}\smash{\begin{tabular}[t]{l}$m \! = \! 1$\end{tabular}}}}%
    \put(0.22544006,2.72132816){\color[rgb]{0,0.47843137,1}\makebox(0,0)[lt]{\lineheight{1.25}\smash{\begin{tabular}[t]{l}$G_1$\end{tabular}}}}%
    \put(0.5397583,2.54530496){\color[rgb]{0.82352941,0,0}\makebox(0,0)[lt]{\lineheight{1.25}\smash{\begin{tabular}[t]{l}$G_2$\end{tabular}}}}%
    \put(0.76905977,1.61396802){\color[rgb]{0,0,0}\makebox(0,0)[lt]{\lineheight{1.25}\smash{\begin{tabular}[t]{l}$m \! = \! 0$\end{tabular}}}}%
    \put(0.00396064,2.4915944){\color[rgb]{0,0,0}\makebox(0,0)[lt]{\lineheight{1.25}\smash{\begin{tabular}[t]{l}(A)\end{tabular}}}}%
    \put(-0.00118701,1.75743844){\color[rgb]{0,0,0}\makebox(0,0)[lt]{\lineheight{1.25}\smash{\begin{tabular}[t]{l}(B)\end{tabular}}}}%
    \put(0.00151819,1.02181899){\color[rgb]{0,0,0}\makebox(0,0)[lt]{\lineheight{1.25}\smash{\begin{tabular}[t]{l}(C)\end{tabular}}}}%
    \put(-0.00413518,0.28446355){\color[rgb]{0,0,0}\makebox(0,0)[lt]{\lineheight{1.25}\smash{\begin{tabular}[t]{l}(D)\end{tabular}}}}%
    \put(0.18093681,0.76240965){\color[rgb]{0,0,0}\makebox(0,0)[lt]{\lineheight{1.25}\smash{\begin{tabular}[t]{l}$m \! = \! 0$\end{tabular}}}}%
    \put(0.18396098,0.02766352){\color[rgb]{0,0,0}\makebox(0,0)[lt]{\lineheight{1.25}\smash{\begin{tabular}[t]{l}$m \! = \! 0$\end{tabular}}}}%
    \put(0.76923561,0.87159096){\color[rgb]{0,0,0}\makebox(0,0)[lt]{\lineheight{1.25}\smash{\begin{tabular}[t]{l}$m \! = \! 1$\end{tabular}}}}%
    \put(0.16826632,1.27878476){\color[rgb]{0,0,0}\makebox(0,0)[lt]{\lineheight{1.25}\smash{\begin{tabular}[t]{l}$m \! = \! 2$\end{tabular}}}}%
    \put(0.75951352,1.27706513){\color[rgb]{0.97647059,1,1}\makebox(0,0)[lt]{\lineheight{1.25}\smash{\begin{tabular}[t]{l}$m \! = \! 3$\end{tabular}}}}%
    \put(0.16188784,0.54463665){\color[rgb]{0,0,0}\makebox(0,0)[lt]{\lineheight{1.25}\smash{\begin{tabular}[t]{l}$m \! = \! 2$\end{tabular}}}}%
    \put(0.77818085,0.14130018){\color[rgb]{0,0,0}\makebox(0,0)[lt]{\lineheight{1.25}\smash{\begin{tabular}[t]{l}$m \! = \! 1$\end{tabular}}}}%
    \put(0,0){\includegraphics[width=\unitlength,page=2]{subdivision_process.pdf}}%
    \put(0.10993738,2.88555037){\color[rgb]{0,0,0}\makebox(0,0)[lt]{\lineheight{1.25}\smash{\begin{tabular}[t]{l}Regular \\subdivision\end{tabular}}}}%
    \put(0.71713697,2.4547332){\color[rgb]{0.38823529,0.38823529,0.38823529}\makebox(0,0)[lt]{\lineheight{1.25}\smash{\begin{tabular}[t]{l}{\small intersected} \\{\small background}\\{\small elements}\end{tabular}}}}%
    \put(0.11266608,2.15406038){\color[rgb]{0,0,0}\makebox(0,0)[lt]{\lineheight{1.25}\smash{\begin{tabular}[t]{l}Templated\\subdivision $G_1$\end{tabular}}}}%
    \put(0.11232519,1.41664532){\color[rgb]{0,0,0}\makebox(0,0)[lt]{\lineheight{1.25}\smash{\begin{tabular}[t]{l}Templated\\subdivision $G_2$\end{tabular}}}}%
    \put(0.11137009,0.67922963){\color[rgb]{0,0,0}\makebox(0,0)[lt]{\lineheight{1.25}\smash{\begin{tabular}[t]{l}Apply\\material map\end{tabular}}}}%
  \end{picture}%
\endgroup%

%% file: figures/ancestry_and_entities_new.pdf_tex
\begingroup%
  \makeatletter%
  \providecommand\color[2][]{%
    \errmessage{(Inkscape) Color is used for the text in Inkscape, but the package 'color.sty' is not loaded}%
    \renewcommand\color[2][]{}%
  }%
  \providecommand\transparent[1]{%
    \errmessage{(Inkscape) Transparency is used (non-zero) for the text in Inkscape, but the package 'transparent.sty' is not loaded}%
    \renewcommand\transparent[1]{}%
  }%
  \providecommand\rotatebox[2]{#2}%
  \newcommand*\fsize{\dimexpr\f@size pt\relax}%
  \newcommand*\lineheight[1]{\fontsize{\fsize}{#1\fsize}\selectfont}%
  \ifx\svgwidth\undefined%
    \setlength{\unitlength}{225.03058006bp}%
    \ifx\svgscale\undefined%
      \relax%
    \else%
      \setlength{\unitlength}{\unitlength * \real{\svgscale}}%
    \fi%
  \else%
    \setlength{\unitlength}{\svgwidth}%
  \fi%
  \global\let\svgwidth\undefined%
  \global\let\svgscale\undefined%
  \makeatother%
  \begin{picture}(1,1.6189971)%
    \lineheight{1}%
    \setlength\tabcolsep{0pt}%
    \put(0.16315377,0.561479){\color[rgb]{0.82352941,0,0.00392157}\transparent{0.89999998}\makebox(0,0)[lt]{\lineheight{1.25}\smash{\begin{tabular}[t]{l}$\mathrm{v}_3\!:\!(1,12)$\end{tabular}}}}%
    \put(0,0){\includegraphics[width=\unitlength,page=1]{ancestry_and_entities_new.pdf}}%
    \put(0.45511608,1.1288072){\color[rgb]{0,0.47843137,1}\makebox(0,0)[lt]{\lineheight{1.25}\smash{\begin{tabular}[t]{l}$9\!:\!(2,3)$\end{tabular}}}}%
    \put(0.07014826,0.80904047){\color[rgb]{0,0.47843137,1}\makebox(0,0)[lt]{\lineheight{1.25}\smash{\begin{tabular}[t]{l}$4\!:\!(0,4)$\end{tabular}}}}%
    \put(0.51654543,0.5274731){\color[rgb]{0,0.47843137,1}\transparent{0.69999999}\makebox(0,0)[lt]{\lineheight{1.25}\smash{\begin{tabular}[t]{l}$1$\end{tabular}}}}%
    \put(0.96164711,0.74272705){\color[rgb]{0,0.47843137,1}\transparent{0.69999999}\makebox(0,0)[lt]{\lineheight{1.25}\smash{\begin{tabular}[t]{l}$2$\end{tabular}}}}%
    \put(0.96054382,1.18698319){\color[rgb]{0,0.47843137,1}\transparent{0.69999999}\makebox(0,0)[lt]{\lineheight{1.25}\smash{\begin{tabular}[t]{l}$6$\end{tabular}}}}%
    \put(0.64527751,1.45476628){\color[rgb]{0,0.47843137,1}\transparent{0.69999999}\makebox(0,0)[lt]{\lineheight{1.25}\smash{\begin{tabular}[t]{l}$7$\end{tabular}}}}%
    \put(0.60792946,0.95109525){\color[rgb]{0,0.47843137,1}\transparent{0.69999999}\makebox(0,0)[lt]{\lineheight{1.25}\smash{\begin{tabular}[t]{l}$3$\end{tabular}}}}%
    \put(0.54624447,1.03934259){\color[rgb]{0,0.47843137,1}\transparent{0.69999999}\makebox(0,0)[lt]{\lineheight{1.25}\smash{\begin{tabular}[t]{l}$5$\end{tabular}}}}%
    \put(0,0){\includegraphics[width=\unitlength,page=2]{ancestry_and_entities_new.pdf}}%
    \put(0.33524722,1.09564703){\color[rgb]{0.82352941,0,0}\transparent{0.90196103}\makebox(0,0)[lt]{\lineheight{1.25}\smash{\begin{tabular}[t]{l}$\mathrm{v}_2$\end{tabular}}}}%
    \put(0.34300868,1.20726968){\color[rgb]{0.82352941,0,0}\transparent{0.90196103}\makebox(0,0)[lt]{\lineheight{1.25}\smash{\begin{tabular}[t]{l}$\mathrm{v}_1$\end{tabular}}}}%
    \put(0.80789171,0.31533777){\color[rgb]{0.82352941,0,0.00392157}\transparent{0.89999998}\makebox(0,0)[lt]{\lineheight{1.25}\smash{\begin{tabular}[t]{l}$2$\end{tabular}}}}%
    \put(0.60927339,0.33463514){\color[rgb]{0.82352941,0,0.00392157}\transparent{0.89999998}\makebox(0,0)[lt]{\lineheight{1.25}\smash{\begin{tabular}[t]{l}$3$\end{tabular}}}}%
    \put(0.28465258,0.36796273){\color[rgb]{0,0,0.00392157}\transparent{0.89999998}\makebox(0,0)[lt]{\lineheight{1.25}\smash{\begin{tabular}[t]{l}$x$\end{tabular}}}}%
    \put(0.12479546,0.39538041){\color[rgb]{0,0,0.00392157}\transparent{0.89999998}\makebox(0,0)[lt]{\lineheight{1.25}\smash{\begin{tabular}[t]{l}$y$\end{tabular}}}}%
    \put(0.1966924,0.43056495){\color[rgb]{0,0,0.00392157}\transparent{0.89999998}\makebox(0,0)[lt]{\lineheight{1.25}\smash{\begin{tabular}[t]{l}$z$\end{tabular}}}}%
    \put(0.90051953,1.45635586){\color[rgb]{0,0,0.00392157}\transparent{0.89999998}\makebox(0,0)[lt]{\lineheight{1.25}\smash{\begin{tabular}[t]{l}(A)\end{tabular}}}}%
    \put(0.91193934,0.5288551){\color[rgb]{0,0,0.00392157}\transparent{0.89999998}\makebox(0,0)[lt]{\lineheight{1.25}\smash{\begin{tabular}[t]{l}(B)\end{tabular}}}}%
    \put(0.58032475,0.21694075){\color[rgb]{0.82352941,0,0.00392157}\transparent{0.89999998}\makebox(0,0)[lt]{\lineheight{1.25}\smash{\begin{tabular}[t]{l}$4$\end{tabular}}}}%
    \put(0.69686491,0.13206728){\color[rgb]{0.82352941,0,0.00392157}\transparent{0.89999998}\makebox(0,0)[lt]{\lineheight{1.25}\smash{\begin{tabular}[t]{l}$5$\end{tabular}}}}%
    \put(0.6948863,0.3878325){\color[rgb]{0.82352941,0,0.00392157}\transparent{0.89999998}\makebox(0,0)[lt]{\lineheight{1.25}\smash{\begin{tabular}[t]{l}$6$\end{tabular}}}}%
    \put(0.81286343,0.01791899){\color[rgb]{0.33333333,0.70980392,0}\makebox(0,0)[lt]{\lineheight{1.25}\smash{\begin{tabular}[t]{l}$1$\end{tabular}}}}%
    \put(0.8308134,0.11995475){\color[rgb]{0.33333333,0.70980392,0}\makebox(0,0)[lt]{\lineheight{1.25}\smash{\begin{tabular}[t]{l}$2$\end{tabular}}}}%
    \put(0.55271082,0.14083189){\color[rgb]{0.33333333,0.70980392,0}\makebox(0,0)[lt]{\lineheight{1.25}\smash{\begin{tabular}[t]{l}$3$\end{tabular}}}}%
    \put(0.54706386,0.05102111){\color[rgb]{0.33333333,0.70980392,0}\makebox(0,0)[lt]{\lineheight{1.25}\smash{\begin{tabular}[t]{l}$4$\end{tabular}}}}%
    \put(0.74242344,0.32870352){\color[rgb]{0.33333333,0.70980392,0}\makebox(0,0)[lt]{\lineheight{1.25}\smash{\begin{tabular}[t]{l}$5$\end{tabular}}}}%
    \put(0.8258115,0.47699822){\color[rgb]{0.33333333,0.70980392,0}\makebox(0,0)[lt]{\lineheight{1.25}\smash{\begin{tabular}[t]{l}$6$\end{tabular}}}}%
    \put(0.61389529,0.51047361){\color[rgb]{0.33333333,0.70980392,0}\makebox(0,0)[lt]{\lineheight{1.25}\smash{\begin{tabular}[t]{l}$7$\end{tabular}}}}%
    \put(0.5258319,0.33477161){\color[rgb]{0.33333333,0.70980392,0}\makebox(0,0)[lt]{\lineheight{1.25}\smash{\begin{tabular}[t]{l}$8$\end{tabular}}}}%
    \put(0.70616714,0.06396688){\color[rgb]{0.33333333,0.70980392,0}\makebox(0,0)[lt]{\lineheight{1.25}\smash{\begin{tabular}[t]{l}$9$\end{tabular}}}}%
    \put(0.93198459,0.23528879){\color[rgb]{0.33333333,0.70980392,0}\makebox(0,0)[lt]{\lineheight{1.25}\smash{\begin{tabular}[t]{l}$10$\end{tabular}}}}%
    \put(0.66053615,0.45966351){\color[rgb]{0.33333333,0.70980392,0}\makebox(0,0)[lt]{\lineheight{1.25}\smash{\begin{tabular}[t]{l}$11$\end{tabular}}}}%
    \put(0.43430967,0.26867577){\color[rgb]{0.33333333,0.70980392,0}\makebox(0,0)[lt]{\lineheight{1.25}\smash{\begin{tabular}[t]{l}$12$\end{tabular}}}}%
    \put(0.18094978,1.07626124){\color[rgb]{0.82352941,0,0.00392157}\transparent{0.89999998}\makebox(0,0)[lt]{\lineheight{1.25}\smash{\begin{tabular}[t]{l}$\mathrm{v}_3$\end{tabular}}}}%
    \put(0.42397893,0.89432892){\color[rgb]{0,0.47843137,1}\makebox(0,0)[lt]{\lineheight{1.25}\smash{\begin{tabular}[t]{l}$15\!:\!(2,4)$\end{tabular}}}}%
    \put(0,0){\includegraphics[width=\unitlength,page=3]{ancestry_and_entities_new.pdf}}%
    \put(0.24648898,1.58138058){\color[rgb]{0.60784314,0.60784314,0.60784314}\makebox(0,0)[lt]{\lineheight{1.25}\smash{\begin{tabular}[t]{l}$F_1\!:\!(3,1)$\end{tabular}}}}%
    \put(0.27473106,0.74067174){\color[rgb]{0.82352941,0,0}\transparent{0.90196103}\makebox(0,0)[lt]{\lineheight{1.25}\smash{\begin{tabular}[t]{l}$\mathrm{v}_1\!:\!(2,3)$\end{tabular}}}}%
    \put(-0.00369963,0.67042877){\color[rgb]{0.60784314,0.60784314,0.60784314}\makebox(0,0)[lt]{\lineheight{1.25}\smash{\begin{tabular}[t]{l}$F_4\!:\!(3,1)$\end{tabular}}}}%
    \put(0.27457556,0.62373851){\color[rgb]{0.82352941,0,0}\transparent{0.90196103}\makebox(0,0)[lt]{\lineheight{1.25}\smash{\begin{tabular}[t]{l}$\mathrm{v}_2\!:\!(2,4)$\end{tabular}}}}%
    \put(0.21009198,1.29187367){\color[rgb]{0,0.47843137,1}\makebox(0,0)[lt]{\lineheight{1.25}\smash{\begin{tabular}[t]{l}$8\!:\!(0,8)$\end{tabular}}}}%
    \put(0,0){\includegraphics[width=\unitlength,page=4]{ancestry_and_entities_new.pdf}}%
    \put(0.05142637,1.42942206){\color[rgb]{0.60784314,0.60784314,0.60784314}\makebox(0,0)[lt]{\lineheight{1.25}\smash{\begin{tabular}[t]{l}$F_2\!:\!(2,3)$\end{tabular}}}}%
    \put(0.02269954,0.9137545){\color[rgb]{0.60784314,0.60784314,0.60784314}\makebox(0,0)[lt]{\lineheight{1.25}\smash{\begin{tabular}[t]{l}$F_3\!:\!(2,4)$\end{tabular}}}}%
    \put(0,0){\includegraphics[width=\unitlength,page=5]{ancestry_and_entities_new.pdf}}%
    \put(0.79512028,0.1923822){\color[rgb]{0.82352941,0,0.00392157}\transparent{0.89999998}\makebox(0,0)[lt]{\lineheight{1.25}\smash{\begin{tabular}[t]{l}$1$\end{tabular}}}}%
  \end{picture}%
\endgroup%

%% file: figures/fg_mesh_abrev.pdf_tex
\begingroup%
  \makeatletter%
  \providecommand\color[2][]{%
    \errmessage{(Inkscape) Color is used for the text in Inkscape, but the package 'color.sty' is not loaded}%
    \renewcommand\color[2][]{}%
  }%
  \providecommand\transparent[1]{%
    \errmessage{(Inkscape) Transparency is used (non-zero) for the text in Inkscape, but the package 'transparent.sty' is not loaded}%
    \renewcommand\transparent[1]{}%
  }%
  \providecommand\rotatebox[2]{#2}%
  \newcommand*\fsize{\dimexpr\f@size pt\relax}%
  \newcommand*\lineheight[1]{\fontsize{\fsize}{#1\fsize}\selectfont}%
  \ifx\svgwidth\undefined%
    \setlength{\unitlength}{210.75bp}%
    \ifx\svgscale\undefined%
      \relax%
    \else%
      \setlength{\unitlength}{\unitlength * \real{\svgscale}}%
    \fi%
  \else%
    \setlength{\unitlength}{\svgwidth}%
  \fi%
  \global\let\svgwidth\undefined%
  \global\let\svgscale\undefined%
  \makeatother%
  \begin{picture}(1,1.99644128)%
    \lineheight{1}%
    \setlength\tabcolsep{0pt}%
    \put(0,0){\includegraphics[width=\unitlength,page=1]{fg_mesh_abrev.pdf}}%
    \put(0.49822064,1.2455516){\makebox(0,0)[t]{\lineheight{1.25}\smash{\begin{tabular}[t]{c}\textbf{Foreground Mesh }$\FgMesh$\end{tabular}}}}%
    \put(0,0){\includegraphics[width=\unitlength,page=2]{fg_mesh_abrev.pdf}}%
    \put(0.12811388,1.15302491){\makebox(0,0)[lt]{\lineheight{1.25}\smash{\begin{tabular}[t]{l}+ Background mesh $\mathcal{H}$\end{tabular}}}}%
    \put(0,0){\includegraphics[width=\unitlength,page=3]{fg_mesh_abrev.pdf}}%
    \put(0.12811388,1.06049822){\makebox(0,0)[lt]{\lineheight{1.25}\smash{\begin{tabular}[t]{l}+ Array of fg. vertices $\mathcal{V}$\end{tabular}}}}%
    \put(0,0){\includegraphics[width=\unitlength,page=4]{fg_mesh_abrev.pdf}}%
    \put(0.12811388,0.96797153){\makebox(0,0)[lt]{\lineheight{1.25}\smash{\begin{tabular}[t]{l}+ Array of fg. cells $\mathcal{C}$\end{tabular}}}}%
    \put(0,0){\includegraphics[width=\unitlength,page=5]{fg_mesh_abrev.pdf}}%
    \put(0.12811388,0.87544484){\makebox(0,0)[lt]{\lineheight{1.25}\smash{\begin{tabular}[t]{l}+ Array of child meshes $\mathrm{\ListOfCMs}$\end{tabular}}}}%
    \put(0,0){\includegraphics[width=\unitlength,page=6]{fg_mesh_abrev.pdf}}%
    \put(0.12811388,0.78291815){\makebox(0,0)[lt]{\lineheight{1.25}\smash{\begin{tabular}[t]{l}+ Array of Subphases $\mathcal{S}$\end{tabular}}}}%
    \put(0,0){\includegraphics[width=\unitlength,page=7]{fg_mesh_abrev.pdf}}%
    \put(0.12811388,0.69039146){\makebox(0,0)[lt]{\lineheight{1.25}\smash{\begin{tabular}[t]{l}+ Subphase graph $\mathcal{G}_S$\end{tabular}}}}%
    \put(0,0){\includegraphics[width=\unitlength,page=8]{fg_mesh_abrev.pdf}}%
    \put(0.12811388,0.59786477){\color[rgb]{0,0,0}\makebox(0,0)[lt]{\lineheight{1.25}\smash{\begin{tabular}[t]{l}+ Graph of disconnected\end{tabular}}}}%
    \put(0.13776166,0.54472074){\color[rgb]{0,0,0}\makebox(0,0)[lt]{\lineheight{1.25}\smash{\begin{tabular}[t]{l}   subphase neighbors $\mathcal{G}_I$\end{tabular}}}}%
    \put(0,0){\includegraphics[width=\unitlength,page=9]{fg_mesh_abrev.pdf}}%
    \put(0.12811388,0.46263345){\makebox(0,0)[lt]{\lineheight{1.25}\smash{\begin{tabular}[t]{l}+ Facet connectivity $\mathcal{F}$\end{tabular}}}}%
    \put(0,0){\includegraphics[width=\unitlength,page=10]{fg_mesh_abrev.pdf}}%
    \put(0.22878638,1.75942453){\makebox(0,0)[t]{\lineheight{1.25}\smash{\begin{tabular}[t]{c}\textbf{Foreground}\\\textbf{Vertex }$\Vertex$\end{tabular}}}}%
    \put(0,0){\includegraphics[width=\unitlength,page=11]{fg_mesh_abrev.pdf}}%
    \put(0.02135231,1.61565836){\makebox(0,0)[lt]{\lineheight{1.25}\smash{\begin{tabular}[t]{l}+ Coordinates $\pos$\end{tabular}}}}%
    \put(0,0){\includegraphics[width=\unitlength,page=12]{fg_mesh_abrev.pdf}}%
    \put(0.02135231,1.54092527){\color[rgb]{0,0,0}\makebox(0,0)[lt]{\lineheight{1.25}\smash{\begin{tabular}[t]{l}+ Bg. ancestor\end{tabular}}}}%
    \put(0.02403615,1.4911032){\color[rgb]{0,0,0}\makebox(0,0)[lt]{\lineheight{1.25}\smash{\begin{tabular}[t]{l}    rank $r$ \end{tabular}}}}%
    \put(0,0){\includegraphics[width=\unitlength,page=13]{fg_mesh_abrev.pdf}}%
    \put(0.02135231,1.42348754){\color[rgb]{0,0,0}\makebox(0,0)[lt]{\lineheight{1.25}\smash{\begin{tabular}[t]{l}+ Bg. ancestor\end{tabular}}}}%
    \put(0.02403615,1.37366548){\color[rgb]{0,0,0}\makebox(0,0)[lt]{\lineheight{1.25}\smash{\begin{tabular}[t]{l}    index $a$\end{tabular}}}}%
    \put(0,0){\includegraphics[width=\unitlength,page=14]{fg_mesh_abrev.pdf}}%
    \put(0.74733096,1.93238434){\makebox(0,0)[t]{\lineheight{1.25}\smash{\begin{tabular}[t]{c}\textbf{Foreground Cell }$\Cell$\end{tabular}}}}%
    \put(0,0){\includegraphics[width=\unitlength,page=15]{fg_mesh_abrev.pdf}}%
    \put(0.51957295,1.83985765){\color[rgb]{0,0,0}\makebox(0,0)[lt]{\lineheight{1.25}\smash{\begin{tabular}[t]{l}+ Array of vertices $\mathrm{V}$\end{tabular}}}}%
    \put(0.57651246,1.79003559){\color[rgb]{0,0,0}\makebox(0,0)[lt]{\lineheight{1.25}\smash{\begin{tabular}[t]{l}(sorted by ordinal)\end{tabular}}}}%
    \put(0,0){\includegraphics[width=\unitlength,page=16]{fg_mesh_abrev.pdf}}%
    \put(0.51957295,1.71886121){\makebox(0,0)[lt]{\lineheight{1.25}\smash{\begin{tabular}[t]{l}+ Material $m$\end{tabular}}}}%
    \put(0,0){\includegraphics[width=\unitlength,page=17]{fg_mesh_abrev.pdf}}%
    \put(0.51957295,1.63701068){\color[rgb]{0,0,0}\makebox(0,0)[lt]{\lineheight{1.25}\smash{\begin{tabular}[t]{l}+ Bg. element \end{tabular}}}}%
    \put(0.57826222,1.58718861){\color[rgb]{0,0,0}\makebox(0,0)[lt]{\lineheight{1.25}\smash{\begin{tabular}[t]{l}membership $E$\end{tabular}}}}%
    \put(0,0){\includegraphics[width=\unitlength,page=18]{fg_mesh_abrev.pdf}}%
    \put(0.51957295,1.50533808){\color[rgb]{0,0,0}\makebox(0,0)[lt]{\lineheight{1.25}\smash{\begin{tabular}[t]{l}+ Subphase \end{tabular}}}}%
    \put(0.57714502,1.45551601){\color[rgb]{0,0,0}\makebox(0,0)[lt]{\lineheight{1.25}\smash{\begin{tabular}[t]{l}membership $S$\end{tabular}}}}%
    \put(0,0){\includegraphics[width=\unitlength,page=19]{fg_mesh_abrev.pdf}}%
    \put(0.51957295,1.37366548){\makebox(0,0)[lt]{\lineheight{1.25}\smash{\begin{tabular}[t]{l}+ Parallel ID $I$\end{tabular}}}}%
    \put(0,0){\includegraphics[width=\unitlength,page=20]{fg_mesh_abrev.pdf}}%
    \put(0.49822064,0.3202847){\makebox(0,0)[t]{\lineheight{1.25}\smash{\begin{tabular}[t]{c}\textbf{Child Mesh }$\ChildMesh$\end{tabular}}}}%
    \put(0,0){\includegraphics[width=\unitlength,page=21]{fg_mesh_abrev.pdf}}%
    \put(0.12811388,0.22775801){\makebox(0,0)[lt]{\lineheight{1.25}\smash{\begin{tabular}[t]{l}+ List of fg. cells $\mathrm{C}$\end{tabular}}}}%
    \put(0,0){\includegraphics[width=\unitlength,page=22]{fg_mesh_abrev.pdf}}%
    \put(0.12811388,0.13523132){\makebox(0,0)[lt]{\lineheight{1.25}\smash{\begin{tabular}[t]{l}+ List of fg. vertices $\mathrm{V}$\end{tabular}}}}%
    \put(0,0){\includegraphics[width=\unitlength,page=23]{fg_mesh_abrev.pdf}}%
    \put(0.12811388,0.04982206){\makebox(0,0)[lt]{\lineheight{1.25}\smash{\begin{tabular}[t]{l}+ Vertex param. coords. $\SetDef{\xi_{\Vertex}}_{\Vertex=1}^{\mathtt{size}(\mathrm{V})}$\end{tabular}}}}%
  \end{picture}%
\endgroup%

%% file: figures/reg_sub_template.pdf_tex
\begingroup%
  \makeatletter%
  \providecommand\color[2][]{%
    \errmessage{(Inkscape) Color is used for the text in Inkscape, but the package 'color.sty' is not loaded}%
    \renewcommand\color[2][]{}%
  }%
  \providecommand\transparent[1]{%
    \errmessage{(Inkscape) Transparency is used (non-zero) for the text in Inkscape, but the package 'transparent.sty' is not loaded}%
    \renewcommand\transparent[1]{}%
  }%
  \providecommand\rotatebox[2]{#2}%
  \newcommand*\fsize{\dimexpr\f@size pt\relax}%
  \newcommand*\lineheight[1]{\fontsize{\fsize}{#1\fsize}\selectfont}%
  \ifx\svgwidth\undefined%
    \setlength{\unitlength}{270.60377226bp}%
    \ifx\svgscale\undefined%
      \relax%
    \else%
      \setlength{\unitlength}{\unitlength * \real{\svgscale}}%
    \fi%
  \else%
    \setlength{\unitlength}{\svgwidth}%
  \fi%
  \global\let\svgwidth\undefined%
  \global\let\svgscale\undefined%
  \makeatother%
  \begin{picture}(1,1.95086939)%
    \lineheight{1}%
    \setlength\tabcolsep{0pt}%
    \put(0,0){\includegraphics[width=\unitlength,page=1]{reg_sub_template.pdf}}%
    \put(0.17699172,0.81462224){\color[rgb]{0.82352941,0,0.00392157}\makebox(0,0)[lt]{\lineheight{1.25}\smash{\begin{tabular}[t]{l}\small$(2,\!3)$\end{tabular}}}}%
    \put(0,0){\includegraphics[width=\unitlength,page=2]{reg_sub_template.pdf}}%
    \put(0.09182874,1.90795368){\color[rgb]{0,0,0}\makebox(0,0)[lt]{\lineheight{1.25}\smash{\begin{tabular}[t]{l}(A)\end{tabular}}}}%
    \put(0.09008489,1.34128904){\color[rgb]{0,0,0}\makebox(0,0)[lt]{\lineheight{1.25}\smash{\begin{tabular}[t]{l}(B)\end{tabular}}}}%
    \put(0,0){\includegraphics[width=\unitlength,page=3]{reg_sub_template.pdf}}%
    \put(0.50380279,1.71877669){\color[rgb]{0.82352941,0,0.00392157}\makebox(0,0)[lt]{\lineheight{1.25}\smash{\begin{tabular}[t]{l}\small$(2,\!0)$\end{tabular}}}}%
    \put(0.49150023,0.69334185){\color[rgb]{0.82352941,0,0.00392157}\makebox(0,0)[lt]{\lineheight{1.25}\smash{\begin{tabular}[t]{l}\small$(3,\!0)$\end{tabular}}}}%
    \put(0.87667202,0.404073){\color[rgb]{0.82352941,0,0.00392157}\makebox(0,0)[lt]{\lineheight{1.25}\smash{\begin{tabular}[t]{l}\small$(2,\!1)$\end{tabular}}}}%
    \put(0.37192662,0.18448729){\color[rgb]{0.82352941,0,0.00392157}\makebox(0,0)[lt]{\lineheight{1.25}\smash{\begin{tabular}[t]{l}\small$(2,\!5)$\end{tabular}}}}%
    \put(0.49657227,1.20130787){\color[rgb]{0.82352941,0,0.00392157}\makebox(0,0)[lt]{\lineheight{1.25}\smash{\begin{tabular}[t]{l}\small$(2,\!6)$\end{tabular}}}}%
    \put(0.02436035,0.27882667){\color[rgb]{0.82352941,0,0.00392157}\makebox(0,0)[lt]{\lineheight{1.25}\smash{\begin{tabular}[t]{l}\small$(2,\!4)$\end{tabular}}}}%
    \put(0.71427859,0.92634213){\color[rgb]{0.82352941,0,0.00392157}\makebox(0,0)[lt]{\lineheight{1.25}\smash{\begin{tabular}[t]{l}$(2,\!2)$\end{tabular}}}}%
    \put(0,0){\includegraphics[width=\unitlength,page=4]{reg_sub_template.pdf}}%
    \put(0.32180852,1.47259027){\color[rgb]{0.17254902,0.17254902,0.17254902}\makebox(0,0)[lt]{\lineheight{1.25}\smash{\begin{tabular}[t]{l}\small$x$\end{tabular}}}}%
    \put(0.2186252,1.58373489){\color[rgb]{0.17254902,0.17254902,0.17254902}\makebox(0,0)[lt]{\lineheight{1.25}\smash{\begin{tabular}[t]{l}\small$y$\end{tabular}}}}%
    \put(0,0){\includegraphics[width=\unitlength,page=5]{reg_sub_template.pdf}}%
    \put(0.94278114,0.0876007){\color[rgb]{0.17254902,0.17254902,0.17254902}\makebox(0,0)[lt]{\lineheight{1.25}\smash{\begin{tabular}[t]{l}\small$x$\end{tabular}}}}%
    \put(0.87127327,0.11528911){\color[rgb]{0.17254902,0.17254902,0.17254902}\makebox(0,0)[lt]{\lineheight{1.25}\smash{\begin{tabular}[t]{l}\small$z$\end{tabular}}}}%
    \put(0.78922365,0.06492166){\color[rgb]{0.17254902,0.17254902,0.17254902}\makebox(0,0)[lt]{\lineheight{1.25}\smash{\begin{tabular}[t]{l}\small$y$\end{tabular}}}}%
  \end{picture}%
\endgroup%

%% file: figures/subdivision_templates.pdf_tex
\begingroup%
  \makeatletter%
  \providecommand\color[2][]{%
    \errmessage{(Inkscape) Color is used for the text in Inkscape, but the package 'color.sty' is not loaded}%
    \renewcommand\color[2][]{}%
  }%
  \providecommand\transparent[1]{%
    \errmessage{(Inkscape) Transparency is used (non-zero) for the text in Inkscape, but the package 'transparent.sty' is not loaded}%
    \renewcommand\transparent[1]{}%
  }%
  \providecommand\rotatebox[2]{#2}%
  \newcommand*\fsize{\dimexpr\f@size pt\relax}%
  \newcommand*\lineheight[1]{\fontsize{\fsize}{#1\fsize}\selectfont}%
  \ifx\svgwidth\undefined%
    \setlength{\unitlength}{243.64492341bp}%
    \ifx\svgscale\undefined%
      \relax%
    \else%
      \setlength{\unitlength}{\unitlength * \real{\svgscale}}%
    \fi%
  \else%
    \setlength{\unitlength}{\svgwidth}%
  \fi%
  \global\let\svgwidth\undefined%
  \global\let\svgscale\undefined%
  \makeatother%
  \begin{picture}(1,1.69123346)%
    \lineheight{1}%
    \setlength\tabcolsep{0pt}%
    \put(0,0){\includegraphics[width=\unitlength,page=1]{subdivision_templates.pdf}}%
    \put(-0.00071969,1.64513782){\color[rgb]{0,0,0.00392157}\makebox(0,0)[lt]{\lineheight{1.25}\smash{\begin{tabular}[t]{l}(A)\end{tabular}}}}%
    \put(0.00069182,1.10143989){\color[rgb]{0,0,0.00392157}\makebox(0,0)[lt]{\lineheight{1.25}\smash{\begin{tabular}[t]{l}(B)\end{tabular}}}}%
    \put(0.12606835,1.27681186){\color[rgb]{0,0,0.00392157}\makebox(0,0)[lt]{\lineheight{1.25}\smash{\begin{tabular}[t]{l}2 edges\end{tabular}}}}%
    \put(0.02527826,0.63703701){\color[rgb]{0,0,0.00392157}\makebox(0,0)[lt]{\lineheight{1.25}\smash{\begin{tabular}[t]{l}4 edges\end{tabular}}}}%
    \put(0.56641631,0.63358183){\color[rgb]{0,0,0.00392157}\makebox(0,0)[lt]{\lineheight{1.25}\smash{\begin{tabular}[t]{l}3 edges\end{tabular}}}}%
    \put(0.52192426,1.27744181){\color[rgb]{0,0,0.00392157}\makebox(0,0)[lt]{\lineheight{1.25}\smash{\begin{tabular}[t]{l}1 edge + 1 vertex\end{tabular}}}}%
    \put(0.00793061,0.06875366){\color[rgb]{0,0,0.00392157}\makebox(0,0)[lt]{\lineheight{1.25}\smash{\begin{tabular}[t]{l}2 edges\end{tabular}}}}%
    \put(0.01077582,0.00063488){\color[rgb]{0,0,0.00392157}\makebox(0,0)[lt]{\lineheight{1.25}\smash{\begin{tabular}[t]{l}+ 1 vertex\end{tabular}}}}%
    \put(0.5311991,0.06875366){\color[rgb]{0,0,0.00392157}\makebox(0,0)[lt]{\lineheight{1.25}\smash{\begin{tabular}[t]{l}1 edge\end{tabular}}}}%
    \put(0.53432071,0.00063488){\color[rgb]{0,0,0.00392157}\makebox(0,0)[lt]{\lineheight{1.25}\smash{\begin{tabular}[t]{l}+ 2 vertices\end{tabular}}}}%
  \end{picture}%
\endgroup%

%% file: figures/proximity_during_subdivision.pdf_tex
\begingroup%
  \makeatletter%
  \providecommand\color[2][]{%
    \errmessage{(Inkscape) Color is used for the text in Inkscape, but the package 'color.sty' is not loaded}%
    \renewcommand\color[2][]{}%
  }%
  \providecommand\transparent[1]{%
    \errmessage{(Inkscape) Transparency is used (non-zero) for the text in Inkscape, but the package 'transparent.sty' is not loaded}%
    \renewcommand\transparent[1]{}%
  }%
  \providecommand\rotatebox[2]{#2}%
  \newcommand*\fsize{\dimexpr\f@size pt\relax}%
  \newcommand*\lineheight[1]{\fontsize{\fsize}{#1\fsize}\selectfont}%
  \ifx\svgwidth\undefined%
    \setlength{\unitlength}{219.93240477bp}%
    \ifx\svgscale\undefined%
      \relax%
    \else%
      \setlength{\unitlength}{\unitlength * \real{\svgscale}}%
    \fi%
  \else%
    \setlength{\unitlength}{\svgwidth}%
  \fi%
  \global\let\svgwidth\undefined%
  \global\let\svgscale\undefined%
  \makeatother%
  \begin{picture}(1,1.61983134)%
    \lineheight{1}%
    \setlength\tabcolsep{0pt}%
    \put(0.47559625,1.58134295){\color[rgb]{0,0.47843137,1}\makebox(0,0)[lt]{\lineheight{1.25}\smash{\begin{tabular}[t]{l}$+$\end{tabular}}}}%
    \put(0.10445732,1.106181){\color[rgb]{0,0.47843137,1}\makebox(0,0)[lt]{\lineheight{1.25}\smash{\begin{tabular}[t]{l}$+$\end{tabular}}}}%
    \put(0.81029136,1.05197666){\color[rgb]{0,0.47843137,1}\makebox(0,0)[lt]{\lineheight{1.25}\smash{\begin{tabular}[t]{l}$-$\end{tabular}}}}%
    \put(0,0){\includegraphics[width=\unitlength,page=1]{proximity_during_subdivision.pdf}}%
    \put(0.63335515,1.34048158){\color[rgb]{0.82352941,0,0.00392157}\makebox(0,0)[lt]{\lineheight{1.25}\smash{\begin{tabular}[t]{l}$0$\end{tabular}}}}%
    \put(0.23264235,0.9826896){\color[rgb]{0.82352941,0,0.00392157}\makebox(0,0)[lt]{\lineheight{1.25}\smash{\begin{tabular}[t]{l}$m= \langle 1,0 \rangle =2$\end{tabular}}}}%
    \put(0.53271223,0.88799131){\color[rgb]{0,0,0.00392157}\makebox(0,0)[lt]{\lineheight{1.25}\smash{\begin{tabular}[t]{l}apply subdivision \\template\end{tabular}}}}%
    \put(0.05846182,0.01033044){\color[rgb]{0.82352941,0,0.00392157}\makebox(0,0)[lt]{\lineheight{1.25}\smash{\begin{tabular}[t]{l}$\langle 1,0,1 \rangle$\end{tabular}}}}%
    \put(0.06310826,0.08882659){\color[rgb]{0.82352941,0,0.00392157}\makebox(0,0)[lt]{\lineheight{1.25}\smash{\begin{tabular}[t]{l}$2m+1 =5$\end{tabular}}}}%
    \put(0.58624967,0.02396806){\color[rgb]{0.82352941,0,0.00392157}\makebox(0,0)[lt]{\lineheight{1.25}\smash{\begin{tabular}[t]{l}$\langle 1,0,0 \rangle$\end{tabular}}}}%
    \put(0.58563017,0.09962317){\color[rgb]{0.82352941,0,0.00392157}\makebox(0,0)[lt]{\lineheight{1.25}\smash{\begin{tabular}[t]{l}$2m+0 =4$\end{tabular}}}}%
    \put(0.36766389,1.10022005){\color[rgb]{0.82352941,0,0.00392157}\makebox(0,0)[lt]{\lineheight{1.25}\smash{\begin{tabular}[t]{l}$0$\end{tabular}}}}%
    \put(0.5319116,1.15137114){\color[rgb]{0.82352941,0,0.00392157}\makebox(0,0)[lt]{\lineheight{1.25}\smash{\begin{tabular}[t]{l}$\Gamma$\end{tabular}}}}%
    \put(0.28267029,0.6949973){\color[rgb]{0,0.47843137,1}\makebox(0,0)[lt]{\lineheight{1.25}\smash{\begin{tabular}[t]{l}$+$\end{tabular}}}}%
    \put(-0.00378537,0.22316334){\color[rgb]{0,0.47843137,1}\makebox(0,0)[lt]{\lineheight{1.25}\smash{\begin{tabular}[t]{l}$+$\end{tabular}}}}%
    \put(0,0){\includegraphics[width=\unitlength,page=2]{proximity_during_subdivision.pdf}}%
    \put(0.23508494,0.2030971){\color[rgb]{0,0.47843137,1}\makebox(0,0)[lt]{\lineheight{1.25}\smash{\begin{tabular}[t]{l}$0$\end{tabular}}}}%
    \put(0.45809217,0.73055376){\color[rgb]{0,0.47843137,1}\makebox(0,0)[lt]{\lineheight{1.25}\smash{\begin{tabular}[t]{l}$+$\end{tabular}}}}%
    \put(0,0){\includegraphics[width=\unitlength,page=3]{proximity_during_subdivision.pdf}}%
    \put(0.61058475,0.48369919){\color[rgb]{0,0.47843137,1}\makebox(0,0)[lt]{\lineheight{1.25}\smash{\begin{tabular}[t]{l}$0$\end{tabular}}}}%
    \put(0.37770252,0.15722786){\color[rgb]{0,0.47843137,1}\makebox(0,0)[lt]{\lineheight{1.25}\smash{\begin{tabular}[t]{l}$0$\end{tabular}}}}%
    \put(0.91716191,0.19893957){\color[rgb]{0,0.47843137,1}\makebox(0,0)[lt]{\lineheight{1.25}\smash{\begin{tabular}[t]{l}$-$\end{tabular}}}}%
    \put(0,0){\includegraphics[width=\unitlength,page=4]{proximity_during_subdivision.pdf}}%
    \put(0.77641444,0.43531302){\color[rgb]{0,0.47843137,1}\makebox(0,0)[lt]{\lineheight{1.25}\smash{\begin{tabular}[t]{l}$0$\end{tabular}}}}%
    \put(0.65496601,0.22238572){\color[rgb]{0,0.47843137,1}\makebox(0,0)[lt]{\lineheight{1.25}\smash{\begin{tabular}[t]{l}$P\!=\!-$\end{tabular}}}}%
    \put(0.4152276,0.44911505){\color[rgb]{0,0.47843137,1}\makebox(0,0)[lt]{\lineheight{1.25}\smash{\begin{tabular}[t]{l}$P\!=\!+$\end{tabular}}}}%
    \put(0.12294079,0.29288464){\color[rgb]{0,0.47843137,1}\makebox(0,0)[lt]{\lineheight{1.25}\smash{\begin{tabular}[t]{l}$P\!=\!+$\end{tabular}}}}%
    \put(0.50016044,0.20385072){\color[rgb]{0,0.47843137,1}\makebox(0,0)[lt]{\lineheight{1.25}\smash{\begin{tabular}[t]{l}$0$\end{tabular}}}}%
    \put(0,0){\includegraphics[width=\unitlength,page=5]{proximity_during_subdivision.pdf}}%
  \end{picture}%
\endgroup%

%% file: figures/subphase.pdf_tex
\begingroup%
  \makeatletter%
  \providecommand\color[2][]{%
    \errmessage{(Inkscape) Color is used for the text in Inkscape, but the package 'color.sty' is not loaded}%
    \renewcommand\color[2][]{}%
  }%
  \providecommand\transparent[1]{%
    \errmessage{(Inkscape) Transparency is used (non-zero) for the text in Inkscape, but the package 'transparent.sty' is not loaded}%
    \renewcommand\transparent[1]{}%
  }%
  \providecommand\rotatebox[2]{#2}%
  \newcommand*\fsize{\dimexpr\f@size pt\relax}%
  \newcommand*\lineheight[1]{\fontsize{\fsize}{#1\fsize}\selectfont}%
  \ifx\svgwidth\undefined%
    \setlength{\unitlength}{144bp}%
    \ifx\svgscale\undefined%
      \relax%
    \else%
      \setlength{\unitlength}{\unitlength * \real{\svgscale}}%
    \fi%
  \else%
    \setlength{\unitlength}{\svgwidth}%
  \fi%
  \global\let\svgwidth\undefined%
  \global\let\svgscale\undefined%
  \makeatother%
  \begin{picture}(1,0.63541667)%
    \lineheight{1}%
    \setlength\tabcolsep{0pt}%
    \put(0,0){\includegraphics[width=\unitlength,page=1]{subphase.pdf}}%
    \put(0.49479167,0.54166667){\makebox(0,0)[t]{\lineheight{1.25}\smash{\begin{tabular}[t]{c}\textbf{Subphase }$\Subphase$\end{tabular}}}}%
    \put(0,0){\includegraphics[width=\unitlength,page=2]{subphase.pdf}}%
    \put(0.03125,0.40625){\makebox(0,0)[lt]{\lineheight{1.25}\smash{\begin{tabular}[t]{l}+ Bg. element membership $\BgElemInd$\end{tabular}}}}%
    \put(0,0){\includegraphics[width=\unitlength,page=3]{subphase.pdf}}%
    \put(0.03125,0.29166667){\makebox(0,0)[lt]{\lineheight{1.25}\smash{\begin{tabular}[t]{l}+ Ordinal $\Unzipping$\end{tabular}}}}%
    \put(0,0){\includegraphics[width=\unitlength,page=4]{subphase.pdf}}%
    \put(0.03125,0.17708333){\makebox(0,0)[lt]{\lineheight{1.25}\smash{\begin{tabular}[t]{l}+ Array of fg. cells $\mathrm{C}$\end{tabular}}}}%
    \put(0,0){\includegraphics[width=\unitlength,page=5]{subphase.pdf}}%
    \put(0.03125,0.05729167){\makebox(0,0)[lt]{\lineheight{1.25}\smash{\begin{tabular}[t]{l}{\color{gray}+ Parallel ID $I$}\end{tabular}}}}%
  \end{picture}%
\endgroup%

%% file: figures/bg_facet_descendants.pdf_tex
\begingroup%
  \makeatletter%
  \providecommand\color[2][]{%
    \errmessage{(Inkscape) Color is used for the text in Inkscape, but the package 'color.sty' is not loaded}%
    \renewcommand\color[2][]{}%
  }%
  \providecommand\transparent[1]{%
    \errmessage{(Inkscape) Transparency is used (non-zero) for the text in Inkscape, but the package 'transparent.sty' is not loaded}%
    \renewcommand\transparent[1]{}%
  }%
  \providecommand\rotatebox[2]{#2}%
  \newcommand*\fsize{\dimexpr\f@size pt\relax}%
  \newcommand*\lineheight[1]{\fontsize{\fsize}{#1\fsize}\selectfont}%
  \ifx\svgwidth\undefined%
    \setlength{\unitlength}{401.00322573bp}%
    \ifx\svgscale\undefined%
      \relax%
    \else%
      \setlength{\unitlength}{\unitlength * \real{\svgscale}}%
    \fi%
  \else%
    \setlength{\unitlength}{\svgwidth}%
  \fi%
  \global\let\svgwidth\undefined%
  \global\let\svgscale\undefined%
  \makeatother%
  \begin{picture}(1,0.32138699)%
    \lineheight{1}%
    \setlength\tabcolsep{0pt}%
    \put(0,0){\includegraphics[width=\unitlength,page=1]{./figures/bg_facet_descendants.pdf}}%
    \put(0.91018779,0.27592794){\color[rgb]{0.82352941,0,0.00392157}\makebox(0,0)[lt]{\lineheight{1.25}\smash{\begin{tabular}[t]{l}$\Gamma$\end{tabular}}}}%
  \end{picture}%
\endgroup%

%% file: figures/subphase_generation.pdf_tex
\begingroup%
  \makeatletter%
  \providecommand\color[2][]{%
    \errmessage{(Inkscape) Color is used for the text in Inkscape, but the package 'color.sty' is not loaded}%
    \renewcommand\color[2][]{}%
  }%
  \providecommand\transparent[1]{%
    \errmessage{(Inkscape) Transparency is used (non-zero) for the text in Inkscape, but the package 'transparent.sty' is not loaded}%
    \renewcommand\transparent[1]{}%
  }%
  \providecommand\rotatebox[2]{#2}%
  \newcommand*\fsize{\dimexpr\f@size pt\relax}%
  \newcommand*\lineheight[1]{\fontsize{\fsize}{#1\fsize}\selectfont}%
  \ifx\svgwidth\undefined%
    \setlength{\unitlength}{191.23728354bp}%
    \ifx\svgscale\undefined%
      \relax%
    \else%
      \setlength{\unitlength}{\unitlength * \real{\svgscale}}%
    \fi%
  \else%
    \setlength{\unitlength}{\svgwidth}%
  \fi%
  \global\let\svgwidth\undefined%
  \global\let\svgscale\undefined%
  \makeatother%
  \begin{picture}(1,2.19541958)%
    \lineheight{1}%
    \setlength\tabcolsep{0pt}%
    \put(0,0){\includegraphics[width=\unitlength,page=1]{subphase_generation.pdf}}%
    \put(0.38923801,0.26795591){\color[rgb]{0,0,0}\makebox(0,0)[lt]{\lineheight{1.25}\smash{\begin{tabular}[t]{l}$7$\end{tabular}}}}%
    \put(0,0){\includegraphics[width=\unitlength,page=2]{subphase_generation.pdf}}%
    \put(0.1953419,0.07371417){\color[rgb]{0,0,0}\makebox(0,0)[lt]{\lineheight{1.25}\smash{\begin{tabular}[t]{l}$1$\end{tabular}}}}%
    \put(0,0){\includegraphics[width=\unitlength,page=3]{subphase_generation.pdf}}%
    \put(0.19670615,0.26916935){\color[rgb]{0,0,0}\makebox(0,0)[lt]{\lineheight{1.25}\smash{\begin{tabular}[t]{l}$4$\end{tabular}}}}%
    \put(0.19605644,0.41711947){\color[rgb]{0,0,0}\makebox(0,0)[lt]{\lineheight{1.25}\smash{\begin{tabular}[t]{l}$5$\end{tabular}}}}%
    \put(0.38815338,0.41704189){\color[rgb]{0,0,0}\makebox(0,0)[lt]{\lineheight{1.25}\smash{\begin{tabular}[t]{l}$8$\end{tabular}}}}%
    \put(0.53529464,0.23709721){\color[rgb]{0,0,0}\makebox(0,0)[lt]{\lineheight{1.25}\smash{\begin{tabular}[t]{l}$6$\end{tabular}}}}%
    \put(0.4632241,0.07571193){\color[rgb]{0,0,0}\makebox(0,0)[lt]{\lineheight{1.25}\smash{\begin{tabular}[t]{l}$2$\end{tabular}}}}%
    \put(0.78938351,0.07567783){\color[rgb]{0,0,0}\makebox(0,0)[lt]{\lineheight{1.25}\smash{\begin{tabular}[t]{l}$3$\end{tabular}}}}%
    \put(0.70168121,0.41775672){\color[rgb]{0,0,0}\makebox(0,0)[lt]{\lineheight{1.25}\smash{\begin{tabular}[t]{l}$11$\end{tabular}}}}%
    \put(0.8460417,0.34170384){\color[rgb]{0,0,0}\makebox(0,0)[lt]{\lineheight{1.25}\smash{\begin{tabular}[t]{l}$12$\end{tabular}}}}%
    \put(0.70127653,0.26581363){\color[rgb]{0,0,0}\makebox(0,0)[lt]{\lineheight{1.25}\smash{\begin{tabular}[t]{l}$10$\end{tabular}}}}%
    \put(0.53816719,0.44627145){\color[rgb]{0,0,0}\makebox(0,0)[lt]{\lineheight{1.25}\smash{\begin{tabular}[t]{l}$9$\end{tabular}}}}%
    \put(3.45007457,-0.09791914){\color[rgb]{0,0,0}\makebox(0,0)[lt]{\lineheight{1.25}\smash{\begin{tabular}[t]{l}$7$\end{tabular}}}}%
    \put(0,0){\includegraphics[width=\unitlength,page=4]{subphase_generation.pdf}}%
    \put(3.44898983,0.05116687){\color[rgb]{0,0,0}\makebox(0,0)[lt]{\lineheight{1.25}\smash{\begin{tabular}[t]{l}$8$\end{tabular}}}}%
    \put(3.59613126,-0.12877784){\color[rgb]{0,0,0}\makebox(0,0)[lt]{\lineheight{1.25}\smash{\begin{tabular}[t]{l}$6$\end{tabular}}}}%
    \put(3.52406054,-0.29016309){\color[rgb]{0,0,0}\makebox(0,0)[lt]{\lineheight{1.25}\smash{\begin{tabular}[t]{l}$2$\end{tabular}}}}%
    \put(3.85021996,-0.29019719){\color[rgb]{0,0,0}\makebox(0,0)[lt]{\lineheight{1.25}\smash{\begin{tabular}[t]{l}$3$\end{tabular}}}}%
    \put(3.76251777,0.0518817){\color[rgb]{0,0,0}\makebox(0,0)[lt]{\lineheight{1.25}\smash{\begin{tabular}[t]{l}$11$\end{tabular}}}}%
    \put(3.90687814,-0.02417118){\color[rgb]{0,0,0}\makebox(0,0)[lt]{\lineheight{1.25}\smash{\begin{tabular}[t]{l}$12$\end{tabular}}}}%
    \put(3.76211314,-0.10006142){\color[rgb]{0,0,0}\makebox(0,0)[lt]{\lineheight{1.25}\smash{\begin{tabular}[t]{l}$10$\end{tabular}}}}%
    \put(3.59900369,0.08039643){\color[rgb]{0,0,0}\makebox(0,0)[lt]{\lineheight{1.25}\smash{\begin{tabular}[t]{l}$9$\end{tabular}}}}%
    \put(0.38923795,0.96462098){\color[rgb]{0,0,0}\makebox(0,0)[lt]{\lineheight{1.25}\smash{\begin{tabular}[t]{l}$7$\end{tabular}}}}%
    \put(0,0){\includegraphics[width=\unitlength,page=5]{subphase_generation.pdf}}%
    \put(0.19534183,0.77037933){\color[rgb]{0,0,0}\makebox(0,0)[lt]{\lineheight{1.25}\smash{\begin{tabular}[t]{l}$1$\end{tabular}}}}%
    \put(0,0){\includegraphics[width=\unitlength,page=6]{subphase_generation.pdf}}%
    \put(0.19670608,0.96583442){\color[rgb]{0,0,0}\makebox(0,0)[lt]{\lineheight{1.25}\smash{\begin{tabular}[t]{l}$4$\end{tabular}}}}%
    \put(0.19605637,1.11378465){\color[rgb]{0,0,0}\makebox(0,0)[lt]{\lineheight{1.25}\smash{\begin{tabular}[t]{l}$5$\end{tabular}}}}%
    \put(0.38815333,1.11370707){\color[rgb]{0,0,0}\makebox(0,0)[lt]{\lineheight{1.25}\smash{\begin{tabular}[t]{l}$8$\end{tabular}}}}%
    \put(0.53529459,0.93376234){\color[rgb]{0,0,0}\makebox(0,0)[lt]{\lineheight{1.25}\smash{\begin{tabular}[t]{l}$6$\end{tabular}}}}%
    \put(0.46322404,0.77237703){\color[rgb]{0,0,0}\makebox(0,0)[lt]{\lineheight{1.25}\smash{\begin{tabular}[t]{l}$2$\end{tabular}}}}%
    \put(0.78938351,0.77234299){\color[rgb]{0,0,0}\makebox(0,0)[lt]{\lineheight{1.25}\smash{\begin{tabular}[t]{l}$3$\end{tabular}}}}%
    \put(0.70168121,1.1144219){\color[rgb]{0,0,0}\makebox(0,0)[lt]{\lineheight{1.25}\smash{\begin{tabular}[t]{l}$11$\end{tabular}}}}%
    \put(0.8460417,1.03836891){\color[rgb]{0,0,0}\makebox(0,0)[lt]{\lineheight{1.25}\smash{\begin{tabular}[t]{l}$12$\end{tabular}}}}%
    \put(0.70127653,0.96247876){\color[rgb]{0,0,0}\makebox(0,0)[lt]{\lineheight{1.25}\smash{\begin{tabular}[t]{l}$10$\end{tabular}}}}%
    \put(0.53816713,1.14293657){\color[rgb]{0,0,0}\makebox(0,0)[lt]{\lineheight{1.25}\smash{\begin{tabular}[t]{l}$9$\end{tabular}}}}%
    \put(3.44270666,0.51994111){\color[rgb]{0,0,0}\makebox(0,0)[lt]{\lineheight{1.25}\smash{\begin{tabular}[t]{l}$7$\end{tabular}}}}%
    \put(0,0){\includegraphics[width=\unitlength,page=7]{subphase_generation.pdf}}%
    \put(3.24881038,0.32569946){\color[rgb]{0,0,0}\makebox(0,0)[lt]{\lineheight{1.25}\smash{\begin{tabular}[t]{l}$1$\end{tabular}}}}%
    \put(0,0){\includegraphics[width=\unitlength,page=8]{subphase_generation.pdf}}%
    \put(3.25017468,0.5211546){\color[rgb]{0,0,0}\makebox(0,0)[lt]{\lineheight{1.25}\smash{\begin{tabular}[t]{l}$4$\end{tabular}}}}%
    \put(3.2495251,0.66910473){\color[rgb]{0,0,0}\makebox(0,0)[lt]{\lineheight{1.25}\smash{\begin{tabular}[t]{l}$5$\end{tabular}}}}%
    \put(3.44162193,0.66902715){\color[rgb]{0,0,0}\makebox(0,0)[lt]{\lineheight{1.25}\smash{\begin{tabular}[t]{l}$8$\end{tabular}}}}%
    \put(3.58876336,0.48908241){\color[rgb]{0,0,0}\makebox(0,0)[lt]{\lineheight{1.25}\smash{\begin{tabular}[t]{l}$6$\end{tabular}}}}%
    \put(3.51669264,0.32769716){\color[rgb]{0,0,0}\makebox(0,0)[lt]{\lineheight{1.25}\smash{\begin{tabular}[t]{l}$2$\end{tabular}}}}%
    \put(3.59163579,0.69825671){\color[rgb]{0,0,0}\makebox(0,0)[lt]{\lineheight{1.25}\smash{\begin{tabular}[t]{l}$9$\end{tabular}}}}%
    \put(0.00256798,2.07827165){\color[rgb]{0,0,0}\makebox(0,0)[lt]{\lineheight{1.25}\smash{\begin{tabular}[t]{l}$n_u(4) \! = \! 2$\end{tabular}}}}%
    \put(0.01366277,1.63476823){\color[rgb]{0,0,0}\makebox(0,0)[lt]{\lineheight{1.25}\smash{\begin{tabular}[t]{l}$n_u(1) \! = \! 1$\end{tabular}}}}%
    \put(0.34360725,1.63082226){\color[rgb]{0,0,0}\makebox(0,0)[lt]{\lineheight{1.25}\smash{\begin{tabular}[t]{l}$n_u(2) \! = \! 1$\end{tabular}}}}%
    \put(0.71670053,1.63092765){\color[rgb]{0,0,0}\makebox(0,0)[lt]{\lineheight{1.25}\smash{\begin{tabular}[t]{l}$n_u(3) \! = \! 1$\end{tabular}}}}%
    \put(0.33084675,2.09185215){\color[rgb]{0,0,0}\makebox(0,0)[lt]{\lineheight{1.25}\smash{\begin{tabular}[t]{l}$n_u(5) \! = \! 4$\end{tabular}}}}%
    \put(0.70503235,2.04670586){\color[rgb]{0,0,0}\makebox(0,0)[lt]{\lineheight{1.25}\smash{\begin{tabular}[t]{l}$n_u(6) \! = \! 3$\end{tabular}}}}%
    \put(0.00876989,1.18813296){\color[rgb]{0,0,0}\makebox(0,0)[lt]{\lineheight{1.25}\smash{\begin{tabular}[t]{l}(B)\end{tabular}}}}%
    \put(0.00671394,0.49992619){\color[rgb]{0,0,0}\makebox(0,0)[lt]{\lineheight{1.25}\smash{\begin{tabular}[t]{l}(C)\end{tabular}}}}%
    \put(0.00104557,2.15966352){\color[rgb]{0,0,0}\makebox(0,0)[lt]{\lineheight{1.25}\smash{\begin{tabular}[t]{l}(A)\end{tabular}}}}%
    \put(0,0){\includegraphics[width=\unitlength,page=9]{subphase_generation.pdf}}%
    \put(0.10066938,0.6361111){\color[rgb]{0,0,0}\makebox(0,0)[lt]{\lineheight{1.25}\smash{\begin{tabular}[t]{l}Generate\\subphase graph\end{tabular}}}}%
    \put(0.33355902,0.11766289){\color[rgb]{0.82352941,0,0.00392157}\makebox(0,0)[lt]{\lineheight{1.25}\smash{\begin{tabular}[t]{l}$\mathcal{G}_I$\end{tabular}}}}%
    \put(0.80863283,1.13886134){\color[rgb]{0.82352941,0,0.00392157}\makebox(0,0)[lt]{\lineheight{1.25}\smash{\begin{tabular}[t]{l}$\mathcal{G}_I$\end{tabular}}}}%
    \put(0.86643507,0.21630469){\color[rgb]{0,0.49019608,0.97254902}\makebox(0,0)[lt]{\lineheight{1.25}\smash{\begin{tabular}[t]{l}$\mathcal{G}_S$\end{tabular}}}}%
  \end{picture}%
\endgroup%

%% file: figures/flood_fill.pdf_tex
\begingroup%
  \makeatletter%
  \providecommand\color[2][]{%
    \errmessage{(Inkscape) Color is used for the text in Inkscape, but the package 'color.sty' is not loaded}%
    \renewcommand\color[2][]{}%
  }%
  \providecommand\transparent[1]{%
    \errmessage{(Inkscape) Transparency is used (non-zero) for the text in Inkscape, but the package 'transparent.sty' is not loaded}%
    \renewcommand\transparent[1]{}%
  }%
  \providecommand\rotatebox[2]{#2}%
  \newcommand*\fsize{\dimexpr\f@size pt\relax}%
  \newcommand*\lineheight[1]{\fontsize{\fsize}{#1\fsize}\selectfont}%
  \ifx\svgwidth\undefined%
    \setlength{\unitlength}{184.77864183bp}%
    \ifx\svgscale\undefined%
      \relax%
    \else%
      \setlength{\unitlength}{\unitlength * \real{\svgscale}}%
    \fi%
  \else%
    \setlength{\unitlength}{\svgwidth}%
  \fi%
  \global\let\svgwidth\undefined%
  \global\let\svgscale\undefined%
  \makeatother%
  \begin{picture}(1,1.67011082)%
    \lineheight{1}%
    \setlength\tabcolsep{0pt}%
    \put(0,0){\includegraphics[width=\unitlength,page=1]{flood_fill.pdf}}%
    \put(0.39990894,1.26561164){\color[rgb]{0,0,0}\makebox(0,0)[lt]{\lineheight{1.25}\smash{\begin{tabular}[t]{l}$7$\end{tabular}}}}%
    \put(0,0){\includegraphics[width=\unitlength,page=2]{flood_fill.pdf}}%
    \put(0.19923533,1.06458035){\color[rgb]{0,0,0}\makebox(0,0)[lt]{\lineheight{1.25}\smash{\begin{tabular}[t]{l}$1$\end{tabular}}}}%
    \put(0,0){\includegraphics[width=\unitlength,page=3]{flood_fill.pdf}}%
    \put(0.20064731,1.26686726){\color[rgb]{0,0,0}\makebox(0,0)[lt]{\lineheight{1.25}\smash{\begin{tabular}[t]{l}$4$\end{tabular}}}}%
    \put(0.19997503,1.41998886){\color[rgb]{0,0,0}\makebox(0,0)[lt]{\lineheight{1.25}\smash{\begin{tabular}[t]{l}$5$\end{tabular}}}}%
    \put(0.41000558,0.93694695){\color[rgb]{0.82352941,0,0}\makebox(0,0)[lt]{\lineheight{1.25}\smash{\begin{tabular}[t]{l}$\mathrm{supp}(N_{\BfIndex})$\end{tabular}}}}%
    \put(0.13763351,0.51951785){\color[rgb]{0.82352941,0,0}\makebox(0,0)[lt]{\lineheight{1.25}\smash{\begin{tabular}[t]{l}$\mathrm{supp} \, \widetilde{N}_{\BfIndex}^{\EnrLvl=1}$\end{tabular}}}}%
    \put(0.39878629,1.41990857){\color[rgb]{0,0,0}\makebox(0,0)[lt]{\lineheight{1.25}\smash{\begin{tabular}[t]{l}$8$\end{tabular}}}}%
    \put(0.55107058,1.23367427){\color[rgb]{0,0,0}\makebox(0,0)[lt]{\lineheight{1.25}\smash{\begin{tabular}[t]{l}$6$\end{tabular}}}}%
    \put(0.47648098,1.06664799){\color[rgb]{0,0,0}\makebox(0,0)[lt]{\lineheight{1.25}\smash{\begin{tabular}[t]{l}$2$\end{tabular}}}}%
    \put(0.81404101,1.06661265){\color[rgb]{0,0,0}\makebox(0,0)[lt]{\lineheight{1.25}\smash{\begin{tabular}[t]{l}$3$\end{tabular}}}}%
    \put(0.7232731,1.42064826){\color[rgb]{0,0,0}\makebox(0,0)[lt]{\lineheight{1.25}\smash{\begin{tabular}[t]{l}$11$\end{tabular}}}}%
    \put(0.87267936,1.34193719){\color[rgb]{0,0,0}\makebox(0,0)[lt]{\lineheight{1.25}\smash{\begin{tabular}[t]{l}$12$\end{tabular}}}}%
    \put(0.72285409,1.26339442){\color[rgb]{0,0,0}\makebox(0,0)[lt]{\lineheight{1.25}\smash{\begin{tabular}[t]{l}$10$\end{tabular}}}}%
    \put(0.55404388,1.45015985){\color[rgb]{0,0,0}\makebox(0,0)[lt]{\lineheight{1.25}\smash{\begin{tabular}[t]{l}$9$\end{tabular}}}}%
    \put(0.25902349,1.19530293){\color[rgb]{0,0.49019608,0.97254902}\makebox(0,0)[lt]{\lineheight{1.25}\smash{\begin{tabular}[t]{l}$\mathcal{G}_{p}$\end{tabular}}}}%
    \put(0.21596396,1.61759741){\color[rgb]{0,0,0}\makebox(0,0)[lt]{\lineheight{1.25}\smash{\begin{tabular}[t]{l}$\mathtt{prune}(\mathcal{G}_S, \mathrm{supp}(N_{\BfIndex}))$\end{tabular}}}}%
    \put(0,0){\includegraphics[width=\unitlength,page=4]{flood_fill.pdf}}%
    \put(0.22250771,0.841292){\color[rgb]{0,0,0}\makebox(0,0)[lt]{\lineheight{1.25}\smash{\begin{tabular}[t]{l}$\mathtt{flood\_fill}(\mathcal{G}_P)$\end{tabular}}}}%
    \put(0.14436218,0.71521407){\color[rgb]{0,0,0}\makebox(0,0)[lt]{\lineheight{1.25}\smash{\begin{tabular}[t]{l}$R_{\BfIndex} = \{\,\{1,4,7\},\,\{2,6\},\,\{5,8\},\,\{9\}\,\}$\end{tabular}}}}%
    \put(0.35199223,0.64657295){\color[rgb]{0,0,0}\makebox(0,0)[lt]{\lineheight{1.25}\smash{\begin{tabular}[t]{l}${\scriptsize ~ \, \EnrLvl\!=\!1 ~~~~~ \EnrLvl\!=\!2 ~~~\, \EnrLvl\!=\!3 ~~ \, \EnrLvl\!=\!4}$\end{tabular}}}}%
    \put(0,0){\includegraphics[width=\unitlength,page=5]{flood_fill.pdf}}%
    \put(-0.00430462,1.50688389){\color[rgb]{0,0,0}\makebox(0,0)[lt]{\lineheight{1.25}\smash{\begin{tabular}[t]{l}(A)\end{tabular}}}}%
    \put(-0.00410775,0.71786938){\color[rgb]{0,0,0}\makebox(0,0)[lt]{\lineheight{1.25}\smash{\begin{tabular}[t]{l}(B)\end{tabular}}}}%
    \put(-0.00322293,0.55469807){\color[rgb]{0,0,0}\makebox(0,0)[lt]{\lineheight{1.25}\smash{\begin{tabular}[t]{l}(C)\end{tabular}}}}%
    \put(0.43558086,0.41319745){\color[rgb]{0,0,0}\makebox(0,0)[lt]{\lineheight{1.25}\smash{\begin{tabular}[t]{l}$7$\end{tabular}}}}%
    \put(0,0){\includegraphics[width=\unitlength,page=6]{flood_fill.pdf}}%
    \put(0.2445102,0.2385297){\color[rgb]{0,0,0}\makebox(0,0)[lt]{\lineheight{1.25}\smash{\begin{tabular}[t]{l}$1$\end{tabular}}}}%
    \put(0.24629718,0.35350714){\color[rgb]{0,0,0}\makebox(0,0)[lt]{\lineheight{1.25}\smash{\begin{tabular}[t]{l}$4$\end{tabular}}}}%
  \end{picture}%
\endgroup%

%% file: figures/unzipping.pdf_tex
\begingroup%
  \makeatletter%
  \providecommand\color[2][]{%
    \errmessage{(Inkscape) Color is used for the text in Inkscape, but the package 'color.sty' is not loaded}%
    \renewcommand\color[2][]{}%
  }%
  \providecommand\transparent[1]{%
    \errmessage{(Inkscape) Transparency is used (non-zero) for the text in Inkscape, but the package 'transparent.sty' is not loaded}%
    \renewcommand\transparent[1]{}%
  }%
  \providecommand\rotatebox[2]{#2}%
  \newcommand*\fsize{\dimexpr\f@size pt\relax}%
  \newcommand*\lineheight[1]{\fontsize{\fsize}{#1\fsize}\selectfont}%
  \ifx\svgwidth\undefined%
    \setlength{\unitlength}{283.53712427bp}%
    \ifx\svgscale\undefined%
      \relax%
    \else%
      \setlength{\unitlength}{\unitlength * \real{\svgscale}}%
    \fi%
  \else%
    \setlength{\unitlength}{\svgwidth}%
  \fi%
  \global\let\svgwidth\undefined%
  \global\let\svgscale\undefined%
  \makeatother%
  \begin{picture}(1,1.66727221)%
    \lineheight{1}%
    \setlength\tabcolsep{0pt}%
    \put(0,0){\includegraphics[width=\unitlength,page=1]{unzipping.pdf}}%
    \put(0.17629649,0.2499838){\color[rgb]{0,0,0}\makebox(0,0)[lt]{\lineheight{1.25}\smash{\begin{tabular}[t]{l}$\scriptstyle \BgElemInd=1$\end{tabular}}}}%
    \put(0.3621022,0.14772009){\color[rgb]{0,0,0}\makebox(0,0)[lt]{\lineheight{1.25}\smash{\begin{tabular}[t]{l}$\scriptstyle \BgElemInd=2$\end{tabular}}}}%
    \put(0.53396983,0.04334861){\color[rgb]{0,0,0}\makebox(0,0)[lt]{\lineheight{1.25}\smash{\begin{tabular}[t]{l}$\scriptstyle \BgElemInd=3$\end{tabular}}}}%
    \put(0.35213191,0.35149102){\color[rgb]{0,0,0}\makebox(0,0)[lt]{\lineheight{1.25}\smash{\begin{tabular}[t]{l}$\scriptstyle \BgElemInd=4$\end{tabular}}}}%
    \put(0.52278525,0.25496101){\color[rgb]{0,0,0}\makebox(0,0)[lt]{\lineheight{1.25}\smash{\begin{tabular}[t]{l}$\scriptstyle \BgElemInd=5$\end{tabular}}}}%
    \put(0.81823911,0.08418883){\color[rgb]{0,0,0}\makebox(0,0)[lt]{\lineheight{1.25}\smash{\begin{tabular}[t]{l}Background \\mesh $\BgMesh$\end{tabular}}}}%
    \put(0.04432557,0.86904561){\color[rgb]{0,0,0}\makebox(0,0)[lt]{\lineheight{1.25}\smash{\begin{tabular}[t]{l}$\Xi_{\BgElemInd}^1$\end{tabular}}}}%
    \put(0.84695093,0.49942022){\color[rgb]{0.82352941,0,0}\makebox(0,0)[lt]{\lineheight{1.25}\smash{\begin{tabular}[t]{l}copy \\elements \\$n_{\Unzipping}(\BgElemInd)$ \\times\end{tabular}}}}%
    \put(0.4507822,1.38001773){\color[rgb]{0,0,0}\makebox(0,0)[lt]{\lineheight{1.25}\smash{\begin{tabular}[t]{l}$\Xi_{\BgElemInd}^3$\end{tabular}}}}%
    \put(0.20634874,1.2319279){\color[rgb]{0,0,0}\makebox(0,0)[lt]{\lineheight{1.25}\smash{\begin{tabular}[t]{l}$\Xi_{\BgElemInd}^2$\end{tabular}}}}%
    \put(0.44502378,1.63317406){\color[rgb]{0,0,0}\makebox(0,0)[lt]{\lineheight{1.25}\smash{\begin{tabular}[t]{l}$\Xi_{\BgElemInd}^4$\end{tabular}}}}%
    \put(0.69884553,0.15695844){\color[rgb]{0,0,0}\makebox(0,0)[lt]{\lineheight{1.25}\smash{\begin{tabular}[t]{l}$\scriptstyle \BgElemInd=6$\end{tabular}}}}%
    \put(0.27996189,0.86933115){\color[rgb]{0.98823529,1,1}\rotatebox{-30}{\makebox(0,0)[lt]{\lineheight{1.25}\smash{\begin{tabular}[t]{l}$\scriptstyle \Subphase=4$\end{tabular}}}}}%
    \put(0.61261646,1.56504242){\color[rgb]{0.38431373,0.38431373,0.38431373}\rotatebox{-30}{\makebox(0,0)[lt]{\lineheight{1.25}\smash{\begin{tabular}[t]{l}$\scriptstyle \Subphase=9$\end{tabular}}}}}%
    \put(0.70148606,0.64889081){\color[rgb]{0.38431373,0.38431373,0.38431373}\rotatebox{-30}{\makebox(0,0)[lt]{\lineheight{1.25}\smash{\begin{tabular}[t]{l}$\scriptstyle \Subphase=10$\end{tabular}}}}}%
    \put(0.73046665,0.9262496){\color[rgb]{0.38431373,0.38431373,0.38431373}\rotatebox{-30}{\makebox(0,0)[lt]{\lineheight{1.25}\smash{\begin{tabular}[t]{l}$\scriptstyle \Subphase=11$\end{tabular}}}}}%
    \put(0.50076939,1.05517295){\color[rgb]{1,1,1}\rotatebox{-30}{\makebox(0,0)[lt]{\lineheight{1.25}\smash{\begin{tabular}[t]{l}$\scriptstyle \Subphase=7$\end{tabular}}}}}%
    \put(0.48445859,0.77032481){\color[rgb]{0.38431373,0.38431373,0.38431373}\rotatebox{-9.9624513}{\makebox(0,0)[lt]{\lineheight{1.25}\smash{\begin{tabular}[t]{l}$\scriptstyle \Subphase=6$\end{tabular}}}}}%
    \put(0.39547536,1.20329064){\color[rgb]{0.98823529,1,1}\rotatebox{-30.60312}{\makebox(0,0)[lt]{\lineheight{1.25}\smash{\begin{tabular}[t]{l}$\scriptstyle \Subphase=5$\end{tabular}}}}}%
    \put(0.61098643,1.33908596){\color[rgb]{0.98823529,1,1}\rotatebox{-30.60312}{\makebox(0,0)[lt]{\lineheight{1.25}\smash{\begin{tabular}[t]{l}$\scriptstyle \Subphase=8$\end{tabular}}}}}%
    \put(0,0){\includegraphics[width=\unitlength,page=2]{unzipping.pdf}}%
    \put(0.57342201,0.54769042){\color[rgb]{0.98823529,1,1}\rotatebox{-30}{\makebox(0,0)[lt]{\lineheight{1.25}\smash{\begin{tabular}[t]{l}$\scriptstyle \Subphase=3$\end{tabular}}}}}%
    \put(0,0){\includegraphics[width=\unitlength,page=3]{unzipping.pdf}}%
    \put(0.7140067,1.1403733){\color[rgb]{0.98823529,1,1}\rotatebox{-30.60312}{\makebox(0,0)[lt]{\lineheight{1.25}\smash{\begin{tabular}[t]{l}$\scriptstyle \Subphase=12$\end{tabular}}}}}%
    \put(0,0){\includegraphics[width=\unitlength,page=4]{unzipping.pdf}}%
    \put(0.2827011,0.61773191){\color[rgb]{1,1,1}\rotatebox{-30}{\makebox(0,0)[lt]{\lineheight{1.25}\smash{\begin{tabular}[t]{l}$\scriptstyle \Subphase=2$\end{tabular}}}}}%
    \put(0.06901494,0.7437924){\color[rgb]{0.98823529,1,1}\rotatebox{-30}{\makebox(0,0)[lt]{\lineheight{1.25}\smash{\begin{tabular}[t]{l}$\scriptstyle \Subphase=1$\end{tabular}}}}}%
  \end{picture}%
\endgroup%

%% file: figures/cluster_generation.pdf_tex
\begingroup%
  \makeatletter%
  \providecommand\color[2][]{%
    \errmessage{(Inkscape) Color is used for the text in Inkscape, but the package 'color.sty' is not loaded}%
    \renewcommand\color[2][]{}%
  }%
  \providecommand\transparent[1]{%
    \errmessage{(Inkscape) Transparency is used (non-zero) for the text in Inkscape, but the package 'transparent.sty' is not loaded}%
    \renewcommand\transparent[1]{}%
  }%
  \providecommand\rotatebox[2]{#2}%
  \newcommand*\fsize{\dimexpr\f@size pt\relax}%
  \newcommand*\lineheight[1]{\fontsize{\fsize}{#1\fsize}\selectfont}%
  \ifx\svgwidth\undefined%
    \setlength{\unitlength}{177.25200142bp}%
    \ifx\svgscale\undefined%
      \relax%
    \else%
      \setlength{\unitlength}{\unitlength * \real{\svgscale}}%
    \fi%
  \else%
    \setlength{\unitlength}{\svgwidth}%
  \fi%
  \global\let\svgwidth\undefined%
  \global\let\svgscale\undefined%
  \makeatother%
  \begin{picture}(1,0.81635642)%
    \lineheight{1}%
    \setlength\tabcolsep{0pt}%
    \put(0,0){\includegraphics[width=\unitlength,page=1]{cluster_generation.pdf}}%
    \put(0.03118865,0.61367809){\color[rgb]{0,0,0}\makebox(0,0)[lt]{\lineheight{1.25}\smash{\begin{tabular}[t]{l}$E=4$\end{tabular}}}}%
    \put(0.36489641,0.76416309){\color[rgb]{0,0,0}\makebox(0,0)[lt]{\lineheight{1.25}\smash{\begin{tabular}[t]{l}$m=2$\end{tabular}}}}%
    \put(0.3630776,0.01483864){\color[rgb]{0,0,0}\makebox(0,0)[lt]{\lineheight{1.25}\smash{\begin{tabular}[t]{l}$m=0$\end{tabular}}}}%
    \put(0.63915193,0.00895017){\color[rgb]{0,0,0}\makebox(0,0)[lt]{\lineheight{1.25}\smash{\begin{tabular}[t]{l}$\Xi_4^1$\end{tabular}}}}%
    \put(0.63427846,0.77728338){\color[rgb]{0,0,0}\makebox(0,0)[lt]{\lineheight{1.25}\smash{\begin{tabular}[t]{l}$\Xi_4^2$\end{tabular}}}}%
    \put(0,0){\includegraphics[width=\unitlength,page=2]{cluster_generation.pdf}}%
    \put(0.09868347,0.30683854){\color[rgb]{0,0,0}\makebox(0,0)[lt]{\lineheight{1.25}\smash{\begin{tabular}[t]{l}$4$\end{tabular}}}}%
    \put(0.09798252,0.46646201){\color[rgb]{0,0,0}\makebox(0,0)[lt]{\lineheight{1.25}\smash{\begin{tabular}[t]{l}$5$\end{tabular}}}}%
    \put(0,0){\includegraphics[width=\unitlength,page=3]{cluster_generation.pdf}}%
    \put(0.78650881,0.15425483){\color[rgb]{0.60784314,0.60784314,0.60784314}\makebox(0,0)[lt]{\lineheight{1.25}\smash{\begin{tabular}[t]{l}$\Gamma_{0,2}^I$\end{tabular}}}}%
    \put(0.78115938,0.52570378){\color[rgb]{0.60784314,0.60784314,0.60784314}\makebox(0,0)[lt]{\lineheight{1.25}\smash{\begin{tabular}[t]{l}$\Gamma_{2,0}^I$\end{tabular}}}}%
    \put(0,0){\includegraphics[width=\unitlength,page=4]{cluster_generation.pdf}}%
    \put(0.15415329,0.3878874){\color[rgb]{0.82352941,0,0.00392157}\makebox(0,0)[lt]{\lineheight{1.25}\smash{\begin{tabular}[t]{l}$\mathcal{G}_I$\end{tabular}}}}%
    \put(0.92100334,0.55485457){\color[rgb]{0.82352941,0,0.00392157}\makebox(0,0)[lt]{\lineheight{1.25}\smash{\begin{tabular}[t]{l}$\normal$\end{tabular}}}}%
    \put(0.92046405,0.33055789){\color[rgb]{0.82352941,0,0.00392157}\makebox(0,0)[lt]{\lineheight{1.25}\smash{\begin{tabular}[t]{l}$\normal$\end{tabular}}}}%
    \put(0,0){\includegraphics[width=\unitlength,page=5]{cluster_generation.pdf}}%
  \end{picture}%
\endgroup%

%% file: figures/ghost_cluster_generation.pdf_tex
\begingroup%
  \makeatletter%
  \providecommand\color[2][]{%
    \errmessage{(Inkscape) Color is used for the text in Inkscape, but the package 'color.sty' is not loaded}%
    \renewcommand\color[2][]{}%
  }%
  \providecommand\transparent[1]{%
    \errmessage{(Inkscape) Transparency is used (non-zero) for the text in Inkscape, but the package 'transparent.sty' is not loaded}%
    \renewcommand\transparent[1]{}%
  }%
  \providecommand\rotatebox[2]{#2}%
  \newcommand*\fsize{\dimexpr\f@size pt\relax}%
  \newcommand*\lineheight[1]{\fontsize{\fsize}{#1\fsize}\selectfont}%
  \ifx\svgwidth\undefined%
    \setlength{\unitlength}{387.43597556bp}%
    \ifx\svgscale\undefined%
      \relax%
    \else%
      \setlength{\unitlength}{\unitlength * \real{\svgscale}}%
    \fi%
  \else%
    \setlength{\unitlength}{\svgwidth}%
  \fi%
  \global\let\svgwidth\undefined%
  \global\let\svgscale\undefined%
  \makeatother%
  \begin{picture}(1,0.30204395)%
    \lineheight{1}%
    \setlength\tabcolsep{0pt}%
    \put(0,0){\includegraphics[width=\unitlength,page=1]{ghost_cluster_generation.pdf}}%
    \put(0.04962477,0.23908667){\color[rgb]{0,0,0}\makebox(0,0)[lt]{\lineheight{1.25}\smash{\begin{tabular}[t]{l}$E=5$\end{tabular}}}}%
    \put(0.00065242,0.2791126){\color[rgb]{0,0,0}\makebox(0,0)[lt]{\lineheight{1.25}\smash{\begin{tabular}[t]{l}(A)\end{tabular}}}}%
    \put(-0.79239827,-0.64655053){\color[rgb]{0,0,0}\makebox(0,0)[lt]{\lineheight{1.25}\smash{\begin{tabular}[t]{l}(B)\end{tabular}}}}%
    \put(0.20546188,0.23873277){\color[rgb]{0,0,0}\makebox(0,0)[lt]{\lineheight{1.25}\smash{\begin{tabular}[t]{l}$E=6$\end{tabular}}}}%
    \put(0.03775026,0.09662297){\color[rgb]{0,0,0}\makebox(0,0)[lt]{\lineheight{1.25}\smash{\begin{tabular}[t]{l}$7$\end{tabular}}}}%
    \put(0,0){\includegraphics[width=\unitlength,page=2]{ghost_cluster_generation.pdf}}%
    \put(0.03721489,0.17021138){\color[rgb]{0,0,0}\makebox(0,0)[lt]{\lineheight{1.25}\smash{\begin{tabular}[t]{l}$8$\end{tabular}}}}%
    \put(0.1098434,0.08139118){\color[rgb]{0,0,0}\makebox(0,0)[lt]{\lineheight{1.25}\smash{\begin{tabular}[t]{l}$6$\end{tabular}}}}%
    \put(0.19197138,0.17056422){\color[rgb]{0,0,0}\makebox(0,0)[lt]{\lineheight{1.25}\smash{\begin{tabular}[t]{l}$11$\end{tabular}}}}%
    \put(0.2632273,0.13302479){\color[rgb]{0,0,0}\makebox(0,0)[lt]{\lineheight{1.25}\smash{\begin{tabular}[t]{l}$12$\end{tabular}}}}%
    \put(0.19177166,0.09556552){\color[rgb]{0,0,0}\makebox(0,0)[lt]{\lineheight{1.25}\smash{\begin{tabular}[t]{l}$10$\end{tabular}}}}%
    \put(0.11126133,0.18463898){\color[rgb]{0,0,0}\makebox(0,0)[lt]{\lineheight{1.25}\smash{\begin{tabular}[t]{l}$9$\end{tabular}}}}%
    \put(0.1472423,0.02423338){\color[rgb]{0,0.49019608,0.97254902}\makebox(0,0)[lt]{\lineheight{1.25}\smash{\begin{tabular}[t]{l}$\mathcal{G}_S$\end{tabular}}}}%
    \put(0,0){\includegraphics[width=\unitlength,page=3]{ghost_cluster_generation.pdf}}%
    \put(-0.77684207,-0.97757125){\color[rgb]{0,0,0}\rotatebox{30}{\makebox(0,0)[lt]{\lineheight{1.25}\smash{\begin{tabular}[t]{l}$m \! = \! 0$\end{tabular}}}}}%
    \put(-1.0068951,-0.70080523){\color[rgb]{0,0,0}\rotatebox{-10.278728}{\makebox(0,0)[lt]{\lineheight{1.25}\smash{\begin{tabular}[t]{l}$m \! = \! 0$\end{tabular}}}}}%
    \put(-0.64946764,-0.68695434){\color[rgb]{0,0,0}\rotatebox{-30}{\makebox(0,0)[lt]{\lineheight{1.25}\smash{\begin{tabular}[t]{l}$m \! = \! 2$\end{tabular}}}}}%
    \put(-0.63244371,-0.80858428){\color[rgb]{0,0,0}\rotatebox{-30}{\makebox(0,0)[lt]{\lineheight{1.25}\smash{\begin{tabular}[t]{l}$m \! = \! 1$\end{tabular}}}}}%
    \put(-1.02201216,-0.81839486){\color[rgb]{0,0,0}\rotatebox{30}{\makebox(0,0)[lt]{\lineheight{1.25}\smash{\begin{tabular}[t]{l}$\Xi[5][4]$\end{tabular}}}}}%
    \put(-1.02311817,-0.67725893){\color[rgb]{0,0,0}\rotatebox{30}{\makebox(0,0)[lt]{\lineheight{1.25}\smash{\begin{tabular}[t]{l}$\Xi[5][1]$\end{tabular}}}}}%
    \put(-0.69871696,-0.81799473){\color[rgb]{0,0,0}\rotatebox{30}{\makebox(0,0)[lt]{\lineheight{1.25}\smash{\begin{tabular}[t]{l}$\Xi[5][2]$\end{tabular}}}}}%
    \put(-0.69863682,-0.67601177){\color[rgb]{0,0,0}\rotatebox{30}{\makebox(0,0)[lt]{\lineheight{1.25}\smash{\begin{tabular}[t]{l}$\Xi[5][3]$\end{tabular}}}}}%
    \put(0.76672531,0.28864914){\color[rgb]{0,0,0}\makebox(0,0)[lt]{\lineheight{1.25}\smash{\begin{tabular}[t]{l}$\Xi_5^4$\end{tabular}}}}%
    \put(0.61779307,0.2887197){\color[rgb]{0,0,0}\makebox(0,0)[lt]{\lineheight{1.25}\smash{\begin{tabular}[t]{l}$\Xi_6^2$\end{tabular}}}}%
    \put(0.61880452,0.12630441){\color[rgb]{0,0,0}\makebox(0,0)[lt]{\lineheight{1.25}\smash{\begin{tabular}[t]{l}$\Xi_6^1$\end{tabular}}}}%
    \put(0.91760879,0.28931937){\color[rgb]{0,0,0}\makebox(0,0)[lt]{\lineheight{1.25}\smash{\begin{tabular}[t]{l}$\Xi_6^3$\end{tabular}}}}%
    \put(0,0){\includegraphics[width=\unitlength,page=4]{ghost_cluster_generation.pdf}}%
    \put(-0.78439995,-0.86431512){\color[rgb]{0,0,0}\rotatebox{-30}{\makebox(0,0)[lt]{\lineheight{1.25}\smash{\begin{tabular}[t]{l}$\Xi[5][3]$\end{tabular}}}}}%
    \put(-0.78376848,-0.72364801){\color[rgb]{0,0,0}\rotatebox{-30}{\makebox(0,0)[lt]{\lineheight{1.25}\smash{\begin{tabular}[t]{l}$\Xi[5][3]$\end{tabular}}}}}%
    \put(-0.4608827,-0.72364815){\color[rgb]{0,0,0}\rotatebox{-30}{\makebox(0,0)[lt]{\lineheight{1.25}\smash{\begin{tabular}[t]{l}$\Xi[5][2]$\end{tabular}}}}}%
    \put(-0.46081852,-0.86494772){\color[rgb]{0,0,0}\rotatebox{-30}{\makebox(0,0)[lt]{\lineheight{1.25}\smash{\begin{tabular}[t]{l}$\Xi[5][1]$\end{tabular}}}}}%
    \put(0.4482322,0.21560526){\color[rgb]{0,0,0}\makebox(0,0)[lt]{\lineheight{1.25}\smash{\begin{tabular}[t]{l}$m \! = \! 2$\end{tabular}}}}%
    \put(0.7512843,0.053289){\color[rgb]{0,0,0}\makebox(0,0)[lt]{\lineheight{1.25}\smash{\begin{tabular}[t]{l}$m \! = \! 0$\end{tabular}}}}%
    \put(0.75551366,0.24555826){\color[rgb]{0,0,0}\makebox(0,0)[lt]{\lineheight{1.25}\smash{\begin{tabular}[t]{l}$m \! = \! 0$\end{tabular}}}}%
    \put(0.4436515,0.08559654){\color[rgb]{0,0,0}\makebox(0,0)[lt]{\lineheight{1.25}\smash{\begin{tabular}[t]{l}$m \! = \! 1$\end{tabular}}}}%
    \put(0.46758601,0.28882213){\color[rgb]{0,0,0}\makebox(0,0)[lt]{\lineheight{1.25}\smash{\begin{tabular}[t]{l}$\Xi_5^3$\end{tabular}}}}%
    \put(0.91942205,0.12725282){\color[rgb]{0,0,0}\makebox(0,0)[lt]{\lineheight{1.25}\smash{\begin{tabular}[t]{l}$\Xi_6^3$\end{tabular}}}}%
    \put(0.46757931,0.12631122){\color[rgb]{0,0,0}\makebox(0,0)[lt]{\lineheight{1.25}\smash{\begin{tabular}[t]{l}$\Xi_5^2$\end{tabular}}}}%
    \put(0.76770964,0.12631122){\color[rgb]{0,0,0}\makebox(0,0)[lt]{\lineheight{1.25}\smash{\begin{tabular}[t]{l}$\Xi_5^1$\end{tabular}}}}%
    \put(0.25375037,0.07547383){\color[rgb]{0,0.49019608,0.97254902}\makebox(0,0)[lt]{\lineheight{1.25}\smash{\begin{tabular}[t]{l}(4)\end{tabular}}}}%
    \put(0.13120977,0.13003103){\color[rgb]{0,0.49019608,0.97254902}\makebox(0,0)[lt]{\lineheight{1.25}\smash{\begin{tabular}[t]{l}(2)\end{tabular}}}}%
    \put(0.06903479,0.13975826){\color[rgb]{0,0.49019608,0.97254902}\makebox(0,0)[lt]{\lineheight{1.25}\smash{\begin{tabular}[t]{l}(1)\end{tabular}}}}%
    \put(0.25565101,0.19418172){\color[rgb]{0,0.49019608,0.97254902}\makebox(0,0)[lt]{\lineheight{1.25}\smash{\begin{tabular}[t]{l}(3)\end{tabular}}}}%
    \put(-0.83272155,-0.7025957){\color[rgb]{0,0.49019608,0.97254902}\makebox(0,0)[lt]{\lineheight{1.25}\smash{\begin{tabular}[t]{l}(1)\end{tabular}}}}%
    \put(-0.62508577,-0.9195498){\color[rgb]{0,0.49019608,0.97254902}\makebox(0,0)[lt]{\lineheight{1.25}\smash{\begin{tabular}[t]{l}(2)\end{tabular}}}}%
    \put(-0.51656522,-0.70307827){\color[rgb]{0,0.49019608,0.97254902}\makebox(0,0)[lt]{\lineheight{1.25}\smash{\begin{tabular}[t]{l}(3)\end{tabular}}}}%
    \put(-0.949578,-0.92091876){\color[rgb]{0,0.49019608,0.97254902}\makebox(0,0)[lt]{\lineheight{1.25}\smash{\begin{tabular}[t]{l}(4)\end{tabular}}}}%
    \put(0,0){\includegraphics[width=\unitlength,page=5]{ghost_cluster_generation.pdf}}%
    \put(0.36713244,0.28109571){\color[rgb]{0,0,0}\makebox(0,0)[lt]{\lineheight{1.25}\smash{\begin{tabular}[t]{l}(B)\end{tabular}}}}%
    \put(0.54121612,0.28213586){\color[rgb]{0,0.49019608,0.97254902}\makebox(0,0)[lt]{\lineheight{1.25}\smash{\begin{tabular}[t]{l}(1)\end{tabular}}}}%
    \put(0.54128137,0.12023361){\color[rgb]{0,0.49019608,0.97254902}\makebox(0,0)[lt]{\lineheight{1.25}\smash{\begin{tabular}[t]{l}(2)\end{tabular}}}}%
    \put(0.84484166,0.28313045){\color[rgb]{0,0.49019608,0.97254902}\makebox(0,0)[lt]{\lineheight{1.25}\smash{\begin{tabular}[t]{l}(3)\end{tabular}}}}%
    \put(0.84484166,0.12110125){\color[rgb]{0,0.49019608,0.97254902}\makebox(0,0)[lt]{\lineheight{1.25}\smash{\begin{tabular}[t]{l}(4)\end{tabular}}}}%
  \end{picture}%
\endgroup%

%% file: figures/aura.pdf_tex
\begingroup%
  \makeatletter%
  \providecommand\color[2][]{%
    \errmessage{(Inkscape) Color is used for the text in Inkscape, but the package 'color.sty' is not loaded}%
    \renewcommand\color[2][]{}%
  }%
  \providecommand\transparent[1]{%
    \errmessage{(Inkscape) Transparency is used (non-zero) for the text in Inkscape, but the package 'transparent.sty' is not loaded}%
    \renewcommand\transparent[1]{}%
  }%
  \providecommand\rotatebox[2]{#2}%
  \newcommand*\fsize{\dimexpr\f@size pt\relax}%
  \newcommand*\lineheight[1]{\fontsize{\fsize}{#1\fsize}\selectfont}%
  \ifx\svgwidth\undefined%
    \setlength{\unitlength}{170.44008871bp}%
    \ifx\svgscale\undefined%
      \relax%
    \else%
      \setlength{\unitlength}{\unitlength * \real{\svgscale}}%
    \fi%
  \else%
    \setlength{\unitlength}{\svgwidth}%
  \fi%
  \global\let\svgwidth\undefined%
  \global\let\svgscale\undefined%
  \makeatother%
  \begin{picture}(1,2.1770458)%
    \lineheight{1}%
    \setlength\tabcolsep{0pt}%
    \put(0,0){\includegraphics[width=\unitlength,page=1]{aura.pdf}}%
    \put(0.29025219,2.12008129){\color[rgb]{0.82352941,0,0}\makebox(0,0)[lt]{\lineheight{1.25}\smash{\begin{tabular}[t]{l}$p = 3$\end{tabular}}}}%
    \put(0.18568919,2.03605114){\color[rgb]{0.33333333,0.70980392,0}\makebox(0,0)[lt]{\lineheight{1.25}\smash{\begin{tabular}[t]{l}$\Gamma$\end{tabular}}}}%
    \put(0.11589779,0.70278049){\color[rgb]{0.33333333,0.70980392,0}\makebox(0,0)[lt]{\lineheight{1.25}\smash{\begin{tabular}[t]{l}$\Gamma$\end{tabular}}}}%
    \put(0.66306648,2.12108217){\color[rgb]{0.82352941,0,0}\makebox(0,0)[lt]{\lineheight{1.25}\smash{\begin{tabular}[t]{l}$p = 4$\end{tabular}}}}%
    \put(0.28797152,1.28474755){\color[rgb]{0.82352941,0,0}\makebox(0,0)[lt]{\lineheight{1.25}\smash{\begin{tabular}[t]{l}$p = 1$\end{tabular}}}}%
    \put(0.66145095,1.28404434){\color[rgb]{0.82352941,0,0}\makebox(0,0)[lt]{\lineheight{1.25}\smash{\begin{tabular}[t]{l}$p = 2$\end{tabular}}}}%
    \put(0.22790881,1.12399277){\color[rgb]{0.82352941,0,0}\makebox(0,0)[lt]{\lineheight{1.25}\smash{\begin{tabular}[t]{l}$p = 3$\end{tabular}}}}%
    \put(0.72463712,1.12415265){\color[rgb]{0.82352941,0,0}\makebox(0,0)[lt]{\lineheight{1.25}\smash{\begin{tabular}[t]{l}$p = 4$\end{tabular}}}}%
    \put(0.2256572,0.12419135){\color[rgb]{0.82352941,0,0}\makebox(0,0)[lt]{\lineheight{1.25}\smash{\begin{tabular}[t]{l}$p = 1$\end{tabular}}}}%
    \put(0.72723272,0.1244264){\color[rgb]{0.82352941,0,0}\makebox(0,0)[lt]{\lineheight{1.25}\smash{\begin{tabular}[t]{l}$p = 2$\end{tabular}}}}%
    \put(0,0){\includegraphics[width=\unitlength,page=2]{aura.pdf}}%
    \put(0.55281285,1.99443779){\color[rgb]{0,0.47843137,1}\makebox(0,0)[lt]{\lineheight{1.25}\smash{\begin{tabular}[t]{l}$\mathrm{supp}(N_{\BfIndex})$\end{tabular}}}}%
    \put(0.47533869,0.07450916){\color[rgb]{0.38823529,0.38823529,0.38823529}\makebox(0,0)[lt]{\lineheight{1.25}\smash{\begin{tabular}[t]{l}Aura\end{tabular}}}}%
    \put(0.46419446,0.00079923){\color[rgb]{0.00392157,0.00392157,0.00392157}\makebox(0,0)[lt]{\lineheight{1.25}\smash{\begin{tabular}[t]{l}Owned domain\end{tabular}}}}%
    \put(-0.00209912,2.1369268){\color[rgb]{0.00392157,0.00392157,0.00392157}\makebox(0,0)[lt]{\lineheight{1.25}\smash{\begin{tabular}[t]{l}(A)\end{tabular}}}}%
    \put(-0.0046668,1.13892711){\color[rgb]{0.00392157,0.00392157,0.00392157}\makebox(0,0)[lt]{\lineheight{1.25}\smash{\begin{tabular}[t]{l}(B)\end{tabular}}}}%
    \put(0,0){\includegraphics[width=\unitlength,page=3]{aura.pdf}}%
  \end{picture}%
\endgroup%

%% file: figures/multi_beam.pdf_tex
\begingroup%
  \makeatletter%
  \providecommand\color[2][]{%
    \errmessage{(Inkscape) Color is used for the text in Inkscape, but the package 'color.sty' is not loaded}%
    \renewcommand\color[2][]{}%
  }%
  \providecommand\transparent[1]{%
    \errmessage{(Inkscape) Transparency is used (non-zero) for the text in Inkscape, but the package 'transparent.sty' is not loaded}%
    \renewcommand\transparent[1]{}%
  }%
  \providecommand\rotatebox[2]{#2}%
  \newcommand*\fsize{\dimexpr\f@size pt\relax}%
  \newcommand*\lineheight[1]{\fontsize{\fsize}{#1\fsize}\selectfont}%
  \ifx\svgwidth\undefined%
    \setlength{\unitlength}{366.03568971bp}%
    \ifx\svgscale\undefined%
      \relax%
    \else%
      \setlength{\unitlength}{\unitlength * \real{\svgscale}}%
    \fi%
  \else%
    \setlength{\unitlength}{\svgwidth}%
  \fi%
  \global\let\svgwidth\undefined%
  \global\let\svgscale\undefined%
  \makeatother%
  \begin{picture}(1,0.40189649)%
    \lineheight{1}%
    \setlength\tabcolsep{0pt}%
    \put(0,0){\includegraphics[width=\unitlength,page=1]{multi_beam.pdf}}%
    \put(0.02848652,0.37618595){\color[rgb]{0.82352941,0,0.00392157}\makebox(0,0)[lt]{\lineheight{1.25}\smash{\begin{tabular}[t]{l}$\traction$\end{tabular}}}}%
    \put(0,0){\includegraphics[width=\unitlength,page=2]{multi_beam.pdf}}%
    \put(0.24867886,0.2519799){\color[rgb]{0.82352941,0,0.00392157}\makebox(0,0)[lt]{\lineheight{1.25}\smash{\begin{tabular}[t]{l}$h_0$\end{tabular}}}}%
    \put(0.09311046,0.19715252){\color[rgb]{0.82352941,0,0.00392157}\makebox(0,0)[lt]{\lineheight{1.25}\smash{\begin{tabular}[t]{l}$\mathrm{supp}(N)$\end{tabular}}}}%
    \put(0.00846505,0.22366489){\color[rgb]{0.82352941,0,0.00392157}\makebox(0,0)[lt]{\lineheight{1.25}\smash{\begin{tabular}[t]{l}$\delta_o$\end{tabular}}}}%
    \put(0.11201313,0.31101614){\color[rgb]{0.82352941,0,0.00392157}\makebox(0,0)[lt]{\lineheight{1.25}\smash{\begin{tabular}[t]{l}$\delta_d$\end{tabular}}}}%
    \put(0.08378556,0.28177424){\color[rgb]{0.82352941,0,0.00392157}\makebox(0,0)[lt]{\lineheight{1.25}\smash{\begin{tabular}[t]{l}$\delta_t$\end{tabular}}}}%
    \put(0,0){\includegraphics[width=\unitlength,page=3]{multi_beam.pdf}}%
    \put(0.57551188,0.20398481){\color[rgb]{0,0,0}\makebox(0,0)[lt]{\lineheight{1.25}\smash{\begin{tabular}[t]{l}$\sigma_{yy}$\end{tabular}}}}%
    \put(0.10064625,0.38321558){\color[rgb]{0,0,0}\makebox(0,0)[lt]{\lineheight{1.25}\smash{\begin{tabular}[t]{l}(A)\end{tabular}}}}%
    \put(0.42004185,0.38314716){\color[rgb]{0,0,0}\makebox(0,0)[lt]{\lineheight{1.25}\smash{\begin{tabular}[t]{l}(B)\end{tabular}}}}%
    \put(0.78665194,0.382402){\color[rgb]{0,0,0}\makebox(0,0)[lt]{\lineheight{1.25}\smash{\begin{tabular}[t]{l}(C)\end{tabular}}}}%
    \put(0.91853729,0.20468672){\color[rgb]{0,0,0}\makebox(0,0)[lt]{\lineheight{1.25}\smash{\begin{tabular}[t]{l}$\sigma_{yy}$\end{tabular}}}}%
    \put(0.57547844,0.17739763){\color[rgb]{0,0,0}\makebox(0,0)[lt]{\lineheight{1.25}\smash{\begin{tabular}[t]{l}$1$\end{tabular}}}}%
    \put(0.91822805,0.09490682){\color[rgb]{0,0,0}\makebox(0,0)[lt]{\lineheight{1.25}\smash{\begin{tabular}[t]{l}$0$\end{tabular}}}}%
    \put(0.5751113,0.09571508){\color[rgb]{0,0,0}\makebox(0,0)[lt]{\lineheight{1.25}\smash{\begin{tabular}[t]{l}$0$\end{tabular}}}}%
    \put(0.57562355,0.01240726){\color[rgb]{0,0,0}\makebox(0,0)[lt]{\lineheight{1.25}\smash{\begin{tabular}[t]{l}$-1$\end{tabular}}}}%
    \put(0.91854521,0.17696491){\color[rgb]{0,0,0}\makebox(0,0)[lt]{\lineheight{1.25}\smash{\begin{tabular}[t]{l}$7$\end{tabular}}}}%
    \put(0.91812607,0.01173005){\color[rgb]{0,0,0}\makebox(0,0)[lt]{\lineheight{1.25}\smash{\begin{tabular}[t]{l}$-7$\end{tabular}}}}%
    \put(0,0){\includegraphics[width=\unitlength,page=4]{multi_beam.pdf}}%
    \put(0.07197993,0.0218684){\color[rgb]{0.17254902,0.17254902,0.17254902}\makebox(0,0)[lt]{\lineheight{1.25}\smash{\begin{tabular}[t]{l}$x$\end{tabular}}}}%
    \put(0.00799481,0.08866667){\color[rgb]{0.17254902,0.17254902,0.17254902}\makebox(0,0)[lt]{\lineheight{1.25}\smash{\begin{tabular}[t]{l}$y$\end{tabular}}}}%
  \end{picture}%
\endgroup%

%% file: figures/multi_beam_L2_error.pdf_tex
\begingroup
  \makeatletter
  \providecommand\color[2][]{%
    \GenericError{(gnuplot) \space\space\space\@spaces}{%
      Package color not loaded in conjunction with
      terminal option `colourtext'%
    }{See the gnuplot documentation for explanation.%
    }{Either use 'blacktext' in gnuplot or load the package
      color.sty in LaTeX.}%
    \renewcommand\color[2][]{}%
  }%
  \providecommand\includegraphics[2][]{%
    \GenericError{(gnuplot) \space\space\space\@spaces}{%
      Package graphicx or graphics not loaded%
    }{See the gnuplot documentation for explanation.%
    }{The gnuplot epslatex terminal needs graphicx.sty or graphics.sty.}%
    \renewcommand\includegraphics[2][]{}%
  }%
  \providecommand\rotatebox[2]{#2}%
  \@ifundefined{ifGPcolor}{%
    \newif\ifGPcolor
    \GPcolorfalse
  }{}%
  \@ifundefined{ifGPblacktext}{%
    \newif\ifGPblacktext
    \GPblacktexttrue
  }{}%
  \let\gplgaddtomacro\g@addto@macro
  \gdef\gplbacktext{}%
  \gdef\gplfronttext{}%
  \makeatother
  \ifGPblacktext
    \def\colorrgb#1{}%
    \def\colorgray#1{}%
  \else
    \ifGPcolor
      \def\colorrgb#1{\color[rgb]{#1}}%
      \def\colorgray#1{\color[gray]{#1}}%
      \expandafter\def\csname LTw\endcsname{\color{white}}%
      \expandafter\def\csname LTb\endcsname{\color{black}}%
      \expandafter\def\csname LTa\endcsname{\color{black}}%
      \expandafter\def\csname LT0\endcsname{\color[rgb]{1,0,0}}%
      \expandafter\def\csname LT1\endcsname{\color[rgb]{0,1,0}}%
      \expandafter\def\csname LT2\endcsname{\color[rgb]{0,0,1}}%
      \expandafter\def\csname LT3\endcsname{\color[rgb]{1,0,1}}%
      \expandafter\def\csname LT4\endcsname{\color[rgb]{0,1,1}}%
      \expandafter\def\csname LT5\endcsname{\color[rgb]{1,1,0}}%
      \expandafter\def\csname LT6\endcsname{\color[rgb]{0,0,0}}%
      \expandafter\def\csname LT7\endcsname{\color[rgb]{1,0.3,0}}%
      \expandafter\def\csname LT8\endcsname{\color[rgb]{0.5,0.5,0.5}}%
    \else
      \def\colorrgb#1{\color{black}}%
      \def\colorgray#1{\color[gray]{#1}}%
      \expandafter\def\csname LTw\endcsname{\color{white}}%
      \expandafter\def\csname LTb\endcsname{\color{black}}%
      \expandafter\def\csname LTa\endcsname{\color{black}}%
      \expandafter\def\csname LT0\endcsname{\color{black}}%
      \expandafter\def\csname LT1\endcsname{\color{black}}%
      \expandafter\def\csname LT2\endcsname{\color{black}}%
      \expandafter\def\csname LT3\endcsname{\color{black}}%
      \expandafter\def\csname LT4\endcsname{\color{black}}%
      \expandafter\def\csname LT5\endcsname{\color{black}}%
      \expandafter\def\csname LT6\endcsname{\color{black}}%
      \expandafter\def\csname LT7\endcsname{\color{black}}%
      \expandafter\def\csname LT8\endcsname{\color{black}}%
    \fi
  \fi
    \setlength{\unitlength}{0.0500bp}%
    \ifx\gptboxheight\undefined%
      \newlength{\gptboxheight}%
      \newlength{\gptboxwidth}%
      \newsavebox{\gptboxtext}%
    \fi%
    \setlength{\fboxrule}{0.5pt}%
    \setlength{\fboxsep}{1pt}%
    \definecolor{tbcol}{rgb}{1,1,1}%
\begin{picture}(3514.00,2664.00)%
    \gplgaddtomacro\gplbacktext{%
      \csname LTb\endcsname
      \put(486,572){\makebox(0,0)[r]{\strut{}$10^{-6}$}}%
      \put(486,1251){\makebox(0,0)[r]{\strut{}$10^{-4}$}}%
      \put(486,1930){\makebox(0,0)[r]{\strut{}$10^{-2}$}}%
      \put(486,2609){\makebox(0,0)[r]{\strut{}$10^{0}$}}%
      \put(561,396){\makebox(0,0){\strut{}$10^{2}$}}%
      \put(1516,396){\makebox(0,0){\strut{}$10^{3}$}}%
      \put(2471,396){\makebox(0,0){\strut{}$10^{4}$}}%
      \put(3425,396){\makebox(0,0){\strut{}$10^{5}$}}%
    }%
    \gplgaddtomacro\gplfronttext{%
      \csname LTb\endcsname
      \put(3026,2433){\makebox(0,0)[r]{\strut{}Enrich.: Off }}%
      \csname LTb\endcsname
      \put(3026,2207){\makebox(0,0)[r]{\strut{}Enrich.: On }}%
      \csname LTb\endcsname
      \put(3026,1981){\makebox(0,0)[r]{\strut{}$r=4$ }}%
      \csname LTb\endcsname
      \put(-136,1590){\rotatebox{-270.00}{\makebox(0,0){\footnotesize$\norm{\DiscDispl-\displ}{\mathrm{L2}}/\norm{\displ}{\mathrm{L2}}$}}}%
      \put(1993,154){\makebox(0,0){\small Number of DOFs}}%
    }%
    \gplbacktext
    \put(0,0){\includegraphics[width={175.70bp},height={133.20bp}]{multi_beam_L2_error.pdf}}%
    \gplfronttext
  \end{picture}%
\endgroup

%% file: figures/multi_beam_H1s_error.pdf_tex
\begingroup
  \makeatletter
  \providecommand\color[2][]{%
    \GenericError{(gnuplot) \space\space\space\@spaces}{%
      Package color not loaded in conjunction with
      terminal option `colourtext'%
    }{See the gnuplot documentation for explanation.%
    }{Either use 'blacktext' in gnuplot or load the package
      color.sty in LaTeX.}%
    \renewcommand\color[2][]{}%
  }%
  \providecommand\includegraphics[2][]{%
    \GenericError{(gnuplot) \space\space\space\@spaces}{%
      Package graphicx or graphics not loaded%
    }{See the gnuplot documentation for explanation.%
    }{The gnuplot epslatex terminal needs graphicx.sty or graphics.sty.}%
    \renewcommand\includegraphics[2][]{}%
  }%
  \providecommand\rotatebox[2]{#2}%
  \@ifundefined{ifGPcolor}{%
    \newif\ifGPcolor
    \GPcolorfalse
  }{}%
  \@ifundefined{ifGPblacktext}{%
    \newif\ifGPblacktext
    \GPblacktexttrue
  }{}%
  \let\gplgaddtomacro\g@addto@macro
  \gdef\gplbacktext{}%
  \gdef\gplfronttext{}%
  \makeatother
  \ifGPblacktext
    \def\colorrgb#1{}%
    \def\colorgray#1{}%
  \else
    \ifGPcolor
      \def\colorrgb#1{\color[rgb]{#1}}%
      \def\colorgray#1{\color[gray]{#1}}%
      \expandafter\def\csname LTw\endcsname{\color{white}}%
      \expandafter\def\csname LTb\endcsname{\color{black}}%
      \expandafter\def\csname LTa\endcsname{\color{black}}%
      \expandafter\def\csname LT0\endcsname{\color[rgb]{1,0,0}}%
      \expandafter\def\csname LT1\endcsname{\color[rgb]{0,1,0}}%
      \expandafter\def\csname LT2\endcsname{\color[rgb]{0,0,1}}%
      \expandafter\def\csname LT3\endcsname{\color[rgb]{1,0,1}}%
      \expandafter\def\csname LT4\endcsname{\color[rgb]{0,1,1}}%
      \expandafter\def\csname LT5\endcsname{\color[rgb]{1,1,0}}%
      \expandafter\def\csname LT6\endcsname{\color[rgb]{0,0,0}}%
      \expandafter\def\csname LT7\endcsname{\color[rgb]{1,0.3,0}}%
      \expandafter\def\csname LT8\endcsname{\color[rgb]{0.5,0.5,0.5}}%
    \else
      \def\colorrgb#1{\color{black}}%
      \def\colorgray#1{\color[gray]{#1}}%
      \expandafter\def\csname LTw\endcsname{\color{white}}%
      \expandafter\def\csname LTb\endcsname{\color{black}}%
      \expandafter\def\csname LTa\endcsname{\color{black}}%
      \expandafter\def\csname LT0\endcsname{\color{black}}%
      \expandafter\def\csname LT1\endcsname{\color{black}}%
      \expandafter\def\csname LT2\endcsname{\color{black}}%
      \expandafter\def\csname LT3\endcsname{\color{black}}%
      \expandafter\def\csname LT4\endcsname{\color{black}}%
      \expandafter\def\csname LT5\endcsname{\color{black}}%
      \expandafter\def\csname LT6\endcsname{\color{black}}%
      \expandafter\def\csname LT7\endcsname{\color{black}}%
      \expandafter\def\csname LT8\endcsname{\color{black}}%
    \fi
  \fi
    \setlength{\unitlength}{0.0500bp}%
    \ifx\gptboxheight\undefined%
      \newlength{\gptboxheight}%
      \newlength{\gptboxwidth}%
      \newsavebox{\gptboxtext}%
    \fi%
    \setlength{\fboxrule}{0.5pt}%
    \setlength{\fboxsep}{1pt}%
    \definecolor{tbcol}{rgb}{1,1,1}%
\begin{picture}(3514.00,2664.00)%
    \gplgaddtomacro\gplbacktext{%
      \csname LTb\endcsname
      \put(486,572){\makebox(0,0)[r]{\strut{}$10^{-5}$}}%
      \put(486,979){\makebox(0,0)[r]{\strut{}$10^{-4}$}}%
      \put(486,1387){\makebox(0,0)[r]{\strut{}$10^{-3}$}}%
      \put(486,1794){\makebox(0,0)[r]{\strut{}$10^{-2}$}}%
      \put(486,2202){\makebox(0,0)[r]{\strut{}$10^{-1}$}}%
      \put(486,2609){\makebox(0,0)[r]{\strut{}$10^{0}$}}%
      \put(561,396){\makebox(0,0){\strut{}$10^{2}$}}%
      \put(1516,396){\makebox(0,0){\strut{}$10^{3}$}}%
      \put(2471,396){\makebox(0,0){\strut{}$10^{4}$}}%
      \put(3425,396){\makebox(0,0){\strut{}$10^{5}$}}%
    }%
    \gplgaddtomacro\gplfronttext{%
      \csname LTb\endcsname
      \put(3026,2433){\makebox(0,0)[r]{\strut{}Enrich.: Off }}%
      \csname LTb\endcsname
      \put(3026,2207){\makebox(0,0)[r]{\strut{}Enrich.: On }}%
      \csname LTb\endcsname
      \put(3026,1981){\makebox(0,0)[r]{\strut{}$r=2$ }}%
      \csname LTb\endcsname
      \put(-136,1590){\rotatebox{-270.00}{\makebox(0,0){\footnotesize$\abs{\DiscDispl-\displ}_{\mathrm{H1}}/\abs{\displ}_{\mathrm{H1}}$}}}%
      \put(1993,154){\makebox(0,0){\small Number of DOFs}}%
    }%
    \gplbacktext
    \put(0,0){\includegraphics[width={175.70bp},height={133.20bp}]{multi_beam_H1s_error.pdf}}%
    \gplfronttext
  \end{picture}%
\endgroup

%% file: figures/multi_beam_cond.pdf_tex
\begingroup
  \makeatletter
  \providecommand\color[2][]{%
    \GenericError{(gnuplot) \space\space\space\@spaces}{%
      Package color not loaded in conjunction with
      terminal option `colourtext'%
    }{See the gnuplot documentation for explanation.%
    }{Either use 'blacktext' in gnuplot or load the package
      color.sty in LaTeX.}%
    \renewcommand\color[2][]{}%
  }%
  \providecommand\includegraphics[2][]{%
    \GenericError{(gnuplot) \space\space\space\@spaces}{%
      Package graphicx or graphics not loaded%
    }{See the gnuplot documentation for explanation.%
    }{The gnuplot epslatex terminal needs graphicx.sty or graphics.sty.}%
    \renewcommand\includegraphics[2][]{}%
  }%
  \providecommand\rotatebox[2]{#2}%
  \@ifundefined{ifGPcolor}{%
    \newif\ifGPcolor
    \GPcolorfalse
  }{}%
  \@ifundefined{ifGPblacktext}{%
    \newif\ifGPblacktext
    \GPblacktexttrue
  }{}%
  \let\gplgaddtomacro\g@addto@macro
  \gdef\gplbacktext{}%
  \gdef\gplfronttext{}%
  \makeatother
  \ifGPblacktext
    \def\colorrgb#1{}%
    \def\colorgray#1{}%
  \else
    \ifGPcolor
      \def\colorrgb#1{\color[rgb]{#1}}%
      \def\colorgray#1{\color[gray]{#1}}%
      \expandafter\def\csname LTw\endcsname{\color{white}}%
      \expandafter\def\csname LTb\endcsname{\color{black}}%
      \expandafter\def\csname LTa\endcsname{\color{black}}%
      \expandafter\def\csname LT0\endcsname{\color[rgb]{1,0,0}}%
      \expandafter\def\csname LT1\endcsname{\color[rgb]{0,1,0}}%
      \expandafter\def\csname LT2\endcsname{\color[rgb]{0,0,1}}%
      \expandafter\def\csname LT3\endcsname{\color[rgb]{1,0,1}}%
      \expandafter\def\csname LT4\endcsname{\color[rgb]{0,1,1}}%
      \expandafter\def\csname LT5\endcsname{\color[rgb]{1,1,0}}%
      \expandafter\def\csname LT6\endcsname{\color[rgb]{0,0,0}}%
      \expandafter\def\csname LT7\endcsname{\color[rgb]{1,0.3,0}}%
      \expandafter\def\csname LT8\endcsname{\color[rgb]{0.5,0.5,0.5}}%
    \else
      \def\colorrgb#1{\color{black}}%
      \def\colorgray#1{\color[gray]{#1}}%
      \expandafter\def\csname LTw\endcsname{\color{white}}%
      \expandafter\def\csname LTb\endcsname{\color{black}}%
      \expandafter\def\csname LTa\endcsname{\color{black}}%
      \expandafter\def\csname LT0\endcsname{\color{black}}%
      \expandafter\def\csname LT1\endcsname{\color{black}}%
      \expandafter\def\csname LT2\endcsname{\color{black}}%
      \expandafter\def\csname LT3\endcsname{\color{black}}%
      \expandafter\def\csname LT4\endcsname{\color{black}}%
      \expandafter\def\csname LT5\endcsname{\color{black}}%
      \expandafter\def\csname LT6\endcsname{\color{black}}%
      \expandafter\def\csname LT7\endcsname{\color{black}}%
      \expandafter\def\csname LT8\endcsname{\color{black}}%
    \fi
  \fi
    \setlength{\unitlength}{0.0500bp}%
    \ifx\gptboxheight\undefined%
      \newlength{\gptboxheight}%
      \newlength{\gptboxwidth}%
      \newsavebox{\gptboxtext}%
    \fi%
    \setlength{\fboxrule}{0.5pt}%
    \setlength{\fboxsep}{1pt}%
    \definecolor{tbcol}{rgb}{1,1,1}%
\begin{picture}(3968.00,2834.00)%
    \gplgaddtomacro\gplbacktext{%
      \csname LTb\endcsname
      \put(571,594){\makebox(0,0)[r]{\strut{}$10^{4}$}}%
      \put(571,998){\makebox(0,0)[r]{\strut{}$10^{5}$}}%
      \put(571,1402){\makebox(0,0)[r]{\strut{}$10^{6}$}}%
      \put(571,1805){\makebox(0,0)[r]{\strut{}$10^{7}$}}%
      \put(571,2209){\makebox(0,0)[r]{\strut{}$10^{8}$}}%
      \put(571,2613){\makebox(0,0)[r]{\strut{}$10^{9}$}}%
      \put(634,418){\makebox(0,0){\strut{}$0.4$}}%
      \put(1442,418){\makebox(0,0){\strut{}$0.6$}}%
      \put(2251,418){\makebox(0,0){\strut{}$0.8$}}%
      \put(3059,418){\makebox(0,0){\strut{}$1.0$}}%
      \put(3867,418){\makebox(0,0){\strut{}$1.2$}}%
    }%
    \gplgaddtomacro\gplfronttext{%
      \csname LTb\endcsname
      \put(2626,2463){\makebox(0,0)[r]{\strut{}Ghost: Off}}%
      \csname LTb\endcsname
      \put(2626,2263){\makebox(0,0)[r]{\strut{}Ghost: On}}%
      \csname LTb\endcsname
      \put(95,1603){\rotatebox{-270.00}{\makebox(0,0){\strut{}$\kappa$}}}%
      \put(2250,154){\makebox(0,0){\strut{}$o/h_0 \ \text{(x-offset)}$}}%
    }%
    \gplbacktext
    \put(0,0){\includegraphics[width={198.40bp},height={141.70bp}]{multi_beam_cond.pdf}}%
    \gplfronttext
  \end{picture}%
\endgroup

%% file: figures/brick_wall_setup.pdf_tex
\begingroup%
  \makeatletter%
  \providecommand\color[2][]{%
    \errmessage{(Inkscape) Color is used for the text in Inkscape, but the package 'color.sty' is not loaded}%
    \renewcommand\color[2][]{}%
  }%
  \providecommand\transparent[1]{%
    \errmessage{(Inkscape) Transparency is used (non-zero) for the text in Inkscape, but the package 'transparent.sty' is not loaded}%
    \renewcommand\transparent[1]{}%
  }%
  \providecommand\rotatebox[2]{#2}%
  \newcommand*\fsize{\dimexpr\f@size pt\relax}%
  \newcommand*\lineheight[1]{\fontsize{\fsize}{#1\fsize}\selectfont}%
  \ifx\svgwidth\undefined%
    \setlength{\unitlength}{204.86033042bp}%
    \ifx\svgscale\undefined%
      \relax%
    \else%
      \setlength{\unitlength}{\unitlength * \real{\svgscale}}%
    \fi%
  \else%
    \setlength{\unitlength}{\svgwidth}%
  \fi%
  \global\let\svgwidth\undefined%
  \global\let\svgscale\undefined%
  \makeatother%
  \begin{picture}(1,0.82257388)%
    \lineheight{1}%
    \setlength\tabcolsep{0pt}%
    \put(0,0){\includegraphics[width=\unitlength,page=1]{brick_wall_setup.pdf}}%
    \put(0.25851407,0.15593779){\color[rgb]{0,0,0}\makebox(0,0)[lt]{\lineheight{1.25}\smash{\begin{tabular}[t]{l}$\delta$\end{tabular}}}}%
    \put(0.42282901,0.06911977){\color[rgb]{0,0,0}\makebox(0,0)[lt]{\lineheight{1.25}\smash{\begin{tabular}[t]{l}$\delta$\end{tabular}}}}%
    \put(0.00791556,0.50880167){\color[rgb]{0,0,0}\makebox(0,0)[lt]{\lineheight{1.25}\smash{\begin{tabular}[t]{l}$1$\end{tabular}}}}%
    \put(0.14975529,0.77799836){\color[rgb]{0,0,0}\makebox(0,0)[lt]{\lineheight{1.25}\smash{\begin{tabular}[t]{l}$1$\end{tabular}}}}%
    \put(0.62910583,0.71272328){\color[rgb]{0,0,0}\makebox(0,0)[lt]{\lineheight{1.25}\smash{\begin{tabular}[t]{l}$2$\end{tabular}}}}%
    \put(0.58800506,0.55252448){\color[rgb]{0.17254902,0.17254902,0.17254902}\makebox(0,0)[lt]{\lineheight{1.25}\smash{\begin{tabular}[t]{l}$x$\end{tabular}}}}%
    \put(0.60094808,0.64706409){\color[rgb]{0.17254902,0.17254902,0.17254902}\makebox(0,0)[lt]{\lineheight{1.25}\smash{\begin{tabular}[t]{l}$y$\end{tabular}}}}%
    \put(0.50519777,0.68440037){\color[rgb]{0.17254902,0.17254902,0.17254902}\makebox(0,0)[lt]{\lineheight{1.25}\smash{\begin{tabular}[t]{l}$z$\end{tabular}}}}%
  \end{picture}%
\endgroup%

%% file: figures/brick_wall_cut_configs.pdf_tex
\begingroup%
  \makeatletter%
  \providecommand\color[2][]{%
    \errmessage{(Inkscape) Color is used for the text in Inkscape, but the package 'color.sty' is not loaded}%
    \renewcommand\color[2][]{}%
  }%
  \providecommand\transparent[1]{%
    \errmessage{(Inkscape) Transparency is used (non-zero) for the text in Inkscape, but the package 'transparent.sty' is not loaded}%
    \renewcommand\transparent[1]{}%
  }%
  \providecommand\rotatebox[2]{#2}%
  \newcommand*\fsize{\dimexpr\f@size pt\relax}%
  \newcommand*\lineheight[1]{\fontsize{\fsize}{#1\fsize}\selectfont}%
  \ifx\svgwidth\undefined%
    \setlength{\unitlength}{142.93268039bp}%
    \ifx\svgscale\undefined%
      \relax%
    \else%
      \setlength{\unitlength}{\unitlength * \real{\svgscale}}%
    \fi%
  \else%
    \setlength{\unitlength}{\svgwidth}%
  \fi%
  \global\let\svgwidth\undefined%
  \global\let\svgscale\undefined%
  \makeatother%
  \begin{picture}(1,4.29253792)%
    \lineheight{1}%
    \setlength\tabcolsep{0pt}%
    \put(0,0){\includegraphics[width=\unitlength,page=1]{brick_wall_cut_configs.pdf}}%
    \put(0.17293884,3.315772){\color[rgb]{0,0,0}\makebox(0,0)[lt]{\lineheight{1.25}\smash{\begin{tabular}[t]{l}(A)\end{tabular}}}}%
    \put(0.17293884,2.25557709){\color[rgb]{0,0,0}\makebox(0,0)[lt]{\lineheight{1.25}\smash{\begin{tabular}[t]{l}(B)\end{tabular}}}}%
    \put(0.17293884,1.1846473){\color[rgb]{0,0,0}\makebox(0,0)[lt]{\lineheight{1.25}\smash{\begin{tabular}[t]{l}(C)\end{tabular}}}}%
    \put(0.17293884,0.10531909){\color[rgb]{0,0,0}\makebox(0,0)[lt]{\lineheight{1.25}\smash{\begin{tabular}[t]{l}(D)\end{tabular}}}}%
    \put(0,0){\includegraphics[width=\unitlength,page=2]{brick_wall_cut_configs.pdf}}%
  \end{picture}%
\endgroup%

%% file: figures/cond_brick_wall.pdf_tex
\begingroup
  \makeatletter
  \providecommand\color[2][]{%
    \GenericError{(gnuplot) \space\space\space\@spaces}{%
      Package color not loaded in conjunction with
      terminal option `colourtext'%
    }{See the gnuplot documentation for explanation.%
    }{Either use 'blacktext' in gnuplot or load the package
      color.sty in LaTeX.}%
    \renewcommand\color[2][]{}%
  }%
  \providecommand\includegraphics[2][]{%
    \GenericError{(gnuplot) \space\space\space\@spaces}{%
      Package graphicx or graphics not loaded%
    }{See the gnuplot documentation for explanation.%
    }{The gnuplot epslatex terminal needs graphicx.sty or graphics.sty.}%
    \renewcommand\includegraphics[2][]{}%
  }%
  \providecommand\rotatebox[2]{#2}%
  \@ifundefined{ifGPcolor}{%
    \newif\ifGPcolor
    \GPcolorfalse
  }{}%
  \@ifundefined{ifGPblacktext}{%
    \newif\ifGPblacktext
    \GPblacktexttrue
  }{}%
  \let\gplgaddtomacro\g@addto@macro
  \gdef\gplbacktext{}%
  \gdef\gplfronttext{}%
  \makeatother
  \ifGPblacktext
    \def\colorrgb#1{}%
    \def\colorgray#1{}%
  \else
    \ifGPcolor
      \def\colorrgb#1{\color[rgb]{#1}}%
      \def\colorgray#1{\color[gray]{#1}}%
      \expandafter\def\csname LTw\endcsname{\color{white}}%
      \expandafter\def\csname LTb\endcsname{\color{black}}%
      \expandafter\def\csname LTa\endcsname{\color{black}}%
      \expandafter\def\csname LT0\endcsname{\color[rgb]{1,0,0}}%
      \expandafter\def\csname LT1\endcsname{\color[rgb]{0,1,0}}%
      \expandafter\def\csname LT2\endcsname{\color[rgb]{0,0,1}}%
      \expandafter\def\csname LT3\endcsname{\color[rgb]{1,0,1}}%
      \expandafter\def\csname LT4\endcsname{\color[rgb]{0,1,1}}%
      \expandafter\def\csname LT5\endcsname{\color[rgb]{1,1,0}}%
      \expandafter\def\csname LT6\endcsname{\color[rgb]{0,0,0}}%
      \expandafter\def\csname LT7\endcsname{\color[rgb]{1,0.3,0}}%
      \expandafter\def\csname LT8\endcsname{\color[rgb]{0.5,0.5,0.5}}%
    \else
      \def\colorrgb#1{\color{black}}%
      \def\colorgray#1{\color[gray]{#1}}%
      \expandafter\def\csname LTw\endcsname{\color{white}}%
      \expandafter\def\csname LTb\endcsname{\color{black}}%
      \expandafter\def\csname LTa\endcsname{\color{black}}%
      \expandafter\def\csname LT0\endcsname{\color{black}}%
      \expandafter\def\csname LT1\endcsname{\color{black}}%
      \expandafter\def\csname LT2\endcsname{\color{black}}%
      \expandafter\def\csname LT3\endcsname{\color{black}}%
      \expandafter\def\csname LT4\endcsname{\color{black}}%
      \expandafter\def\csname LT5\endcsname{\color{black}}%
      \expandafter\def\csname LT6\endcsname{\color{black}}%
      \expandafter\def\csname LT7\endcsname{\color{black}}%
      \expandafter\def\csname LT8\endcsname{\color{black}}%
    \fi
  \fi
    \setlength{\unitlength}{0.0500bp}%
    \ifx\gptboxheight\undefined%
      \newlength{\gptboxheight}%
      \newlength{\gptboxwidth}%
      \newsavebox{\gptboxtext}%
    \fi%
    \setlength{\fboxrule}{0.5pt}%
    \setlength{\fboxsep}{1pt}%
    \definecolor{tbcol}{rgb}{1,1,1}%
\begin{picture}(4250.00,2834.00)%
    \gplgaddtomacro\gplbacktext{%
      \csname LTb\endcsname
      \put(600,616){\makebox(0,0)[r]{\strut{}$0$}}%
      \csname LTb\endcsname
      \put(600,1115){\makebox(0,0)[r]{\strut{}$10$}}%
      \csname LTb\endcsname
      \put(600,1615){\makebox(0,0)[r]{\strut{}$20$}}%
      \csname LTb\endcsname
      \put(600,2114){\makebox(0,0)[r]{\strut{}$30$}}%
      \csname LTb\endcsname
      \put(600,2613){\makebox(0,0)[r]{\strut{}$40$}}%
      \csname LTb\endcsname
      \put(680,440){\makebox(0,0){\strut{}$0$}}%
      \csname LTb\endcsname
      \put(1368,440){\makebox(0,0){\strut{}$0.001$}}%
      \csname LTb\endcsname
      \put(2056,440){\makebox(0,0){\strut{}$0.01$}}%
      \csname LTb\endcsname
      \put(2745,440){\makebox(0,0){\strut{}$0.1$}}%
      \csname LTb\endcsname
      \put(3433,440){\makebox(0,0){\strut{}$1$}}%
      \csname LTb\endcsname
      \put(4121,440){\makebox(0,0){\strut{}$10$}}%
    }%
    \gplgaddtomacro\gplfronttext{%
      \csname LTb\endcsname
      \put(1115,979){\makebox(0,0)[l]{\strut{}No Ghost}}%
      \csname LTb\endcsname
      \put(1115,779){\makebox(0,0)[l]{\strut{}With Ghost}}%
      \csname LTb\endcsname
      \put(207,1614){\rotatebox{-270.00}{\makebox(0,0){\strut{}$\mathrm{log}_{10}(\kappa)$}}}%
      \put(2400,154){\makebox(0,0){\strut{}$\delta/h$}}%
    }%
    \gplbacktext
    \put(0,0){\includegraphics[width={212.50bp},height={141.70bp}]{cond_brick_wall.pdf}}%
    \gplfronttext
  \end{picture}%
\endgroup

%% file: figures/nTop_sandwich_setup.pdf_tex
\begingroup%
  \makeatletter%
  \providecommand\color[2][]{%
    \errmessage{(Inkscape) Color is used for the text in Inkscape, but the package 'color.sty' is not loaded}%
    \renewcommand\color[2][]{}%
  }%
  \providecommand\transparent[1]{%
    \errmessage{(Inkscape) Transparency is used (non-zero) for the text in Inkscape, but the package 'transparent.sty' is not loaded}%
    \renewcommand\transparent[1]{}%
  }%
  \providecommand\rotatebox[2]{#2}%
  \newcommand*\fsize{\dimexpr\f@size pt\relax}%
  \newcommand*\lineheight[1]{\fontsize{\fsize}{#1\fsize}\selectfont}%
  \ifx\svgwidth\undefined%
    \setlength{\unitlength}{448.0043335bp}%
    \ifx\svgscale\undefined%
      \relax%
    \else%
      \setlength{\unitlength}{\unitlength * \real{\svgscale}}%
    \fi%
  \else%
    \setlength{\unitlength}{\svgwidth}%
  \fi%
  \global\let\svgwidth\undefined%
  \global\let\svgscale\undefined%
  \makeatother%
  \begin{picture}(1,0.38560363)%
    \lineheight{1}%
    \setlength\tabcolsep{0pt}%
    \put(0,0){\includegraphics[width=\unitlength,page=1]{nTop_sandwich_setup.pdf}}%
    \put(0.64423748,0.03544593){\color[rgb]{0,0,0}\makebox(0,0)[lt]{\lineheight{1.25}\smash{\begin{tabular}[t]{l}(A)\end{tabular}}}}%
    \put(0.42993254,0.01771356){\color[rgb]{0,0,0}\makebox(0,0)[lt]{\lineheight{1.25}\smash{\begin{tabular}[t]{l}(B)\end{tabular}}}}%
    \put(0.2433862,0.01729249){\color[rgb]{0,0,0}\makebox(0,0)[lt]{\lineheight{1.25}\smash{\begin{tabular}[t]{l}(C)\end{tabular}}}}%
  \end{picture}%
\endgroup%

%% file: figures/size_scaling_eta.pdf_tex
\begingroup
  \makeatletter
  \providecommand\color[2][]{%
    \GenericError{(gnuplot) \space\space\space\@spaces}{%
      Package color not loaded in conjunction with
      terminal option `colourtext'%
    }{See the gnuplot documentation for explanation.%
    }{Either use 'blacktext' in gnuplot or load the package
      color.sty in LaTeX.}%
    \renewcommand\color[2][]{}%
  }%
  \providecommand\includegraphics[2][]{%
    \GenericError{(gnuplot) \space\space\space\@spaces}{%
      Package graphicx or graphics not loaded%
    }{See the gnuplot documentation for explanation.%
    }{The gnuplot epslatex terminal needs graphicx.sty or graphics.sty.}%
    \renewcommand\includegraphics[2][]{}%
  }%
  \providecommand\rotatebox[2]{#2}%
  \@ifundefined{ifGPcolor}{%
    \newif\ifGPcolor
    \GPcolorfalse
  }{}%
  \@ifundefined{ifGPblacktext}{%
    \newif\ifGPblacktext
    \GPblacktexttrue
  }{}%
  \let\gplgaddtomacro\g@addto@macro
  \gdef\gplbacktext{}%
  \gdef\gplfronttext{}%
  \makeatother
  \ifGPblacktext
    \def\colorrgb#1{}%
    \def\colorgray#1{}%
  \else
    \ifGPcolor
      \def\colorrgb#1{\color[rgb]{#1}}%
      \def\colorgray#1{\color[gray]{#1}}%
      \expandafter\def\csname LTw\endcsname{\color{white}}%
      \expandafter\def\csname LTb\endcsname{\color{black}}%
      \expandafter\def\csname LTa\endcsname{\color{black}}%
      \expandafter\def\csname LT0\endcsname{\color[rgb]{1,0,0}}%
      \expandafter\def\csname LT1\endcsname{\color[rgb]{0,1,0}}%
      \expandafter\def\csname LT2\endcsname{\color[rgb]{0,0,1}}%
      \expandafter\def\csname LT3\endcsname{\color[rgb]{1,0,1}}%
      \expandafter\def\csname LT4\endcsname{\color[rgb]{0,1,1}}%
      \expandafter\def\csname LT5\endcsname{\color[rgb]{1,1,0}}%
      \expandafter\def\csname LT6\endcsname{\color[rgb]{0,0,0}}%
      \expandafter\def\csname LT7\endcsname{\color[rgb]{1,0.3,0}}%
      \expandafter\def\csname LT8\endcsname{\color[rgb]{0.5,0.5,0.5}}%
    \else
      \def\colorrgb#1{\color{black}}%
      \def\colorgray#1{\color[gray]{#1}}%
      \expandafter\def\csname LTw\endcsname{\color{white}}%
      \expandafter\def\csname LTb\endcsname{\color{black}}%
      \expandafter\def\csname LTa\endcsname{\color{black}}%
      \expandafter\def\csname LT0\endcsname{\color{black}}%
      \expandafter\def\csname LT1\endcsname{\color{black}}%
      \expandafter\def\csname LT2\endcsname{\color{black}}%
      \expandafter\def\csname LT3\endcsname{\color{black}}%
      \expandafter\def\csname LT4\endcsname{\color{black}}%
      \expandafter\def\csname LT5\endcsname{\color{black}}%
      \expandafter\def\csname LT6\endcsname{\color{black}}%
      \expandafter\def\csname LT7\endcsname{\color{black}}%
      \expandafter\def\csname LT8\endcsname{\color{black}}%
    \fi
  \fi
    \setlength{\unitlength}{0.0500bp}%
    \ifx\gptboxheight\undefined%
      \newlength{\gptboxheight}%
      \newlength{\gptboxwidth}%
      \newsavebox{\gptboxtext}%
    \fi%
    \setlength{\fboxrule}{0.5pt}%
    \setlength{\fboxsep}{1pt}%
    \definecolor{tbcol}{rgb}{1,1,1}%
\begin{picture}(4250.00,2720.00)%
    \gplgaddtomacro\gplbacktext{%
      \csname LTb\endcsname
      \put(600,616){\makebox(0,0)[r]{\strut{}$0.6$}}%
      \csname LTb\endcsname
      \put(600,1087){\makebox(0,0)[r]{\strut{}$0.8$}}%
      \csname LTb\endcsname
      \put(600,1558){\makebox(0,0)[r]{\strut{}$1.0$}}%
      \csname LTb\endcsname
      \put(600,2028){\makebox(0,0)[r]{\strut{}$1.2$}}%
      \csname LTb\endcsname
      \put(600,2499){\makebox(0,0)[r]{\strut{}$1.4$}}%
      \csname LTb\endcsname
      \put(680,440){\makebox(0,0){\strut{}$1$}}%
      \csname LTb\endcsname
      \put(1254,440){\makebox(0,0){\strut{}$2$}}%
      \csname LTb\endcsname
      \put(1827,440){\makebox(0,0){\strut{}$4$}}%
      \csname LTb\endcsname
      \put(2401,440){\makebox(0,0){\strut{}$8$}}%
      \csname LTb\endcsname
      \put(2974,440){\makebox(0,0){\strut{}$16$}}%
      \csname LTb\endcsname
      \put(3548,440){\makebox(0,0){\strut{}$32$}}%
      \csname LTb\endcsname
      \put(4122,440){\makebox(0,0){\strut{}$64$}}%
    }%
    \gplgaddtomacro\gplfronttext{%
      \csname LTb\endcsname
      \put(3686,1379){\makebox(0,0)[r]{\strut{}tessellation}}%
      \csname LTb\endcsname
      \put(3686,1179){\makebox(0,0)[r]{\strut{}enrichment}}%
      \csname LTb\endcsname
      \put(3686,979){\makebox(0,0)[r]{\strut{}ghost}}%
      \csname LTb\endcsname
      \put(3686,779){\makebox(0,0)[r]{\strut{}$\lambda'_{loc}$}}%
      \csname LTb\endcsname
      \put(141,1557){\rotatebox{-270.00}{\makebox(0,0){\strut{}$\eta$}}}%
      \put(2400,154){\makebox(0,0){\strut{}$s/s_0$}}%
    }%
    \gplbacktext
    \put(0,0){\includegraphics[width={212.50bp},height={136.00bp}]{size_scaling_eta.pdf}}%
    \gplfronttext
  \end{picture}%
\endgroup

%% file: figures/size_scaling_nu.pdf_tex
\begingroup
  \makeatletter
  \providecommand\color[2][]{%
    \GenericError{(gnuplot) \space\space\space\@spaces}{%
      Package color not loaded in conjunction with
      terminal option `colourtext'%
    }{See the gnuplot documentation for explanation.%
    }{Either use 'blacktext' in gnuplot or load the package
      color.sty in LaTeX.}%
    \renewcommand\color[2][]{}%
  }%
  \providecommand\includegraphics[2][]{%
    \GenericError{(gnuplot) \space\space\space\@spaces}{%
      Package graphicx or graphics not loaded%
    }{See the gnuplot documentation for explanation.%
    }{The gnuplot epslatex terminal needs graphicx.sty or graphics.sty.}%
    \renewcommand\includegraphics[2][]{}%
  }%
  \providecommand\rotatebox[2]{#2}%
  \@ifundefined{ifGPcolor}{%
    \newif\ifGPcolor
    \GPcolorfalse
  }{}%
  \@ifundefined{ifGPblacktext}{%
    \newif\ifGPblacktext
    \GPblacktexttrue
  }{}%
  \let\gplgaddtomacro\g@addto@macro
  \gdef\gplbacktext{}%
  \gdef\gplfronttext{}%
  \makeatother
  \ifGPblacktext
    \def\colorrgb#1{}%
    \def\colorgray#1{}%
  \else
    \ifGPcolor
      \def\colorrgb#1{\color[rgb]{#1}}%
      \def\colorgray#1{\color[gray]{#1}}%
      \expandafter\def\csname LTw\endcsname{\color{white}}%
      \expandafter\def\csname LTb\endcsname{\color{black}}%
      \expandafter\def\csname LTa\endcsname{\color{black}}%
      \expandafter\def\csname LT0\endcsname{\color[rgb]{1,0,0}}%
      \expandafter\def\csname LT1\endcsname{\color[rgb]{0,1,0}}%
      \expandafter\def\csname LT2\endcsname{\color[rgb]{0,0,1}}%
      \expandafter\def\csname LT3\endcsname{\color[rgb]{1,0,1}}%
      \expandafter\def\csname LT4\endcsname{\color[rgb]{0,1,1}}%
      \expandafter\def\csname LT5\endcsname{\color[rgb]{1,1,0}}%
      \expandafter\def\csname LT6\endcsname{\color[rgb]{0,0,0}}%
      \expandafter\def\csname LT7\endcsname{\color[rgb]{1,0.3,0}}%
      \expandafter\def\csname LT8\endcsname{\color[rgb]{0.5,0.5,0.5}}%
    \else
      \def\colorrgb#1{\color{black}}%
      \def\colorgray#1{\color[gray]{#1}}%
      \expandafter\def\csname LTw\endcsname{\color{white}}%
      \expandafter\def\csname LTb\endcsname{\color{black}}%
      \expandafter\def\csname LTa\endcsname{\color{black}}%
      \expandafter\def\csname LT0\endcsname{\color{black}}%
      \expandafter\def\csname LT1\endcsname{\color{black}}%
      \expandafter\def\csname LT2\endcsname{\color{black}}%
      \expandafter\def\csname LT3\endcsname{\color{black}}%
      \expandafter\def\csname LT4\endcsname{\color{black}}%
      \expandafter\def\csname LT5\endcsname{\color{black}}%
      \expandafter\def\csname LT6\endcsname{\color{black}}%
      \expandafter\def\csname LT7\endcsname{\color{black}}%
      \expandafter\def\csname LT8\endcsname{\color{black}}%
    \fi
  \fi
    \setlength{\unitlength}{0.0500bp}%
    \ifx\gptboxheight\undefined%
      \newlength{\gptboxheight}%
      \newlength{\gptboxwidth}%
      \newsavebox{\gptboxtext}%
    \fi%
    \setlength{\fboxrule}{0.5pt}%
    \setlength{\fboxsep}{1pt}%
    \definecolor{tbcol}{rgb}{1,1,1}%
\begin{picture}(4250.00,2720.00)%
    \gplgaddtomacro\gplbacktext{%
      \csname LTb\endcsname
      \put(600,616){\makebox(0,0)[r]{\strut{}$0.6$}}%
      \csname LTb\endcsname
      \put(600,1087){\makebox(0,0)[r]{\strut{}$0.8$}}%
      \csname LTb\endcsname
      \put(600,1558){\makebox(0,0)[r]{\strut{}$1.0$}}%
      \csname LTb\endcsname
      \put(600,2028){\makebox(0,0)[r]{\strut{}$1.2$}}%
      \csname LTb\endcsname
      \put(600,2499){\makebox(0,0)[r]{\strut{}$1.4$}}%
      \csname LTb\endcsname
      \put(680,440){\makebox(0,0){\strut{}$1$}}%
      \csname LTb\endcsname
      \put(1254,440){\makebox(0,0){\strut{}$2$}}%
      \csname LTb\endcsname
      \put(1827,440){\makebox(0,0){\strut{}$4$}}%
      \csname LTb\endcsname
      \put(2401,440){\makebox(0,0){\strut{}$8$}}%
      \csname LTb\endcsname
      \put(2974,440){\makebox(0,0){\strut{}$16$}}%
      \csname LTb\endcsname
      \put(3548,440){\makebox(0,0){\strut{}$32$}}%
      \csname LTb\endcsname
      \put(4122,440){\makebox(0,0){\strut{}$64$}}%
    }%
    \gplgaddtomacro\gplfronttext{%
      \csname LTb\endcsname
      \put(3686,1179){\makebox(0,0)[r]{\strut{}tessellation}}%
      \csname LTb\endcsname
      \put(3686,979){\makebox(0,0)[r]{\strut{}enrichment}}%
      \csname LTb\endcsname
      \put(3686,779){\makebox(0,0)[r]{\strut{}$\lambda'_{glob}$}}%
      \csname LTb\endcsname
      \put(141,1557){\rotatebox{-270.00}{\makebox(0,0){\strut{}$\mu$}}}%
      \put(2400,154){\makebox(0,0){\strut{}$s/s_0$}}%
    }%
    \gplbacktext
    \put(0,0){\includegraphics[width={212.50bp},height={136.00bp}]{size_scaling_nu.pdf}}%
    \gplfronttext
  \end{picture}%
\endgroup

%% file: figures/weak_scaling_eta.pdf_tex
\begingroup
  \makeatletter
  \providecommand\color[2][]{%
    \GenericError{(gnuplot) \space\space\space\@spaces}{%
      Package color not loaded in conjunction with
      terminal option `colourtext'%
    }{See the gnuplot documentation for explanation.%
    }{Either use 'blacktext' in gnuplot or load the package
      color.sty in LaTeX.}%
    \renewcommand\color[2][]{}%
  }%
  \providecommand\includegraphics[2][]{%
    \GenericError{(gnuplot) \space\space\space\@spaces}{%
      Package graphicx or graphics not loaded%
    }{See the gnuplot documentation for explanation.%
    }{The gnuplot epslatex terminal needs graphicx.sty or graphics.sty.}%
    \renewcommand\includegraphics[2][]{}%
  }%
  \providecommand\rotatebox[2]{#2}%
  \@ifundefined{ifGPcolor}{%
    \newif\ifGPcolor
    \GPcolorfalse
  }{}%
  \@ifundefined{ifGPblacktext}{%
    \newif\ifGPblacktext
    \GPblacktexttrue
  }{}%
  \let\gplgaddtomacro\g@addto@macro
  \gdef\gplbacktext{}%
  \gdef\gplfronttext{}%
  \makeatother
  \ifGPblacktext
    \def\colorrgb#1{}%
    \def\colorgray#1{}%
  \else
    \ifGPcolor
      \def\colorrgb#1{\color[rgb]{#1}}%
      \def\colorgray#1{\color[gray]{#1}}%
      \expandafter\def\csname LTw\endcsname{\color{white}}%
      \expandafter\def\csname LTb\endcsname{\color{black}}%
      \expandafter\def\csname LTa\endcsname{\color{black}}%
      \expandafter\def\csname LT0\endcsname{\color[rgb]{1,0,0}}%
      \expandafter\def\csname LT1\endcsname{\color[rgb]{0,1,0}}%
      \expandafter\def\csname LT2\endcsname{\color[rgb]{0,0,1}}%
      \expandafter\def\csname LT3\endcsname{\color[rgb]{1,0,1}}%
      \expandafter\def\csname LT4\endcsname{\color[rgb]{0,1,1}}%
      \expandafter\def\csname LT5\endcsname{\color[rgb]{1,1,0}}%
      \expandafter\def\csname LT6\endcsname{\color[rgb]{0,0,0}}%
      \expandafter\def\csname LT7\endcsname{\color[rgb]{1,0.3,0}}%
      \expandafter\def\csname LT8\endcsname{\color[rgb]{0.5,0.5,0.5}}%
    \else
      \def\colorrgb#1{\color{black}}%
      \def\colorgray#1{\color[gray]{#1}}%
      \expandafter\def\csname LTw\endcsname{\color{white}}%
      \expandafter\def\csname LTb\endcsname{\color{black}}%
      \expandafter\def\csname LTa\endcsname{\color{black}}%
      \expandafter\def\csname LT0\endcsname{\color{black}}%
      \expandafter\def\csname LT1\endcsname{\color{black}}%
      \expandafter\def\csname LT2\endcsname{\color{black}}%
      \expandafter\def\csname LT3\endcsname{\color{black}}%
      \expandafter\def\csname LT4\endcsname{\color{black}}%
      \expandafter\def\csname LT5\endcsname{\color{black}}%
      \expandafter\def\csname LT6\endcsname{\color{black}}%
      \expandafter\def\csname LT7\endcsname{\color{black}}%
      \expandafter\def\csname LT8\endcsname{\color{black}}%
    \fi
  \fi
    \setlength{\unitlength}{0.0500bp}%
    \ifx\gptboxheight\undefined%
      \newlength{\gptboxheight}%
      \newlength{\gptboxwidth}%
      \newsavebox{\gptboxtext}%
    \fi%
    \setlength{\fboxrule}{0.5pt}%
    \setlength{\fboxsep}{1pt}%
    \definecolor{tbcol}{rgb}{1,1,1}%
\begin{picture}(4250.00,2834.00)%
    \gplgaddtomacro\gplbacktext{%
      \csname LTb\endcsname
      \put(600,616){\makebox(0,0)[r]{\strut{}$0.0$}}%
      \csname LTb\endcsname
      \put(600,1015){\makebox(0,0)[r]{\strut{}$0.2$}}%
      \csname LTb\endcsname
      \put(600,1415){\makebox(0,0)[r]{\strut{}$0.4$}}%
      \csname LTb\endcsname
      \put(600,1814){\makebox(0,0)[r]{\strut{}$0.6$}}%
      \csname LTb\endcsname
      \put(600,2214){\makebox(0,0)[r]{\strut{}$0.8$}}%
      \csname LTb\endcsname
      \put(600,2613){\makebox(0,0)[r]{\strut{}$1.0$}}%
      \csname LTb\endcsname
      \put(680,440){\makebox(0,0){\strut{}$1$}}%
      \csname LTb\endcsname
      \put(1254,440){\makebox(0,0){\strut{}$2$}}%
      \csname LTb\endcsname
      \put(1827,440){\makebox(0,0){\strut{}$4$}}%
      \csname LTb\endcsname
      \put(2401,440){\makebox(0,0){\strut{}$8$}}%
      \csname LTb\endcsname
      \put(2974,440){\makebox(0,0){\strut{}$16$}}%
      \csname LTb\endcsname
      \put(3548,440){\makebox(0,0){\strut{}$32$}}%
      \csname LTb\endcsname
      \put(4122,440){\makebox(0,0){\strut{}$64$}}%
    }%
    \gplgaddtomacro\gplfronttext{%
      \csname LTb\endcsname
      \put(1115,1379){\makebox(0,0)[l]{\strut{}tessellation}}%
      \csname LTb\endcsname
      \put(1115,1179){\makebox(0,0)[l]{\strut{}enrichment}}%
      \csname LTb\endcsname
      \put(1115,979){\makebox(0,0)[l]{\strut{}ghost}}%
      \csname LTb\endcsname
      \put(1115,779){\makebox(0,0)[l]{\strut{}$\lambda'_{loc}$}}%
      \csname LTb\endcsname
      \put(141,1614){\rotatebox{-270.00}{\makebox(0,0){\strut{}$\eta$}}}%
      \put(2400,154){\makebox(0,0){\strut{}$s/s_0$}}%
    }%
    \gplbacktext
    \put(0,0){\includegraphics[width={212.50bp},height={141.70bp}]{weak_scaling_eta.pdf}}%
    \gplfronttext
  \end{picture}%
\endgroup

%% file: figures/weak_scaling_nu.pdf_tex
\begingroup
  \makeatletter
  \providecommand\color[2][]{%
    \GenericError{(gnuplot) \space\space\space\@spaces}{%
      Package color not loaded in conjunction with
      terminal option `colourtext'%
    }{See the gnuplot documentation for explanation.%
    }{Either use 'blacktext' in gnuplot or load the package
      color.sty in LaTeX.}%
    \renewcommand\color[2][]{}%
  }%
  \providecommand\includegraphics[2][]{%
    \GenericError{(gnuplot) \space\space\space\@spaces}{%
      Package graphicx or graphics not loaded%
    }{See the gnuplot documentation for explanation.%
    }{The gnuplot epslatex terminal needs graphicx.sty or graphics.sty.}%
    \renewcommand\includegraphics[2][]{}%
  }%
  \providecommand\rotatebox[2]{#2}%
  \@ifundefined{ifGPcolor}{%
    \newif\ifGPcolor
    \GPcolorfalse
  }{}%
  \@ifundefined{ifGPblacktext}{%
    \newif\ifGPblacktext
    \GPblacktexttrue
  }{}%
  \let\gplgaddtomacro\g@addto@macro
  \gdef\gplbacktext{}%
  \gdef\gplfronttext{}%
  \makeatother
  \ifGPblacktext
    \def\colorrgb#1{}%
    \def\colorgray#1{}%
  \else
    \ifGPcolor
      \def\colorrgb#1{\color[rgb]{#1}}%
      \def\colorgray#1{\color[gray]{#1}}%
      \expandafter\def\csname LTw\endcsname{\color{white}}%
      \expandafter\def\csname LTb\endcsname{\color{black}}%
      \expandafter\def\csname LTa\endcsname{\color{black}}%
      \expandafter\def\csname LT0\endcsname{\color[rgb]{1,0,0}}%
      \expandafter\def\csname LT1\endcsname{\color[rgb]{0,1,0}}%
      \expandafter\def\csname LT2\endcsname{\color[rgb]{0,0,1}}%
      \expandafter\def\csname LT3\endcsname{\color[rgb]{1,0,1}}%
      \expandafter\def\csname LT4\endcsname{\color[rgb]{0,1,1}}%
      \expandafter\def\csname LT5\endcsname{\color[rgb]{1,1,0}}%
      \expandafter\def\csname LT6\endcsname{\color[rgb]{0,0,0}}%
      \expandafter\def\csname LT7\endcsname{\color[rgb]{1,0.3,0}}%
      \expandafter\def\csname LT8\endcsname{\color[rgb]{0.5,0.5,0.5}}%
    \else
      \def\colorrgb#1{\color{black}}%
      \def\colorgray#1{\color[gray]{#1}}%
      \expandafter\def\csname LTw\endcsname{\color{white}}%
      \expandafter\def\csname LTb\endcsname{\color{black}}%
      \expandafter\def\csname LTa\endcsname{\color{black}}%
      \expandafter\def\csname LT0\endcsname{\color{black}}%
      \expandafter\def\csname LT1\endcsname{\color{black}}%
      \expandafter\def\csname LT2\endcsname{\color{black}}%
      \expandafter\def\csname LT3\endcsname{\color{black}}%
      \expandafter\def\csname LT4\endcsname{\color{black}}%
      \expandafter\def\csname LT5\endcsname{\color{black}}%
      \expandafter\def\csname LT6\endcsname{\color{black}}%
      \expandafter\def\csname LT7\endcsname{\color{black}}%
      \expandafter\def\csname LT8\endcsname{\color{black}}%
    \fi
  \fi
    \setlength{\unitlength}{0.0500bp}%
    \ifx\gptboxheight\undefined%
      \newlength{\gptboxheight}%
      \newlength{\gptboxwidth}%
      \newsavebox{\gptboxtext}%
    \fi%
    \setlength{\fboxrule}{0.5pt}%
    \setlength{\fboxsep}{1pt}%
    \definecolor{tbcol}{rgb}{1,1,1}%
\begin{picture}(4250.00,2550.00)%
    \gplgaddtomacro\gplbacktext{%
      \csname LTb\endcsname
      \put(600,616){\makebox(0,0)[r]{\strut{}$0.6$}}%
      \csname LTb\endcsname
      \put(600,1473){\makebox(0,0)[r]{\strut{}$0.8$}}%
      \csname LTb\endcsname
      \put(600,2329){\makebox(0,0)[r]{\strut{}$1.0$}}%
      \csname LTb\endcsname
      \put(680,440){\makebox(0,0){\strut{}$1$}}%
      \csname LTb\endcsname
      \put(1254,440){\makebox(0,0){\strut{}$2$}}%
      \csname LTb\endcsname
      \put(1827,440){\makebox(0,0){\strut{}$4$}}%
      \csname LTb\endcsname
      \put(2401,440){\makebox(0,0){\strut{}$8$}}%
      \csname LTb\endcsname
      \put(2974,440){\makebox(0,0){\strut{}$16$}}%
      \csname LTb\endcsname
      \put(3548,440){\makebox(0,0){\strut{}$32$}}%
      \csname LTb\endcsname
      \put(4122,440){\makebox(0,0){\strut{}$64$}}%
    }%
    \gplgaddtomacro\gplfronttext{%
      \csname LTb\endcsname
      \put(1115,1179){\makebox(0,0)[l]{\strut{}tessellation}}%
      \csname LTb\endcsname
      \put(1115,979){\makebox(0,0)[l]{\strut{}enrichment}}%
      \csname LTb\endcsname
      \put(1115,779){\makebox(0,0)[l]{\strut{}$\lambda'_{glob}$}}%
      \csname LTb\endcsname
      \put(141,1472){\rotatebox{-270.00}{\makebox(0,0){\strut{}$\mu$}}}%
      \put(2400,154){\makebox(0,0){\strut{}$s/s_0$}}%
    }%
    \gplbacktext
    \put(0,0){\includegraphics[width={212.50bp},height={127.50bp}]{weak_scaling_nu.pdf}}%
    \gplfronttext
  \end{picture}%
\endgroup

%% file: figures/strong_scaling_eta.pdf_tex
\begingroup
  \makeatletter
  \providecommand\color[2][]{%
    \GenericError{(gnuplot) \space\space\space\@spaces}{%
      Package color not loaded in conjunction with
      terminal option `colourtext'%
    }{See the gnuplot documentation for explanation.%
    }{Either use 'blacktext' in gnuplot or load the package
      color.sty in LaTeX.}%
    \renewcommand\color[2][]{}%
  }%
  \providecommand\includegraphics[2][]{%
    \GenericError{(gnuplot) \space\space\space\@spaces}{%
      Package graphicx or graphics not loaded%
    }{See the gnuplot documentation for explanation.%
    }{The gnuplot epslatex terminal needs graphicx.sty or graphics.sty.}%
    \renewcommand\includegraphics[2][]{}%
  }%
  \providecommand\rotatebox[2]{#2}%
  \@ifundefined{ifGPcolor}{%
    \newif\ifGPcolor
    \GPcolorfalse
  }{}%
  \@ifundefined{ifGPblacktext}{%
    \newif\ifGPblacktext
    \GPblacktexttrue
  }{}%
  \let\gplgaddtomacro\g@addto@macro
  \gdef\gplbacktext{}%
  \gdef\gplfronttext{}%
  \makeatother
  \ifGPblacktext
    \def\colorrgb#1{}%
    \def\colorgray#1{}%
  \else
    \ifGPcolor
      \def\colorrgb#1{\color[rgb]{#1}}%
      \def\colorgray#1{\color[gray]{#1}}%
      \expandafter\def\csname LTw\endcsname{\color{white}}%
      \expandafter\def\csname LTb\endcsname{\color{black}}%
      \expandafter\def\csname LTa\endcsname{\color{black}}%
      \expandafter\def\csname LT0\endcsname{\color[rgb]{1,0,0}}%
      \expandafter\def\csname LT1\endcsname{\color[rgb]{0,1,0}}%
      \expandafter\def\csname LT2\endcsname{\color[rgb]{0,0,1}}%
      \expandafter\def\csname LT3\endcsname{\color[rgb]{1,0,1}}%
      \expandafter\def\csname LT4\endcsname{\color[rgb]{0,1,1}}%
      \expandafter\def\csname LT5\endcsname{\color[rgb]{1,1,0}}%
      \expandafter\def\csname LT6\endcsname{\color[rgb]{0,0,0}}%
      \expandafter\def\csname LT7\endcsname{\color[rgb]{1,0.3,0}}%
      \expandafter\def\csname LT8\endcsname{\color[rgb]{0.5,0.5,0.5}}%
    \else
      \def\colorrgb#1{\color{black}}%
      \def\colorgray#1{\color[gray]{#1}}%
      \expandafter\def\csname LTw\endcsname{\color{white}}%
      \expandafter\def\csname LTb\endcsname{\color{black}}%
      \expandafter\def\csname LTa\endcsname{\color{black}}%
      \expandafter\def\csname LT0\endcsname{\color{black}}%
      \expandafter\def\csname LT1\endcsname{\color{black}}%
      \expandafter\def\csname LT2\endcsname{\color{black}}%
      \expandafter\def\csname LT3\endcsname{\color{black}}%
      \expandafter\def\csname LT4\endcsname{\color{black}}%
      \expandafter\def\csname LT5\endcsname{\color{black}}%
      \expandafter\def\csname LT6\endcsname{\color{black}}%
      \expandafter\def\csname LT7\endcsname{\color{black}}%
      \expandafter\def\csname LT8\endcsname{\color{black}}%
    \fi
  \fi
    \setlength{\unitlength}{0.0500bp}%
    \ifx\gptboxheight\undefined%
      \newlength{\gptboxheight}%
      \newlength{\gptboxwidth}%
      \newsavebox{\gptboxtext}%
    \fi%
    \setlength{\fboxrule}{0.5pt}%
    \setlength{\fboxsep}{1pt}%
    \definecolor{tbcol}{rgb}{1,1,1}%
\begin{picture}(4250.00,2834.00)%
    \gplgaddtomacro\gplbacktext{%
      \csname LTb\endcsname
      \put(600,616){\makebox(0,0)[r]{\strut{}$0.2$}}%
      \csname LTb\endcsname
      \put(600,1015){\makebox(0,0)[r]{\strut{}$0.4$}}%
      \csname LTb\endcsname
      \put(600,1415){\makebox(0,0)[r]{\strut{}$0.6$}}%
      \csname LTb\endcsname
      \put(600,1814){\makebox(0,0)[r]{\strut{}$0.8$}}%
      \csname LTb\endcsname
      \put(600,2214){\makebox(0,0)[r]{\strut{}$1.0$}}%
      \csname LTb\endcsname
      \put(600,2613){\makebox(0,0)[r]{\strut{}$1.2$}}%
      \csname LTb\endcsname
      \put(680,440){\makebox(0,0){\strut{}$1$}}%
      \csname LTb\endcsname
      \put(1254,440){\makebox(0,0){\strut{}$2$}}%
      \csname LTb\endcsname
      \put(1827,440){\makebox(0,0){\strut{}$4$}}%
      \csname LTb\endcsname
      \put(2401,440){\makebox(0,0){\strut{}$8$}}%
      \csname LTb\endcsname
      \put(2974,440){\makebox(0,0){\strut{}$16$}}%
      \csname LTb\endcsname
      \put(3548,440){\makebox(0,0){\strut{}$32$}}%
      \csname LTb\endcsname
      \put(4122,440){\makebox(0,0){\strut{}$64$}}%
    }%
    \gplgaddtomacro\gplfronttext{%
      \csname LTb\endcsname
      \put(1115,1379){\makebox(0,0)[l]{\strut{}tessellation}}%
      \csname LTb\endcsname
      \put(1115,1179){\makebox(0,0)[l]{\strut{}enrichment}}%
      \csname LTb\endcsname
      \put(1115,979){\makebox(0,0)[l]{\strut{}ghost}}%
      \csname LTb\endcsname
      \put(1115,779){\makebox(0,0)[l]{\strut{}$\lambda'_{loc}$}}%
      \csname LTb\endcsname
      \put(141,1614){\rotatebox{-270.00}{\makebox(0,0){\strut{}$\eta$}}}%
      \put(2400,154){\makebox(0,0){\strut{}$c/c_0$}}%
    }%
    \gplbacktext
    \put(0,0){\includegraphics[width={212.50bp},height={141.70bp}]{strong_scaling_eta.pdf}}%
    \gplfronttext
  \end{picture}%
\endgroup

%% file: figures/strong_scaling_eta_small.pdf_tex
\begingroup
  \makeatletter
  \providecommand\color[2][]{%
    \GenericError{(gnuplot) \space\space\space\@spaces}{%
      Package color not loaded in conjunction with
      terminal option `colourtext'%
    }{See the gnuplot documentation for explanation.%
    }{Either use 'blacktext' in gnuplot or load the package
      color.sty in LaTeX.}%
    \renewcommand\color[2][]{}%
  }%
  \providecommand\includegraphics[2][]{%
    \GenericError{(gnuplot) \space\space\space\@spaces}{%
      Package graphicx or graphics not loaded%
    }{See the gnuplot documentation for explanation.%
    }{The gnuplot epslatex terminal needs graphicx.sty or graphics.sty.}%
    \renewcommand\includegraphics[2][]{}%
  }%
  \providecommand\rotatebox[2]{#2}%
  \@ifundefined{ifGPcolor}{%
    \newif\ifGPcolor
    \GPcolorfalse
  }{}%
  \@ifundefined{ifGPblacktext}{%
    \newif\ifGPblacktext
    \GPblacktexttrue
  }{}%
  \let\gplgaddtomacro\g@addto@macro
  \gdef\gplbacktext{}%
  \gdef\gplfronttext{}%
  \makeatother
  \ifGPblacktext
    \def\colorrgb#1{}%
    \def\colorgray#1{}%
  \else
    \ifGPcolor
      \def\colorrgb#1{\color[rgb]{#1}}%
      \def\colorgray#1{\color[gray]{#1}}%
      \expandafter\def\csname LTw\endcsname{\color{white}}%
      \expandafter\def\csname LTb\endcsname{\color{black}}%
      \expandafter\def\csname LTa\endcsname{\color{black}}%
      \expandafter\def\csname LT0\endcsname{\color[rgb]{1,0,0}}%
      \expandafter\def\csname LT1\endcsname{\color[rgb]{0,1,0}}%
      \expandafter\def\csname LT2\endcsname{\color[rgb]{0,0,1}}%
      \expandafter\def\csname LT3\endcsname{\color[rgb]{1,0,1}}%
      \expandafter\def\csname LT4\endcsname{\color[rgb]{0,1,1}}%
      \expandafter\def\csname LT5\endcsname{\color[rgb]{1,1,0}}%
      \expandafter\def\csname LT6\endcsname{\color[rgb]{0,0,0}}%
      \expandafter\def\csname LT7\endcsname{\color[rgb]{1,0.3,0}}%
      \expandafter\def\csname LT8\endcsname{\color[rgb]{0.5,0.5,0.5}}%
    \else
      \def\colorrgb#1{\color{black}}%
      \def\colorgray#1{\color[gray]{#1}}%
      \expandafter\def\csname LTw\endcsname{\color{white}}%
      \expandafter\def\csname LTb\endcsname{\color{black}}%
      \expandafter\def\csname LTa\endcsname{\color{black}}%
      \expandafter\def\csname LT0\endcsname{\color{black}}%
      \expandafter\def\csname LT1\endcsname{\color{black}}%
      \expandafter\def\csname LT2\endcsname{\color{black}}%
      \expandafter\def\csname LT3\endcsname{\color{black}}%
      \expandafter\def\csname LT4\endcsname{\color{black}}%
      \expandafter\def\csname LT5\endcsname{\color{black}}%
      \expandafter\def\csname LT6\endcsname{\color{black}}%
      \expandafter\def\csname LT7\endcsname{\color{black}}%
      \expandafter\def\csname LT8\endcsname{\color{black}}%
    \fi
  \fi
    \setlength{\unitlength}{0.0500bp}%
    \ifx\gptboxheight\undefined%
      \newlength{\gptboxheight}%
      \newlength{\gptboxwidth}%
      \newsavebox{\gptboxtext}%
    \fi%
    \setlength{\fboxrule}{0.5pt}%
    \setlength{\fboxsep}{1pt}%
    \definecolor{tbcol}{rgb}{1,1,1}%
\begin{picture}(4250.00,2834.00)%
    \gplgaddtomacro\gplbacktext{%
      \csname LTb\endcsname
      \put(600,616){\makebox(0,0)[r]{\strut{}$0.0$}}%
      \csname LTb\endcsname
      \put(600,949){\makebox(0,0)[r]{\strut{}$0.2$}}%
      \csname LTb\endcsname
      \put(600,1282){\makebox(0,0)[r]{\strut{}$0.4$}}%
      \csname LTb\endcsname
      \put(600,1615){\makebox(0,0)[r]{\strut{}$0.6$}}%
      \csname LTb\endcsname
      \put(600,1947){\makebox(0,0)[r]{\strut{}$0.8$}}%
      \csname LTb\endcsname
      \put(600,2280){\makebox(0,0)[r]{\strut{}$1.0$}}%
      \csname LTb\endcsname
      \put(600,2613){\makebox(0,0)[r]{\strut{}$1.2$}}%
      \csname LTb\endcsname
      \put(680,440){\makebox(0,0){\strut{}$1$}}%
      \csname LTb\endcsname
      \put(1254,440){\makebox(0,0){\strut{}$2$}}%
      \csname LTb\endcsname
      \put(1827,440){\makebox(0,0){\strut{}$4$}}%
      \csname LTb\endcsname
      \put(2401,440){\makebox(0,0){\strut{}$8$}}%
      \csname LTb\endcsname
      \put(2974,440){\makebox(0,0){\strut{}$16$}}%
      \csname LTb\endcsname
      \put(3548,440){\makebox(0,0){\strut{}$32$}}%
      \csname LTb\endcsname
      \put(4122,440){\makebox(0,0){\strut{}$64$}}%
    }%
    \gplgaddtomacro\gplfronttext{%
      \csname LTb\endcsname
      \put(1115,1379){\makebox(0,0)[l]{\strut{}tessellation}}%
      \csname LTb\endcsname
      \put(1115,1179){\makebox(0,0)[l]{\strut{}enrichment}}%
      \csname LTb\endcsname
      \put(1115,979){\makebox(0,0)[l]{\strut{}ghost}}%
      \csname LTb\endcsname
      \put(1115,779){\makebox(0,0)[l]{\strut{}$\lambda^{\prime}_{loc}$}}%
      \csname LTb\endcsname
      \put(141,1614){\rotatebox{-270.00}{\makebox(0,0){\strut{}$\eta$}}}%
      \put(2400,154){\makebox(0,0){\strut{}$c/c_0$}}%
    }%
    \gplbacktext
    \put(0,0){\includegraphics[width={212.50bp},height={141.70bp}]{strong_scaling_eta_small.pdf}}%
    \gplfronttext
  \end{picture}%
\endgroup

%% file: figures/fg_mesh.pdf_tex
\begingroup%
  \makeatletter%
  \providecommand\color[2][]{%
    \errmessage{(Inkscape) Color is used for the text in Inkscape, but the package 'color.sty' is not loaded}%
    \renewcommand\color[2][]{}%
  }%
  \providecommand\transparent[1]{%
    \errmessage{(Inkscape) Transparency is used (non-zero) for the text in Inkscape, but the package 'transparent.sty' is not loaded}%
    \renewcommand\transparent[1]{}%
  }%
  \providecommand\rotatebox[2]{#2}%
  \newcommand*\fsize{\dimexpr\f@size pt\relax}%
  \newcommand*\lineheight[1]{\fontsize{\fsize}{#1\fsize}\selectfont}%
  \ifx\svgwidth\undefined%
    \setlength{\unitlength}{444bp}%
    \ifx\svgscale\undefined%
      \relax%
    \else%
      \setlength{\unitlength}{\unitlength * \real{\svgscale}}%
    \fi%
  \else%
    \setlength{\unitlength}{\svgwidth}%
  \fi%
  \global\let\svgwidth\undefined%
  \global\let\svgscale\undefined%
  \makeatother%
  \begin{picture}(1,1.31587838)%
    \lineheight{1}%
    \setlength\tabcolsep{0pt}%
    \put(0,0){\includegraphics[width=\unitlength,page=1]{fg_mesh.pdf}}%
    \put(0.24493243,1.28547297){\makebox(0,0)[t]{\lineheight{1.25}\smash{\begin{tabular}[t]{c}\textbf{Background Element }$\BgElemInd$\end{tabular}}}}%
    \put(0,0){\includegraphics[width=\unitlength,page=2]{fg_mesh.pdf}}%
    \put(0.01013514,1.24155405){\makebox(0,0)[lt]{\lineheight{1.25}\smash{\begin{tabular}[t]{l}+ Array of vertices $\mathrm{V}$ (sorted by ordinal)\end{tabular}}}}%
    \put(0,0){\includegraphics[width=\unitlength,page=3]{fg_mesh.pdf}}%
    \put(0.01013514,1.19763514){\makebox(0,0)[lt]{\lineheight{1.25}\smash{\begin{tabular}[t]{l}+ Basis $\mathcal{B}_E = \SetDef{N_{\LocalBfIndex}(\xi)}_{\LocalBfIndex=1}^{n_{\LocalBfIndex}(E)}$\end{tabular}}}}%
    \put(0,0){\includegraphics[width=\unitlength,page=4]{fg_mesh.pdf}}%
    \put(0.01013514,1.15371622){\makebox(0,0)[lt]{\lineheight{1.25}\smash{\begin{tabular}[t]{l}+ Local to global basis fnct. map $\mathrm{IEN}: \LocalBfIndex\mapsto\BfIndex$\end{tabular}}}}%
    \put(0,0){\includegraphics[width=\unitlength,page=5]{fg_mesh.pdf}}%
    \put(0.01013514,1.11148649){\makebox(0,0)[lt]{\lineheight{1.25}\smash{\begin{tabular}[t]{l}+ Parallel ID $I$, owning processor $p$\end{tabular}}}}%
    \put(0,0){\includegraphics[width=\unitlength,page=6]{fg_mesh.pdf}}%
    \put(0.01013514,1.05743243){\makebox(0,0)[lt]{\lineheight{1.25}\smash{\begin{tabular}[t]{l}+ $\pos,\pfrac{\pos}{\xi},\cdots \leftarrow \mathtt{interpolate\_space}(\xi)$\end{tabular}}}}%
    \put(0,0){\includegraphics[width=\unitlength,page=7]{fg_mesh.pdf}}%
    \put(0.75,1.28547297){\makebox(0,0)[t]{\lineheight{1.25}\smash{\begin{tabular}[t]{c}\textbf{Cluster }$\Cluster$\end{tabular}}}}%
    \put(0,0){\includegraphics[width=\unitlength,page=8]{fg_mesh.pdf}}%
    \put(0.55743243,1.24155405){\makebox(0,0)[lt]{\lineheight{1.25}\smash{\begin{tabular}[t]{l}+ Background element $E$\end{tabular}}}}%
    \put(0,0){\includegraphics[width=\unitlength,page=9]{fg_mesh.pdf}}%
    \put(0.55743243,1.19763514){\makebox(0,0)[lt]{\lineheight{1.25}\smash{\begin{tabular}[t]{l}+ Quadrature points $\SetDef{(\xi_q,\omega_q)}_{q=1}^{n_q(E)}$\end{tabular}}}}%
    \put(0,0){\includegraphics[width=\unitlength,page=10]{fg_mesh.pdf}}%
    \put(0.55743243,1.15033784){\color[rgb]{0,0,0}\makebox(0,0)[lt]{\lineheight{1.25}\smash{\begin{tabular}[t]{l}+ For side clusters: normals associated \end{tabular}}}}%
    \put(0.55507382,1.12241423){\color[rgb]{0,0,0}\makebox(0,0)[lt]{\lineheight{1.25}\smash{\begin{tabular}[t]{l}    with quadrature points $\SetDef{\normal_q}_{q=1}^{n_q(E)}$\end{tabular}}}}%
    \put(0,0){\includegraphics[width=\unitlength,page=11]{fg_mesh.pdf}}%
    \put(0.21452703,0.98648649){\makebox(0,0)[t]{\lineheight{1.25}\smash{\begin{tabular}[t]{c}\textbf{Background Mesh }$\BgMesh$\end{tabular}}}}%
    \put(0,0){\includegraphics[width=\unitlength,page=12]{fg_mesh.pdf}}%
    \put(0.04391892,0.94087838){\makebox(0,0)[lt]{\lineheight{1.25}\smash{\begin{tabular}[t]{l}+ Array of elements $\mathcal{E} = \SetDef{E}_{E=1}^{n_E}$ \end{tabular}}}}%
    \put(0,0){\includegraphics[width=\unitlength,page=13]{fg_mesh.pdf}}%
    \put(0.04391892,0.88513514){\makebox(0,0)[lt]{\lineheight{1.25}\smash{\begin{tabular}[t]{l}+ $\mathtt{get\_number\_of\_basis\_fncts}()$\end{tabular}}}}%
    \put(0,0){\includegraphics[width=\unitlength,page=14]{fg_mesh.pdf}}%
    \put(0.04391892,0.84121622){\makebox(0,0)[lt]{\lineheight{1.25}\smash{\begin{tabular}[t]{l}+ $\mathtt{get\_elements\_in\_support}(\BfIndex)$\end{tabular}}}}%
    \put(0,0){\includegraphics[width=\unitlength,page=15]{fg_mesh.pdf}}%
    \put(0.04391892,0.7972973){\makebox(0,0)[lt]{\lineheight{1.25}\smash{\begin{tabular}[t]{l}+ $\mathtt{get\_entities\_on\_element}(r,E)$\end{tabular}}}}%
    \put(0,0){\includegraphics[width=\unitlength,page=16]{fg_mesh.pdf}}%
    \put(0.04391892,0.75337838){\makebox(0,0)[lt]{\lineheight{1.25}\smash{\begin{tabular}[t]{l}+ $\mathtt{get\_elements\_on\_entity}(r,e)$\end{tabular}}}}%
    \put(0,0){\includegraphics[width=\unitlength,page=17]{fg_mesh.pdf}}%
    \put(0.04391892,0.70945946){\makebox(0,0)[lt]{\lineheight{1.25}\smash{\begin{tabular}[t]{l}+ $\mathtt{get\_basis\_fnct\_ID}(\BfIndex)$\end{tabular}}}}%
    \put(0,0){\includegraphics[width=\unitlength,page=18]{fg_mesh.pdf}}%
    \put(0.21114865,0.24324324){\makebox(0,0)[t]{\lineheight{1.25}\smash{\begin{tabular}[t]{c}\textbf{Entity Connectivity }$\FacetConn$\end{tabular}}}}%
    \put(0,0){\includegraphics[width=\unitlength,page=19]{fg_mesh.pdf}}%
    \put(0.01013514,0.19932432){\makebox(0,0)[lt]{\lineheight{1.25}\smash{\begin{tabular}[t]{l}+ Entity to Cell $\mathrm{EtC} = \SetDef{\SetDef{\cdots}}_{e=1}^{\mathtt{NumEntities}}$\end{tabular}}}}%
    \put(0,0){\includegraphics[width=\unitlength,page=20]{fg_mesh.pdf}}%
    \put(0.01013514,0.15540541){\color[rgb]{0,0,0}\makebox(0,0)[lt]{\lineheight{1.25}\smash{\begin{tabular}[t]{l}+ Cell to Entity $\mathrm{CtE} = \SetDef{\SetDef{\cdots}}_{c=1}^{\mathtt{NumCells}}$\end{tabular}}}}%
    \put(0.01076055,0.12862963){\color[rgb]{0,0,0}\makebox(0,0)[lt]{\lineheight{1.25}\smash{\begin{tabular}[t]{l}   (entities sorted by ordinal)\end{tabular}}}}%
    \put(0,0){\includegraphics[width=\unitlength,page=21]{fg_mesh.pdf}}%
    \put(0.69256757,1.03547297){\makebox(0,0)[t]{\lineheight{1.25}\smash{\begin{tabular}[t]{c}\textbf{Foreground Mesh }$\FgMesh$\end{tabular}}}}%
    \put(0,0){\includegraphics[width=\unitlength,page=22]{fg_mesh.pdf}}%
    \put(0.53378378,0.99155405){\makebox(0,0)[lt]{\lineheight{1.25}\smash{\begin{tabular}[t]{l}+ Background mesh $\mathcal{H}$\end{tabular}}}}%
    \put(0,0){\includegraphics[width=\unitlength,page=23]{fg_mesh.pdf}}%
    \put(0.53378378,0.94763514){\makebox(0,0)[lt]{\lineheight{1.25}\smash{\begin{tabular}[t]{l}+ Array of fg. vertices $\mathcal{V}$\end{tabular}}}}%
    \put(0,0){\includegraphics[width=\unitlength,page=24]{fg_mesh.pdf}}%
    \put(0.53378378,0.90371622){\makebox(0,0)[lt]{\lineheight{1.25}\smash{\begin{tabular}[t]{l}+ Array of fg. cells $\mathcal{C}$\end{tabular}}}}%
    \put(0,0){\includegraphics[width=\unitlength,page=25]{fg_mesh.pdf}}%
    \put(0.53378378,0.8597973){\makebox(0,0)[lt]{\lineheight{1.25}\smash{\begin{tabular}[t]{l}+ Array of child meshes $\mathrm{CMS}$\end{tabular}}}}%
    \put(0,0){\includegraphics[width=\unitlength,page=26]{fg_mesh.pdf}}%
    \put(0.53378378,0.81587838){\makebox(0,0)[lt]{\lineheight{1.25}\smash{\begin{tabular}[t]{l}+ Array of Subphases $\mathcal{S}$\end{tabular}}}}%
    \put(0,0){\includegraphics[width=\unitlength,page=27]{fg_mesh.pdf}}%
    \put(0.53378378,0.77195946){\makebox(0,0)[lt]{\lineheight{1.25}\smash{\begin{tabular}[t]{l}+ Subphase graph $\mathcal{G}_S$\end{tabular}}}}%
    \put(0,0){\includegraphics[width=\unitlength,page=28]{fg_mesh.pdf}}%
    \put(0.53378378,0.72804054){\color[rgb]{0,0,0}\makebox(0,0)[lt]{\lineheight{1.25}\smash{\begin{tabular}[t]{l}+ Graph of disconnected\end{tabular}}}}%
    \put(0.53461766,0.70105635){\color[rgb]{0,0,0}\makebox(0,0)[lt]{\lineheight{1.25}\smash{\begin{tabular}[t]{l}   subphase neighbors $\mathcal{G}_I$\end{tabular}}}}%
    \put(0,0){\includegraphics[width=\unitlength,page=29]{fg_mesh.pdf}}%
    \put(0.53378378,0.66385135){\makebox(0,0)[lt]{\lineheight{1.25}\smash{\begin{tabular}[t]{l}+ Facet connectivity $\mathcal{F}$\end{tabular}}}}%
    \put(0,0){\includegraphics[width=\unitlength,page=30]{fg_mesh.pdf}}%
    \put(0.90540541,0.55912162){\makebox(0,0)[t]{\lineheight{1.25}\smash{\begin{tabular}[t]{c}\textbf{Foreground}\end{tabular}}}}%
    \put(0.90540541,0.53547297){\makebox(0,0)[t]{\lineheight{1.25}\smash{\begin{tabular}[t]{c}\textbf{Vertex }$\Vertex$\end{tabular}}}}%
    \put(0,0){\includegraphics[width=\unitlength,page=31]{fg_mesh.pdf}}%
    \put(0.82263514,0.48986486){\makebox(0,0)[lt]{\lineheight{1.25}\smash{\begin{tabular}[t]{l}+ Coordinates $\pos$\end{tabular}}}}%
    \put(0,0){\includegraphics[width=\unitlength,page=32]{fg_mesh.pdf}}%
    \put(0.82263514,0.45439189){\color[rgb]{0,0,0}\makebox(0,0)[lt]{\lineheight{1.25}\smash{\begin{tabular}[t]{l}+ Bg. ancestor\end{tabular}}}}%
    \put(0.81968689,0.42897425){\color[rgb]{0,0,0}\makebox(0,0)[lt]{\lineheight{1.25}\smash{\begin{tabular}[t]{l}    rank $r$\end{tabular}}}}%
    \put(0,0){\includegraphics[width=\unitlength,page=33]{fg_mesh.pdf}}%
    \put(0.82263514,0.39864865){\color[rgb]{0,0,0}\makebox(0,0)[lt]{\lineheight{1.25}\smash{\begin{tabular}[t]{l}+ Bg. ancestor\end{tabular}}}}%
    \put(0.81909726,0.37234662){\color[rgb]{0,0,0}\makebox(0,0)[lt]{\lineheight{1.25}\smash{\begin{tabular}[t]{l}    index $a$\end{tabular}}}}%
    \put(0,0){\includegraphics[width=\unitlength,page=34]{fg_mesh.pdf}}%
    \put(0.62922297,0.59797297){\makebox(0,0)[t]{\lineheight{1.25}\smash{\begin{tabular}[t]{c}\textbf{Foreground Cell }$\Cell$\end{tabular}}}}%
    \put(0,0){\includegraphics[width=\unitlength,page=35]{fg_mesh.pdf}}%
    \put(0.48310811,0.55405405){\makebox(0,0)[lt]{\lineheight{1.25}\smash{\begin{tabular}[t]{l}+ Array of vertices $\mathrm{V}$\end{tabular}}}}%
    \put(0,0){\includegraphics[width=\unitlength,page=36]{fg_mesh.pdf}}%
    \put(0.48310811,0.51351351){\makebox(0,0)[lt]{\lineheight{1.25}\smash{\begin{tabular}[t]{l}+ Material $m$\end{tabular}}}}%
    \put(0,0){\includegraphics[width=\unitlength,page=37]{fg_mesh.pdf}}%
    \put(0.48310811,0.47466216){\makebox(0,0)[lt]{\lineheight{1.25}\smash{\begin{tabular}[t]{l}+ Bg. element membership $E$\end{tabular}}}}%
    \put(0,0){\includegraphics[width=\unitlength,page=38]{fg_mesh.pdf}}%
    \put(0.48310811,0.43581081){\makebox(0,0)[lt]{\lineheight{1.25}\smash{\begin{tabular}[t]{l}+ Subphase membership $S$\end{tabular}}}}%
    \put(0,0){\includegraphics[width=\unitlength,page=39]{fg_mesh.pdf}}%
    \put(0.48310811,0.39527027){\makebox(0,0)[lt]{\lineheight{1.25}\smash{\begin{tabular}[t]{l}+ Parallel ID $I$\end{tabular}}}}%
    \put(0,0){\includegraphics[width=\unitlength,page=40]{fg_mesh.pdf}}%
    \put(0.18581081,0.64864865){\makebox(0,0)[t]{\lineheight{1.25}\smash{\begin{tabular}[t]{c}\textbf{Child Mesh }$\ChildMesh$\end{tabular}}}}%
    \put(0,0){\includegraphics[width=\unitlength,page=41]{fg_mesh.pdf}}%
    \put(0.01013514,0.60472973){\makebox(0,0)[lt]{\lineheight{1.25}\smash{\begin{tabular}[t]{l}+ List of fg. cells $\mathrm{C}$\end{tabular}}}}%
    \put(0,0){\includegraphics[width=\unitlength,page=42]{fg_mesh.pdf}}%
    \put(0.01013514,0.56081081){\makebox(0,0)[lt]{\lineheight{1.25}\smash{\begin{tabular}[t]{l}+ List of fg. vertices $\mathrm{V}$\end{tabular}}}}%
    \put(0,0){\includegraphics[width=\unitlength,page=43]{fg_mesh.pdf}}%
    \put(0.01013514,0.52027027){\makebox(0,0)[lt]{\lineheight{1.25}\smash{\begin{tabular}[t]{l}+ Vertex param. coords. $\SetDef{\xi_v}_{v=1}^{\mathtt{size}(\mathrm{V})}$\end{tabular}}}}%
    \put(0,0){\includegraphics[width=\unitlength,page=44]{fg_mesh.pdf}}%
    \put(0.16047297,0.46283784){\makebox(0,0)[t]{\lineheight{1.25}\smash{\begin{tabular}[t]{c}\textbf{Subphase }$\Subphase$\end{tabular}}}}%
    \put(0,0){\includegraphics[width=\unitlength,page=45]{fg_mesh.pdf}}%
    \put(0.01013514,0.41891892){\makebox(0,0)[lt]{\lineheight{1.25}\smash{\begin{tabular}[t]{l}+ Bg. element membership $E$\end{tabular}}}}%
    \put(0,0){\includegraphics[width=\unitlength,page=46]{fg_mesh.pdf}}%
    \put(0.01013514,0.38175676){\makebox(0,0)[lt]{\lineheight{1.25}\smash{\begin{tabular}[t]{l}+ Ordinal $u$\end{tabular}}}}%
    \put(0,0){\includegraphics[width=\unitlength,page=47]{fg_mesh.pdf}}%
    \put(0.01013514,0.34459459){\makebox(0,0)[lt]{\lineheight{1.25}\smash{\begin{tabular}[t]{l}+ Array of fg. cells $\mathrm{C}$\end{tabular}}}}%
    \put(0,0){\includegraphics[width=\unitlength,page=48]{fg_mesh.pdf}}%
    \put(0.01013514,0.30574324){\makebox(0,0)[lt]{\lineheight{1.25}\smash{\begin{tabular}[t]{l}+ Parallel ID $I$\end{tabular}}}}%
    \put(0,0){\includegraphics[width=\unitlength,page=49]{fg_mesh.pdf}}%
    \put(0.73648649,0.31081081){\makebox(0,0)[t]{\lineheight{1.25}\smash{\begin{tabular}[t]{c}\textbf{Queue }$\Queue$\end{tabular}}}}%
    \put(0,0){\includegraphics[width=\unitlength,page=50]{fg_mesh.pdf}}%
    \put(0.4847973,0.26689189){\color[rgb]{0,0,0}\makebox(0,0)[lt]{\lineheight{1.25}\smash{\begin{tabular}[t]{l}+ Map $\mathrm{M}: I_r \mapsto i$ \end{tabular}}}}%
    \put(0.48450248,0.23734675){\color[rgb]{0,0,0}\makebox(0,0)[lt]{\lineheight{1.25}\smash{\begin{tabular}[t]{l}   (maps request identifier $I_r$ to position $i$ in queue)\end{tabular}}}}%
    \put(0,0){\includegraphics[width=\unitlength,page=51]{fg_mesh.pdf}}%
    \put(0.4847973,0.20101351){\makebox(0,0)[lt]{\lineheight{1.25}\smash{\begin{tabular}[t]{l}+ Request queue: array of initializer lists\end{tabular}}}}%
    \put(0,0){\includegraphics[width=\unitlength,page=52]{fg_mesh.pdf}}%
    \put(0.4847973,0.15878378){\makebox(0,0)[lt]{\lineheight{1.25}\smash{\begin{tabular}[t]{l}+ Maximum existing entity index $e_{max}$\end{tabular}}}}%
    \put(0,0){\includegraphics[width=\unitlength,page=53]{fg_mesh.pdf}}%
    \put(0.4847973,0.09966216){\makebox(0,0)[lt]{\lineheight{1.25}\smash{\begin{tabular}[t]{l}+ $\mathtt{request\_exists}(I_r)$\end{tabular}}}}%
    \put(0,0){\includegraphics[width=\unitlength,page=54]{fg_mesh.pdf}}%
    \put(0.4847973,0.0625){\makebox(0,0)[lt]{\lineheight{1.25}\smash{\begin{tabular}[t]{l}+ $\mathtt{get\_index}(I_r)$\end{tabular}}}}%
    \put(0,0){\includegraphics[width=\unitlength,page=55]{fg_mesh.pdf}}%
    \put(0.4847973,0.02533784){\makebox(0,0)[lt]{\lineheight{1.25}\smash{\begin{tabular}[t]{l}+ $\mathtt{queue}(I_r, \{\text{Initializer List}\} )$\end{tabular}}}}%
    \put(0,0){\includegraphics[width=\unitlength,page=56]{fg_mesh.pdf}}%
  \end{picture}%
\endgroup%